\newif\iffigs\figstrue
\documentclass[10pt,a4paper]{article}
\usepackage{latexsym,amssymb,lscape,graphics,setspace}
\usepackage{amsmath}
\usepackage{graphicx}        % pacchetto grafico standard LaTeX  per l'inclusione di file immagine
\usepackage{longtable}
\usepackage[vcentermath]{youngtab}
\usepackage{multirow}
\usepackage{booktabs}
\usepackage{multirow}
\usepackage{color}
\usepackage{slashed,epsfig}
\usepackage{amsfonts}
\usepackage{cite}
\usepackage{verbatim}
\usepackage{colortbl}
\usepackage[table]{xcolor}
\usepackage{hyperref}       % questo pacchetto e'  quello che produce un ipertesto
\usepackage{array}
\usepackage{fancyhdr}
\usepackage{tikz}
\usepackage{multicol}
\usepackage{todonotes}
\makeatletter
\renewcommand{\paragraph}{%
  \@startsection{paragraph}{4}%
  {\z@}{1ex \@plus 1ex \@minus .2ex}{-1em}%
  {\normalfont\normalsize\bfseries}%
}
\makeatother

\newcommand{\mathsym}[1]{{}}

\newtheorem{definizione}{Definition}[section]
\newtheorem{teorema}{Theorem}[section]
\newtheorem{proposizione}{Proposition}[section]

\newtheorem{statement}{Statement}[section]

\newcommand{\bd}{\begin{definizione}}
\newcommand{\ed}{\end{definizione}}

\def\IC{\relax\,\hbox{$\inbar\kern-.3em{\rm C}$}}
\def\IG{\relax\,\hbox{$\inbar\kern-.3em{\rm G}$}}
\def\IB{\relax{\rm I\kern-.18em B}}
\def\ID{\relax{\rm I\kern-.18em D}}
\def\IL{\relax{\rm I\kern-.18em L}}
\def\IF{\relax{\rm I\kern-.18em F}}
\def\IH{\relax{\rm I\kern-.18em H}}
\def\II{\relax{\rm I\kern-.17em I}}
\def\IN{\relax{\rm I\kern-.18em N}}
\def\IP{\relax{\rm I\kern-.18em P}}
\def\IQ{\relax\,\hbox{$\inbar\kern-.3em{\rm Q}$}}
\def\bfzero{\relax\,\hbox{$\inbar\kern-.3em{\rm 0}$}}
\def\IK{\relax{\rm I\kern-.18em K}}
\def\IG{\relax\,\hbox{$\inbar\kern-.3em{\rm G}$}}
 \font\cmss=cmss10 \font\cmsss=cmss10 at 7pt
\def\IR{\relax{\rm I\kern-.18em R}}
\def\ZZ{\relax\ifmmode\mathchoice
{\hbox{\cmss Z\kern-.4em Z}}{\hbox{\cmss Z\kern-.4em Z}}
{\lower.9pt\hbox{\cmsss Z\kern-.4em Z}} {\lower1.2pt\hbox{\cmsss
Z\kern-.4em Z}}\else{\cmss Z\kern-.4em Z}\fi}
\def\bfone{\relax{\rm 1\kern-.35em 1}}

\def\diag{{\rm diag}}

\def\Solv{\mathop{\rm Solv}\nolimits}

\def\inbar{\vrule height1.5ex width.4pt depth0pt}
\def\bfzero{\relax{\rm I\kern-.18em 0}}
\def\bfone{\relax{\rm 1\kern-.35em 1}}

\def\twovec#1#2{\left(\begin{array}{c}
{#1}\nonumber \\ {#2}\nonumber \\
\end{array}
\right)}
\def\o#1#2{{{#1}\over{#2}}}
\DeclareFontFamily{U}{rsf}{} \DeclareFontShape{U}{rsf}{m}{n}{
  <5> <6> rsfs5 <7> <8> <9> rsfs7 <10-> rsfs10}{}
\DeclareMathAlphabet\Scr{U}{rsf}{m}{n}

\def\T{T}

\newcommand{\Sp}{\mathop{\rm {}Sp}}

\newcommand{\ft}[2]{{\textstyle\frac{#1}{#2}}}
\def\tilde{\widetilde}

\def\1bar{1\hskip -.275cm -}
\def\2bar{2\hskip -.275cm -}
\def\3bar{3\hskip -.275cm -}

\newsavebox{\uuunit}
\sbox{\uuunit}
                 {\setlength{\unitlength}{0.825em}
                      \begin{picture}(0.6,0.7)
                                      \thinlines
                                      \put(0,0){\line(1,0){0.5}}
                                      \put(0.15,0){\line(0,1){0.7}}
                                      \put(0.35,0){\line(0,1){0.8}}
                                     \multiput(0.3,0.8)(-0.04,-0.02){10}{\rule{0.5pt}{0.5pt}}
                      \end {picture}}

\makeatletter \@addtoreset{equation}{section} \makeatother

\def\bfone{\relax{\rm 1\kern-.35em 1}}

\def\bfone{\relax{\rm 1\kern-.35em 1}}
\font\cmss=cmss10 \font\cmsss=cmss10 at 7pt

\newcommand{\so}{\mathfrak{so}}
\newcommand{\su}{\mathfrak{su}}
\newcommand{\spalg}{\mathfrak{sp}}
\newcommand{\usp}{\mathfrak{usp}}

\newcommand{\sym}{\mathfrak{sp}}
\newcommand{\slal}{\mathfrak{sl}}

\begin{document}
\begin{titlepage}
\begin{center}
{{\large {\sc  The PAINT GROUP TITS-SATAKE THEORY\\[5pt] of HYPERBOLIC SYMMETRIC SPACES:}}\\[5pt]
   DISTANCE FUNCTION, PAINT INVARIANTS and DISCRETE SUBGROUPS${}^\dagger$}
 \par\bigskip{\sc Ugo Bruzzo$^{a,b}$, Pietro Fr\'e\,$^{c,d}$,  and Mario Trigiante$^{e,f}$}  \par
\smallskip
{\sl \small
${}^a\,$
Departamento de Matem\'atica, Instituto de Ci\^encias Exatas,  Universidade Federal \\ de Minas Gerais,
Av.~Ant\^onio Carlos 6627,   Belo Horizonte, MG  31270-901, Brazil
 \\[2pt]
${}^b\,$IGAP (Institute for Geometry and Physics), Trieste, Italy \\[2pt]
${}^c\,$ {Emeritus at} Dipartimento di Fisica, Universit\`a di Torino, Via P. Giuria 1, I-10125 Torino, Italy \\[2pt]
${}^{d}\,${Consultant at} Additati\&Partners Consulting s.r.l, Via Filippo Pacini 36, I-51100 Pistoia, Italy \\[2pt]
${}^e\,$Dipartimento DISAT, Politecnico di Torino,
C.so Duca degli Abruzzi 24, I-10129 Torino, Italy\\[2pt]
${}^f\,$INFN, Sezione di Torino\\[2pt]
E-mail:  {\tt ubruzzo@ufmg.br, pietro.fre@unito.it,
mario.trigiante@polito.it} }
\begin{abstract}
The present paper, which is partly a review and turned out to be the corner stone of an entire new development in geometrical deep learning, aims at presenting, in a unified
mathematical framework, a complex and articulated lore regarding noncompact symmetric spaces $\mathrm{U/H}$, where $\mathrm{U}$ is a simple Lie
group, whose Lie algebra $\mathbb{U}$, is a noncompact real section of a complex simple Lie algebra and $\mathrm{H}\subset\mathrm{U}$ is the
maximal compact Lie subgroup. All such manifolds are \emph{Riemannian normal manifolds}, according to Alekseevsky's definition, in the sense
that they are metrically equivalent to a solvable Lie group manifold $\mathcal{S}_{\mathrm{U/H}}$. This identification provides a vision in
which, on one side one can derive quite explicit and challenging formulas for the unique distance function between points of the manifold, on
the other one, one is able to organize the entire set of the available manifolds in universality classes distinguished by their common
{\em Tits-Satake submanifold,} and, correspondingly, by their noncompact rank. The members of the class are distinguished by their different {\em paint group,}
this notion having been established by two of the present authors and their collaborators in the years 2007-2009. In relation with the
construction of neural networks, these mathematical structures offer unique possibilities for replacing \emph{ad hoc activation functions} with
the naturally defined nonlinear operations that relate Lie algebras to Lie Groups and viceversa. The paint group invariants offer new tokens
both to construct algorithms and inspect (hopefully to control) their working. A conspicuous part of the paper is devoted to the study and
systematic construction of parabolic/elliptic  discrete subgroups of the Lie groups $\mathrm{SO(r,r+q)}$, in view of discretization and/or
tessellations of the space to which data are to be mapped. Furthermore it is shown how the ingredients of Special K\"ahler Geometry and the
$c$-map, well known in the supergravity literature, provide a unified classification scheme of the relevant park of Tits-Satake universality
classes with noncompact rank $r\leq 4$.
\end{abstract}
\vfill
\end{center}
\par
\vfill\begin{minipage}{\textwidth} \small
\parbox{\textwidth}{\hrulefill} \\
Date: \today \\
% MSC 2020: 14H60, 14J60, 14N15 \\
% Keywords: Higgs bundles, Grassmann bundles, semistability, spectral varieties\\
U.B.'s   research is partly supported by Bolsa de Produtividade 305343/2025-4 of Brazilian CNPq, PRIN project 2022BTA242  ``Geometry of algebraic structures: moduli, invariants, deformations'' and INdAM-GNSAGA. P.G.F.~acknowledges  support  by the Company \emph{Additati\&Partners Consulting s.r.l} during all the
time of development of the present research.
\end{minipage}
\end{titlepage}
{\small \tableofcontents} \noindent {}
\newpage
\section{Introduction}
The aim of the present paper, as stated in its version of March 2025 when it was firstly posted on ArXiv, was to present, in a self-contained and   well-organized way, a
mathematical lore that was  developed in the context of Supergravity theory, and was here enriched with new strategic
ingredients in the perspective of applying such structures to an upgrading and improvement of the mathematical framework of Neural Networks, not
excluding other possible applications. A little more than one year after its first posting, by the time we prepare it for submission to a journal, what in March 2025 was a speculative suggestion  has been transformed to a large extent  into reality through the developments described in the following 
list of papers of which the first is the present one:
\cite{pgtstheory,TSnaviga,naviga,tassellandum,axialbeltra,geotermico,secondtemperature,terzatemperatura}.
As we illustrate in some detail in the conclusion section \ref{zakluchenye}, the \textit{supervised} \textbf{Cartan Neural Networks} do by now exist in the formulation presented in \cite{TSnaviga} and the road to \textit{unsupervised ones with reinforcement} has been opened by the results that two of us with A. Sorin have obtained since December 2005 in the three papers \cite{geotermico,secondtemperature,terzatemperatura} concerned with Souriau Thermodynamics \cite{souriaub1,souriaub2,barbarpapad4,
marlentropia,caldobarbaresco,barbaresco2,barbaresco3,marlegibbs,barbarpapad1,barbarpapad2,barbarpapad3,nebbo} and its extension. At the root of all such promising research lines stand the \textbf{Paint Group / Tits Satake  theory} (PGTS) of non-compact symmetric spaces systematically exposed in this paper and the focal intuition that the geometrical lore established within the 50 year long development of Supergravity Theory, or inspired by it, is strategic for an efficient geometrical refoundation of Artificial Intelligence.  
\par The  unconventional setup of this article, its content and the motivations for its writing can be
understood only if we start from the general considerations exposed in the following subsection.
\par  Indeed the main ambition of, the present
paper is that of providing a sound mathematical basis for \emph{a New Paradigm in the construction of  Geometric Neural Networks Architectures}
that will be developed in a series of subsquent papers partially already \emph{in fieri}.  For this new paradigm we have devised the  acronym
  \emph{PGTS Theory} whose rationale will be clear after reading the next subsection.
\subsection{Mathematical Theories nursed by Supergravity/Supersymmetry}
While Mathematics and Physics have enjoyed a fruitful, two-thousand-year interaction, the last forty years saw a particularly fertile manifestation of this relationship: the intense cross-fertilization of Supergravity theory with algebra and Geometry across their many facets.  The inventory of examples is very long, but for the purposes of this paper, we present only a short illustrative list, focusing on those items most relevant to our issues and goals.
\begin{description}
\item[1)] The discovery under the  name of  \emph{Cartan Integrable Systems} \cite{DAuria:1982ada,DAuria:1982uck} of the  supersymmetric
extension of Sullivan's bosonic \emph{Free Differential Algebras} \cite{sullivano}.
\item[2)] The discovery and formalization of \emph{Special
K\"ahler Geometry} \cite{skagazino_1,skgmat_1,skagazino2,specHomgeoA1,pupallo1,pupallo2,De_Jaegher_1998} in its two versions, \emph{local} and
\emph{rigid}.
\item[3)] The \emph{Solvable Group}  representation of \emph{noncompact homogeneous spaces}
    \cite{Cremmer:1979up},\cite{Fre:1996js,Andrianopoli:1996bq,Andrianopoli:1997wi,Trigiante:1997zba,pancetta2}.
\item[4)] The discovery of the \emph{$c$-map} from \emph{Special
K\"ahler Manifolds} to \emph{Quaternionic K\"ahler manifolds}, and the discovery of the \emph{$c^\star$-map} from \emph{Special K\"ahler
Manifolds} to \emph{Pseudo-Quaternionic K\"ahler manifolds}
\cite{Cecotti:1988ad,titsusataku,Bergshoeff:2008be,nostrog22,bacchihollititti,Fr__2013}.
\item[5)] The discovery of the concepts of \emph{paint
and subpaint groups}, the formalization of the \emph{Tits-Satake projection} and of its commutativity with the $c$-map
\cite{noipainted,titsusataku,rocek2006quaternion}.
\item[6)] \emph{Mirror Symmetry} in the Geometry of Calabi-Yau moduli spaces
\cite{picfucperiod_1,mirsym_1,mirsym_3,mirsym_4, mirsym_5,mirsym_6,mirsym_7,mirsym_8,mirsym_9,mirsym_10,mirsym_11, picfucperiod_3,picfucperiod_5}
\item[7)] The discovery of \emph{Englert equation} that generalizes to $7$-dimensions the time-honoured \emph{Beltrami equation} operating
instead in $3$-dimensions\cite{Englert:1982vs,Fre:2016gat,Cerchiai:2018shs}.
\item[8)] The $T$-tensor and the general \emph{embedding tensor} formulation of gauging of extended supergravities that, once extracted from the
    supergravity context, is a geometrical theory per se \cite{deWit:1982bul},
    \cite{Cordaro:1998tx},\cite{Nicolai:2000sc},\cite{deWit:2002vt} (see
    \cite{Samtleben:2008pe},\cite{Trigiante:2016mnt} for reviews).
\item[9)] The cross influence on the theory of \emph{restricted holonomy manifolds} as in
particular $\mathrm{G_2}$-manifolds and $\mathrm{Spin(7)}$-manifolds \cite{spinettamarengo,Fre:2006zxd}.
\end{description}
Why has the special relationship between Supersymmetry (SUSY)/Supergravity and Mathematics been so much more creative and productive than the general cross-relation between Physics and Mathematics?
The answer is simple: to a large extent, this success is not a simple cross-link between the two sciences. Instead, it uncovers internal connections between different mathematical fields and various structure groupings of mathematical objects.
Although Supergravity and supersymmetric field theories were initially conceived as physical theories, their foundation rests on a purely mathematical principle --- not an experimental one. This principle is the realization of Lagrangian field theories that must incorporate supersymmetry (represented by a super Lie algebra) as an infinitesimal invariance. This requirement proved to be both extremely stringent and deeply creative.
Throughout its almost fifty-year history, supergravity Lagrangian field theory has acted as an intermediate step connecting structures that belong to pure mathematics. It unveils hidden, purely mathematical connections that would not have been observed or utilized as organizing principles without the guiding hand of supersymmetry. This general phenomenon, evident in all the examples we have reviewed, we term the
\emph{precognition of supersymmetry}.

Ultimately, we can dispense with the supergravity or supersymmetric gauge theory Lagrangians that led to the discovery of new mathematical structures or new organizational patterns, and consider these structures for their own sake within pure mathematics. Mathematicians are accustomed to organizing their discipline into a well-established hierarchy that has its own principles, basic definitions, and logical order, often derived from historical development. The {\em precognition of supersymmetry}, which indirectly stems from a {\em higher symmetry principle,} provides structural viewpoints that organize objects from different classical mathematical families into unsuspected --- and at first glance bizarre --- but well-defined chains based on a unique logic. The $c$-map, one of our main concerns in the present paper, will be recognized by the reader as an outstanding example of this phenomenon.

\subsection{Application perspectives of a newly  reassembled mathematical theory and Data Science}
Summarizing, there is a rich set of mathematical structures connected with supergravity that constitute chapters of pure mathematics, partially
preexisting in  supergravity, partially newly developed in the context of supergravity, which are worth analysing and grouping according to the
mathematical logic suggested by  \emph{the precognition of supersymmetry}, rather than dismembered according to classical mathematical ordering.
The organization of such structures according to the  \emph{the precognition of supersymmetry} builds up a systematic theory of its
own which happens to be solid, rich, useful and liable to further developments within its own scope. Furthermore, by cutting the historical
relation between such  theory and its supergravity mother, together with the physics-inspired applications that the latter
has suggested, the aforementioned   theory \emph{opens up to applications in other scientific fields} besides physics. One
application field that will be preeminent in our minds while exposing our scheme and results is \emph{Data Science}, with particular
reference to the developments put forward by the authors of \cite{francesi1,francesi2,francesi3}, which we might collectively describe as
\emph{Hyperbolic Learning}. Yet, although extremely interesting, those developments are not the primary reasons for considering the \emph{supergravity nursed mathematical theory} we are going to expose as relevant to Data Science.
The fundamental reason is deeper, and makes \emph{Hyperbolic Learning} a confirmation of the need to adopt a different scheme in the construction of Geometric Neural Networks, rather than a primary motivation.
  The main reason is the  discomfort   of any mathematically minded
scientist with the use of the so-called \emph{point-wise activation functions} (\emph{sigmoid, arc tangent, Heaveside   function, etc.}) on
individual components of a pretended vector, which { constitutes a \emph{mathematical heresy} and} makes the whole architecture of neural
networks dependent on the choice of arbitrarily chosen bases, this happening repeatedly in every map from one layer to the next.  Such
discomfort has emerged also in the community of Data Scientists and the vision of a truly geometrical reformulation of neural networks has been
advanced by several authors \cite{4_GoverH_CNN1,5_GoverH_CNN2,6_G_CNN,Bronstein_2017,7_Graph_Man_DL,19_flowsmanif,librobronstiano}. In particular,
according to such visions, a sound  \emph{Geometrical Deep Learning} should be based on the \emph{principle of covariance, or equivariance},
with respect to some group of symmetry, of the transformations from one layer to the next, typically thought of as \emph{convolutions}.
Consequently,  in view of the need for non-linear maps, the layers cannot be thought of as  vector spaces; rather, they must be
characterized as vector bundles or \emph{"spaces of functions"} on differentiable manifolds $\mathcal{M}_n$.  At the same time the very core of
learning algorithms, the \emph{loss-function}, requires the \emph{notion of distance} between the points representing data;  this brings us to
those differentiable manifolds of negative curvature (Hyperbolic Spaces) that admit uniquely defined geodesics connecting any two points, providing
a notion of distance in terms of length of such geodesic arcs.

The future applications of the theory  we are discussing here are therefore to be looked for in all those areas of research of whatever
Science where the same geometrical needs as in Data Science are preeminent, the distance function being the first on the priority list.

\subsection{The Supergravity nursed Mathematical Theory exposed in this paper}
 The theory we aim to present here in a systematic and logically structured way is tentatively described by this article title. For
convenience we might shorten it into the acronym \emph{PGTS Theory}\footnote{\emph{P}aint \emph{G}roup \emph{T}its-\emph{S}atake}
\emph{ of  Hyperbolic Homogeneous Spaces}.  It is that one which emerges from a   fusion of items 2),3),4),5) of the previous list with
the addition of new ingredients introduced in this article. Purified of the original supergravity concepts that favoured the discovery and the
framing of its various components, the   lore presented here acquires a full mathematical stature of its own. We deem it worth being considered for  future developments in a purely mathematical perspective and for the possible establishment of new theorems and proofs.  This viewpoint
was already assumed by one of the present authors in his book \cite{advancio}, written some years ago, which contains the formulation in the
that perspective of several of the ingredients which are here utilized and reframed in a new hopefully tighter and better shaped logical
ordering. Indeed the new ingredients that we introduce into the texture provide a new quality and a new overall logical ordering of all the items.
We describe such new elements below, after summarizing, within a historical perspective, the older elements of the new reconstruction. It is
appropriate to stress that the  systematic theory assembled and described in this article is especially liable to allow for   applications to
Data Science, as we just advocated above, but not only there. Indeed, possibly also in connection with the Data Science applications, yet independently
interesting in its own mathematical sake, there is an ensuing \emph{new constructive strategy} for the determination of discrete subgroups of the
noncompact isometry groups that include infinite solvable subgroups, substituting abelian lattices of Euclidean or Minkowskian spaces. An example
of this mechanism is provided in this paper by the explicit construction of a subgroup of  $\mathrm{Sp(4,\mathbb{Z})}$, up to our knowledge
previously unknown in the mathematical literature, which is of the type described above.  In the subsequent section \ref{generaloneSOpq},  we were
able to generalize in a fully systematic form the first example to all the symmetric spaces $\mathrm{SO(p,q)/SO(p)\times SO(q)}$.  Such
constructions might be useful in the discretization of manifolds both for Data Science and for Lattice Quantum Field Theory. Tessellation schemes
and Tessellation Coxeter Group constructions can take advantage of the systematic theory  here exposed, and similarly, we expect applications in
Markov Diffusion Processes both in the continuum and on discretized spaces and graphs.
\subsection{Historical outline of the already existing ingredients of the PGTS Theory}
\label{leluka}
%%%%%%%%%%%%%%%%%
In this subsection we preliminarily present in a more or less extended manner, depending on the case,  the various concepts and already
existing mathematical structures that enter, as a building blocks, the new theoretical setup we plan to advocate in the next subsection of this
introduction. All the mathematical ingredients mentioned in the present subsection have been initiated or substantially developed within
supergravity theory.  In our exposition of the needed conceptual mathematical structures, we briefly sketch their history and the supergravity
motivations for their development.
\subsubsection{The solvable Lie group representation}
The first and most basic ingredient of the whole theory is given by following
\begin{teorema}\label{teoUHsolvrep}
Every noncompact symmetric space $\mathrm{U/H}$ where $\mathrm{U}$ is a finite-dimensional simple Lie group whose Lie algebra
$\mathbb{U}$ is a real section $\mathbb{G}_R$, different from the unique maximally compact one, of any of the simple complex Lie algebras $\mathfrak{a}_\ell,\mathfrak{b}_\ell,\mathfrak{c}_\ell,\mathfrak{d}_\ell,\mathfrak{f}_4,
\mathfrak{g}_2,\mathfrak{e}_{6,7,8}$ and $H \subset U$ is its maximal compact subgroup (with Lie algebra $\mathbb{H}$), is \emph{metrically equivalent to a solvable group manifold}:
\begin{equation}\label{solvablemma1}
  \mathcal{S}_{\mathrm{U/H}} \,\sim \, \exp\left[ Solv(\mathrm{U/H})\right]
\end{equation}
whose Lie algebra $Solv(\mathrm{U/H}) \subset \mathbb{U}$ is an appropriate subalgebra of the full isometry algebra $\mathbb{U}$.
\end{teorema}
Metrically equivalent means that:
\begin{enumerate}
  \item As differentiable manifolds the coset manifold $\mathrm{U/H}$ and the group manifold $\mathcal{S}_{\mathrm{U/H}}$ are diffeomorphic.
  \item The $\mathrm{U}$ invariant metric $g$  on $\mathrm{U/H}$, which is unique up to an overall constant, coincides with the metric on
$\mathcal{S}_{\mathrm{U/H}}$ obtained through left  transport to any point of $\mathcal{S}_{\mathrm{U/H}}$ of a suitable positive definite, symmetric  quadratic form $\langle \, , \, \rangle_g $ defined on $ Solv(\mathrm{U/H})$ (which is the tangent space to the Identity $\mathbf{e}\in\mathcal{S}_{\mathrm{U/H}}$)
\begin{equation}\label{solvablemma2}
\langle \, , \, \rangle_g  : \quad  Solv(\mathrm{U/H}) \otimes Solv(\mathrm{U/H}) \to  \mathbb{R}
\end{equation}
satisfying
$$\langle X,Y\rangle_g = \langle Y,X \rangle_g, \quad
 \langle X\, , \, X \rangle_g\,  \geq \, 0 , \quad
\langle X\, , \, X \rangle_g\,  = \, 0 \Rightarrow  X = 0
\quad \text{for all}\ X,Y \in Solv(\mathrm{U/H}).
$$
\end{enumerate} The proof of the above theorem is constructive and we shall provide it later, by showing how $Solv(\mathrm{U/H})$ is explicitly built,
starting from the root system of the corresponding complex Lie algebra $\mathbb{G}$ and the Cartan-Weyl basis of generators for the latter.
\par
What is important to note is that one might invert the terms of the argument and consider instead the notion of \emph{normed solvable Lie algebras}
\begin{definizione}\label{definormosolvalg}
\emph{A normed solvable Lie algebra} is a pair $\left(Solv,\langle \, , \, \rangle_n\right)$  where $Solv$ is a finite-dimensional
Lie solvable algebra and $\langle \, , \, \rangle_n$ is a  positive definite non degenerate quadratic  form on it:
which is requested to be invariant with respect to the adjoint action of $Solv$:
\begin{equation}\label{invarianorma}
  \forall X,Y,Z \in Solv  \quad : \quad \langle \left[Z\, ,\, X\right]\, , \, Y \rangle_n\, + \, \langle X\, , \, \left[Z\, ,\, Y\right] \rangle_n \, = \, 0
\end{equation}
\begin{equation} \label{normedsolva}
Solv\otimes Solv \to \mathbb{R}
\end{equation}
satisfying
$$\langle X,Y\rangle_n = \langle Y,X \rangle_n, \quad
 \langle X\, , \, X \rangle_n,  \geq \, 0 , \quad
\langle X\, , \, X \rangle_n\,  = \, 0 \Rightarrow  X = 0
\quad \text{for all}\ X,Y \in Solv
$$
which is requested to be invariant with respect to the adjoint action of $Solv$:
$$\langle \left[Z\, ,\, X\right]\, , \, Y \rangle_n\, + \, \langle X\, , \, \left[Z\, ,\, Y\right] \rangle_n \, = \, 0 \qquad \text{for all}\ X,Y,Z \in Solv .$$

\end{definizione}
The solvable group manifold
  $\mathcal{S}  = \exp\left[ Solv\right]$
 inherits
from the invariant quadratic form $\langle \, , \,
\rangle_n$ defined on its Lie algebra a Riemannian
metric (the left transport of $\langle \, , \,
\rangle_n$ from the identity to any other point),
for which $\mathcal{S}$ is a group of isometries.
\par
In 1975 Alekseveesky  considered the problem of
constructing all quaternionic K\"ahler manifolds
with a transitive solvable group of isometries
\cite{Alekseevsky1975} and obtained the following result.
\begin{definizione}\label{Alexnormalmanif}
A Riemannian manifold $\left ( \mathcal{M},g \right)$ is   a \emph{normal Riemannian homogeneous space} if it admits a  solvable Lie group
$\mathcal{S}_{\mathcal{M}}$ of isometries that acts on $\mathcal{M}$  freely and transitively.
\end{definizione}

\begin{proposizione} \cite{Alekseevsky1975}
There is a one-to-one correspondence between normed solvable Lie algebras (up to isomorphism) and normal Riemannian homogeneous spaces (up to isometry).
\end{proposizione}
\noindent{\bf Proof.}   Given a normed solvable Lie algebra $Solv$, the corresponding normal Riemannian homogeneous space is the simply connected, connected Lie group $\exp[Solv]$ equipped with the Riemannian metric given by the left transport of the non-degenerate symmetric form on $Solv$.
Vice versa, given a normal Riemannian homogeneous space $\mathcal M$ with solvable isometry group $\mathcal S_{\mathcal M}$, one establishes a diffeomorphism $f\colon\mathcal S_{\mathcal M}\to\mathcal M$ by choosing a point $p\in\mathcal M$ and letting $f(g) = gp$. The normed solvable Lie algebra $Solv$ is the Lie algebra of $\mathcal S_{\mathcal M}$ equipped with the symmetric form given by the pullback Riemannian metric evaluated at the identity. Different choices of $p$ produce isomorphic normed solvable Lie algebras.

%The obvious theorem behind the definition \ref{Alexnormalmanif} is %discussed in \cite{BorelTits}, \cite{Helgasonobook}  and states that if a
%Riemannian manifold $\left ( \mathcal{M},g \right)$ admits a simply %transitive solvable group of isometries $\exp[\Solv_{\mathcal{M}}]$, then %it
%is metrically equivalent to this solvable group manifold the quadratic %form being obtained as in eq.~\ref{identifico}.
%\par

\medskip
Using definition \ref{Alexnormalmanif}, theorem \ref{teoUHsolvrep} can be reformulated by saying that all noncompact symmetric spaces
$\mathrm{U/H}$ are \emph{normal homogeneous spaces}, the converse  however is not true. There are examples of normal homogeneous spaces that are
not symmetric spaces.  In the case of noncompact symmetric spaces
$\mathrm{U/H}$, the solvable group $\mathcal{S}_{\mathcal{M}}$ is defined by the Iwasawa decomposition
$U=\mathcal{S}_{\mathrm{U/H}}\cdot H$, so that
\begin{equation}
   \mathcal{M}= \frac{\mathrm{U}}{\mathrm{H}} \simeq \mathcal{S}_{\mathrm{U/H}}.\label{Iwasawa}
\end{equation}\par
Although Theorem \ref{teoUHsolvrep} was more or less present in the mathematical literature in implicit form, its
constructive proof which provided a tool of great value for explicit constructions of Supergravity  Lagrangians,  classification of their
gaugings, analysis of cosmic billiards, black hole solutions and scalar potentials came in the years 1997-2004 through the work of
supergravity experts in the already quoted papers \cite{Andrianopoli:1996bq,Andrianopoli:1997wi,Trigiante:1997zba}.
\subsubsection{Special K\"ahler Geometry}
The second fundamental pre-existing ingredient of the \emph{PGTS theory}   is \emph{Special K\"ahler Geometry}.
%%%%%%%%%%%%%%%%%%%%%%%%%%%%%%%%%%%%%%%%%%%%%%%%%
The discovery of this class of complex   manifolds and the ample development of their geometry in its two versions, local and rigid,  is undebatably due to Supergravity/Supersymmetry authors, much before that either one of the latter attracted the attention of
mathematicians.  Furthermore neither one of these  K\"ahler geometries might be regarded as a physicists' construction that, at a later stage,
was reformulated in mathematical terms by mathematicians.  Actually, the precise mathematical formulation and the very mathematical definition of
both types of Special K\"ahler geometry was worked out and published by authors of the Theoretical Physics community, who also established the
nomenclature for this chapter of complex differential geometry and the precise relations between the two classes of manifolds, the local and rigid
ones. Such relation is a limiting procedure, somewhat reminescent of group contraction,  from the local case to the rigid one. Historically the
notion of \emph{local Special K\"ahler Geometry}, formulated in \emph{special coordinates}  was introduced in 1984 by B. de Wit et al. in
\cite{SKGaggio1,SKGaggio2}. More or less at the same time, in 1984, what shortly later became known as \emph{rigid Special K\"ahler Geometry}
was derived from the superspace approach to supersymmetric lagrangian field theories in \cite{sierratown,Gates:1983py}. The precise mathematical
definition of \emph{(local) Special K\"ahler Geometry} was introduced in 1990 in two formulations that were later proved to be equivalent. On
one side  Strominger in \cite{skgmat_1},  linked \emph{(local) Special K\"ahler Geometry}   to the theory of complex structure and  K\"ahler
structure deformations of Calabi-Yau threefolds and formulated its definition in terms of a flat symplectic bundle whose sections are related to
the periods of the unique $\Omega^{(3,0)}$-form on the basis of homology $3$-cycles. In this way, however, he overlooked the matter of fact that
Homogeneous Special K\"ahler Manifolds (of the local type) do exist which cannot be identified with moduli spaces of any Calabi-Yau threefold.
On the other side,   in \cite{skagazino2} Castellani, D'Auria and Ferrara based the definition of  \emph{(local) Special K\"ahler Geometry} on
a special identity satisfied by the Riemann tensor of the corresponding K\"ahler metric which follows from $\mathcal{N}=2$ supersymmetry when the
complex coordinates of the manifold are identified with the scalar fields of $\mathcal{N}=2$ vector multiplets coupled to supergravity. The two
mathematical definitions, both motivated by the use of such manifolds in the context of supersymmetric field theories,  are formulated in a
mathematical intrinsic language, and they were almost immediately shown to be equivalent/ Indeed both
definitions can be summarized by means of a set of \emph{formal Picard Fuchs equations} disconnected from specific varieties and their
deformations (for an early review see the appropriate chapters of \cite{Fre:1995bc}). In 1994 the famous paper on the non-perturbative structure
of $\mathcal{N}=2$ Supersymmetrric Gauge  Theories by Seiberg-Witten \cite{Seiberg:1994rs} brought attention  to the underlying   \emph{Special
K\"ahler Geometry of rigid type}, that was extensively used by the authors.
%without acknowledging its preexisting mathematical general definition
%and naming.
For the rigid case, there is also a connection with the deformation theory of algebraic variety complex structures and with the
corresponding Picard-Fuchs equations. The relevant varieties, in this case, are Riemann surfaces rather than Calabi-Yau threefolds, yet, just as
in the local case, the mathematical definition of rigid Special K\"ahler Geometry is autonomous, precise and fully general. In the 1996 paper
\cite{pupallo2}, which is also rather famous for the number of citations, the mathematical definitions  of both \emph{Special K\"ahler geometries}
was clearly spelled out and the contraction procedure, named \emph{ rigid limit} from the first to the second was also sketched. The other 1996
paper \cite{Bill__1996} was devoted to the issue of   Picard Fuchs equations of rigid special geometry  applied to the case of suitably chosen
Riemann surfaces, in particular hyperelliptic.  In the later 1998 paper \cite{Bill__1998} the issue of deriving the rigid limit from the
degeneration to Riemann surfaces of a class of Calabi-Yau threefolds that have a fibered structure was extensively investigated.
%%%%%%%%%%%%%%%%%%%%%%%%%%%%%%%%%%%%%%%%%%%%%%%%%%
\par
In the present context the relevant type of Special K\"ahler Geometry is the local one\footnote{In the supergravity literature the current use is
that when Special K\"ahler Geometry is mentioned without further qualifiers one means that of local type. If one wants to refer to rigid Special
K\"ahler Geometry then one adds the qualifier rigid. Such convention we adopt also in this paper. In order to avoid confusions it is useful to
stress that in the mathematical literature developed after the two versions of Special K\"ahler Geometries had been discovered in the supergravity
community, the current use is to name Special K\"ahler Geometry without qualifier the rigid one and to name Projective K\"ahler Geometry the local
one. Irrespectively from the name local versus rigid, that it is reminiscent of the supersymmetry origin of both (local or rigid) and
irrespectively of the priority in the discovery of these mathematical structures it remains the intrinsic fact that the rigid case is a limiting
case of the local one, where the fundamental flat symplectic bundle providing the definition, degenerates, via an Inonu Wigner contraction of the
structure group, in the same way as the Lorentz group contracts to the Galilei group, or the de-Sitter group contracts to the Poincar\'e one.}
and the  formulation of its definition,  that we present in section \ref{speckalgeosum}, was  first published in \cite{FreLect} and then
included in the already quoted complete exposition of matter coupled $\mathcal{N}=2,D=4$ supergravity  theory in \cite{pupallo2}. The exposition
of section \ref{speckalgeosum} closely follows  that given in the book \cite{advancio}.
\par The set of Special K\"ahler Manifolds (of local type) is vast and it contains  manifolds with no continuous isometries, as the already
mentioned moduli spaces of \emph{K\"ahler structure deformations} (\emph{i.e.} $(1,1)$ Dolbeault cohomology classes) and of \emph{Complex
structure deformations} (\emph{i.e.} $(2,1)$ Dolbeault cohomology classes) of \emph{smooth Calabi-Yau threefolds} but also \emph{homogeneous
spaces} including \emph{symmetric spaces}. The classification of \emph{homogeneous special geometries} was achieved in
\cite{vandersuppa,specHomgeoA1,specHomgeoA2,deWit:1995tf,Cortes}.
\subsubsection{The $c$-map}
The third pre-existing ingredient which plays a relevant role in the new and enriched setup presented in this article, just as it already did in
the systematics of $\mathcal{N}$-extended and matter coupled supergravity models and in their filiations,  is the \emph{$c$-map}. It consists of
a quite non trivial and universal construction that allows to build a quaternionic K\"ahler manifold  $\mathcal{QM}_{4n+4}$ of real dimension
$\mathrm{dim}_{\mathbb{R}} \mathcal{QM}_{4n+4} \, = \, 4n+4$ starting from any Special K\"ahler manifold $\mathcal{SK}_{n}$ of complex dimension
$\mathrm{dim}_{\mathbb{C}} \mathcal{SK}_{n} \, = \, n$.  The embryo of the $c$-map construction is due to Cecotti in \cite{Cecotti:1988ad},
followed by the seminal paper \cite{catlantide1} by Ferrara and Sabharwal. It was utilized and developped to a large extent by de Wit, Van
Proeyen, Vanderseypen, Alekseevsky and Cort\'es in their, in depth, classification study of Homogeneous Special Geometries, already quoted above
\cite{vandersuppa,specHomgeoA1,specHomgeoA2,deWit:1995tf,Cortes}. Further perfectioned in a more mathematical phrasing   in \cite{titsusataku} and
\cite{catlantide2}, $c$-map received its final, general formulation by Fr\'e, Sorin and Trigiante in \cite{cmappotto}. In particular in
\cite{cmappotto}  the crucial general expression was derived of the $\mathrm{SU(2)}$-connection of the Quaternionic K\"ahler Manifolds in terms of
generic Special K\"ahler data. This final general form of the $c$-map construction was reported in \cite{advancio} and we present it here in
section \ref{cmappusquat}.
\subsubsection{The Tits-Satake projection and the paint group.}
The fourth preexisting ingredient of PGTS Theory, which is inextricably intertwined with the previous three  and, once combined with the new
elements to be described below, constitutes the very backbone of the whole theory, is the Tits-Satake projection with the associated notion of
paint group, introduced for the first time  in \cite{noipainted} and then fully formalized in \cite{titsusataku}. In view of  the central
relevance of this two-component ingredient, we anticipate a shortened presentation of its definition and we also
briefly summarize the historical developments that led to its conceptualization.
\par
One of the monumental achievements of \'Elie Cartan was the classification of all Riemannian symmetric spaces, which has the classification of
real sections of complex simple Lie algebras as prerequisite. Everty complex simple Lie algebra $\mathfrak{g}$  corresponds
to a root system $\Phi_{\mathfrak{g}}$ which is codified and fully determined by a corresponding Dynkin diagram encoding the structure of the Cartan matrix. The root system can be partitioned in two disjoint subsets
$\Phi_{\mathfrak{g}} \, = \, \Phi_{\mathfrak{g}}^+ \bigcup \Phi_{\mathfrak{g}}^-$ of equal cardinality, containing respectively the positive and negative roots. A convenient starting point for the discussion of real sections is provided by the introduction of the Cartan-Weyl basis of generators:
\begin{equation}\label{CartanWeylbasis}
  T_A \, = \, \left\{ \underbrace{H_i }_{\text{Cartan generators}}\, , \, \underbrace{E^{\alpha} }_{\alpha \in \Phi_{\mathfrak{g}}^+}\, , \, \underbrace{E^{-\beta}}_{\beta \in\Phi_{\mathfrak{g}}^+} \right\}\quad ; \quad i\, =\, 1,\dots, r,
\end{equation}
where $r$ is the rank of $\mathfrak{g}$ and hence the dimension of the chosen Cartan subalgebra $\mathcal{H}\subset \mathfrak{g}$
(they are all equivalent up to conjugation). Furthermore one can normalize the Cartan generators $H_i $  in such a way that the Killing Cartan metric restricted to $\mathcal{H}$ is
\begin{equation}\label{KronKilmet}
  \kappa (H_i , H_j) \, = \, \delta_{ij} \, = \, \text{Kronecker delta.}
\end{equation}
Correspondingly the elements of the root system become real valued vectors in $\mathbb{R}^r$:
\begin{equation}\label{gommapiuma}
  \alpha(H_i) \, = \, \alpha^i \in \mathbb{R}
\end{equation}
In the so defined Cartan-Weyl basis, the commutation relations of the complex Lie algebra $\mathfrak{g}$ which, as a vector space,  is the complex span of  $T_A$ :
\begin{equation}\label{gcomplex}
  \mathfrak{g} \, = \, \mathrm{span}_{\mathbb{C}} \left( T_A \right),
\end{equation}
become
\begin{eqnarray}\label{CarWeylform}
  \left[ H_i\, ,\, H_j \right ] &=& 0~, \nonumber\\
  \left[ H_i\, ,\, E^\alpha \right ] &=& \alpha_i \, E^\alpha ~,\nonumber\\
  \left[ E^\alpha \, ,\, E^\beta \right ]&=& N(\alpha,\beta) \,
  E^{\alpha + \beta} \quad \mbox{if $(\alpha+\beta) \in \Phi_{g}$} \quad ; \quad N(\alpha,\beta) = \pm 1 \nonumber \\
\left[ E^\alpha \, ,\, E^\beta \right ]&=& 0\, \quad \mbox{if $(\alpha+\beta) \notin \Phi_{g}$} \nonumber\\
  \left[ E^\alpha \, ,\, E^{-\alpha} \right ] &=& 2 \, \alpha \, \cdot \, \mathcal{H}~, \nonumber\\
\end{eqnarray}
A basic theorem whose proof by direct construction is very simple  (see for instance \cite{pietronewbook} ) is the following?
\begin{teorema}\label{teotriangolo}For all finite dimensional complex simple Lie algebras $\mathfrak{g}$,
the Cartan Weyl generators (\ref{CartanWeylbasis}) can be realized, in the fundamental defining representation and in all the other linear
representations, in such a way that all the $T_A$ are real valued matrices, and moreover, the Cartan generators $H_i$ are diagonal matrices,  the
step operators $E^{\alpha}$ for $\alpha >0$ are strictly upper triangular matrices and the step operators $E^{-\alpha}=(E^\alpha)^T $ are strictly
lower triangular matrices obtained as the transposed of the corresponding matrices $E^\alpha$.
\end{teorema}
Using the basis  $T_A$ with the properties stated in theorem \ref{teotriangolo}, it is immediate to see that there exist two universal real
sections of any complex simple Lie algebra:
\begin{description}
  \item[a)]\emph{The maximally split real section $\mathbb{G}_{\mathbf{max}} \subset \mathfrak{g}$}.
  This is defined by assuming that in the linear combinations $c^A \, T_A$ the allowed
  coefficients $c^A$ are all real. In any linear representation of $\mathbb{G}_{\mathbf{max}}$ the matrices
  representing
  \begin{equation}\label{bagarraT}
    T_A \, \equiv \, \left\{ H_i \, , \,  E^{\alpha}\, , \, E^{-\alpha}\right\}
  \end{equation}
 are all \emph{real}. From the representations of
  $\mathbb{G}_{\mathbf{max}}$, by taking linear combinations of the generators with complex coefficients one
  obtains all the linear representations of the complex Lie algebra $\mathfrak{g}$.
  \item[b)]\emph{The maximally compact real section $\mathbb{G}_{\mathbf{comp}}$}. This real section, whose
  exponentiation produces a compact Lie group, is obtained by allowing linear combinations with real coefficients of the
following set of generators:
  \begin{equation}\label{emorraT}
    {T}_A^c \, \equiv \, \left\{ {\it i}\,H_i \, , \, {\it i} \left( E^{\alpha}+ E^{-\alpha}\right) \, , \,
    \left( E^{\alpha}- E^{-\alpha}\right)\right\}
  \end{equation}
In all linear representations of
$\mathbb{G}_{\mathbf{comp}}$ the matrices representing the generators
${T}_A^c$ are \emph{anti-hermitian}.
\end{description}
\par
All other   real sections of the complex Lie algebra $\mathfrak{g}$ are obtained by   the   Cartan involutions
of the maximally compact real section  $\mathbb{G}_{\mathbf{comp}}$.
\begin{definizione}
\label{cartainvolva} Let:
\begin{equation}\label{cavacchia}
    \theta \quad : \quad \mathbb{G}_{\mathbf{comp}} \, \rightarrow \, \mathbb{G}_{\mathbf{comp}}
\end{equation}
be a linear automorphism  of the compact Lie algebra $\mathbb{G}_{\mathbf{comp}}$ . By definition we have:
\begin{equation}\label{caorsoset}
    \forall \alpha,\beta \in \mathbb{R} \quad , \quad \forall \mathbf{X},\mathbf{Y} \in \mathbb{G}_{\mathbf{comp}}\quad :
    \quad\left\{\begin{array}{rcl}
\theta\left (\alpha\, \mathbf{X} + \beta \, \mathbf{Y}\right)  & =  & \alpha \, \theta(\mathbf{X}) \,
+ \, \beta \, \theta(\mathbf{Y}) \\
 \theta\left(\left[ \mathbf{X}\, , \, \mathbf{Y}\right]\right) & = & \left[ \theta(\mathbf{X}) \, ,
 \, \theta(\mathbf{Y})\right].
 \end{array}
     \right.
\end{equation}
If $\theta^2 \, =\,\mathbf{Id}$ then $\theta$ is called a Cartan
involution of the Lie algebra
$\mathbb{G}_{\mathbf{comp}}$.
\end{definizione}
For any Cartan involution $\theta$ the possible eigenvalues are $\pm 1$. This allows us to split the entire Lie algebra
$\mathbb{G}_{\mathbf{comp}}$ in two subspaces corresponding to the eigenvalues $1$ and $-1$, respectively:
\begin{eqnarray}\label{cariniello}
    \mathbb{G}_{\mathbf{comp}}\, = \, \mathbb{H}_\theta\oplus \mathbb{P}_\theta
\end{eqnarray}
One immediately realizes that:
\begin{equation}
  \left[ \mathbb{H}_\theta\, , \,\mathbb{H}_\theta \right ] \, \subset \,  \mathbb{H}_\theta \quad ; \quad
  \left[ \mathbb{H}_\theta\, , \,\mathbb{P}_\theta \right ] \, \subset \, \mathbb{P}_\theta \quad ; \quad
  \left[ \mathbb{P}_\theta\, , \,\mathbb{P}_\theta \right ]
  \,\subset\,\mathbb{H}_\theta\label{catibaleno}
\end{equation}
The first of eqs.~(\ref{catibaleno}) tells us that $\mathbb{H}_\theta$ is a subalgebra, the second states that the decomposition (\ref{cariniello}) is reductive,\footnote{Note that this notion of ``reductive" is not the standard one in mathematics; see \cite{Humphreys-LieAlg}.} the third shows that it is not only reductive but also symmetric. Hence  any Cartan involution
singles out a compact coset manifold which is a symmetric Riemannian space.
\begin{equation}\label{carlicatus}
    \mathfrak{M}^{c}_\theta \, = \, \frac{\mathrm{G_{c}}}{\mathrm{H_\theta}} \quad \text{where} \quad \mathrm{H_\theta} \equiv
    \exp\left[\mathbb{H}_\theta\right] \quad ; \quad \quad \mathrm{G_c} \equiv
    \exp\left[\mathbb{G}_{\mathbf{comp}}\right]
\end{equation}
The structure \ref{catibaleno} has also another important consequence. If we define the vector space
\begin{equation}\label{wiccus}
   \mathbb{G}_\theta \,=\,  \mathbb{H}_\theta\oplus
   \mathbb{K}_\theta\quad ; \quad \mathbb{K}_\theta\,\equiv \,  {\it i} \,\mathbb{P}_\theta,
\end{equation}
we see that $\mathbb{G}_\theta $ is closed under the Lie bracket and hence it is a Lie algebra. It is some real
section  of the complex Lie algebra $\mathfrak{g}$ from which we started and we can consider the new, generally non
compact, symmetric space
\begin{equation}\label{giruclando}
    \mathcal{M}_\theta \, = \, \frac{\mathrm{G_\theta}}{\mathrm{H_\theta}} \quad ;
    \quad \mathrm{H_\theta} \equiv
    \exp\left[\mathbb{H}_\theta\right] \quad ; \quad \quad \mathrm{G_\theta} \equiv
    \exp\left[\mathbb{G}_\theta\right]
\end{equation}
An important theorem for which we refer the reader to  classical
textbooks \cite{gilmorro,Helgasonobook,knappo}\footnote{The proof is
also summarized in appendix B of\cite{ennalunga}.} states that, up to isomorphism, all
real sections $\mathbb{G}_{\mathrm{R}} \subset \mathfrak{g} $ of a simple complex Lie algebra $\mathfrak{g}$
are obtained in this way as a $\mathbb{G}_\theta$ for a convenient choice of the Cartan involution $\theta$. Furthermore as part of the same theorem one has that $\theta$ can always be chosen so that it
maps the compact Cartan subalgebra into itself:
\begin{equation}\label{caleppiotto}
    \theta \, : \,\mathbb{\mathcal{H}}_c \, \rightarrow \, \mathbb{\mathcal{H}}_c
\end{equation}
\par
This short discussion shows that the classification of real sections
of a complex Lie algebra $\mathfrak{g}$ is in one-to-one
correspondence with the classification of symmetric spaces, the
complexification of whose Lie algebra of isometries is
$\mathfrak{g}$.
\par
Let us now consider the action of the Cartan involution $\theta$ on the Cartan subalgebra:
\begin{equation}\label{compatalga}
  {\mathcal{H}}_c \equiv \text{span} \{{\mathit{i}} H_j\}
\end{equation}
of the maximal compact section $\mathbb{G}_{\mathbf{comp}}$. The subspace of $ {\mathcal{H}}_c$ belonging to the eigenspace of $\theta$ corresponding to eigenvalue $1$ is, according to eq.~(\ref{cariniello}), part of the subalgebra $\mathbb{H}_\theta = \mathbb{H}_c$,
while the subspace of $ {\mathcal{H}}_c$ belonging to the eigenspace corresponding to the eigenvalue $-1$ is, according once again to
eq.~ (\ref{cariniello}), part of the subspace $\mathbb{P}$.
 It follows that the Cartan subalgebra of the real section $\mathbb{G}_\theta$ can be split in two subalgebras as follows:
\begin{eqnarray}
\label{liutprandodacremona}
 \mathbb{G}_\theta \, \supset\, \underbrace{\mathcal{H}_\theta}_{\text{CSA of $\mathbb{G}_\theta$} }   &=& \mathcal{H}_\theta^{comp} \, \oplus \, \mathcal{H}_\theta^{n.c.} \nonumber \\
  \mathcal{H}_\theta^{comp} &\equiv & {\mathcal{H}}_c \bigcap \mathbb{H}_\theta \nonumber\\
  \mathcal{H}_\theta^{n.c.}&=& \mathit{i} \, {\mathcal{H}}_c \bigcap \mathbb{P}_\theta   \subset  \mathbb{\mathbb{K}}_\theta
\end{eqnarray}
Since, as stated in the quoted theorems, any real section $\mathbb{G}_\mathrm{R}$ of a complex $\mathfrak{g}$ is isomorphic to $\mathbb{G}_\theta$ for a suitable $\theta$,  it follows that the splitting of the Cartan subalgebra of any $\mathbb{G}_\mathrm{R}$ into a compact and a noncompact subalgebra has an intrinsic meaning. Hence for every real section $\mathbb{G}_{\mathrm{R}}$ we can define the compact and noncompact rank as follows:
\begin{equation}\label{palombodifiume}
  r_{n.c.} \, = \, \dim  \mathcal{H}_\theta^{n.c.} \leq r \quad ; \quad r_c \, = \, \dim  \mathcal{H}_\theta^{comp}  \, = \, r-r_{n.c.}
\end{equation}
Having clarified the relation between real sections $\mathbb{G}_{\mathrm{R}}$ of a complex Lie algebra $\mathfrak{g}$, Cartan involutions $\theta$ and noncompact Riemannian symmetric spaces $\frac{\mathrm{G_\theta}}{\mathrm{H_\theta}}$,
let us then reconsider a symmetric coset manifold:
\begin{equation}\label{simmetrus}
    \mathcal{M}_{\mathbb{G}_\mathrm{R}} \, = \, \frac{\mathrm{G}_\mathrm{R}}{\mathrm{H}_c}
\end{equation}
where $\mathrm{G}_\mathrm{R}$ is a simple Lie group and $\mathrm{H}_c\subset \mathrm{G}_\mathrm{R}$ is its maximal compact subgroup. The Lie algebra $\mathbb{H}_c$ of the denominator
$\mathrm{H}_c$ is the maximal compact subalgebra $\mathbb{H}_c \subset \mathbb{G}_\mathbb{R}$. Denoting as usual by $\mathbb{K}$ the orthogonal complement of $\mathbb{H}_c$ in $\mathbb{G}_\mathrm{R}$
\begin{equation}
  \mathbb{G}_\mathrm{R} = \mathbb{H}_c \, \oplus \,\mathbb{K}
\label{Grdecompo}
\end{equation}
we define  as \emph{rank of the coset} $\mathrm{G_\mathrm{R}/H_c}$ the dimension of the noncompact
Cartan subalgebra:
\begin{equation}
\mbox{rank} \left( \mathrm{G_\mathrm{R}/H_c}\right) \, = \,   r_{nc} \, \equiv \, \mbox{dim} \,
  \mathcal{H}^{n.c.} \quad ; \quad \mathcal{H}^{n.c.} \, \equiv \,
  \mbox{CSA}_{\mathbb{G}_\mathrm{R}} \, \bigcap \, \mathbb{K} \quad ; \quad \mathcal{H}^{comp} \, \equiv \,
  \mbox{CSA}_{\mathbb{G}_\mathrm{R}} \, \bigcap \, \mathbb{H}_c
\label{rncdefi}
\end{equation}
Next, choosing a basis of $\mathbb{\mathcal{H}}_c$ aligned
with the simple roots:
\begin{equation}\label{gominalo}
   \mathbb{\mathcal{H}}_c \, = \, \text{span} \left\{ {\it i}\,H_{\alpha_i}\right\},
\end{equation}
we see that by  vector space  duality  the action of the Cartan involution $\theta$ is  transferred from the Cartan subalgebra to the simple roots $\alpha_i$ and hence to the entire root
lattice. As a consequence we can introduce the notion of real and imaginary roots. One argues as follows.
Every vector in the dual of the full Cartan subalgebra\index{Cartan subalgebra}, in particular every root
$\alpha$, can be decomposed into its parallel and its transverse
part to $\mathcal{H}^{n.c.}$:
\begin{equation}
  \alpha = \alpha_{||} \, \oplus \, \alpha_{\bot}.
\label{splittus}
\end{equation}
\par
A root $\alpha$  is said to be  \emph{imaginary} if $\alpha_{||} \, = \, 0$,  \emph{real} if $\alpha_{\bot}\,=\, 0$.
%Generically a root is complex.

Given the  Dynkin diagram\index{Dynkin diagram} of a complex
Lie algebra, one can characterize a real section\index{real section}
by mentioning which of the simple root\index{simple root}s are
imaginary. We do this by painting black the imaginary roots. The
result is a Tits-Satake diagram\index{Tits-Satake diagram} like that
in Figure \ref{careliusA} which corresponds to the real Lie algebra
$\so(p,2\ell-p+1)$ for $p>2$, $\,\,\ell>2$.
\begin{figure}
\centering \vskip -3cm
\begin{picture}(100,250)
%%%%%%%%%%%%%%%%%%%%%%%%%%
\put (-70,155){$\mathfrak{b}_\ell$} \put (-30,160){\circle {10}} \put (-33,145){$\alpha_1$} \put
(-25,160){\line (1,0){20}} \put (0,160){\circle{10}} \put (-3,145){$\alpha_2$} \put (5,160){\line (1,0){20}}
\put (30,160){\circle {10}} \put (27,145){$\alpha_3$} \put (37,160){$\dots$}
%\put (45,260){\line (1,0){20}}
\put (52,160){\circle {10}} \put (49,145){$\alpha_{p}$} \put(57,160){\line(1,0){20}}
 \put (110,160){\circle* {10}}
 \put(82,160){\circle*{10}}\put (89,160){$\dots$}
\put(107,145){$\alpha_{\ell-2}$} \put (115,160){\line (1,0){20}} \put (140,160){\circle* {10}} \put
(137,145){$\alpha_{\ell-1}$} \put (145,163){\line (1,0){20}} \put (145,158){\line (1,0){20}} \put
(148,157){{\LARGE$>$}} \put (170,160){\circle* {10}} \put (167,145){$\alpha_\ell$}
\end{picture}
\vskip -4.5cm \caption{\label{careliusA} The Tits-Satake
diagram\index{Tits-Satake diagram} representing the real form
$\so(p,2\ell-p+1)$ of the complex $\so(2\ell+1)$ Lie algebra}
\end{figure}
\paragraph{\sc The algebraic Tits-Satake Projection.}
What has been reviewed above is sufficient to introduce the original notion of   Tits-Satake projection, which may be given at the
level of root systems. Let us consider a simple real Lie algebra $\mathbb{G}_\mathrm{R}$ which in general is neither maximally split nor maximally compact;  its noncompact rank satisfies the condition $0<r_{n.c}<r$, having called $r=\text{rank}(\mathfrak{g})$ and $\mathfrak{g}=
\text{complexification}(\mathbb{G}_\mathrm{R})$. Let us recall that the root system $\Phi_{\mathfrak{g}}$  is a set of
  vectors in $\mathcal{H}^\star$,    the dual of the   Cartan subalgebra. Next consider the splitting of   the
roots $\alpha \in \Phi_{\mathfrak{g}}$ according to eq.~(\ref{splittus}) and define the \emph{algebro-geometric TS projection} as
\begin{eqnarray}\label{algeoTS}
  \pi_{\mathrm{agTS}} & : & \mathcal{H}^\star \to (\mathcal{H}^{r_{n.c.}})^\star \nonumber\\
\forall \mathbf{v}\in \mathcal{H}^\star \, :\,  \pi_{\mathrm{agTS}}(\mathbf{v}) & =&  \mathbf{v}_{||}
\end{eqnarray}
The image by $\pi_{\mathrm{agTS}}$ of the root system $\Phi_{\mathfrak{g}}$ is a new system of vectors in $(\mathcal{H}^{r_{n.c.}})^\star$
\begin{equation}\label{tallero}
  \pi_{\mathrm{agTS}}\left(\Phi_{\mathfrak{g}}\right) \, = \, \Phi_{\mathfrak{g}}^{\mathrm{TS}} \, \subset \, (\mathcal{H}^{r_{n.c.}})^\star.
\end{equation}
A priori there is no reason to expect that $\Phi_{\mathfrak{g}}^{\mathrm{TS}} $ should be a root system, yet there is the following:
\begin{teorema}\label{ruttosistnuovo}
If $\mathfrak{g}=\text{\rm complexification}(\mathbb{G}_\mathrm{R})$ is any of the simple Lie algebras except $\mathfrak{a}_\ell$, then
$\Phi_{\mathfrak{g}}^{\mathrm{TS}}$ is the root system of another simple complex Lie algebra $\mathfrak{g}_{TS}$ of rank $r=r_{n.c.}$. Furthermore,
$\mathfrak{g}_{TS}$ is the complexification of a \emph{maximally split} simple real Lie algebra $\mathbb{G}_{\mathrm{TS}}$  which is a
subalgebra of the original non-maximally split real Lie algebra $\mathbb{G}_{\mathrm{TS}} \subset \mathbb{G}_\mathrm{R}$. The Tits-Satake
subalgebra $\mathbb{G}_{\mathrm{TS}}$ is the unique, up to conjugation, maximally split subalgebra of rank equal to the noncompact rank.
\end{teorema}
\paragraph{\sc The paint group and Cosmic Billiards.}
The notion of paint group emerged in the context of the application of $\mathrm{U/H}$ geometry to cosmological solutions of higher dimensional
supergravity without scalar potential (ungauged supergravity) displaying the very interesting universal mechanism named \emph{cosmic billiard}.
Cosmic billiards were proposed in the years  1999-2005 by several authors
\cite{cosmicbilliardliterature1A,cosmicbilliardliterature1B,cosmicbilliardliterature1C,cosmicbilliardliterature1D,cosmicbilliardliterature2A,
cosmicbilliardliterature2B,cosmicbilliardliterature3A,cosmicbilliardliterature3B,cosmicbilliardliterature3C,cosmicbilliardliterature3D,
cosmicbilliardliterature3E,cosmicbilliardliterature4A,cosmicbilliardliterature4B,cosmicbilliardliterature4C} as an elaboration of the much earlier
seminal ideas of V.A. Belinsky, I.M. Khalatnikov, E.M. Lifshitz \cite{Kassner1,Kassner2,Kassner3}. According to that conception, which is better
described as that of \emph{rigid cosmic billiard},  the various dimensions of a higher dimensional gravitational theory are identified with the
generators of the Cartan subalgebra $\mathcal{H}$ of a Lie algebra motivated by supergravity, and cosmic evolution takes place in a Weyl chamber of
$\mathcal{H}$. Considering the  Cartan scalar fields as the coordinate of a fictitious ball, during cosmic evolution such a ball scatters on the
walls of the Weyl chambers and this pictorial image of the phenomenon is at the origin of its denomination \emph{cosmic billiard}. The reason
why we described the original conception as \emph{rigid billiards} is that the walls of the Weyl chamber were initially looked at as static
walls. In the years 2003-2007, two of us (P.F. and M.T.), in collaboration with our students and also with A.S. Sorin, elaborated the conception
of \emph{soft cosmic billiard} that relies on the full integrability of geodesic equations for $\mathrm{U/H}$ manifolds and on the
identification of such geodesic equations with the complete set of field equations of higher dimensional supergravity, while looking at solutions
that depend only on one evolutionary parameter, the cosmic time $t$ \cite{noiconsasha,sashaebog,arrowtime}. In the soft cosmic billiard
picture the bosonic non-gravitational fields of the supergravity theory (scalars and vector fields) are in a one-to-one correspondence with the
positive roots of the solvable Lie algebra (see theorem \ref{teoUHsolvrep}) and the simple ones with walls of the Weyl Chamber. Since these fields
evolve in time as the Cartan ones that describe the space-time metric components, the walls causing the scattering of the cosmic ball rise and
decay, namely  they are soft.
\par
In this context the distinction between compact and noncompact directions of the Cartan subalgebra appeared essential, and this brought the Tits-Satake projection into the game. Indeed in 2005, Fr\'e, Gargiulo and Rulik constructed explicit examples of soft cosmic billiards in the case of a
\emph{non-maximally split symmetric manifolds} $\mathrm{U/H}$ and analyzed the role of the Tits-Satake projection, which led to the introduction
of the new mathematical concept of \emph{paint group} \cite{noipainted}. The tale is simple: considering the algebraic Tits-Satake projection
briefly described in the previous lines, it turns out that from the point of view of the Tits-Satake root system $\Phi_{TS}$, some of the positive
roots have a unique preimage in the original full algebra root system $\Phi_{\mathfrak{g}}$, while others have a multiplet of preimages. Since the
walls of the relevant Weyl chamber are in correspondence with the $\Phi_{TS}$ roots, we were obviously interested in the multiplicity of such walls, which in a way acquired a color, namely, they were \emph{painted}. We observed that the multiplets of preimages of each Tits-Satake positive root made up
linear representations of a compact subgroup $\mathrm{G_{Paint}\subset \mathrm{U}}$, and that such a phenomenon is universal. Different coset
manifolds $\mathrm{U/H}$ whose Lie algebra $\mathbb{U}$ projects to the same Tits-Satake subalgebra $\mathbb{G}_{TS}$ are distinguished by  a
different paint group $\mathrm{G_{Paint}}$, which is their intrinsic characterization within the universality class labeled by $\mathbb{G}_{TS}$.

In 2007 P. Fr\'e, F. Gargiulo, J. Rosseel, K. Rulik, M. Trigiante and A. Van Proeyen \cite{titsusataku} axiomatized the Tits-Satake projection
for all homogeneous special geometries classified, as we said, in \cite{vandersuppa,specHomgeoA1,specHomgeoA2,deWit:1995tf,Cortes}. The authors of
\cite{titsusataku} based their formulation of the projection on the intrinsic definition of the \emph{paint group} as the group of outer
automorphisms of the solvable transitive group of motion of the homogeneous manifold. This is the theory that will be explained in this article
and that constitutes an essential part of the PGTS backbone. Up to our knowledge, the very notion of paint group as a main component of the Tits
Satake projection theory  was never developed in the mathematical literature prior to the papers \cite{noipainted,titsusataku}.

\paragraph{\sc The TS projection in PGTS theory.}
In the PGTS theory resulting from  assembling the so far mentioned ingredients, the Tits-Satake subalgebra and the paint group are going to play
an essential role. As we explain in the sequel, the noncompact symmetric spaces $\mathcal{M}\, =\,\mathrm{U/H}$ have a   well defined
distance function,   which is the arc  length  of the unique  always
existing geodesic that passes through two given points. The distance function is related with the fundamental harmonic $\mathfrak{harm}_0$,
corresponding to the eigenfunctions of the quadratic \emph{Laplace Beltrami operator} $\Delta_{LB}$ belonging to the lowest eigenvalue
$\lambda_0$ of its (discrete) spectrum, and it is a function only of the coordinates of the Tits-Satake submanifold $\mathcal{M}_{\mathrm{TS}}\,
\equiv \, \mathrm{U_{TS}/H_{TS}} \subset \mathcal{M}$ and of a  number $n(r_{n.c.})$, depending on the chosen no-compact rank $r_{n.c.}$ of the
paint group invariants $\mathbf{U}_i$, ($i=1,\dots,n(r_{n.c.})$.  For $r_{n.c.}\leq 4$, the number  $n(r_{n.c.})$ is of the order of
unity or few tens, even if the dimension of $\mathcal{M}$ can be of the order of several thousands. Hence the PGTS theory has a built-in mechanism
of dimensional reduction in the representation of Big Data that are grouped into \emph{similarity classes} corresponding to level set
hypersurfaces defined by\footnote{Later on we will be more precise about the complementary notions of
$\mathrm{G_{subpaint}}\subset \mathrm{G_{Paint}}$ and $\mathrm{G_{subTS}}\subset \mathrm{G_{TS}}$ groups.}.
\begin{equation}\label{lelivelle}
 \mathcal{ M} \, \supset \, \Sigma_{\lambda_1\dots\lambda_{n(r_{n.c.})}} \, \equiv \, \left\{ p \in  \mathcal{M} \, \mid \,
 \mathbf{U}_i(p) \, = \, \lambda_i, \quad i= 1,\dots, n(r_{n.c.})\right\}.
\end{equation}

\subsection{The new ingredients of the PGTS theory introduced in this paper}
\label{mittel} The essential new ingredient introduced in this paper, which is based on previous knowledge, yet, once formulated in an intrinsic
way, leads to a new vision and to   new results, is the distance function expressed by means of a square norm of solvable group elements.
As remarked above, the Definition \ref{Alexnormalmanif} of \emph{normal homogeneous Riemannian spaces}, due to Alekseevsky, arises from the
notion of normed solvable Lie algebra as given in definition \ref{definormosolvalg}. The norm on the solvable Lie algebra induces a metric on the
group manifold $\mathcal{S}\, =  \, \exp[Solv]$ and,  if $\mathcal{S}$ is metrically equivalent to a symmetric space $\mathrm{U/H}$\footnote{The
normal Riemannian homogenous spaces that are not symmetric spaces are a few families that for the moment we disregard.}, which is the case of
our interest, the geodesic equations form an integrable system yielding, as stated above, a unique,  explicitly calculable, geodesic passing
through any two points. This allows us to go one step beyond Alekseevskii by introducing the more restrictive notion of normed solvable Lie
group.

\begin{definizione}\label{definormedsolgroupB}
A solvable Lie group $\mathcal{S}$ is a \emph{normed solvable Lie group} if it admits a function $\mathrm{N^2} \colon  \mathcal{S} \to \mathbb{R}_+$
 such that $\mathrm{N^2}(s)=0$ if and only if $s=e$, the identity element of $\mathcal S$,
and there exists a Riemannian metric $g$ on $\mathcal{S}$ such that the associated geodesic distance satisfies
 \begin{equation}\label{rovescio}
   d^2_g(s_1,s_2) \, = \, \mathrm{N^2}(s_1^{-1} \cdot s_2)
 \end{equation}
\end{definizione}
The function $\mathrm{N^2}$ is said to be a \emph{norm function.}
All noncompact symmetric spaces $\mathrm{U/H}$,  where the Lie algebra $\mathbb{U}$ of $\mathrm{U}$ is a noncompact real section of a complex
simple Lie algebra $\mathbb{G}$,  are \emph{not only normal homogeneous spaces} but also \emph{normed solvable Lie groups}.

The elaboration of algorithmic procedures for the \emph{explicit calculation} of the norm function in terms of \emph{Paint or subpaint group
invariants} is one of the main goals of PGTS in view of applications to Data Science. As we are going to see, the basic calculation expresses the
norm function in terms of $\mathrm{arccosh}(\Lambda_i)$ where $\Lambda_i$ are the $r_{(n.c.)}$ (all non negative) roots of an algebraic equation
of order $r_{(n.c.)}$,   the noncompact rank. This is the secular equation for the eigenvalues of an $r_{(n.c.)}\times r_{(n.c.)}$
matrix $\mathfrak{m}(\Upsilon)$ depending on the solvable coordinates that parameterize the solvable group elements.  In view of Galois' theorem we
are guaranteed to have explicit algebraic expressions for the norm function up to $r_{(n.c.)}\, = \,4 $. Surprizingly  the upper limit
$r_{(n.c.)}\leq 4$ is also the limit that is automatically satisfied by all $D=4$ supergravity relevant $\mathrm{U/H}$ manifolds if the number of
supersymmetries is $\mathcal{N}=2$ (see Table 5.4 at page 235 of \cite{advancio}). For $\mathcal{N}=3,4,5,6$  one goes up to $r_{(n.c.)} \, = \,
6$ and for the unique theory $\mathcal{N}=7,8$ one has
 $r_{(n.c.)} \, = \, 7$. In his original classification of quaternionic K\"ahler manifolds, also the   Alekseevsky  implicitly chose
 $r_{(n.c.)} \, = \, 4$.
\par
These remarks all conspire to suggest that in a first application to Data Science of PGTS theory one should better   consider only
the cases $r_{(n.c.)} \, = \,1,2,3,4$, an upper bound on $r_{(n.c.)}$ meeting also the main desire of mapping Big Dimensional Data to spaces of
much smaller dimensions (the Tits-Satake submanifold) plus a small bunch of subpaint group invariant qualifiers.
\paragraph{An inspiring chinese box scheme for $r_{(n.c.).} \leq 4$.} \label{cinesino} The     rewarding result one obtains adopting the aforementioned upper bound on
$r_{(n.c.)}$ is that, apart from the easily computable Paint and subpaint invariants, the geometry of the four relevant Tits-Satake mnifolds can
be completely analyzed, relying on special K\"ahler geometry and on the $c$-map, in terms of only two building blocks, namely, the Poincar\'e-Lobachevsky hyperbolic plane $\mathbb{H}_2 = \mathrm{PSL}(2,\mathbb{R})/\mathrm{SO(2)} $ and the upper Siegel plane of order 2
$\mathbb{IS}_2 = \mathrm{Sp(4,\mathbb{R})}/\mathrm{SU(2)} $, as we shall explain in later sections in full detail.
\paragraph{Three viewpoints on the Tits-Satake clustering mechanism.}
Another   new ingredient summarized in section \ref{balubbus} is provided by the three different and complementary viewpoints on the
Tits-Satake grouping of points that correspond either to the original formulation in terms of a projection  and on the fibers,   or in the Grassmannian leaves, or still in normal subgroups of the solvable group.
\paragraph{An innovative group-theoretical scheme of discretization.}
A third new ingredient, illustrated in section \ref{plantageneti}, is an innovative method for determining discrete subgroups, usable for
tessellations inspired by the solvable group structure of the manifold. As an example of the method a  subgroup of $\mathrm{Sp(4,\mathbb{Z})}$ is constructed in section \ref{plantageneti}. Then in section \ref{generaloneSOpq} the example is extended to a fully general and exhaustive construction of parabolic subgroups of all groups $\mathrm{SO(p,q,\mathbb{Z})}$.
\paragraph{Summary.}
Altogether, apart from the inclusion of the new above quoted ingredients, what is     new and innovative is the conceptual re-organization of
the whole lore about the geometry of noncompact symmetric spaces $\mathrm{U/H}$, which is purged of much of its Supergravity roots and
reassembled in the perspective of possible applications to Data Science.

 Let us then turn to a thorough presentation of the PGTS theory and of its developments; the latter were not mentioned in the
introduction  since they will be more easily explained after all the building blocks have been laid down.
\section{Special K\"ahler and Special Quaternionic K\"ahler Geometry}
In view of its importance to our programme, we are going to provide a concise introduction to
\emph{Special K\"ahler Geometries} and   the \emph{$c$-map}; see \cite{FreLect,CRTVP,Lledo-et-al} for an introduction, \cite{Freed} for an exposition of Special K\"ahler Geometry for mathematicians and \cite{Cortes-Tulli} for an explanation of the $c$-map.
\subsection{Special K\"ahler Geometry (of the local type)}
\label{speckalgeosum}  Let us begin by
summarizing some relevant concepts and definitions that are propaedeutical to the main definition.
\subsubsection{Hodge--K\"ahler manifolds}
\def\mom{{M(k, \IC)}}
Consider a {\sl line bundle} ${\cal L} {\stackrel{\pi}{\longrightarrow}} {\cal M}$ over a K\"ahler manifold ${\cal M}$, i.e., a
holomorphic vector bundle of rank $r=1$.
Its only topological invariant is its first Chern class, which can be represented by the differential 2-form
\begin{equation}
c_1 ( {\cal L} ) \, =\, \o{i}{2} \, {\bar \partial} \, \left ( \, h^{-1} \, \partial \, h \, \right )\, = \,
\o{i}{2} \, {\bar \partial} \,\partial \, \mbox{log} \,  h \label{chernclass23}
\end{equation}
where   $h(z,{\bar z})$, a real function, is a local representative of    an hermitian fiber metric on ${\cal L}$.
eq.~(\ref{chernclass23}) is the starting point for the definition of Hodge--K\"ahler manifolds. A K\"ahler manifold ${\cal M}$ is a Hodge manifold
if and only if there exists a line bundle ${\cal L} {\stackrel{\pi}{\longrightarrow}} {\cal M}$ such that its first Chern class equals the
cohomology class of the K\"ahler two-form $\mathrm{K}$:\footnote{An equivalent condition is that the cohomology class of the K\"ahler form is integral.}
\begin{equation}
c_1({\cal L} )~=~\left [ \, \mathrm{K} \, \right ] \label{chernclass25}
\end{equation}
\par
In local terms,   if $\xi (z)$ is local section of $\mathcal L$, one can write
\begin{equation}
\mathrm{K}\, =\, \o{i}{2} \, g_{ij^{\star}} \, dz^{i} \, \wedge \, d{\bar z}^{j^{\star}} \, = \, \o{i}{2} \,
{\bar \partial} \,\partial \, \mbox{log} \,\parallel \, \xi (z) \,
\parallel^2
\label{chernclass26}
\end{equation}
Recalling the   expression of the K\"ahler metric in terms of the K\"ahler potential $ g_{ij^{\star}}\, =\, {\partial}_i \,
{\partial}_{j^{\star}} {\mathcal{K}} (z,{\bar z})$, it follows from eq.~(\ref{chernclass26}) that, if the manifold ${\cal M}$ is a Hodge manifold,
  the exponential of the K\"ahler potential can be interpreted as a metric $h(z,{\bar z}) \, = \, \exp \left ( {\cal K} (z,{\bar z})\right )$
on an appropriate line bundle ${\cal L}$.
\par
\subsubsection{Connection on the line bundle}
On any holomorphic line bundle ${\cal L}$ equipped with an hermitian metric there is a canonical hermitian connection, whose local connection forms are
\begin{equation}
\begin{array}{ccccccc}
{\theta}& \equiv & h^{-1} \, \partial  \, h = {\o{1}{h}}\, \partial_i h \, dz^{i} &; & {\bar \theta}& \equiv
& h^{-1} \, {\bar \partial}  \, h = {\o{1}{h}} \,
\partial_{i^\star} h  \,
d{\bar z}^{i^\star} . \cr
\end{array}
\label{canconline}
\end{equation}
For the line-bundle advocated by the Hodge-K\"ahler structure we have
\begin{equation}
\left  [ \, {\bar \partial}\,\theta \,  \right ] \, = \, c_1({\cal L}) \, = \, [\mathrm{K}] \label{curvc1}
\end{equation}
and since the fiber metric $h$ can be identified with the exponential of the K\"ahler potential we obtain
\begin{equation}
\begin{array}{ccccccc}
{\theta}& = &  \partial  \,{\cal K} =  \partial_i {\cal K} dz^{i} & ; & {\bar \theta}& = &   {\bar \partial}
\, {\cal K} =
\partial_{i^\star} {\cal K}
d{\bar z}^{i^\star}\cr
\end{array}
\label{curvconline}
\end{equation}
To define special K\"ahler geometry,  in addition to the  line--bundle ${\cal L}$, one needs a flat holomorphic vector bundle ${\cal
SV} \, \longrightarrow \, {\cal M}_n$ of rank  $2\, n_V$,\footnote{In Supergravity $n_V$ is the total number of vector fields in the
theory.}  where $n_V = n+1$,  with $n$ the complex dimension of the base
manifold ${\cal M}_n$  whose coordinates  we denote $z^i$. As we shall leter show,
the flat vector bundle $\mathcal{SV}$ will be symplectic.

We shall make extensive use of covariant derivatives with respect to the canonical connection of the line bundle ${\cal L}$. Let us
review its normalization. If
$\mbox{exp}[f_{\alpha\beta}(z)]$ is the transition function between two local trivializations of the hermitian line bundle ${\cal L}
{\stackrel{\pi}{\longrightarrow}} {\cal M}$, the transition function in the corresponding principal $\mathrm{U(1)}$--bundle ${\cal U} \,
\longrightarrow {\cal M}$ is   $\mbox{exp}[{\rm i}{\rm Im}f_{\alpha\beta}(z)]$ and the K\"ahler potentials in two different charts are related
by ${\cal K}_\beta = {\cal K}_\alpha + f_{\alpha\beta}   + {\bar {f}}_{\alpha\beta}$. At the level of connections this correspondence is
formulated by setting: $\mbox{ $\mathrm{U(1)}$--connection}   \equiv   {\cal Q} \,  = \, \mbox{Im} \theta = -{\o{\rm i}{2}}   \left ( \theta -
{\bar \theta} \right)$. If we apply this formula to the case of the $\mathrm{U(1)}$--bundle ${\cal U} \, \longrightarrow \, {\cal M}$ associated
with the line--bundle ${\cal L}$ whose first Chern class equals the K\"ahler class, we get:
\begin{equation}
{\cal Q}  =    {\o{\rm i}{2}} \left ( \partial_i {\cal K} dz^{i} -
\partial_{i^\star} {\cal K}
d{\bar z}^{i^\star} \right ) \label{u1conect}
\end{equation}
Let ${\cal U}^p$ be the principal bundle associated with the hermitian line bundle $\mathcal L^p$, and let $\Phi(z,\bar z)$ be a section of it. Its
covariant derivative is $ \nabla \Phi = (d - i p {\cal Q}) \Phi $ or, in components,
\begin{equation}
\begin{array}{ccccccc}
\nabla_i \Phi &=&
 (\partial_i + {1\over 2} p \partial_i {\cal K}) \Phi &; &
\nabla_{i^*}\Phi &=&(\partial_{i^*}-{1\over 2} p \partial_{i^*} {\cal K}) \Phi \cr
\end{array}
\label{scrivo2}
\end{equation}
A covariantly holomorphic section of ${\cal U}$ is defined by the equation: $ \nabla_{i^*} \Phi = 0  $. We can easily map each  section $\Phi (z,
\bar z)$ of ${\cal U}^p$ into a  section of the line--bundle ${\cal L}$ by setting:
\begin{equation}
\tilde{\Phi} = e^{-p {\cal K}/2} \Phi  \,   . \label{mappuccia}
\end{equation}
  With this position we obtain:
\begin{equation}
\begin{array}{ccccccc}
\nabla_i    \tilde{\Phi}&    =& (\partial_i   +   p   \partial_i  {\cal K}) \tilde{\Phi}& ; &
\nabla_{i^*}\tilde{\Phi}&=& \partial_{i^*} \tilde{\Phi}\cr
\end{array}
\end{equation}
Under the map of eq.~(\ref{mappuccia}) covariantly holomorphic sections of ${\cal U}$ flow into holomorphic sections of ${\cal L}$ and viceversa.
\subsubsection{Special K\"ahler manifolds (of the local type)}
We are now ready to give the first of two equivalent definitions of special K\"ahler manifolds.
\begin{definizione}
\label{defspecial}
\label{defispeckelA} A Hodge-K\"ahler manifold is {\bf Special K\"ahler (of the local type)} if there exists a completely symmetric holomorphic
 $W_{i j k}$ of $(T^\star{\cal M})^3 \otimes {\cal L}^2$ such that the
following identities are satisfied:
\begin{eqnarray}
\partial_{m^*}   W_{ijk}& =& 0   \quad   \partial_m  W_{i^*  j^*  k^*}
=0 \nonumber \\
\nabla_{[m}      W_{i]jk}& =&  0
\quad \nabla_{[m}W_{i^*]j^*k^*}= 0 \nonumber \\
{\cal R}_{i^*j\ell^*k}& =&  g_{\ell^*j}g_{ki^*} +g_{\ell^*k}g_{j i^*} - e^{2 {\cal K}} W_{i^* \ell^* s^*}
W_{t k j} g^{s^*t} \label{specialone}
\end{eqnarray}
where $\mathcal R$ is the Riemann tensor of the Levi--Civita connection.
\end{definizione}
In the above equations $\nabla$ denotes the covariant derivative with respect to both the Levi--Civita and the $\mathrm{U(1)}$ holomorphic
connection of eq.~(\ref{u1conect}). In the case of $W_{ijk}$, the $\mathrm{U(1)}$ weight is $p = 2$.
\par
Out of the $W_{ijk}$ we can construct covariantly holomorphic sections of weight 2 and - 2 by setting:
\begin{equation}
C_{ijk}\,=\,W_{ijk}\,e^{  {\cal K}}  \quad ; \quad C_{i^\star j^\star k^\star}\,=\,W_{i^\star j^\star
k^\star}\,e^{  {\cal K}} \label{specialissimo}
\end{equation}
The flat bundle mentioned in the previous subsection apparently does not appear in this definition of special geometry. Yet it is there. It is
indeed the essential ingredient in the second definition whose equivalence to the first we shall shortly provide.
\par
%%%%%%%%%%%%%%%%%%%%%%%%%%%%%%%%%%%%%%%%%%%%%%%%%%
Let ${\cal L} {\stackrel{\pi}{\longrightarrow}} {\cal M}$ denote the complex line bundle whose first Chern class equals the cohomology class of
the K\"ahler form $\mathrm{K}$ of an $n$-dimensional Hodge--K\"ahler manifold ${\cal M}$. Let ${\cal SV} \, \longrightarrow \,{\cal M}$ denote a
holomorphic flat vector bundle of rank $2n+2$ with structure group $\mathrm{Sp(2n+2,\mathbb{R})}$. Consider   tensor bundles of the type ${\cal
H}\,=\,{\cal SV} \otimes {\cal L}$. A typical holomorphic section of such a bundle will be denoted by ${\Omega}$ and will have the following
structure:
\begin{equation}\label{ololo}
  \Omega \, = \,\left( \begin{array}{c}
                   X^\Lambda \\
                   F_ \Sigma
                 \end{array}\right)
   \quad \Lambda,\Sigma =0,1,\dots,n
\end{equation}
By definition the transition functions between two local trivializations $U_i \subset {\cal M}$ and $U_j \subset {\cal M}$ of the bundle ${\cal
H}$ have the following form:
\begin{equation}\label{cambogia}
\left( \begin{array}{c}
                   X^\Lambda \\
                   F_ \Sigma
                 \end{array}\right)_i \, = \, e^{f_{ij}} M_{ij} \left( \begin{array}{c}
                   X^\Lambda \\
                   F_ \Sigma
                 \end{array}\right)_j
\end{equation}
where   $f_{ij}$ are holomorphic maps $U_i \cap U_j \, \rightarrow \,\IC $ while $M_{ij}$ is a constant $\mathrm{Sp(2n+2,\mathbb{R})}$ matrix. For
a consistent definition of the bundle the transition functions are obviously subject to the cocycle condition on a triple overlap:
$e^{f_{ij}+f_{jk}+f_{ki}} = 1 $ and $ M_{ij} M_{jk} M_{ki} = 1 $.
\par
Let ${\rm i}\langle\ \vert\ \rangle$ be the compatible hermitian metric on $\cal H$
\begin{equation}
{\rm i}\langle \Omega \, \vert \, \bar \Omega \rangle \, \equiv \,
- {\rm i} \Omega^\T \left(
\begin{array}{cc}
0 & \mathbf{1} \\
-\mathbf{1} & 0 \\
\end{array}
\right) \,{\bar \Omega} \label{compati}
\end{equation}
\bd We say that a Hodge--K\"ahler manifold ${\cal M}$ is {\bf special K\"ahler} if there exists a bundle ${\cal H}$ of the type described above
such that for some section $\Omega \, \in \, \Gamma({\cal H},{\cal M})$ the K\"ahler two form is given by:
\begin{equation}
\mathrm{K}= \o{\rm i}{2}
 \partial \bar \partial \, \mbox{\rm log} \, \left ({\rm i}\langle \Omega \,
 \vert \, \bar \Omega
\rangle \right )=\frac{i}{2}\,g_{i j^*}\,dz^i\wedge d\bar{z}^{j^*}  \label{compati1}
\end{equation}
\ed From the point of view of local properties, eq.~(\ref{compati1}) implies that we have an expression for the K\"ahler potential in terms of the
holomorphic section $\Omega$:
\begin{equation}
\mathcal{K}\,  = \,  -\mbox{log}\left ({\rm i}\langle \Omega \,
 \vert \, \bar \Omega
\rangle \right )\, =\, -\mbox{log}\left [ {\rm i} \left ({\bar X}^\Lambda F_\Lambda - {\bar F}_\Sigma
X^\Sigma \right ) \right ] \label{specpot}
\end{equation}
The relation between the two definitions of special manifolds is obtained by introducing a non--holomorphic section of the bundle ${\cal H}$
according to:
\begin{equation}\label{covholsec}
V \, = \, \left(
  \begin{array}{c}
    L^{\Lambda} \\
    M_\Sigma \\
  \end{array}
\right) \, \equiv \, e^{\mathcal{K}/2}\Omega \, = \, e^{\mathcal{K}/2}\,
\left(
  \begin{array}{c}
    X^{\Lambda} \\
    F_\Sigma \\
  \end{array}
\right)
\end{equation}
so that eq.~(\ref{specpot}) becomes:
\begin{equation}
1 \, = \,  {\rm i}\langle V  \,
 \vert \, \bar V
\rangle  \, = \,   {\rm i} \left ({\bar L}^\Lambda M_\Lambda - {\bar M}_\Sigma L^\Sigma \right )
\label{specpotuno}
\end{equation}
Since $V$ is related to a holomorphic section by eq.~(\ref{covholsec}) it immediately follows that:
\begin{equation}
\nabla_{i^\star} V \, = \, \left ( \partial_{i^\star} - {\o{1}{2}}
\partial_{i^\star}{\cal K} \right ) \, V \, = \, 0
\label{nonsabeo}
\end{equation}
On the other hand, from eq.~(\ref{specpotuno}), defining:
\begin{eqnarray}
U_i  & = &  \nabla_i V  =   \left ( \partial_{i} + {\o{1}{2}}
\partial_{i}{\cal K} \right ) \, V   \equiv
\left(
\begin{array}{c}
f^{\Lambda}_{i} \\
h_{\Sigma\vert i} \\
\end{array}
\right)
\nonumber\\
{\bar U}_{i^\star}  & = &  \nabla_{i^\star}{\bar V}  =   \left ( \partial_{i^\star} + {\o{1}{2}}
\partial_{i^\star}{\cal K} \right ) \, {\bar V}   \equiv
\left(
\begin{array}{c}
{\bar f}^{\Lambda}_{i^\star} \\
{\bar h}_{\Sigma\vert i^\star} \\
\end{array}
\right)
\label{uvector}
\end{eqnarray}
it follows that:
\begin{equation}
\label{ctensor} \nabla_i U_j  = {\rm i} C_{ijk} \, g^{k\ell^\star} \, {\bar U}_{\ell^\star}
\end{equation}
where $\nabla_i$ denotes the covariant derivative containing both the Levi--Civita connection on the bundle ${\cal TM}$ and the canonical
connection $\theta$ on the line bundle ${\cal L}$. In eq.~(\ref{ctensor}) the symbol $C_{ijk}$ denotes a covariantly holomorphic (
$\nabla_{\ell^\star}C_{ijk}=0$) section of the bundle ${\cal TM}^3\otimes{\cal L}^2$ that is totally symmetric in its indices. This tensor can be
identified with the tensor of eq.~(\ref{specialissimo}) appearing in eq.~(\ref{specialone}). Alternatively, the set of differential equations:
\begin {eqnarray}
&&\nabla _i V  = U_i\nonumber\\
 && \nabla _i U_j = {\rm i} C_{ijk} g^{k \ell^\star} U_{\ell^\star}\nonumber\\
 && \nabla _{i^\star} U_j = g_{{i^\star}j} V\nonumber\\
 &&\nabla _{i^\star} V =0 \label{defaltern}
\end{eqnarray}
with V satisfying eqs.~ (\ref{covholsec}, \ref {specpotuno}) give yet another definition of special geometry. In particular it is easy to find
eq.~(\ref{specialone}) as integrability conditions of(\ref{defaltern})\footnote{We omit the detailed proof that from eqs.~ (\ref{defaltern}) one
obtains eq.~(\ref{specialone}). The essential link between the two formulations resides in the second of eqs.~ (\ref{defaltern}) which identifies the
tensor $C_{ijk}$ with the expression of the derivative of $U_i$ in terms of the same objects $U_k$. }.
\subsubsection{The symmetric matrix $\mathcal{N}_{\Lambda\Sigma}$ in Special K\"ahler Geometry}
\label{scrittaN} Another essential item intrinsically associated with Special K\"ahler Geometry is a complex symmetric matrix
$\mathcal{N}_{\Lambda\Sigma}$  which in Supergravity has its own motivation for its existence, yet from the mathematical viewpoint it is very much
significant that the same $\mathcal{N}_{\Lambda\Sigma}$ constitutes an integral part of the Special Geometry set up and might prove of relevance
also in Data Science applications. In any case it is an essential ingredient in the $c$-map construction. We provide its general definition in the
following lines. Explicitly $\mathcal{N}_{\Lambda\Sigma}$ which, in relation to its interpretation in the case of Calabi-Yau threefolds, is named
the {\it period matrix}, is defined by means of the following relations:
\begin{equation}
{\bar M}_\Lambda = {{\cal N}}_{\Lambda\Sigma}{\bar L}^\Sigma \quad ; \quad h_{\Sigma\vert i} = { {\cal
N}}_{\Lambda\Sigma} f^\Sigma_i \label{etamedia}
\end{equation}
which can be solved introducing the two $(n+1)\times (n+1)$ vectors
\begin{equation}
f^\Lambda_I = \twovec{f^\Lambda_i}{{\bar L}^\Lambda} \quad ; \quad h_{\Lambda \vert I} =  \twovec{h_{\Lambda
\vert i}}{{\bar M}_\Lambda} \label{nuovivec}
\end{equation}
and setting:
\begin{equation}
{{\cal N}}_{\Lambda\Sigma}= h_{\Lambda \vert I} \circ \left ( f^{-1} \right )^I_{\phantom{I} \Sigma}
\label{intriscripen}
\end{equation}
\par
Let us now consider the case where the Special K\"ahler manifold $\mathcal{SK}_n$ of complex dimension $n$ has some isometry group
$\mathrm{U}_{\mathcal{SK}}$. Compatibility with the Special Geometry structure requires the existence of a $2n+2$-dimensional symplectic
representation of such a group that we name the $\mathbf{W}$ representation. In other words,   a symplectic embedding  of the isometry group
$\mathcal{SK}_n$
\begin{equation}
  \mathrm{U}_{\mathcal{SK}} \mapsto \mathrm{Sp(2n+2, \mathbb{R})}
\label{sympembed}
\end{equation}
necessarily exists, such that for each element $\xi \in \mathrm{U}_{\mathcal{SK}}$ we have its representation by means of a suitable real
symplectic matrix:
\begin{equation}
  \xi \mapsto \Lambda_\xi \equiv \left( \begin{array}{cc}
     A_\xi & B_\xi \\
     C_\xi & D_\xi \
  \end{array} \right)
\label{embeddusmatra}
\end{equation}
satisfying the defining relation (in terms of the symplectic skew-symmetric metric $\mathbb{C}$):
\begin{equation}
  \Lambda_\xi ^T \, \underbrace{\left( \begin{array}{cc}
     \mathbf{0}_{n \times n}  & { \mathbf{1}}_{n \times n} \\
     -{ \mathbf{1}}_{n \times n}  & \mathbf{0}_{n \times n}  \
  \end{array} \right)}_{ \equiv \, \mathbb{C}} \, \Lambda_\xi = \underbrace{\left( \begin{array}{cc}
     \mathbf{0}_{n \times n}  & { \mathbf{1}}_{n \times n} \\
     -{ \mathbf{1}}_{n \times n}  & \mathbf{0}_{n \times n}  \
  \end{array} \right)}_{\mathbb{C}}
\label{definingsympe}
\end{equation}
which implies the following relations on the $n \times n$ blocks:
\begin{eqnarray}
A^T_\xi \, C_\xi - C^T_\xi \, A_\xi & = & 0 \nonumber\\
A^T_\xi \, D_\xi - C^T_\xi \, B_\xi& = & \mathbf{1}\nonumber\\
B^T_\xi \, C_\xi - D^T_\xi \, A_\xi& = & - \mathbf{1}\nonumber\\
B^T_\xi \, D_\xi - D^T_\xi \, B_\xi & =  & 0 \label{symplerele}
\end{eqnarray}
Under an element of the isometry group the symplectic section $\Omega$ of Special Geometry transforms as follows:
\begin{equation}
\Omega\left( \xi \, \cdot \, z\right) \, = \, \Lambda_\xi \, \Omega\left ( z \right )
\end{equation}
As a consequence of its definition, under the same isometry the matrix ${\cal N}$ transforms  by means of a generalized linear fractional
transformation:
\begin{equation}
  \mathcal{N}\left(\xi \cdot z,\xi \cdot \bar{z}\right) = \left(  C_\xi + D_\xi \, \mathcal{N}(z,\bar{z})\right)
  \, \left( A_\xi + B_\xi \,\mathcal{N}(z,\bar{z})\right) ^{-1}
\label{Ntransfa}
\end{equation}
%%%%%%%%%%%%%%%%%%%%%%%%%%%%%%%%%%%%%%%%%%%%%%%%%%%%%%%%%%%%%%%%
%%%%%%%%%%%%%%%%%%%%%%%%%%%%%%%%%%%%%%%%%%%%%%%%%%%%%%%%%%%%%%%%
\subsection{Quaternionic K\"ahler versus HyperK\"ahler manifolds}
Next we provide the definition of Quaternionic K\"ahler manifolds that we compare with that of HyperK\"ahler manifolds.
\par
Both a Quaternionic K\"ahler or a HyperK\"ahler manifold $\mathcal{QM}$ is a $4 m$-dimensional real manifold endowed with a metric $h$:
\begin{equation}
d s^2 = h_{u v} (q) d q^u \otimes d q^v   \quad ; \quad u,v=1,\dots, 4  m \label{qmetrica}
\end{equation}
and three complex structures
\begin{equation}
(J^x) \,:~~ T(\mathcal{QM}) \, \longrightarrow \, T(\mathcal{QM}) \qquad \quad (x=1,2,3)
\end{equation}
that satisfy the quaternionic algebra
\begin{equation}
J^x J^y = - \delta^{xy} \, \bfone \,  +  \, \epsilon^{xyz} J^z
%\label{quaternionetta}
\label{quatalgebra}
\end{equation}
and respect to which the metric is hermitian:
\begin{equation}
\forall   \mbox{\bf X} ,\mbox{\bf Y}  \in   T\mathcal{QM}   \,: \quad h \left( J^x \mbox{\bf X}, J^x
\mbox{\bf Y} \right )   = h \left( \mbox{\bf X}, \mbox{\bf Y} \right ) \quad \quad
  (x=1,2,3)
\label{hermit}
\end{equation}
From eq.~ (\ref{hermit}) it follows that one can introduce a triplet of 2-forms
\begin{equation}
\begin{array}{ccccccc}
K^x& = &K^x_{u v} d q^u \wedge d q^v & ; & K^x_{uv} &=&   h_{uw} (J^x)^w_v \cr
\end{array}
\label{iperforme}
\end{equation}
that provide the generalization of the concept of K\"ahler form occurring in  the complex case. The triplet $K^x$ is named the {\it HyperK\"ahler}
form. It is an $\mathrm{SU}(2)$ Lie--algebra valued 2--form  in the same way as the K\"ahler form is a $\mathrm{U(1)}$ Lie--algebra valued
2--form. In the complex case the definition of K\"ahler manifold involves the statement that the K\"ahler 2--form is closed. At the same time in
Hodge--K\"ahler manifolds  the K\"ahler 2--form can be identified with the curvature of a line--bundle which in the case of rigid supersymmetry is
flat. Similar steps can be taken also here and lead to two possibilities: either HyperK\"ahler or Quaternionic K\"ahler manifolds.
\par
Let us  introduce a principal $\mathrm{SU}(2)$--bundle ${\cal SU}$ as follows:
\begin{equation}
{\cal SU} \, \longrightarrow \, \mathcal{QM} \label{su2bundle}
\end{equation}
Let $\omega^x$ denote a connection on such a bundle. To obtain either a HyperK\"ahler or a Quaternionic K\"ahler manifold we must impose the
condition that the HyperK\"ahler 2--form is covariantly closed with respect to the connection $\omega^x$:
\begin{equation}
\nabla K^x \equiv d K^x + \epsilon^{x y z} \omega^y \wedge K^z    \, = \, 0 \label{closkform}
\end{equation}
The only difference between the two kinds of geometries resides in the structure of the ${\cal SU}$--bundle.
\begin{definizione}
\label{defihypkel}
 A
HyperK\"ahler manifold is a $4 m$--dimensional manifold with the structure described above and such that the ${\cal SU}$--bundle is {\bf flat}
\end{definizione}
 Defining the ${\cal SU}$--curvature by:
\begin{equation}
\Omega^x \, \equiv \, d \omega^x + {1\over 2} \epsilon^{x y z} \omega^y \wedge \omega^z \label{su2curv}
\end{equation}
in the HyperK\"ahler case we have:
\begin{equation}
\Omega^x \, = \, 0 \label{piattello}
\end{equation}
Viceversa
\begin{definizione}
\label{defiquatkel}
 A Quaternionic K\"ahler manifold is a $4 m$--dimensional manifold with the structure described above and such that
the curvature of the ${\cal SU}$--bundle is proportional to the HyperK\"ahler 2--form
\end{definizione}
Hence, in the quaternionic case we can
write:
\begin{equation}
\Omega^x \, = \, { {\lambda}}\, K^x \label{piegatello}
\end{equation}
where $\lambda$ is a non vanishing real number.
\par
As a consequence of the above structure the manifold $\mathcal{QM}$ has a holonomy group of the following type:
\begin{eqnarray}
{\rm Hol}(\mathcal{QM})&=& \mathrm{SU}(2)\otimes \mathrm{H} \quad (\mbox{Quaternionic K\"ahler}) \nonumber \\
{\rm Hol}(\mathcal{QM})&=& \bfone \otimes \mathrm{H} \quad (\mbox{HyperK\"ahler}) \nonumber \\ \mathrm{H} &
\subset & \mathrm{Sp (2m,\mathbb{R}) }\label{olonomia}
\end{eqnarray}
In both cases, introducing flat indices $\{A,B,C= 1,2\} \{\alpha,\beta,\gamma = 1,.., 2m\}$  that run, respectively, in the fundamental
representation of $\mathrm{SU}(2)$ and of $\mathrm{Sp}(2m,\mathbb{R})$, we can find a vielbein 1-form
\begin{equation}
{\cal U}^{A\alpha} = {\cal U}^{A\alpha}_u (q) d q^u \label{quatvielbein}
\end{equation}
such that
\begin{equation}
h_{uv} = {\cal U}^{A\alpha}_u {\cal U}^{B\beta}_v \mathbb{C}_{\alpha\beta}\epsilon_{AB} \label{quatmet}
\end{equation}
where $\mathbb{C}_{\alpha \beta} = - \mathbb{C}_{\beta \alpha}$ and $\epsilon_{AB} = - \epsilon_{BA}$ are, respectively, the flat
$\mathrm{Sp}(2m)$ and $\mathrm{Sp}(2) \sim \mathrm{SU}(2)$ invariant metrics. The vielbein ${\cal U}^{A\alpha}$ is covariantly closed with respect
to the $\mathrm{SU}(2)$-connection $\omega^z$ and to some $\mathrm{Sp}(2m,\mathbb{R})$-Lie algebra valued connection $\Delta^{\alpha\beta} =
\Delta^{\beta \alpha}$:
\begin{eqnarray}
\nabla {\cal U}^{A\alpha}& \equiv & d{\cal U}^{A\alpha} +{i\over 2} \omega^x (\epsilon
\sigma_x\epsilon^{-1})^A_{\phantom{A}B} \wedge{\cal U}^{B\alpha} +  \Delta^{\alpha\beta} \wedge
{\cal U}^{A\gamma} \mathbb{C}_{\beta\gamma} =0 \label{quattorsion}
\end{eqnarray}
\noindent where $(\sigma^x)_A^{\phantom{A}B}$ are the standard Pauli matrices. Furthermore ${ \cal U}^{A\alpha}$ satisfies  the reality condition:
\begin{equation}
{\cal U}_{A\alpha} \equiv ({\cal U}^{A\alpha})^* = \epsilon_{AB} \mathbb{C}_{\alpha\beta} {\cal U}^{B\beta}
\label{quatreality}
\end{equation}
eq.~(\ref{quatreality})  defines  the  rule to lower the symplectic indices by means   of  the  flat symplectic   metrics $\epsilon_{AB}$   and
$\mathbb{C}_{\alpha \beta}$. More specifically we can write a stronger version of eq.~ (\ref{quatmet}) \cite{baggherozzo}:
\begin{eqnarray}
({\cal U}^{A\alpha}_u {\cal U}^{B\beta}_v + {\cal U}^{A\alpha}_v {\cal
 U}^{B\beta}_u)\mathbb{C}_{\alpha\beta}&=& h_{uv} \epsilon^{AB}\nonumber\\
 \label{piuforte}
\end{eqnarray}
\noindent We have also the inverse vielbein ${\cal U}^u_{A\alpha}$ defined by the equation
\begin{equation}
{\cal U}^u_{A\alpha} {\cal U}^{A\alpha}_v = \delta^u_v \label{2.64}
\end{equation}
Flattening a pair of indices of the Riemann tensor ${\cal R}^{uv}_{\phantom{uv}{ts}}$ we obtain
\begin{equation}
{\cal R}^{uv}_{\phantom{uv}{ts}} {\cal U}^{\alpha A}_u {\cal U}^{\beta B}_v = -\,{{\rm i}\over 2}
\Omega^x_{ts} \epsilon^{AC}
 (\sigma_x)_C^{\phantom {C}B} \mathbb{C}^{\alpha \beta}+
 \mathbb{R}^{\alpha\beta}_{ts}\epsilon^{AB}
\label{2.65}
\end{equation}
\noindent where $\mathbb{R}^{\alpha\beta}_{ts}$ is the curvature 2-form of the $\mathrm{Sp}(2m) $ connection:
\begin{equation}
d \Delta^{\alpha\beta} + \Delta^{\alpha \gamma} \wedge \Delta^{\delta \beta} \mathbb{C}_{\gamma \delta}
\equiv \mathbb{R}^{\alpha\beta} = \mathbb{R}^{\alpha \beta}_{ts} dq^t \wedge dq^s \label{2.66}
\end{equation}
eq.~ (\ref{2.65}) is the explicit statement that the Levi Civita connection associated with the metric $h$ has a holonomy group contained in
$\mathrm{SU}(2) \otimes \mathrm{Sp}(2m)$. Consider now eqs.~ (\ref{quatalgebra}), (\ref{iperforme}) and (\ref{piegatello}). We easily deduce the
following relation:
\begin{equation}
h^{st} K^x_{us} K^y_{tw} = -   \delta^{xy} h_{uw} +
  \epsilon^{xyz} K^z_{uw}
\label{universala}
\end{equation}
that holds true both in the HyperK\"ahler and in the quaternionic case. In the latter case, using eq.~ (\ref{piegatello}), eq.~ (\ref{universala})
can be rewritten as follows:
\begin{equation}
h^{st} \Omega^x_{us} \Omega^y_{tw} = - \lambda^2 \delta^{xy} h_{uw} + \lambda \epsilon^{xyz} \Omega^z_{uw}
\label{2.67}
\end{equation}
eq.~(\ref{2.67}) implies that the intrinsic components of the curvature
 2-form $\Omega^x$ yield a representation of the quaternion algebra.
In the HyperK\"ahler case such a representation is provided only by the HyperK\"ahler form. In the quaternionic case we can write:
\begin{equation}
\Omega^x_{A\alpha, B \beta} \equiv \Omega^x_{uv} {\cal U}^u_{A\alpha} {\cal U}^v_{B\beta} = - i \lambda
\mathbb{C}_{\alpha\beta} (\sigma_x)_A^{\phantom {A}C}\epsilon _{CB} \label{2.68}
\end{equation}
\noindent Alternatively eq.~(\ref{2.68}) can be rewritten in an intrinsic form as
\begin{equation}
\Omega^x =\,-{\rm i}\, \lambda \mathbb{C}_{\alpha\beta} (\sigma _x)_A^{\phantom {A}C}\epsilon _{CB} {\cal
U}^{\alpha A} \wedge {\cal U}^{\beta B} \label{2.69}
\end{equation}
\noindent whence we also get:
\begin{equation}
{i\over 2} \Omega^x (\sigma_x)_A^{\phantom{A}B} = \lambda{\cal U}_{A\alpha} \wedge {\cal U}^{B\alpha}
\label{2.70}
\end{equation}
%%%%%%%%%%%%%%%%%%%%%%%%%%%%%%%%%%%%%
\subsubsection{The Quaternionic K\"ahler Geometry in the image of the $c$-map} \label{cmappusquat}
%%%%%%%%%%%%%%%%%%%%%%%%%%%%%%%%%%%%%%
Next we consider those Quaternionic K\"ahler manifolds that are in the image of the $c$-map.\footnote{Not all noncompact,  homogeneous
Quaternionic K\"ahler manifolds that are relevant to supergravity (which are \emph{normal}, i.e. exhibiting a solvable group of isometries having
a free and transitive action on it) are in the image of the $c$-map, the only exception being the quaternionic projective spaces
\cite{vandersuppa,Cecotti:1988ad}. In the case of Data Science applications, as we discuss  in  next sections, the Tits-Satake projections of low
$r$ relevant cases are in the image of the $c$-map.} This latter
\begin{equation}\label{cimappo}
    \mbox{$c$-map} \, \, : \,\, \mathcal{SK}_n \, \Longrightarrow \, \mathcal{QM}_{4n+4}
\end{equation}
is a universal construction that, starting from an arbitrary  Special K\"ahler manifold $\mathcal{SK}_n$ of complex dimension $n$, irrespectively
whether it is homogenoeus or not, leads to a unique Quaternionic K\"ahler manifold $\mathcal{QM}_{4n+4}$ of real dimension $4n+4$ which contains
$\mathcal{SK}_n$ as a submanifold. The precise modern definition of the $c$-map, originally introduced in \cite{catlantide1,catlantide2}, is
provided below.
\begin{definizione}
\label{deficimappus} Let $\mathcal{SK}_n$ be a special K\"ahler manifold whose complex coordinates we denote by $z^i$ and whose K\"ahler metric we
denote by $g_{ij^\star}$. Let moreover $\mathcal{N}_{\Lambda\Sigma}(z,{\bar z})$ be the symmetric period matrix defined by
eq.~(\ref{intriscripen}), introduce the following set of $4n+4$ coordinates:
\begin{equation}\label{finnico2}
    \left\{q^u \right\} \, \equiv \, \underbrace{\{U,a\}}_{\mbox{2 real}}\, \bigcup \,
    \underbrace{\underbrace{\{ z^i\}}_{\mbox{n complex}}}_{\mbox{2n real}} \, \bigcup\, \underbrace{\mathbf{Z}
    \, = \, \{ Z^\Lambda \, , \, Z_\Sigma \}}_{\mbox{(2n+2) real}}
\end{equation}
Let us further introduce the following $(\mathrm{2n+2})\times(\mathrm{2n+2}) $ matrix  ${\cal M}_4^{-1}$:
\begin{eqnarray}
\mathcal{M}_4^{-1} & = & \left(\begin{array}{c|c} {\mathrm{Im}}\mathcal{N}\, +\, {\mathrm{Re}}\mathcal{N} \,
{ \mathrm{Im}}\mathcal{N}^{-1}\, {\mathrm{Re}}\mathcal{N} &
\, -{\mathrm{Re}}\mathcal{N}\,{ \mathrm{Im}}\,\mathcal{N}^{-1}\\
\hline -\, { \mathrm{Im}}\mathcal{N}^{-1}\,{\mathrm{Re}}\mathcal{N}  & { \mathrm{Im}}\mathcal{N}^{-1} \
\end{array}\right) \label{inversem4}
\end{eqnarray}
which depends only on the coordinate of the Special K\"ahler manifold. The $c$-map image of $\mathcal{SK}_n$ is the unique Quaternionic K\"ahler
manifold $\mathcal{QM}_{4n+4}$ whose coordinates are the $q^u$ defined in (\ref{finnico2}) and whose metric is given by the following universal
formula
\begin{eqnarray}
ds^2_{\mathcal{QM}} &=&\frac{1}{4} \left(  d{U}^2+\, 4 g_{ij^\star} \,d{z}^j\, d{{\bar z}}^{j^\star} +
e^{-2\,U}\,(d{a}+{\bf Z}^T\mathbb{C}d{{\bf Z}})^2\,-\,2 \, e^{-U}\,d{{\bf Z}}^T\,\mathcal{M}_4^{-1}\,d{{\bf
Z}}\right) \nonumber\\
\label{geodaction}
\end{eqnarray}
\end{definizione}
The metric (\ref{geodaction}) has the following positive definite signature
\begin{equation}
\mbox{sign}\left[ds^2_{\mathcal{QM}}\right] \, = \, \left(\underbrace{+,\dots,+}_{4+4\mathrm{n}}\right)
\end{equation}
since the matrix $\mathcal{M}_4^{-1} $ is negative definite.
\par
In the case the Special K\"ahler pre-image is a symmetric space $\mathrm{U}_{\mathcal{SK}}/\mathrm{H}_{\mathcal{SK}}$, the manifold
$\mathcal{QM}$  turns out to be a symmetric space, $\mathrm{U}_{Q}/\mathrm{H}_{Q}$. We  come back to the issue of symmetric homogeneous
Quaternionic K\"ahler manifolds in section \ref{omosymmetro}
\subsubsection{The HyperK\"ahler two-forms and the $\su(2)$-connection} The reason
why we state that $\mathcal{QM}_{4n+4}$ is Quaternionic K\"ahler is that, by utilizing only the identities of Special K\"ahler Geometry, we can
construct the three complex structures $J_u^{x|v}$ satisfying the quaternionic algebra (\ref{quatalgebra}) the corresponding HyperK\"ahler
two-forms $K^x$  and the $\su(2)$ connection $\omega^x$ with respect to which they are covariantly constant.

The construction is extremely beautiful, it was found in \cite{cmappotto} and it is the following one.

Consider the K\"ahler connection $\mathcal{Q}$ defined by eq.~ (\ref{u1conect}) and furthermore introduce the following differential form:
\begin{equation}
\label{Phidiffe}
    \Phi \, = \, da + \mathbf{Z}^T \, \mathbb{C}\, \mathrm{d}\mathbf{Z}
\end{equation}
Next define the two dimensional representation of both the $\su(2)$ connection and of the HyperK\"ahler $2$-forms as it follows:
\begin{eqnarray}
  \omega  &=& \frac{\rm i}{\sqrt{2}}\,\sum_{x=1}^3 \, \omega^x \, \gamma_x \label{cunnettasu2}\\
  \mathbf{K} &=& \frac{\rm i}{\sqrt{2}}\,\sum_{x=1}^3 \, K^x \, \sigma_x  \label{HypKalmatra}
\end{eqnarray}
where $\gamma_x$ denotes a basis of $2\times 2$ euclidian $\gamma$-matrices for which we utilize the following basis which is  convenient in the
explicit calculations\footnote{The chosen $\gamma$-matrices are a permutation of the standard Pauli matrices divided by $\sqrt{2}$ and multiplied
by $\frac{\rm i}{2}$. They can be used as a basis of anti-hermitian generators for the $\su(2)$ algebra in the fundamental defining
representation.}:
\begin{eqnarray}
  \gamma_1&=& \left(
\begin{array}{ll}
 \frac{1}{\sqrt{2}} & 0 \\
 0 & -\frac{1}{\sqrt{2}}
\end{array}
\right) \nonumber\\
   \gamma_2&=& \left(
\begin{array}{ll}
 0 & -\frac{i}{\sqrt{2}} \\
 \frac{i}{\sqrt{2}} & 0
\end{array}
\right)\nonumber\\
 \gamma_3 &=& \left(
\begin{array}{ll}
 0 & \frac{1}{\sqrt{2}} \\
 \frac{1}{\sqrt{2}} & 0
\end{array}
\right) \label{gamminiB}
\end{eqnarray}
These $\gamma$-matrices satisfy the following Clifford algebra:
\begin{equation}\label{cliffordus}
    \left\{ \gamma_x \, , \, \gamma_y \right \} \, = \, \delta^{xy} \, \mathbf{1}_{2 \times 2}
\end{equation}
and $\frac{\rm i}{2} \, \gamma_x$ provide a basis of generators of the $\su(2)$ algebra.
\par
Having fixed these conventions, the expression of the quaternionic $\su(2)$-connection in terms of Special Geometry structures is encoded in the
following expression for the $2\times 2$-matrix valued $1$-form $\omega$. Explicitly we have:
\begin{equation}\label{omegaSu2}
    \omega \, = \, \left( \begin{array}{cc}
                            -\frac{\rm i}{2} \, \mathcal{Q}  \, - \, \frac{\rm i}{4} \,e^{-U} \, \Phi  &
                            e^{-\frac{U}{2}} \, V^T \, \mathbb{C} \, \mathrm{d}\mathbf{Z}\\
                            - \, e^{-\frac{U}{2}} \, \overline{V}^T \, \mathbb{C} \, \mathrm{d}\mathbf{Z} &
                            \frac{\rm i}{2} \, \mathcal{Q}  \, + \, \frac{\rm i}{4} \,e^{-U} \, \Phi
                          \end{array}
    \right)
\end{equation}
where $V$ and $\overline{V}$ denote the covariantly holomorphic sections of Special geometry defined in eqs.~ (\ref{covholsec}). The curvature of
this connection is obtained from a straight-forward calculation:
\begin{eqnarray}
\label{K2per2}
  \mathbf{K} &\equiv& d\omega \, + \, \omega \, \wedge \, \omega \nonumber \\
  \null &=& \left(\begin{array}{cc}
                 \mathfrak{u}   & \mathfrak{v} \\
                    - \,\overline{\mathfrak{v}}& -\,\mathfrak{u}
                  \end{array}
   \right)
\end{eqnarray}
the independent $2$-form matrix elements being given by the following explicit formulas:
\begin{eqnarray}
\label{uvvb}
 \mathfrak{u} &=& -{\rm i} \frac{1}{2} \, K \, -\frac{1}{8} dS \,\wedge \, d\bar{S}\, - \, e^{-U} \, V^T \,
 \mathbb{C} \,\mathrm{d}\mathbf{Z}\, \wedge\, \bar{V}^T \, \mathbb{C} \,\mathrm{d}\mathbf{Z} \, - \,
 \frac{1}{4} \, e^{-U} \,\mathrm{d}\mathbf{Z}^T \, \wedge \, \mathbb{C} \, \mathrm{d}\mathbf{Z} \nonumber\\
  \mathfrak{v }&=& e^{-\frac{U}{2}} \left( \, DV^T \, \wedge \, \mathbb{C} \, \mathrm{d}\mathbf{Z}\, - \,
  \frac{1}{2} dS \, \wedge \,
  V^T \, \mathbb{C} \, \mathrm{d}\mathbf{Z}\right)\nonumber\\
  \overline{\mathfrak{v }} &=& e^{-\frac{U}{2}} \left( \, D\overline{V}^T \, \wedge \, \mathbb{C} \,
  \mathrm{d}\mathbf{Z}\, - \, \frac{1}{2} d\overline{S} \, \wedge \,
  \overline{V}^T \, \mathbb{C} \, \mathrm{d}\mathbf{Z}\right)
\end{eqnarray}
where
\begin{equation}\label{kalleforma}
    K \, = \, \frac{ {\rm i}}{2} \, g_{ij^\star} \, dz^i \, \wedge \, d\bar{z}^{j^\star}
\end{equation}
is the K\"ahler $2$-form of the Special K\"ahler submanifold and where we have used the following short hand notations:
\begin{eqnarray}
  dS &=& dU \, + \, {\rm i} \, e^{-U}\, \left(da \, + \, \mathbf{Z}^T \, \mathbb{C} \, \mathrm{d}\mathbf{Z}\right)
  \label{firbone1} \\
  d\overline{S} &=& dU \, - \, {\rm i} \, e^{-U}\, \left(da \, + \, \mathbf{Z}^T \, \mathbb{C} \,
  \mathrm{d}\mathbf{Z}\right) \label{firbone2} \\
  DV &=& dz^i \, \nabla_i V \label{firbone3} \\
  D\overline{V} &=& d\bar{z}^{i^\star} \, \nabla_{i^\star} V \label{firbone4}
\end{eqnarray}
The three HyperK\"ahler forms $K^x$ are easily extracted from eqs.~ (\ref{K2per2}-\ref{uvvb}) by collecting the coefficients of the $\gamma$-matrix
expansion and we need not to write their form which is immediately deduced. The relevant thing is that the components of $K^x$ with an index
raised through multiplication with the inverse of the quaternionic metric $h^{uv}$ exactly satisfy the algebra of quaternionic complex structures
(\ref{quatalgebra}). Explicitly we have:
\begin{eqnarray}
  K^x &=& - \, {\rm i} \, 4 \sqrt{2} \, \mbox{Tr} \, \left( \gamma^x \, \mathbf{K}\right) \, \equiv \,
  K^x_{uv} \, dq^u \,\wedge \, dq^v \nonumber\\
  J^{x|s}_u &=& K^x_{uv}\, h^{vs} \nonumber \\
  J^{x|s}_u \, J^{y|v}_s &=& - \delta^{ xy} \, \delta^v_u \, + \, \epsilon^{xyz} \, J^{z|v}_u \label{quatKvera}
\end{eqnarray}
The above formulas are not only the general proof that the Riemaniann manifold $\mathcal{QM}$ defined by the metric (\ref{geodaction}) is indeed a
Quaternionic K\"ahler  manifold, but they  also provide an algorithm to write in terms of Special Geometry structures  the tri-holomorphic moment
map of the principal isometries possessed by $\mathcal{QM}$; this property was of high relevance for gauged supergravity models, but it is not yet
clear whether it is going to play a role in data science applications.
\subsubsection{Homogeneous Symmetric Special Quaternionic
K\"ahler manifolds} \label{omosymmetro} When the Special K\"ahler manifold $\mathcal{SK}_n$ is a symmetric coset space, it turns out that the
metric (\ref{geodaction}) is actually the symmetric metric on an enlarged symmetric coset manifold
\begin{equation}\label{qcosetto}
    \mathcal{QM}_{4n+4} \, = \, \frac{\mathrm{U}_Q}{\mathrm{H}_Q} \, \supset \,
    \frac{\mathrm{U}_{\mathcal{SK}}}{\mathrm{H}_{\mathcal{SK}}}
 \end{equation}
\par
Naming $\Lambda[\mathfrak{g}]$ the $\mathbf{W}$-representation of any finite element of the $\mathfrak{g}\in\mathrm{U}_{\mathcal{SK}}$ group, we
have that the matrix $\mathcal{M}_4(z,\bar{z})$ transforms as follows:
\begin{equation}\label{traduco}
    \mathcal{M}_4\left( \mathfrak{g}\cdot z,\mathfrak{g}\cdot \bar{z}  \right)\, = \,
    \Lambda[\mathfrak{g}] \, \mathcal{M}_4\left(z,\bar{z}\right ) ]\, \Lambda^T[\mathfrak{g}]
\end{equation}
where $\mathfrak{g}\cdot z$ denotes the nonlinear action of $\mathrm{U}_{\mathcal{SK}}$ on the scalar fields. Since the space
$\frac{\mathrm{U}_{\mathcal{SK}}}{\mathrm{H}_{\mathcal{SK}}}$ is homogeneous, after  choosing any reference point $z_0$, all the others can be reached by
a suitable group element $\mathfrak{g}_z$ such that $\mathfrak{g}_z\cdot z_0 \, = \, z$, and we can write
\begin{equation}\label{turnaconto}
 \mathcal{M}_4^{-1}(z,\bar{z}) \,  = \,  \Lambda^T[\mathfrak{g}_z^{-1}] \,
 \mathcal{M}_4^{-1}(z_0,\bar{z}_0) ]\, \Lambda[\mathfrak{g}^{-1}_z] .
\end{equation}
This allows us to introduce a set of $4n+4$ vielbein defined in the following way:
\begin{equation}\label{filibaine}
    E^I_{\mathcal{QM}} \, = \, \frac{1}{2} \, \left\{ dU \, ,
    \, \underbrace{e^i(z)}_{2\,n} \, , \, e^{-U} \,\left(d{a}+{\bf Z}^T\mathbb{C}d{{\bf
Z}}\right) \, , \, \underbrace{e^{-\frac {U} {2}}\, \Lambda[\mathfrak{g}_z^{-1}] \,
\mathrm{d}\mathbf{Z}}_{2n+2} \right\}
\end{equation}
and rewrite the metric (\ref{geodaction}) as it follows:
\begin{equation}\label{cornish}
    ds^2_{\mathcal{QM}} \, = \, E^I_{\mathcal{QM}} \, \mathfrak{q}_{IJ} \, E^J_{\mathcal{QM}}
\end{equation}
where the quadratic symmetric constant tensor  $\mathfrak{q}_{IJ}$ has the following form:
\begin{equation}\label{quadrotta}
    \mathfrak{q}_{IJ} \, = \, \left( \begin{array}{c|c|c|c}
                                       1 & 0 & 0 & 0 \\
                                       \hline
                                       0 & \delta_{ij} & 0 & 0 \\
                                       \hline
                                       0 & 0 & 1 & 0 \\
                                       \hline
                                       0 & 0 & 0 & -\, 2\, \mathcal{M}_4^{-1}(z_0,\bar{z}_0)
                                     \end{array}
    \right)
\end{equation}
The above defined vielbein are endowed with a very special property namely they identically satisfy a set of Maurer Cartan equations:
\begin{equation}\label{MCSolv}
    dE^I_{\mathcal{QM}} \, - \, \frac{1}{2} f^I_{\phantom{I}JK} \, E^J_{\mathcal{QM}} \,
    \wedge \,  E^K_{\mathcal{QM}} \, = \, 0
\end{equation}
where $f^I_{\phantom{I}JK}$ are the structure constants of a solvable Lie algebra $\mathfrak{A}$ which can be identified as follows:
\begin{equation}\label{solvableGHalg}
    \mathfrak{A} \, = \, Solv\left( \frac{\mathrm{U}_\mathcal{Q}}{\mathrm{H}_\mathcal{Q}} \right)
\end{equation}
In the above equation $Solv\left( \frac{\mathrm{U}_\mathcal{Q}}{\mathrm{H}_\mathcal{\mathcal{Q}}} \right)$ denotes the Lie algebra of the solvable
group manifold metrically equivalent to the non-comapact coset manifold $\frac{\mathrm{U}_\mathcal{Q}}{\mathrm{H}_\mathcal{Q}}$ according to what
we explained in section \ref{leluka}. In the case ${\mathrm{U}_{\mathcal{SK}}}$ is a \emph{maximally split} real form of a complex Lie algebra,
then  also ${\mathrm{U}_{\mathcal{Q}}}$ is maximally split and  we have:
\begin{equation}\label{solvableGHalg2}
    Solv\left( \frac{\mathrm{U}_\mathcal{Q}}{\mathrm{H}_\mathcal{Q}} \right) \, =\,
    \mbox{Bor}\left ( \mathbb{U}_\mathcal{Q} \right)
\end{equation}
where $\mbox{Bor}\left ( \mathbb{U}_\mathcal{Q} \right)$ denotes the \emph{Borel subalgebra} of the semi-simple Lie algebra $\mathbb{G}$,
generated by its Cartan generators and by the step operators associated with all positive roots.
\par
According to the  theory summarized in section \ref{leluka} above, the very fact that the vielbein (\ref{filibaine}) satisfies the
Maurer Cartan equations of the Lie algebra  $Solv\left( \frac{\mathrm{U}_\mathcal{Q}}{\mathrm{H}_\mathcal{Q}} \right)$ implies that the metric
(\ref{cornish}) is the symmetric metric on the coset manifold $\frac{\mathrm{U}_\mathcal{Q}}{\mathrm{H}_\mathcal{Q}}$, which therefore admits
continuous isometries associated with all the generators of the Lie algebra $\mathbb{U}_\mathcal{Q}$. For the reader's convenience, the list of
Symmetric Special manifolds and of their Quaternionic K\"ahler  counterparts in the image of the $c$-map is recalled  in Table \ref{homomodelisti},
which reproduces the results of \cite{toineugenio}, according to which there is   a short list of Symmetric Homogeneous Special manifolds
comprising five discrete cases and two infinite series.
\begin{table}
\begin{center}
{\small
\begin{tabular}{||c|c||c||}
  \hline
   $\mathcal{SK}_n$  & $\mathcal{QM}_{4n+4}$ &    $\mbox{dim} \, \mathcal{SK}_n \, = \, $  \\
   Special K\"ahler manifold & Quaternionic K\"ahler manifold & $n$  \\
  \hline
\null & \null &\null \\
 $ \frac{\mathrm{SU(1,1)}}{\mathrm{U(1)}}$ & $ \frac{\mathrm{G_{2(2)}}}{\mathrm{SU(2)\times SU(2)}}$ & $n=1$\\
\null & \null &  \\
\hline
\null & \null &\null \\
  $ \frac{\mathrm{Sp(6,R)}}{\mathrm{SU(3)\times  U(1)}}$ & $ \frac{\mathrm{F_{4(4)}}}{\mathrm{USp(6)\times SU(2)}}$  &$n=6$\\
 \null & \null &  \\
 \null & \null &\null \\
\hline
\null & \null &\null \\
 $ \frac{\mathrm{SU(3,3)}}{\mathrm{SU(3)\times SU(3) \times U(1)}}$ & $
 \frac{\mathrm{E_{6(2)}}}{\mathrm{SU(6)\times SU(2)}}$    &$n=9$\\
 \null & \null &  \\
\null & \null &\null \\
\hline
\null & \null &\null \\
 $ \frac{\mathrm{SO^\star(12)}}{\mathrm{SU(6)\times U(1)}}$ & $
 \frac{\mathrm{E_{7(-5)}}}{\mathrm{SO(12)\times SU(2)}}$  & $n=15$ \\
\null & \null &  \\
\null & \null &\null \\
\hline
\null & \null &\null \\
$ \frac{\mathrm{E_{7(-25)}}}{\mathrm{E_{6(-78)} \times U(1)}}$ & $
\frac{\mathrm{E_{8(-24)}}}{\mathrm{E_{7(-133)}\times SU(2)}}$    &  $n=27$ \\
\null & \null &  \\
\hline
\null & \null &\null \\
 $ \frac{\mathrm{SL(2,\mathbb{R})}}{\mathrm{SO(2)}}\times\frac{\mathrm{SO(2,2+p)}}
 {\mathrm{SO(2)\times SO(2+p)}}$ & $ \frac{\mathrm{SO(4,4+p)}}{\mathrm{SO(4)\times SO(4+p)}}$   & $n=3+p$  \\
  \null & \null &  \\
\hline
\null & \null &\null \\
$ \frac{\mathrm{SU(p+1,1)}}{\mathrm{SU(p+1)\times U(1)}}$ & $
\frac{\mathrm{SU(p+2,2)}}{\mathrm{SU(p+2)\times SU(2)}}$    & $n=p+1$ \\
\null & \null &\null \\
\hline
\end{tabular}
} \caption{List of special K\"ahler symmetric spaces with their Quaternionic K\"ahler  $c$-map images. The
number $n$ denotes the complex dimension  of the Special K\"ahler preimage. On the other hand $4n+4$ is the
real dimension of the Quaternionic K\"ahler $c$-map image. \label{homomodelisti}}
\end{center}
\end{table}
\par
Inspecting the above results we immediately realize that the Lie algebra $\mathbb{U}_{Q}$ contains two universal Heisenberg subalgebras of
dimension $(2n+3)$, namely:
\begin{eqnarray}
\mathbb{U}_{\mathcal{Q}} \, \supset \, \mathbb{H}\mathrm{eis}_1 \, &=& \mbox{span}_\mathbb{R}\,
\left\{\mathbf{W}^{1\alpha} \, , \, \mathbb{Z}_1 \right\} \quad ; \quad  \mathbb{Z}_1  \, =\, L_+  \,
\equiv \, L^1\, + \, L^2 \nonumber\\
  \null &\null & \left[\mathbf{W}^{1\alpha}\, , \, \mathbf{\mathbf{W}}^{1\beta} \right]\, = \, -\,
  \frac{1}{2} \, \mathbb{C}^{\alpha\beta} \, \mathbb{Z}_1  \quad ; \quad \left[\mathbb{Z}_1\, , \,
  \mathbf{W}^{1\beta} \right]\, = \,0 \nonumber\\
  \label{Heinber1}\\
  \mathbb{U}_{\mathcal{Q}} \, \supset \, \mathbb{H}\mathrm{eis}_2 \, &=& \mbox{span}_\mathbb{R}\,
  \left\{\mathbf{W}^{2\alpha} \, , \, \mathbb{Z}_2 \right\} \quad ; \quad  \mathbb{Z}_2  \, =\, L_-  \,
  \equiv \, L^1\, - \, L^2 \nonumber\\
  \null &\null & \left[\mathbf{W}^{2\alpha}\, , \, \mathbf{W}^{2\beta} \right]\, = \, -\, \frac{1}{2} \,
  \mathbb{C}^{\alpha\beta} \, \mathbb{Z}_2  \quad ; \quad \left[\mathbb{Z}_2\, , \, \mathbf{W}^{2\beta} \right]\, =
  \,0 \nonumber\\
  \label{Heinber2}
\end{eqnarray}
The first of these Heisenberg subalgebras of isometries is the universal one, which exists for all Quaternionic K\"ahler manifolds
$\mathcal{QM}_{4n+4}$ lying in the image of the $c$-map, irrespectively  whether the pre-image Special K\"ahler manifold $\mathcal{SK}_n$ is a
symmetric space or not. The second Heisenberg algebra exists only in the case when the Quaternionic K\"ahler manifold $\mathcal{QM}_{4n+4}$ is a
symmetric space. As we discuss in sections \ref{plantageneti} and  \ref{generaloneSOpq}, the Heisenberg subalgebras can play an important role in
the attempts to find discretization schemes via tessellations of the host manifold to which data are mapped, which, as we know,  is anyhow a Lie
group, namely the solvable group manifold. Hence no discretization scheme can avoid the issue of singling out discrete subgroups of such a group.

From this discussion we also realize that the central charge $\mathbb{Z}_1$ is just the $L_+$ generator of a universal $\slal(2,\mathbb{R})_E$ Lie
algebra that exists only in the symmetric space case and which was named the Ehlers algebra in \cite{advancio}.  When $\slal(2,\mathbb{R})_E$ does
exist we can introduce the universal compact generator:
\begin{equation}\label{ruotogrande}
    \mathfrak{S}\, \equiv\, L_+ \, - \, L_- \,= \, 2 \, \lambda^2
\end{equation}
which rotates the two sets of Heisenberg translations one into the other:
\begin{equation}\label{tabarro}
    \left[ \mathfrak{S}\, , \, \mathbf{W}^{i\alpha}\right ] \, = \, \epsilon^{ij} \, \mathbf{W}^{j\alpha}
\end{equation}
%%%%%%%%%%%%%%%%%%%%%%%%%%%%%%
\section{The Solvable Lie algebra of $\mathrm{U/H}$ spaces, the Tits-Satake Projection and the paint group}
We come next to the core of our exposition. We   begin by recalling the list of classical complex Lie algebras. The exceptional Lie algebras also
fall in the general scheme of real sections, Tits-Satake projection, Solvable Subalgebras and paint group characterization, yet they have fixed
dimensionality and, while they happened to be very much relevant in Supergravity, they seem to have less chance of application in Data Science,
although this is not excluded. For this reason they will be mentioned here only occasionally.
\paragraph{\sc The classical series of symmetric noncompact coset manifolds.} Relying on their complete
classification  (see \cite{fre2023book}),  according with standard nomenclature and leaving aside the $5$ exceptional algebras that have not a
variable dimensionality as required by Data Science, the complex simple Lie algebras are distributed into four infinite families  as follows,
where $\ell$ is the rank, namely the dimensionality of the Cartan subalgebra:
\begin{enumerate}
  \item $\mathfrak{a}_\ell$, the complex Lie algebra $\slal(\ell +1,\mathbb{C})$.
  \item $\mathfrak{b}_\ell$, the complex Lie algebra $\so (2\ell +1,\mathbb{C})$.
  \item $\mathfrak{c}_\ell$, the complex Lie algebra $\sym (2\ell,\mathbb{C})$.
  \item $\mathfrak{d}_\ell$, the complex Lie algebra $\so(2\ell,\mathbb{C})$.
\end{enumerate}
The above classical algebras are determined by the corresponding Dynkin diagrams displayed in Figure \ref{pippusdia}.
\begin{figure}
\centering
\begin{picture}(100,200)
\put (-70,185){$\mathfrak{a}_\ell$} \put (-20,190){\circle {10}} \put (-23,175){$\alpha_1$} \put
(-15,190){\line (1,0){20}} \put (10,190){\circle {10}} \put (7,175){$\alpha_2$} \put (15,190){\line
(1,0){20}} \put (40,190){\circle {10}} \put (37,175){$\alpha_3$} \put (47,190){$\dots$}
%\put (45,290){\line (1,0){20}}
\put (70,190){\circle {10}} \put (67,175){$\alpha_{\ell-2}$} \put (75,190){\line (1,0){20}} \put
(100,190){\circle {10}} \put (97,175){$\alpha_{\ell-1}$} \put (105,190){\line (1,0){20}} \put
(130,190){\circle {10}} \put (127,175){$\alpha_\ell$}
%%%%%%%%%%%%%%%%%%%%%%%%%%
\put (-70,155){$\mathfrak{b}_\ell$} \put (-20,160){\circle {10}} \put (-23,145){$\alpha_1$} \put
(-15,160){\line (1,0){20}} \put (10,160){\circle {10}} \put (7,145){$\alpha_2$} \put (15,160){\line
(1,0){20}} \put (40,160){\circle {10}} \put (37,145){$\alpha_3$} \put (47,160){$\dots$}
%\put (45,260){\line (1,0){20}}
\put (70,160){\circle {10}} \put (67,145){$\alpha_{\ell-2}$} \put (75,160){\line (1,0){20}} \put
(100,160){\circle {10}} \put (97,145){$\alpha_{\ell-1}$} \put (105,163){\line (1,0){20}} \put (105,158){\line
(1,0){20}} \put (108,157){{\LARGE$>$}} \put (130,160){\circle {10}} \put (127,145){$\alpha_\ell$}
%%%%%%%%%%%%%%%%%%%%%%%%%%%%%%%%%
\put (-70,125){$\mathfrak{c}_\ell$} \put (-20,130){\circle {10}} \put (-23,115){$\alpha_1$} \put
(-15,130){\line (1,0){20}} \put (10,130){\circle {10}} \put (7,115){$\alpha_2$} \put (15,130){\line
(1,0){20}} \put (40,130){\circle {10}} \put (37,115){$\alpha_3$} \put (47,130){$\dots$}
%\put (45,230){\line (1,0){20}}
\put (70,130){\circle {10}} \put (67,115){$\alpha_{\ell-2}$} \put (75,130){\line (1,0){20}} \put
(100,130){\circle {10}} \put (97,115){$\alpha_{\ell-1}$} \put (105,133){\line (1,0){20}} \put (105,128){\line
(1,0){20}} \put (108,126){{\LARGE$<$}} \put (130,130){\circle {10}} \put (127,115){$\alpha_\ell$}
%%%%%%%%%%%%%%%%%%%%%%%%%%%%%%%%%%%%%%%%%%%%%%%%%%%%%%%%%%%%%%%%%%%%%
\put (-70,65){$\mathfrak{d}_\ell$} \put (-20,70){\circle {10}} \put (-23,55){$\alpha_1$} \put (-15,70){\line
(1,0){20}} \put (10,70){\circle {10}} \put (7,55){$\alpha_2$} \put (15,70){\line (1,0){20}} \put
(40,70){\circle {10}} \put (37,55){$\alpha_3$} \put (47,70){$\dots$}
%\put (45,230){\line (1,0){20}}
\put (70,70){\circle {10}} \put (67,55){$\alpha_{\ell-3}$} \put (75,70){\line (1,0){25}} \put
(105,70){\circle {10}} \put (100,55){$\alpha_{\ell-2}$} \put (110,70){\line (1,1){20}} \put (110,70){\line
(1,-1){20}} \put (133.2,93.2){\circle {10}} \put (133.2,46.8){\circle {10}} \put
(143.2,93.2){$\alpha_{\ell-1}$} \put (143.2,43.8){$\alpha_\ell$}
%%%%%%%%%%%%%%%%%%%%%%%%%%%%%%%%%%%%%%%%%%%%%%%%%%%%%%%%%%%%%%%%%%%%%%
\end{picture}
\vskip -1cm \caption{\label{pippusdia} The Dynkin
diagram\index{Dynkin diagram}s of the four infinite families of
classical simple algebras}
\end{figure}
The unique maximally compact real section of each of these complex Lie algebra is shown below:
\begin{enumerate}
  \item $\mathfrak{a}_\ell \relbar\joinrel\rightsquigarrow  \su(\ell +1)$.
  \item $\mathfrak{b}_\ell \relbar\joinrel\rightsquigarrow  \so(2\ell +1,\mathbb{R})$.
  \item $\mathfrak{c}_\ell \relbar\joinrel\rightsquigarrow  \usp(2\ell)$.
  \item $\mathfrak{d}_\ell \relbar\joinrel\rightsquigarrow  \so(2\ell,\mathbb{R})$.
\end{enumerate}
All other real sections $\mathbb{U}$ are noncompact and therefore lead to symmetric spaces $\mathrm{U/H}$.
\subsection{The maximally split case}
\label{piripicchio} In the case we choose $\mathbb{U}=\mathbb{G}_{\mathrm{max}}$, where $\mathbb{G}_{\mathrm{max}}$ is the maximally split real
section of any of the simple Lie algebras, according  with the definition given in eq.~(\ref{bagarraT}),  the maximally compact subalgebra is just:
\begin{equation}\label{maxcomHms}
  \mathbb{H} \, \equiv \, \text{span}_\mathbb{R}\left\{ E^{\alpha}- E^{-\alpha} \, \mid\, \alpha \in \Phi^+\right\}
\end{equation}
where $\Phi^+$ denotes the subset of all positive roots of the root system. Furthermore  the orthogonal decomposition
\begin{equation}\label{ortofruttamaxsplit}
  \mathbb{U} \, = \, \mathbb{H} \oplus \mathbb{K}
\end{equation}
is achieved by setting:
\begin{equation}\label{kappabasamaxsplit}
   \mathbb{K} \, \equiv \, \text{span}_\mathbb{R}\left\{\mathcal{H}_i ,\,  E^{\alpha}+ E^{-\alpha} \, \mid\,
   \{\mathcal{H}_i\} = \text{basis of CSA}, \,
   \alpha \in \Phi^+\right\}
\end{equation}
In this case the solvable Lie subalgebra $Solv_{\mathrm{U/H}}\subset \mathbb{U}$ whose corresponding solvable Lie group $\mathcal{S}_{U/H}$ is
metrically equivalent to the symmetric space manifold $\mathrm{U/H}$ is just the Borel subalgebra $\mathbb{B}(\mathbb{U})\subset \mathbb{U}$:
\begin{equation}\label{borellino}
  Solv_{\mathrm{U/H}} \, = \, \mathbb{B}(\mathbb{U}) \, \equiv \, \text{span}_\mathbb{R}\left\{\mathcal{H}_i ,\,  E^{\alpha} \, \mid\,
   \{\mathcal{H}_i\} = \text{basis of CSA}, \,
   \alpha \in \Phi^+\right\}
\end{equation}
The equality between the dimensionality of the Borel subalgebra with that of the subspace $\mathbb{K}$ of coset generators and of the symmetric
space $\mathrm{U/H}$ is evident from eqs.~ (\ref{kappabasamaxsplit},\ref{borellino}).
\par
The Borel subalgebra (\ref{borellino}) is turned into a \emph{normed solvable Lie algebra} according with definition \ref{definormosolvalg} by
setting:
\begin{eqnarray}\label{canonicB}
     < \mathcal{H}_i\, ,\, \mathcal{H}_j>_n & = &  2 \, \delta_{ij}~, \nonumber\\
     < \mathcal{H}_i\, ,\, E^\alpha>_n & = &  0 ~,\nonumber\\
     < E^\alpha\, ,\, E^\beta>_n & = &  \delta_{\alpha\beta}
\end{eqnarray}
This quadratic form is invariant according with eq.~(\ref{invarianorma}) and its normalization is absolute  if the generators of the Weyl-Cartan
basis for $\mathbb{U}$ have the standard normalization displayed in eq.~ (\ref{CarWeylform}). This norm on the solvable Lie algebra is precisely
that induced by the unique (up to an overall scale factor) $\mathrm{U}$-invariant Riemannian metric on the symmetric space $\mathrm{U/H}$
according with definition \ref{Alexnormalmanif} and \eqref{solvgroupN}.
\par
In relation with this issue, this is  the appropriate point to mention an algorithmic procedure to find the precise form of the norm $\langle \,
, \, \rangle_n$ induced on the solvable Lie algebra by the symmetric space canonical metric. Such an algorithm  is quite convenient and almost
automatic in the dual basis of $1$-forms, since it utilizes the standard Mathematica code \emph{Vielbgrav23} \cite{codiconi} for the calculation
of the intrinsic components of the curvature 2-form in vielbein Riemannian geometry.
\subsubsection{Algorithm to find the norm  $\langle \, , \, \rangle_n$ on the solvable
Lie algebra $Solv_{\mathrm{U/H}}$ of a symmetric space} Let $\mathcal{T}_A$ be a basis of generators of the solvable Lie algebra
$Solv_{\mathrm{U/H}}$ of a hyperbolic symmetric space, not necessarily maximally split. Let $\mathbb{L}(\boldsymbol{\Upsilon})\in
\mathcal{S}_{\mathrm{U/H}}$ be a generic element of the solvable group parameterized by any convenient set $\Upsilon^I$
($I=1,\dots,\dim (Solv_{\mathrm{U/H}})$) of parameters.  Let
\begin{equation}\label{commodore}
  \Theta \, \equiv \, \mathbb{L}(\boldsymbol{\Upsilon})^{-1} \, \mathrm{d}\left[\mathbb{L}(\boldsymbol{\Upsilon})\right]
\end{equation}
 be a  left-invariant $1$-form and
 \begin{equation}\label{capitanfracasso}
    \Theta \, = \, \sum_{A=1}^{\dim (Solv_{\mathrm{U/H}})} \, \mathbf{e}^A \, \mathcal{T}_A
 \end{equation}
 its expansion in the basis of generators $\mathcal{T}_A$.
 \par
The metric on solvable manifold can be always written as
\begin{equation}\label{metrolinea}
  ds^2 \, = \, \kappa_{AB}\, \mathbf{e}^A \times \mathbf{e}^B,
\end{equation}
where $\kappa_{AB}$ is a symmetric constant matrix. One can use this choice of vielbein to calculate by means of the  \emph{Vielbgrav23 code} \cite{codiconi} the
curvature 2-form $\mathfrak{R}^{AB}$, the Riemann tensor $\mathrm{Rie}^{AB}_{\phantom{AB} CD}$ and the Ricci tensor $\mathrm{Ric}^{A}_{\phantom{A}
C}$ according to the conventions
 \begin{equation}\label{convenzioneginevra}
   \mathfrak{R}^{AB}\, = \, \mathrm{Rie}^{AB}_{\phantom{AB} CD} \, \mathbf{e}^C \wedge\mathbf{e}^D \quad; \quad \mathrm{Ric}^{A}_{\phantom{A} C}
   \, = \, \mathrm{Rie}^{AB}_{\phantom{AB} CB}
 \end{equation}
Since the vielbein $\mathbf{e}^A$  satisfies  the Maurer Cartan equations of the solvable group, we are guaranteed that the contorsion has only
constant components and so does the spin connection and the Riemann tensor. With an arbitrary choice of the matrix $\kappa_{AB}$ the Ricci tensor
is constant but not proportional to the Kronecker Delta. Imposing that $ \mathrm{Ric}^{A}_{\phantom{A} C} \, = \, \lambda \, \delta^A_C$ yields a
number of algebraic equations for $\kappa_{AB}$ which necessarily has one and only one solution $\kappa^{[E]}_{AB}$ up to an overall constant
(depending on $\lambda$) since we know a priori that the $\mathrm{U}$-invariant metric on $\mathrm{U/H}$ is   Einstein. In this way the
precise form of the norm on the solvable Lie algebra induced by the symmetric space can be determined:
\begin{equation}\label{sfogliatinalimone}
  \langle \mathcal{T}_A\, , \, \mathcal{T}_B \rangle_{n} \, = \,  \kappa^{[E]}_{AB}
\end{equation}
\subsubsection{The maximally split $\mathrm{SL(N,\mathbb{R})/SO(N)}$ symmetric space}
\label{realescena}  The simple Lie algebra $\mathfrak{a}_{\mathrm{N-1}}$, identified by the first line  Dynkin diagram in Figure \ref{pippusdia}
(setting $\ell = \mathrm{N-1}$) is the abstract form of the Lie algebra $\slal(\mathrm{N},\mathbb{C})$ of complex traceless matrices in dimension
$\mathrm{N}$. The corresponding Borel subalgebra $\mathbb{B}_\mathrm{N}^{\mathbb{C}} \, \equiv \,\mathbb{B}({\mathfrak{a}}_{{\mathrm{N-1}}})$ is
simply given by  the subset of all upper triangular traceless complex matrices. Hence we define
\begin{equation}
\label{Borelliano} \slal (\mathrm{N},\mathbb{C}) \, \supset \, {\mathbb{B}}_\mathrm{N}^{\mathbb{C}} \, \ni \, \mathfrak{b} \, = \,
\left(
  \begin{array}{cccccc}
    \star & \star & \star & \star & \dots & \star \\
    0 & \star & \star & \star & \dots & \star \\
    0 & 0 & \star & \star & \dots & \star \\
    0 & 0 & 0 & \star & \dots & \star \\
    \dots & \dots & \dots & \dots & \dots & \dots \\
    0 & 0 & 0 & 0 & 0  & \star \\
  \end{array}
\right) \quad ; \quad Tr[\mathfrak{b}]\, = \, 0
\end{equation}
The maximally split real section is obtained by restricting all the traceless $\mathrm{N\times N}$ matrices to have purely real entries. This
obviously applies also to the upper triangular traceless matrices forming the Borel subalgebra.  Hence the real Borellian
${\mathbb{B}}_\mathrm{N}^{\mathbb{R}} \subset \slal(\mathrm{N},\mathbb{R})$  is a solvable subalgebra of $\mathbb{U} \, =
\,\slal(\mathrm{N},\mathbb{R})$ having real dimensions:
\begin{equation}\label{dimenBN}
    \mathrm{dim}_{\mathbb{R}} \,{\mathbb{B}}_\mathrm{N}^\mathbb{R} \, \equiv \, d_\mathrm{N} \, = \,
    \frac{\mathrm{N}(\mathrm{N}+1)}{2} \, - \, 1
\end{equation}
That above in eq.~(\ref{dimenBN}) is also, as it should be, the dimension of the coset manifold $\mathrm{SL(N,\mathbb{R})/SO(N)}$ which, in this
simple case, can also be described  as the manifold of all unimodular symmetric matrices $\mathcal{M}$:
\begin{equation}\label{Mcalli}
  \mathrm{SL(N,\mathbb{R})/SO(N)} \, = \, \left\{ \mathcal{M} \in \mathrm{Hom}(\mathbb{R}^N,\mathbb{R}^N)\, \mid \, \mathcal{M} \, = \, \mathcal{M}^T,
  \quad \mathrm{Det}[\mathcal{M}] \, = \, 1\right\}
\end{equation}
The above concise definition of the symmetric space turns out to be very useful while discussing both the spectrum of harmonic functions on
$\mathrm{U/H}$ and when deriving the integral formulas for the geodesics and, hence, also for the geodesic distance function. Furthermore the
concise definition (\ref{Mcalli}) applies not only to the maximal split $\mathrm{SL(N,\mathbb{R})/SO(N)}$ case, but also to all the other
classical $\mathrm{U/H}$. What happens in those cases is that the symmetric unimodular $\mathcal{M}$ is constrained by additional quadratic
constraints related with the definition of the relevant Lie group.
\subsubsection{The solvable coordinates}
In the maximally split case $\mathrm{U/H} =\mathrm{SL(N,\mathbb{R})/SO(N)}$ we define the mapping from the solvable Lie algebra (in this case the
borellian) to the solvable group introducing an ordered graded basis of generators  $\mathrm{T}^A$ for the Solvable Lie algebra based on the root
height. The generators start at grade $0$, which corresponds to Cartan generators. Next come the grade $1$ generators that are orderly associated
with simple roots. Next we have the grade $2$ generators associated with the roots of hight $2$ and so on up to the unique generators associated
with the unique highest root.
\par
In this way we write:
\begin{equation}\label{carontedimonio}
\mathbb{L}(\Upsilon)=\prod_{A=1}^{n} \,\exp\left[\Upsilon_{A} \mathrm{T}^A\right]
\end{equation} where
$\Upsilon^A$ are the \emph{solvable coordinates} of the symmetric space.
\subsection{The triangular embedding of the $\mathrm{U/H}$ symmetric spaces associated with other real sections}
\label{triangolodellemore} In order to construct the solvable Lie algebra for the other real sections it is very important the following
\begin{statement}\label{statamento}
 Let $\mathrm{N}$ be the real dimension of the fundamental representation of $\mathbb{U}$ and $\mathbb{H}\subset \mathbb{U} $ its maximal
compact subalgebra. There always is  a suitable symmetric matrix $\eta_t$ with the appropriate signature for the real sections of
$\mathfrak{b}_\ell$ and $\mathfrak{d}_\ell$ and a suitable skew-symmetric matrix $\mathbb{C}_t$ for the real sections of $c_\ell$, implying
canonical embeddings
\begin{eqnarray}
 \mathbb{ U} & \hookrightarrow & \slal(\mathrm{N},\mathbb{R})~,\nonumber\\
\mathbb{ U} \, \supset \, \mathbb{H} & \hookrightarrow &
\so\mathrm{(N)} \, \subset \,
\slal(\mathrm{N},\mathbb{R})~. \label{triaembed}
\end{eqnarray}
Such an embedding is determined in all cases by the choice of the basis, where $\Solv\left(\mathrm{U/H} \right)$ is made up  by upper triangular
matrices. Furthermore  in the same basis, in the case of the $\mathfrak{b}_\ell$ and $\mathfrak{d}_\ell$ algebras, all the elements of both
$\mathbb{K}$ and $\mathbb{H}$ are$\eta_t$-skew-symmetric matrices:
\begin{eqnarray}\label{etasym}
\eta_t \, \mathrm{K} \, + \,  \mathrm{K}^T \, \eta_t~ \, = \, 0   &  \forall \, \mathrm{K}\, \in \, \mathbb{K} \nonumber \\
\eta_t \, \mathrm{H} \, + \,   \mathrm{H}^T\,\eta_t~ \, = \,  0   &  \forall \, \mathrm{H}\, \in \, \mathbb{H}
\end{eqnarray}
while, at the same time they must satisfy the conditions
\begin{equation}\label{symeta}
     \mathrm{K} \,- \,  \mathrm{K}^T  \, = \, 0, \qquad   \mathrm{H} \, + \,   \mathrm{H}^T   \, = \, 0.
\end{equation}
The first condition (\ref{etasym}) guarantees that the elements of both $\mathbb{K}$ and $\mathbb{H}$ belong to the $\mathbb{U}$ Lie algebra, while
second condition (\ref{symeta})  guarantees that $\mathbb{H}$ is mapped to $\so(\mathrm{N})$  while $\mathbb{K}$ is mapped instead to the
orthogonal complement of $\so(\mathrm{N})$ in $\slal(\mathrm{N},\mathbb{R})$.
\par
In the case of real sections of the $c_\ell$ Lie algebras, the conditions to be satisfied by the matrix $\mathbb{C}_t$ and the matrices of
$\mathbb{H}$ and $\mathbb{K}$ are identical to   (\ref{etasym}-\ref{symeta}), with $\eta_t$ replaced by $\mathbb{C}_t$. Since the
latter is skew-symmetric, the condition (\ref{etasym}) becomes a condition of $\mathbb{C}_t$ symmetricity.
\end{statement}
In section \ref{sorsgenerale} we shall  show the explicit construction of the $\eta_t$ matrix for the case of the Lie algebras $\so(r,r+2s)$. The choice
of this basis is the strategic instrument to derive the explicit form of all Lie algebra generators and analyse the Tits-Satake projection and
the Paint group. Let us then turn to the latter, already sketched in the introduction, and present it in full-fledged form.
\subsection{Tits-Satake Projection and the paint group}
\label{colorato} In the introductory section we have already anticipated the concept of \emph{Tits-Satake projection}, defined
in eq.~(\ref{algeoTS}) and acting on the root system, according with the decomposition (\ref{splittus}) induced by the Cartan involution that
defines  the considered real section $\mathbb{U}$. We have also presented without proof   Theorem \ref{ruttosistnuovo};   as it will be  the case
 most    theorems in this chapter, the proof will simply a case-by-case analysis.

At first we answer a preliminary question, i.e.,  \emph{What is the solvable Lie
algebra of the solvable group that is metrically equivalent to $\mathrm{U/H}$, when the latter is not maximally split?}.  The example of the maximally
split case clearly suggests what is the required solvable algebra\index{solvable algebra} for other real sections. Indeed we have:
\begin{equation}\label{coreutagreco}
    \Solv\left(\frac{\mathrm{U}}{\mathrm{H}}\right)\, = \, \mathcal{H}^{n.c.}\oplus
    \text{span}\left(\mathcal{E}^\alpha
     \right) \quad \quad \forall \alpha\in\pmb{\Phi}_+ \, / \, \alpha_\parallel \, \neq \, 0,
\end{equation}
where $\mathcal{H}^{n.c.}$ is the noncompact part of the Cartan subalgebra,  and $\mathcal{E}^\alpha$ denotes the combination of step operators
pertaining to the positive roots $\alpha$ that appear in the real section $\mathbb{U}=\mathbb{G}_R$; the sum is extended only to those roots
that are not purely imaginary. Indeed the step operators pertaining to imaginary roots are included into the maximal compact subalgebra, which now
is larger than the number of positive roots. The explicit form of the generators associated with the roots can be worked out only through the use
of the triangular embedding, yet schematically eq.~(\ref{coreutagreco}) is sufficient for our general discussion.
%%%%%%%%%%%%%%%%%%%%%%%%%%%%%%%%%%%%%%%%%%%%%%%%%%%%%%%%%%%%
\subsubsection{The Tits-Satake projection: an  in depth algebraic glance}
\label{Titsatesection} In this subsection we elaborate all the needed details, explaining  how the Tits-Satake projection of a \emph{normed
solvable Lie algebra} is actually performed, and how it is related to the notions of \emph{paint group} $\mathrm{G}_{\mathrm{paint}}$ and
\emph{subpaint group} $\mathrm{G}_{\mathrm{subpaint}} \subset \mathrm{G}_{\mathrm{paint}}$. Such peculiar universal structure of the solvable
algebras $\Solv(\mathrm{U/H})$ had not been observed before the paper \cite{noipainted};  its scope extends beyond symmetric spaces, as it was  demonstrated in
\cite{titsusataku} --- indeed, the Tits-Satake projection can be defined for general normed solvable Lie algebras. Yet, our main interest is in
symmetric spaces, and the  just mentioned notions  were extracted precisely from the case of the Tits-Satake projections of solvable Lie
algebras associated with symmetric spaces $\Solv(\mathrm{G/H})$, so that we shall focus on that case.

The Tits-Satake projection is made in two steps. First one sets all $\alpha_{\bot} = 0 $, projecting the original root system $\Phi_\mathbb{U}$
onto a new system of vectors $\overline\Phi$   in a Euclidean space of dimension equal to the noncompact rank $r_{\rm nc}$,
called a {\em restricted root system}. It is not an ordinary root system as  roots may occur with multiplicities,  and $2\alpha_{||}$ can be a root if $\alpha_{||}$ is one. In the second step, one deletes the multiplicities of the restricted
roots. Thus we have
\begin{equation}
  \pi_{\mathrm{TS}} \quad : \quad \Phi_\mathbb{U} \to  \Phi_{\rm TS}\, ;\qquad
  \Phi_\mathbb{U} \xrightarrow{\alpha_{\bot} = 0}\overline{\Phi} \
\xrightarrow{\mbox{\footnotesize deleting multiplicities}}
  \Phi_{\rm TS}.
 \label{projectoB}
\end{equation}
\emph{If $\overline{\Phi}$ contains no restricted root that is twice another}, then $\Phi_{\rm TS}$ is a root system of simple type.
We will show later that this root subsystem defines a Lie algebra $\mathbb{U}_{\mathrm{TS}}$, \emph{the Tits-Satake subalgebra} of $\mathbb{U}$:
\begin{equation}
  \Phi_{\rm TS} = \mbox{root system of } \mathbb{U}_{\mathrm{TS}}, \qquad
\mathbb{U}_{\mathrm{TS}} \, \subset \,\mathbb{U}. \label{TitsSatakeB}
\end{equation}
The Tits-Satake subalgebra $\mathbb{U}_{\mathrm{TS}}$ is, as a consequence of its own definition, the maximally noncompact real section of its
own complexification. For this reason, considering its maximal compact subalgebra $\mathbb{H}_{\mathrm{TS}} \, \subset \,
\mathbb{U}_{\mathrm{TS}}$ we have a new smaller coset $\mathrm{U_{\mathrm{TS}} / H_{\mathrm{\mathrm{TS}}}}$, which is maximally split, whose
associated solvable algebra $\Solv(\mathrm{U_{\mathrm{TS}} / H_{\mathrm{TS}}})$ is the Borelian of $\mathbb{U}_{\mathrm{TS}}$.
\par
\emph{In the case doubled restricted roots are present in} $\overline{\Phi}$, the projection cannot be expressed in terms of a simple Lie
algebra, but the concept remains the same. The root system is the so-called $bc_r$ system, with $r=r_{\rm nc}$ the noncompact rank of the real
form $\mathbb{U}$. It is the root system of a group $\mathrm{U_{\mathrm{TS}}}$, which is now nonsemisimple. The manifold is similarly defined as
$\mathrm{U_{\mathrm{TS}}/H_{TS}}$, where $\mathrm{H_{TS}}$ is the maximal compact subgroup of $\mathrm{U_{\mathrm{TS}}}$. As we already stressed
this case happens only in the case of real sections of the $\mathfrak{a}_\ell$ Lie algebras.

The next question is: what is the relation between the two solvable Lie algebras $\Solv(\mathrm{U/H})$ and $\Solv(\mathrm{U_{\mathrm{TS}} /
H_{\mathrm{TS}}})$? The answer can be formulated through the following statements A-E.

\paragraph{\sc A)}
In a projection, more than one higher dimensional vector can map to the same lower dimensional one, i.e., xxtin general there will be
several roots of $\Phi_\mathbb{U}$ that have the same image in $\Phi_{\rm TS}$. The imaginary roots vanish under this projection, according to the
definition of \ref{leluka}. Therefore, apart from these imaginary roots, there are two types of roots: those that have a distinct image in the
projected root system, and those that arrange into multiplets with the same projection. We can split the root spaces in subsets according to
whether there is such a degeneracy or not. Calling $\Phi^+_\mathbb{U}$ and $\Phi_{\rm TS}^+$ the sets of positive roots of the two root systems,
we have the following scheme:
\begin{eqnarray}
 && \begin{array}{lclclcc}
 \Phi^+_\mathbb{U} & = & \Phi^\eta& \bigcup &\Phi^\delta& \bigcup & \Phi_{\mathrm{comp}} \\
    \downarrow \pi_{\mathrm{TS}}& & \downarrow \pi_{\mathrm{TS}} &   & \downarrow \pi_{\mathrm{TS}} &   &   \\
\Phi_{\mathrm{TS}}^+ &=& \Phi_{\mathrm{TS}}^\ell & \bigcup & \Phi_{\mathrm{TS}}^s &&  \\
  \end{array} \nonumber\\[3mm]
&&\pi^{-1}_{\mathrm{TS}} \left[ \alpha^\ell \right]  = 1\quad  \forall \alpha^\ell  \in \Phi_{\mathrm{TS}}^\ell ,
\qquad \dim \pi^{-1}_{\mathrm{TS}} \left[ \alpha^s \right]  =   m[\alpha^s]>1 \quad   \forall \alpha^s  \in  \Phi_{\mathrm{TS}}^s .
 \label{SubsetsRootspaces}
\end{eqnarray}
The $\delta  $ part thus contains all the roots that have nontrivial multiplicities under the Tits-Satake projection, while the roots in the $\eta $ part have
no multiplicities; these roots of type $\eta  $ are orthogonal to $\Phi _{\mathrm{comp}}$. This follows from the fact that for any two
root vectors $\alpha $ and $\beta $ where there is no root of the form $\beta +m\alpha $ with $m$ a non-zero integer, the inner product of $\beta
$ and $\alpha $ vanishes. It also follows from this definition that in maximally split symmetric spaces, in which case $\Phi
_{\mathrm{comp}}=\emptyset$, all root vectors are in $\Phi ^\eta $ or $\Phi ^\ell $ (as the Tits-Satake projection is then trivial).
\par
These subsets moreover satisfy the following properties under addition of root vectors:
\begin{equation}
\begin{array}{|l|l|}
\hline \mathbb{U}&  \mathbb{U}_{\mathrm{TS}} \\ \hline \Phi ^\eta +\Phi^\eta \subset\Phi ^\eta &
\Phi_{\mathrm{TS}}^\ell+\Phi_{\mathrm{TS}}^\ell
  \subset\Phi_{\mathrm{TS}}^\ell \\
 \Phi ^\eta +\Phi^\delta  \subset\Phi ^\delta  &
 \Phi_{\mathrm{TS}}^\ell+\Phi_{\mathrm{TS}}^s
  \subset\Phi_{\mathrm{TS}}^s \\
  \Phi ^\delta  +\Phi^\delta  \subset\Phi ^\eta\bigcup\Phi ^\delta
   & \Phi_{\mathrm{TS}}^s+\Phi_{\mathrm{TS}}^s
  \subset\Phi_{\mathrm{TS}}^\ell \bigcup\Phi_{\mathrm{TS}}^s\\
   \Phi_{\mathrm{comp}}+ \Phi ^\eta =\emptyset &   \\
    \Phi_{\mathrm{comp}}+ \Phi ^\delta \subset \Phi ^\delta  &\\
    \hline
\end{array}
\end{equation}
In view of this structure we can enumerate the generators of the solvable algebra $\Solv(\mathrm{U/H})$ in the following way:
\begin{eqnarray}
\Solv(\mathrm{U/H})&=&\left\{H_i,
\Phi_{\alpha^\ell},\Omega_{\alpha^s | I}\right\} \nonumber\\
H_i & \Rightarrow & \mbox{Cartan generators} \nonumber\\
\Phi_{\alpha^\ell} & \Rightarrow & \eta -\mbox{roots} \nonumber\\
\Omega_{\alpha^s | I} &\Rightarrow & \delta -\mbox{roots} \quad ; \quad (I=1,\dots , m[\alpha^s] ).
\label{enumerato}
\end{eqnarray}
The index $I$ enumerating the $m$--roots of $\Phi_{\mathbb{U}}$ that have  the same projection in $\Phi_{\mathrm{TS}}$ is called the \emph{paint
index}.
\paragraph{\sc B)}
There exists a \emph{compact subalgebra} $\mathbb{G}_{\mathrm{paint}} \, \subset
  \, \mathbb{U}$ which acts as  an algebra of outer automorphisms (i.e., outer derivations) of the  solvable
  algebra $\Solv_{\mathbb{U}} \equiv \Solv(\mathrm{U /
  H}) \subset \mathbb{U}$, namely:
\begin{equation}
  \left[ \mathbb{G}_{\mathrm{paint}} \, , \, \Solv_{\mathbb{U}} \right]  \subset
  \Solv_{\mathbb{U}}.
\label{linearrepre}
\end{equation}
\paragraph{\sc C)}
The Cartan generators $H_i$ and the generators $\Phi_{\alpha^\ell} $ are singlets under the action of $\mathbb{G}_{\rm paint}$, i.e.,  each of
them  commutes with the whole of $\mathbb{G}_{\rm paint}$:
\begin{equation}
  \left[ H_i \, , \,\mathbb{G}_{\mathrm{paint}}\right] \, = \, \left[ \Phi_{\alpha^\ell} \, ,
  \,\mathbb{G}_{\rm paint}\right]\, = \, 0
\label{hiphicommuti}
\end{equation}
On the other hand, each of the multiplets of generators $\Omega_{\alpha^s | I}$ constitutes an orbit under the adjoint action of the paint group
${G}_{\mathrm{paint}}$, i.e., a linear representation $\mathbf{D}{[\alpha^s]}$ which, for different roots $\alpha^s$ can be different:
\begin{equation}
  \forall \, X \, \in \, \mathbb{G}_{\mathrm{paint}} \quad : \quad  \left[ X \, , \,
  \Omega_{\alpha^s |
I}\right] \, = \, \left( D^{[\alpha^s]}[X]\right) _I^{\phantom{I}J} \, \Omega_{\alpha^s | J} \label{replicas}
\end{equation}
\paragraph{\sc D)}
The \emph{paint algebra} $\mathbb{G}_{\mathrm{paint}}$ contains a subalgebra
\begin{equation}
  \mathbb{G}^0_{\rm subpaint} \, \subset \, \mathbb{G}_{\mathrm{paint}}
\label{paint0}
\end{equation}
such that with respect to $\mathbb{G}^0_{\mathrm{subpaint}}$,  each $m[\alpha^s]$--dimensional representation $\mathbf{D}{[\alpha^s]}$ branches as
follows:
\begin{equation}
  \mathbf{D}{[\alpha^s]} \,
  \stackrel{\mathbb{G}^0_{\rm subpaint}}{\Longrightarrow} \,
  \underbrace{{\mathbf{1}}}_{\mbox{singlet}}
  \, \oplus \, \underbrace{{\mathbf{J}}}_{(m[\alpha^s]-1)-\mbox{dimensional}
  }
\label{splittatonibus}
\end{equation}
Accordingly, we can split the range of the multiplicity index $I$ as follows:
\begin{equation}
  I = \left\{ 0,x \right\},\qquad  x=1,\dots,m[\alpha^s]-1.
\label{Irango}
\end{equation}
The index $0$ corresponds to the singlet, while $x$ ranges over the representation $\mathbf{J}$.
\paragraph{\sc E)}
The tensor product $\mathbf{J} \otimes \mathbf{J} $ contains both the identity representation $\mathbf{1}$ and the representation $\mathbf{J}$
itself. Furthermore, there is  in the representation $\bigwedge ^3 \mathbf{J}$ a $\mathbb{G}^0_{\mathrm{subpaint}}$-invariant tensor $a^{xyz}$
such that the two solvable Lie algebras $\Solv_{\mathbb{U}}$ and $\Solv_{\mathbb{U}_{\mathrm{TS}}}$ can be written as follows
\begin{equation}
\begin{array}{|l|l|}
\hline
  \Solv_{\mathbb{U}}  & \Solv_{\mathbb{U}_{\mathrm{TS}}} \\
  \hline
 \left[ H_i \, , \, H_j \right] = 0 & \left[ H_i \, , \, H_j \right] =0\\
\left[ H_i \, , \, \Phi_{\alpha^\ell}
  \right] = \alpha^\ell_i \, \Phi_{\alpha^\ell}
  &\left [ H_i \, , \, E^{\alpha^\ell} \right ] = \alpha^\ell_i \\
\left[ H_i \, , \, \Omega_{\alpha^s|I}
  \right] = \alpha^s_i \, \Omega_{\alpha^s|I} &
\left [ H_i \, , \, E^{\alpha^s} \right ] = \alpha^s_i \, E^{\alpha^s}
\\
\left [ \Phi_{\alpha^\ell} \, , \, \Phi_{\beta^\ell} \right ] =N_{\alpha^\ell\beta^\ell} \, \Phi_{\alpha^\ell
+ \beta^\ell}
 &
 \left [ E^{\alpha^\ell} \, , \, E^{\beta^\ell} \right ] = N_{\alpha^\ell\beta^\ell} \,
E^{\alpha^\ell + \beta^\ell}    \\
\left [ \Phi_{\alpha^\ell} \, , \, \Omega_{\beta^s|I} \right ] =N_{\alpha^\ell\beta^s} \, \Omega_{\alpha^\ell
+ \beta^s|I} &
 \left [ E^{\alpha^\ell} \, , \, E^{\beta^s} \right ] = N_{\alpha^\ell\beta^s} \,
E^{\alpha^\ell + \beta^s}   \\
\mbox{If } \alpha^s + \beta^s \in \Phi_{\mathrm{TS}}^\ell: & \\
\phantom{.}\quad\left [ \Omega_{\alpha^s|I} \, , \, \Omega_{\beta^s|J} \right ] =\delta^{IJ} \,
 N_{\alpha^s\beta^s} \Phi_{\alpha^s +\beta^s}
  &
 \left [ E^{\alpha^s} \, , \, E^{\beta^s} \right ] = N_{\alpha^s\beta^s} E^{\alpha^s +\beta^s}
  \\
\mbox{If }
 \alpha^s + \beta^s \in \Phi_{\mathrm{TS}}^s: & \\
\phantom{.}\quad \left\{ \begin{array}{rcl}\left [ \Omega_{\alpha^s|0} \, , \, \Omega_{\beta^s|0} \right ]&
=&
 N_{\alpha^s\beta^s} \Omega_{\alpha^s +\beta^s|0} \\
 %\null \cr
 \left [ \Omega_{\alpha^s|0} \, , \, \Omega_{\beta^s|x} \right ]
&=&
 N_{\alpha^s\beta^s} \Omega_{\alpha^s +\beta^s|x} \\
 %\null \cr
 \left [ \Omega_{\alpha^s|x} \, , \, \Omega_{\beta^s|y} \right ]
&=&
 N_{\alpha^s\beta^s} \left( \delta^{xy} \Omega_{\alpha^s +\beta^s|0}
 + \,  a^{xyz} \, \Omega_{\alpha^s +\beta^s|z} \right) \\
 \end{array}\right. &
 \left [ E^{\alpha^s} \, , \, E^{\beta^s} \right ] = N_{\alpha^s\beta^s} E^{\alpha^s +\beta^s}
  \\
   \hline
\end{array}
 \label{paragonando}
\end{equation}
where $N_{\alpha \beta }=0$ if $\alpha + \beta \notin \Phi_{\mathrm{TS}}$.
\par
\paragraph{\sc Summary.} We can say that for the noncompact, non-maximally split real sections of the simple Lie algebras
(with the only exception of the
$\mathfrak{a}_{\ell}$ series), we have the following four subalgebras:
\begin{equation}\label{cuoredicane}
  \mathbb{U}\supset \mathbb{U}_{\mathrm{TS}}\supset \mathbb{U}_{\mathrm{subTS}}\quad ; \quad \mathbb{U}\supset \mathbb{G}_{\mathrm{Paint}}
  \supset \mathbb{G}_{\mathrm{subpaint}},
\end{equation}
where $\mathbb{U}_{\mathrm{TS}}$ and $\mathbb{U}_{\mathrm{subTS}}$ are both maximally split of total rank $r_{\mathrm{TS}} \, =
\,r_{\mathrm{subTS}} \, = \, r_{n.c.}(\mathbb{U})$, while $\mathbb{G}_{\mathrm{Paint}}$ and $\mathbb{G}_{\mathrm{subpaint}}$ are maximally compact
and have total rank $r_{\mathrm{Paint}}=r_{\mathrm{subpaint}} = r_{\mathrm{comp}}(\mathbb{U})$. Furthermore we have:
\begin{eqnarray}\label{summatheo}
\left[\mathbb{G}_{\mathrm{Paint}}\, , \, \mathbb{U}_{\mathrm{TS}}\right]\, \neq \, 0 \nonumber\\
\left[\mathbb{G}_{\mathrm{subpaint}}\, ,\, \mathbb{U}_{\mathrm{TS}}\right]\, = \, 0\nonumber\\
\left[\mathbb{G}_{\mathrm{Paint}}\, ,\, \mathbb{U}_{\mathrm{subTS}}\right]\, = \, 0
\end{eqnarray}
The Lie subalgebra $\mathbb{U}_{\mathrm{TS}}\subset \mathbb{U}$ has already been described in the previous lines and it is spanned by the
noncompact Cartan generators $H_i$, the step operators $\Phi_{\alpha^\ell}$  associated with the long roots (plus those associated with their
negatives), and  the step operators $\Omega_{\alpha^s\mid 0}$ associated with  the short roots that are in the singlet representation of the
$\mathrm{G}_{subpaint}$ group (plus those associated with their negatives).
\par The subalgebra $\mathbb{U}_{\mathrm{subTS}} \subset\mathbb{U}_{\mathrm{TS}}$ is obtained by dropping also the generators $\Omega_{\alpha^s\mid 0}$
namely restricting oneself only to generators that are paint group singlets.
\par
As we show in detail in  section \ref{sorsgenerale} for the series of Lie algebras $\so(r,r+2s)$, the above algebraic structure reflects into
properties of the symmetric spaces $\mathrm{U/H}$ that admit the following series of projections:
\begin{equation}\label{lollio}
 \mathrm{ \frac{U}{H}}\, \xrightarrow{\pi_{\mathrm{TS}}}\, \mathrm{\frac{U_{TS}}{H_{TS}}}
 \, \xrightarrow{\pi_{\mathrm{subTS}}} \mathrm{\frac{U_{sub}}{H_{subTS}} }
 \end{equation}
 The composition $\pi_{subTS}\, \equiv \, \pi_{\mathrm{sub}} \circ \pi_{\mathrm{TS}}$ allows to write also:
 \begin{equation}\label{crunatore}
    \mathrm{ \frac{U}{H}}\, \xrightarrow{\pi_{\mathrm{subTS}}}  \, \mathrm{\frac{U_{sub}}{H_{subTS}} }.
 \end{equation}
%%%%%%%%%%%%%%%%%%%%%%%%%%%%%%%%%
In the applications to data science it appears of high relevance to investigate the properties of the fibers of these projections, namely:
\begin{equation}\label{fibronilli}
\forall p\in  \mathrm{\frac{U_{TS}}{H_{TS}}} \, : \quad  \mathfrak{F}_{\mathrm{TS}}(p) \, = \, \pi_{\mathrm{TS}}^{-1}\left (p\right) \quad ; \quad
\forall p\in \mathrm{\frac{U_{subTS}}{H_{subTS}}}\, : \quad  \mathfrak{F}_{\mathrm{subTS}}(p) \, = \, \pi^{-1}_{\mathrm{subTS}}\left (p\right).
\end{equation}
These  are  vector spaces with a  linear action of either the subpaint or paint groups. In other words, the manifold
$\mathrm{U/H}$ acquires the structure of the total space of a vector bundle with base manifold either the TS, or the subTS submanifold and
structure group either the subpaint or paint group. In particular  the norm square function on the $\mathrm{U/H}$ symmetric space turns out to
depend only on the points of the base manifold and on a limited number of invariants of that structure group.
\subsubsection{Tits-Satake Universality Classes of Homogeneous Special Geometries}
So far we have analyzed the concepts of Special K\"ahler manifolds and Quaternionic K\"ahler manifolds, and    explained the general construction of
the \emph{$c$-map}, which produces  a Quaternionic K\"ahler manifold starting from a Special K\"ahler  one. The Quaternionic K\"ahler
manifolds that are in the image of the \emph{$c$-map} are called \emph{Special Quaternionic}. Collectively, Special K\"ahler and Special
Quaternionic are called \emph{Special Geometries}. Special Geometries need not  be homogeneous manifolds under the action of a
transitive group of isometries, but some of them are: they are called \emph{Homogeneous Special Geometries}.  Those
discovered in the Supergravity context have been classified, as we already pointed out in the introduction, in
\cite{vandersuppa,specHomgeoA1,specHomgeoA2,deWit:1995tf,Cortes}, and fall into a finite set of series, but some of them are infinite and it might
seem that one ought to examine an infinite number of cases. This is not so because of the Universality classes created by the Tits-Satake
Projection.
\par
What is meant by this   is the following. The Tits-Satake projection   illustrated in the previous subsection  has a number  of
quite  distinctive features:
\begin{enumerate}
  \item $\pi_{\rm TS}$ is a projection operator, so that several different manifolds $\mathcal{SH}_i$ ($i=1,\dots ,r$) may  have the same image
      $\pi_{\rm TS}\left(\mathcal{SH}_i \right) $;
  \item $\pi_{\rm TS}$ preserves the rank of $\mathcal{G}_M$, namely the dimension of the maximal Abelian semisimple subalgebra (Cartan
      subalgebra) of $\mathcal{G}_M$.
  \item $\pi_{\rm TS}$ maps special homogeneous into special homogeneous manifolds: in particular it maps  \emph{special K{\"a}hler} into
      \emph{special K{\"a}hler} and   \emph{Quaternionic} into \emph{Quaternionic}.
  \item $\pi_{\rm TS}$ commutes with  $c$--map, so that we obtain the following commutative diagram:
  \begin{equation}
  \begin{array}{ccc}
   \mbox{Special K{\"a}hler} &
    \stackrel{\mbox{$c$-map}}{\Longrightarrow} & \mbox{Quaternionic-K{\"a}hler} \\
    \begin{array}{cc}
      \pi_{\rm TS} & \Downarrow \\
    \end{array} & \null & \begin{array}{cc}
      \pi_{\rm TS} & \Downarrow \\
    \end{array} \\
    \left( \mbox{Special K{\"a}hler}\right) _{\rm TS} &
    \stackrel{\mbox{$c$-map}}{\Longrightarrow} & \left( \mbox{Quaternionic-K{\"a}hler}\right) _{\rm TS} \
  \end{array}.
\label{diagrammo}
\end{equation}
\end{enumerate}
The main consequence of the above features is that the whole set of special homogeneous manifolds  is distributed into a set of
\emph{universality classes} which turns out to be composed of very few elements.
\par
If we confine ourselves to homogenous symmetric special geometries, then the list of special symmetric manifolds contains only eight items, among
which there are two infinite series. They are displayed in the first column of Table \ref{skTS}.  The $c$-map produces just as many quaternionic (K\"ahler)
manifolds, that  are displayed in the second column of the same table.
\begin{table}[h!]
\begin{center}
$$
\begin{array}{||c|c|c||}
\hline \hline
\null & \null & \null \\
  \mbox{Special K\"ahler}& \mbox{Quaternionic} & \mbox{Tits-Satake projection of Quater.} \\
  \mathcal{SK}_n & \mathcal{QM}_{4n+4} & \mathcal{QM}_{\mathrm{TS}}\\
  \null & \null & \null \\
  \hline
  \null & \null & \null \\
  \frac{\mathrm{U(s+1,1)}}{\mathrm{U(s+1)}\times\mathrm{ U(1)}} &
  \frac{\mathrm{U(s+2,2)}}{\mathrm{U(s+2)}\times\mathrm{ U(2)}}  &
  \frac{\mathrm{U(3,2)}}{\mathrm{U(3)}\times\mathrm{ U(2)}} \\
  \null & \null & \null \\
  \hline
  \null & \null & \null \\
  \frac{\mathrm{SU(1,1)}}{\mathrm{U(1)}}&\frac{\mathrm{G_{(2,2)}}}{\mathrm{SU(2)} \times
  \mathrm{SU(2)}}& \frac{\mathrm{G_{(2,2)}}}{\mathrm{SU(2)} \times \mathrm{SU(2)}} \\
  \null & \null & \null \\
  \hline
  \null & \null & \null \\
  \frac{\mathrm{SU(1,1)}}{\mathrm{U(1)}}\times \frac{\mathrm{SU(1,1)}}{\mathrm{U(1)}}&
  \frac{\mathrm{SO(3,4)}}{\mathrm{SO(3)}\times \mathrm{SO(4)}} &
  \frac{\mathrm{SO(3,4)}}{\mathrm{SO(3)}\times \mathrm{SO(4)}} \\
  \null & \null & \null \\
  \hline
  \null & \null & \null \\
  \frac{\mathrm{SU(1,1)}}{\mathrm{U(1)}}\times \frac{\mathrm{SO(p+2,2)}} {\mathrm{SO(p+2)}\times
  \mathrm{SO(2)}}& \frac{\mathrm{SO(p+4,4)}}{\mathrm{SO(p+4)}\times \mathrm{SO(4)}} &
  \frac{\mathrm{SO(5,4)}}{\mathrm{SO(5)}\times \mathrm{SO(4)}} \\
  \null & \null & \null \\
  \hline
 \begin{array}{c}
 \null\\
     \frac{\mathrm{Sp(6)}}{\mathrm{U(3)}}\\
     \null \\
     \frac{\mathrm{SU(3,3)}}{\mathrm{SU(3)\times SU(3) \times U(1)}} \\
     \null\\
     \frac{\mathrm{SO^\star(12)}}{\mathrm{SU(6)\times  U(1)}} \\
     \null \\
   \frac{\mathrm{E_{(7,-25)}}}{\mathrm{E_{(6,-78)}\times U(1)}}
   \end{array}
 & \begin{array}{c}
 \null\\
     \frac{\mathrm{F_{(4,4)}}}{\mathrm{Usp(6)\times SU(2)}}\\
     \null\\
     \frac{\mathrm{E_{(6,-2)}}}{\mathrm{SU(6)\times SU(2)}}\\
     \null\\
     \frac{\mathrm{E_{(7,-5)}}}{\mathrm{SO(12)\times SU(2)}}\\
     \null \\
      \frac{\mathrm{E_{(8,-24)}}}{\mathrm{E_{(7,-133)}\times SU(2)}}
   \end{array} & \frac{\mathrm{F_{(4,4)}}}{\mathrm{Usp(6)\times SU(2)}} \\
 \null & \null & \null \\
\hline
\end{array}
$$
\caption{The eight series of homogenous symmetric special K\"ahler manifolds (infinite and finite), their
quaternionic counterparts and the grouping of the latter into five Tits-Satake universality
classes.\label{skTS}}
\end{center}
\end{table}
Upon the Tits-Satake projection, this infinite set of models is organized into just five universality classes that are displayed on the third
column of Table \ref{skTS}.
\par
The noncompact symmetric spaces to be utilized in Data Science applications are not necessarily special homogenous spaces, yet in view of the
considerations and plans presented in Section \ref{mittel} all of the spaces that are there in agenda do appear, at the level of
the entire class or of their Tits-Satake projection, in Table \ref{skTS}. This means that the $c$-map has a well founded right to be carefully
considered, and we are going to use it to reduce the menu of building blocks for all cases $r_{n.c.} \leq 4$ to only two items, the Poincar\'e plane and the Siegel plane.
%%%%%%%%%%%%%%%%%%%%%%%%%%%%%%%%%%%%%%%%%%%%%%%%%%%%%%%%%%%%
\subsection{The Symmetric Spaces $\mathrm{SO(r,r+k)/SO(r)\times SO(r+k)}$}
\label{sorrp} Having excluded the symmetric spaces associated with the infinite series of $\mathfrak{a}_\ell$ simple algebras because their Tits-Satake projection yields exotic root systems, we remain with the other three series $\mathfrak{b}_\ell,\mathfrak{c}_\ell,\mathfrak{d}_\ell$.
It appears that, in the context of applications to Data Science a  relevant class of symmetric spaces is the following one:
\begin{equation}\label{parrybar}
\mathcal{M}^{[p,q]}=\frac{\mathrm{G}}{\mathrm{H}}=\frac{{\rm SO}(p,q)}{S[{\rm O}(p)\times {\rm O}(q)]}\,\,,\,\,\,p<q\,.
\end{equation}
These manifolds have rank $p$ and describe the hyperbolic spaces $\mathrm{H}^{(p,q)}$ defined by the equation
\begin{equation}
\sum_{i=0}^{p-1} (x_i)^2-\sum_{j=p}^{n-1} (x_j)^2=1\,,
\end{equation}
where $n=p+q$.
\par
To apply the algebraic machinery described in the previous sections one has to split the family (\ref{parrybar}) in two subfamilies
depending on whether the number $n=p+q$ is even or odd, covering in this way the two classes of real sections associated with the
$\mathfrak{d}_\ell$ and $\mathfrak{b}_\ell$ series of simple Lie algebras, respectively:
\begin{equation}\label{paniepesci}
\mathcal{M}^{[p,q]} \, \rightarrow \,
\begin{array}{cccc}
\null &\mathcal{M}^{[r,s]}_{\mathfrak{d}} &= & \frac{\mathrm{SO(r,r+2s)}}{\mathrm{S[O(r) \times O(r+2s)]}}\\
\nearrow & \null &\null & \null \\
\searrow &\null & \null & \null\\
\null &\mathcal{M}^{[r,s]}_{\mathfrak{b}} &= &\frac{\mathrm{SO(r,r+2s+1)}}{
\mathrm{S[O(r) \times O(r+2s+1)]}}
\end{array}
\end{equation}
\par
The Satake diagrams associated with each of the two families are respectively displayed in Figures \ref{Daldiag} and \ref{Baldiag} where  all the information regarding the final result of the Tits-Satake projection is briefly summarized.
\begin{figure}
\centering
\begin{picture}(100,200)
\put (-120,165){$\mathfrak{d}_{r+s}\left[r,s\right]
\,\, :$}\put(-55,155){$\alpha_1$}\put(-50,170){\circle{10}}\put(-45,170){\line(1,0){20}} \put (-20,170){\circle {10}} \put (-23,155){$\alpha_2$} \put
(-13,170){$\dots$} \put (10,170){\circle {10}} \put (7,155){$\alpha_r$} \put (15,170){\line
(1,0){20}} \put (40,170){\circle*{10}} \put (37,155){$\alpha_{r+1}$} \put (47,170){$\dots$}
%\put (45,230){\line (1,0){20}}
\put (70,170){\circle*{10}} \put (121.2,168){$\alpha_{r+s-2}$} \put (75,170){\line (1,0){25}} \put
(105,170){\circle*{10}} \put (65,155){$\alpha_{r+s-3}$} \put (110,170){\line (1,1){20}} \put
(110,170){\line (1,-1){20}} \put (133.2,193.2){\circle*{10}} \put (133.2,146.8){\circle* {10}} \put
(143.2,193.2){$\alpha_{r+s-1}$} \put (143.2,143.8){$\alpha_{r+s}$}
%%%%%%%%%%%%%%%%%%%%%%%%%%
\end{picture}
\vskip -4cm \caption{\label{Daldiag}
The  Satake diagram of the real section $\mathbb{G}_\mathbb{R}$  of the complex Lie algebra $\mathbb{G}_\mathbb{C} \, = \, \mathfrak{d}_{r+s}$
generating the Lie group $\mathrm{SO(r,r+2s)}$. The total rank of $\mathbb{G}_\mathbb{R}$  is $ \mathrm{rank}= r+s$,
the split rank $ \mathrm{rank}_s = r$, the Tits-Satake subalgebra is $\mathbb{G}_\mathbb{R} \supset\mathbb{G}_{TS} \, = \, \mathfrak{b}_r$,
the paint group is $\mathrm{G_{Paint}}\, = \, \mathrm{SO(2s)}$.}
\end{figure}
%%%%%%%%%%%%%%%%%%%%%%%%%%%%%%%%%%%%%%%%%%%%%%%%%%%%%%%%%%%
\begin{figure}
\centering
\begin{picture}(100,200)
\put (-170,155){$\mathfrak{b}_{r+s}\left[r,s\right]\, \, :$}\put(-100,160){\circle{10}}\put(-103,145){$\alpha_1$}\put(-95,160){\line(1,0){20}}
\put (-70,160){\circle {10}} \put (-73,145){$\alpha_2$}
\put(-60,160){$\dots$}
\put (-40,160){\circle {10}} \put (-43,145){$\alpha_{r-1}$} \put(-35,160){\line(1,0){20}}
%\put (15,160){\line(1,0){20}}
\put(-10,160){\circle{10}}\put(-13,145){$\alpha_{r}$}
\put(-5,160){\line(1,0){20}}
\put (20,160){\circle* {10}} \put (17,145){$\alpha_{r+1}$} \put (29,160){$\dots$}
%\put (45,260){\line (1,0){20}}
\put (50,160){\circle* {10}} \put (47,145){$\alpha_{r+s-2}$} \put (55,160){\line (1,0){40}} \put
(100,160){\circle* {10}} \put (97,145){$\alpha_{r+s-1}$} \put (105,163){\line (1,0){40}} \put (105,158){\line
(1,0){40}} \put (118,157){{\LARGE$>$}} \put (150,160){\circle* {10}} \put (147,145){$\alpha_{r+s}$}
\end{picture}
\vskip -4cm \caption{\label{Baldiag} The Tits-Satake diagram of the real section $\mathbb{G}_\mathbb{R}$  of the complex Lie algebra
$\mathbb{G}_\mathbb{C} \, = \, \mathfrak{b}_{r+s}$  generating the Lie group $\mathrm{SO(r,r+2s+1)}$. The total rank of $\mathbb{G}_\mathbb{R}$
is $ \mathrm{rank}= r+s$, the split rank $ \mathrm{rank}_s = r$, the Tits-Satake subalgebra is $\mathbb{G}_\mathbb{R} \supset\mathbb{G}_{\mathrm{TS}}
\, = \, \mathfrak{b}_r$, the paint group is $\mathrm{G_{Paint}}\, = \, \mathrm{SO(2s+1)}$.}
\end{figure}
\par
Both series of Lie algebras have the same \emph{Tits-Satake subalgebra}, shown in Figure \ref{TSAlgdiag}, and the same \emph{sub Tits
Satake subalgebra}, shown in Figure \ref{subTSAlgdiag}.
\begin{figure}
\centering
\begin{picture}(100,200)
\put (-115,155){$\mathbb{G}_{\mathrm{TS}} \, = \,\mathfrak{b}_{r}\left[r,0\right]\, \,:$}\put(-10,160){\circle{10}}\put(-13,145){$\alpha_{1}$}
\put(-5,160){\line(1,0){20}}
\put (20,160){\circle {10}} \put (17,145){$\alpha_{2}$} \put (29,160){$\dots$}
%\put (45,260){\line (1,0){20}}
\put (50,160){\circle {10}} \put (47,145){$\alpha_{r-2}$} \put (55,160){\line (1,0){40}} \put
(100,160){\circle {10}} \put (97,145){$\alpha_{r-1}$} \put (105,163){\line (1,0){40}} \put (105,158){\line
(1,0){40}} \put (118,157){{\LARGE$>$}} \put (150,160){\circle {10}} \put (147,145){$\alpha_{r}$}
\end{picture}
\vskip -4cm \caption{\label{TSAlgdiag} The Tits-Satake diagram of the Tits-Satake subalgebra $\mathrm{SO(r,r+1)}$.}
\end{figure}
\begin{figure}
\centering
\begin{picture}(100,100)
\put (-125,65){$\mathbb{G}_{\mathrm{subTS}} \, = \,\mathfrak{d}_{r}\left[r,0\right]\, \,:$}
 \put (-20,70){\circle {10}} \put (-23,55){$\alpha_1$} \put (-15,70){\line
(1,0){20}} \put (10,70){\circle {10}} \put (7,55){$\alpha_2$} \put (15,70){\line (1,0){20}} \put
(40,70){\circle {10}} \put (37,55){$\alpha_3$} \put (47,70){$\dots$}
%\put (45,230){\line (1,0){20}}
\put (70,70){\circle {10}} \put (67,55){$\alpha_{r-3}$} \put (75,70){\line (1,0){25}} \put
(105,70){\circle {10}} \put (100,55){$\alpha_{r-2}$} \put (110,70){\line (1,1){20}} \put (110,70){\line
(1,-1){20}} \put (133.2,93.2){\circle {10}} \put (133.2,46.8){\circle {10}} \put
(143.2,93.2){$\alpha_{r-1}$} \put (143.2,43.8){$\alpha_r$}
\end{picture}
\vskip -1cm \caption{\label{subTSAlgdiag} The sub Tits-Satake diagram of the subTS subalgebra $\mathrm{SO(r,r)}$.}
\end{figure}
Hence once the noncompact (or split) rank $r$ is fixed, the two infinite series $\frac{\mathrm{SO(r,r+2s)}}{\mathrm{S[O(r) \times O(r+2s)]}}$ and
$\frac{\mathrm{SO(r,r+2s+1)}}{\mathrm{S[O(r) \times O(r+2s+1)]}}$ project onto the same maximally split  Tits-Satake space
$\mathrm{\frac{SO(r,r+1)}{SO(r)\times SO(r+1)}}$, or onto the same maximally split sub Tits-Satake space $\mathrm{\frac{SO(r,r)}{SO(r)\times
SO(r+1)}}$,  making up a single universality class. What distinguishes the various spaces inside the universality class are the  \emph{Paint
Groups}, either $\mathrm{SO(2s)}$ or $\mathrm{SO(2s+1)}$. The roots of the restricted system $\overline{\Delta}$ split  into two sets
$\Delta^\eta$ and $\Delta^\delta$, which are respectively singlets and either $2s$-tuplets or $(2s+1)$-tuplets under the action of the appropriate
paint group. This means that the solvable coordinates of the coset manifolds split into three classes: the duals of  $r$ Cartan generators,
the duals of $r(r-1)$ long root generators, and the dual of the  $r\times (2s)$ or $r\times (2s+1)$ short roots. The plethora of   short root
solvable coordinates are actually a small set $r$ of solvable coordinates that carry a second hidden index, the paint index $I$ taking a large
number of values.
\par
Because of the above general conclusion, in order to study the full Universality Classes of manifolds (\ref{parrybar}), labeled by the noncompact
rank $r$, it is sufficient to choose one of the two subclasses displayed in eq.~(\ref{paniepesci}). All the final results for the distance
formulas, the Paint and sub Paint invariants and the geometry in the projected space apply to both cases.
We have chosen throughout this paper to deal with the case $\mathcal{M}^{[r,s]}_{\mathfrak{d}}$ related with the $\mathfrak{d}_\ell$ series.
%%%%%%%%%%%%%%%%%%%%%%%%%%%%%%%%%%%%%%%%%%%%%%%%%%%%%%%%%%%%
\subsection{The Lie algebras $\so(r,r+2s)$ and their solvable Lie subalgebra}
\label{sorsgenerale} For this family of noncompact real sections of the complex $\mathfrak{d}_{\ell=r+s}$ Lie algebras, we refer to the Tits-Satake diagram  displayed in Figure \ref{Daldiag}, where all the information regarding the final result of the Tits-Satake
projection is briefly summarized. Recapping what we said in the previous section, the total number of positive roots in the
$\Phi_{\mathfrak{d}_{r+s}}$ root system is
\begin{eqnarray}\label{contidiserva}
  \frac{1}{2} \text{card}[\Phi_{\mathfrak{d}_{r+s}}] & = & (r+s)^2 -(r+s) \nonumber\\
   & = & \underbrace{r^2-r}_{\text{long roots}}\, + \,
  \underbrace{2\, r \, s }_{\text{short roots}}
  \, + \, \underbrace{s^2-s}_{\text{roots of $\mathrm{G_{Paint}}$}}
\end{eqnarray}
where the various addends are interpreted as shown above. Indeed we know that the Tits-Satake root system $\Phi_{TS}$ corresponds to the Dynkin
Diagram $\mathfrak{b}_{r}$ as displayed in Figure \ref{TSAlgdiag}. This is the root system of the complex Lie algebra $\so(2r+1)$, whose maximally
split real section is $\so(r,r+1)$. The number of positive roots in $\Phi_{TS}$ therefore is
\begin{equation}\label{belluccio}
  \frac{1}{2} \text{card}[\Phi^{[r,s]}_{TS}] \, = \, r^2
\end{equation}
As it is evident from the very form of the $\mathfrak{b}_{r}$  Dynkin diagram, in the corresponding root system  $\Phi^{[r,s]}_{TS}$ there are
exactly $r^2-r$ long roots and $r$ short roots. In the root system $\Phi_{\mathfrak{d}_{r+s}}$ we  have one preimage of each of the longroots of
$\Phi^{[r,s]}_{TS}$ and $2s$ pre-images of each of the $r$ short-roots. The remaining $s^2-s$ roots are all imaginary and, in the complete
$\so(r,r+2s)$ Lie algebra, their associated step operators appear only through their maximally compact combinations $E^{\alpha} - E^{-\alpha}$ and
$i (E^{\alpha}+ E^{-\alpha})$, generating, together with noncompact Cartan generators, the paint group $\mathrm{SO(}2s\mathrm{)}$. In the
$\so(r,r+2s)$ Lie algebra, the step operators associated with the $2s$ copies of each short root are linearly transformed one into the other under
the adjoint action of the Paint group $\mathrm{SO(}2s\mathrm{)}$.
\par
The solvable Lie algebra $Solv_{[r,s]}$ associated with this real section is made by the $r$ noncompact Cartan generators $\mathcal{H}_i$, the
$r^2-r$ step operators $E^{\alpha^l}$ associated with positive long roots, and the $2\,rs$ step operators associated with pre-images of the short
roots $E^{\alpha^s_I}$ ($I,1,\dots,2s$). Therefore the total dimension of the solvable Lie algebra is
\begin{equation}\label{gruggo}
 \dim Solv_{[r,s]} \, = r \, + \,r^2-r \, +\,  2\, r \, s \, = r (r+2s) \, = \, \dim  \mathrm{\frac{SO(r,r+2s)}{SO(r)\times SO(r+2s)} }
\end{equation}
as it should be.
\par
A basis of generators for the $\so(r,r+2s)$ Lie algebra, where the above abstract  structure is made explicit, is obtained through the
construction in dimension $\mathrm{N}=2r+2s$ of the matrix $\eta_t$ provided by Statement \ref{statamento}. By the definition of the Lie algebra,
the matrix $\eta_t$ must have signature
\begin{equation}\label{pincuspallus}
 \text{sign}_{r,s} \, = \, \left\{\underbrace{+,\dots, +}_{r+2s \text{ times}},\underbrace{-,\dots,-}_{r \text{ - times}} \right\} .
\end{equation}
The Lie algebra $\so(r,r+2s)$ is then made by the $N\times N$ matrices   fulfilling the   condition
\begin{equation}\label{gornayalavandaia}
  X \in \so(r,r+2s) \quad \Leftrightarrow \quad  X^T \, \eta_t \, + \, \eta_t \, X \, = \, 0.
\end{equation}
The appropriate choice of $\eta_t$ is
\begin{equation}\label{etatdefi}
  \eta \, = \, \eta_t \, = \, \left(
\begin{array}{ccc|c|ccc}
 0 & \cdots & 0 & 0 & 0 & \ldots & 1 \\
 \vdots & \vdots & \vdots & \vdots & \vdots & 1 & \vdots \\
 0 & \ldots & 0 & 0 & 1 & \ldots & 0 \\
 \hline
 0 & \ldots & 0 & \mathbf{1}_{2s \times 2s} & 0 & \ldots & 0  \\
 \hline
 0 & \ldots & 1 & 0 & 0 & \ldots & 0 \\
 \vdots & 1 & \vdots & \vdots & \vdots & \vdots & \vdots \\
 1 & \cdots & 0 & 0 & 0 & \cdots & 0 \\
\end{array}
\right).
\end{equation}
The matrix $\eta_t$ defined above is transformed into the standard one
\begin{equation}\label{etabdefi}
  \eta_b \, \equiv \, \text{diag}\left\{\underbrace{+,\dots, +}_{r+2s \text{ times}},\underbrace{-,\dots,-}_{r \text{ - times}} \right\}
\end{equation}
by the conjugation by a matrix $\Omega$ that can be built for all pairs $r\geq 1,s\geq 1$, and     also has  a general form that  can be easily read off  from its expression in the case $r=3,s=2$:
 \begin{equation}\label{omegagrande}
   \Omega \, = \, \left(
\begin{array}{ccc|cccc|ccc|c}
\hline
\null& r & \null & \null &2 & s&\null & \null & r &\null \\
\hline
 0 & 0 & 0 & 1 & 0 & 0 & 0 & 0 & 0 & 0 &\null \\
 0 & 0 & 0 & 0 & 1 & 0 & 0 & 0 & 0 & 0 & 2 s\\
 0 & 0 & 0 & 0 & 0 & 1 & 0 & 0 & 0 & 0 &\null \\
 0 & 0 & 0 & 0 & 0 & 0 & 1 & 0 & 0 & 0 & \null\\
 \hline
 \frac{1}{\sqrt{2}} & 0 & 0 & 0 & 0 & 0 & 0 & 0 & 0 & \frac{1}{\sqrt{2}} & \null \\
 0 & \frac{1}{\sqrt{2}} & 0 & 0 & 0 & 0 & 0 & 0 & \frac{1}{\sqrt{2}} & 0 & r\\
 0 & 0 & \frac{1}{\sqrt{2}} & 0 & 0 & 0 & 0 & \frac{1}{\sqrt{2}} & 0 & 0 &\null \\
 \hline
 \frac{1}{\sqrt{2}} & 0 & 0 & 0 & 0 & 0 & 0 & 0 & 0 & -\frac{1}{\sqrt{2}}& \null \\
 0 & \frac{1}{\sqrt{2}} & 0 & 0 & 0 & 0 & 0 & 0 & -\frac{1}{\sqrt{2}} & 0 & r \\
 0 & 0 & \frac{1}{\sqrt{2}} & 0 & 0 & 0 & 0 & -\frac{1}{\sqrt{2}} & 0 & 0 & \null \\
\end{array}
\right) .
 \end{equation}
For all values of $r$ and $s$ one has
\begin{equation}\label{etatinetab}
  \Omega \, \eta_t \, \Omega^T \, = \, \eta_b.
\end{equation}
Utilizing the basis $\eta_t$, the solvable subalgebra is easily identified, as, by construction,  it is made by  the matrices that satisfy
the defining condition \eqref{gornayalavandaia} of the Lie algebra $\so(r,r+2s)$,   and  in addition are   upper triangular.
So doing the procedure to explicitly construct the algebra becomes straightforward and best suited to computer algorithms. In that
construction one automatically derives a set of Cartan-Weyl generators that are well well-adapted to the particularly chosen invariant metric
$\eta_t$. In this basis the Cartan involution $\theta$ defining the orthogonal decomposition (\ref{wiccus}) is simply given by matrix
transposition:
\begin{equation}\label{thetatranspo}
 \theta (X) \, = \, X^T \quad  \forall X \in \so(r,r+2s)\,.
\end{equation}
Hence, in accordance   with the notations introduced in the general analysis of the Tits-Satake projection in section \ref{colorato}, the generators of
the solvable Lie algebras are the following:
\begin{equation}\label{generalini}
\boldsymbol{T}\, = \,\left\{\begin{array}{clclcl}
  H_i & =\text{Cartan generators} \quad  &i&=1,\dots, r \quad &\null&\null \\
  \Phi_\alpha^\ell &=  \text{Long root generators} \quad  &\alpha^\ell &=1,\dots, r^2 -r \quad &\null&\null\\
  \Omega_{\alpha^s| I} &= \text{Short root generators} \quad  &\alpha^s &=1,\dots, r \quad  &I&=1,\dots,2s \\
  J_{x} &= \text{generators of the Paint group} \quad &x&=1,\dots, s(2s-1).
  \end{array}\right.
\end{equation}
In the $\eta_t$ basis the matrices representing $H_i$ are diagonal, those representing $\Phi_\alpha^\ell$ and $\Omega_{\alpha^s| I}$ are upper
triangular, while those representing $J_x$ are skew-symmetric.
\par
Correspondingly one introduces the set of solvable coordinates for the coset manifold, collectively denoted $\boldsymbol{\Upsilon}$, with the
  structure
\begin{equation}\label{kyrillov}
  \boldsymbol{\Upsilon} \, =\,  \left\{ \underbrace{\Upsilon_i}_{\text{n.c. Cartan}} \mid
  \underbrace{ \Upsilon_{r+\alpha^\ell}}_{\text{long roots}}\, \mid \underbrace{{\Upsilon}_{r^2+\alpha^s,I}}_{2s \,
  \,\text{preimages  of each short root 5}} \right\},
\end{equation}
and the solvable Lie algebra element singled out by these coordinates is obtained as
\begin{equation}\label{barlumedisperanza}
 \boldsymbol{\mathfrak{s}}(\boldsymbol{\Upsilon}) \, = \, \boldsymbol{\Upsilon} \cdot \boldsymbol{T} \in Solv
\end{equation}
where $\boldsymbol{T}$ denotes the collection of solvable Lie algebra generators in the above described order.
\par
Once the solvable Lie algebra parameterization is established, we need to perform the map of the Lie algebra to the solvable group:
\begin{equation}\label{sigmaexp}
  \Sigma \, : \, Solv \, \stackrel{\Sigma}{\longrightarrow} \, \mathcal{S} \, \equiv \, \exp \left[Solv\right]
\end{equation}
The explicit form of the exponential map $\Sigma$ corresponds to a choice of the solvable coordinates normalization and   is simply a
matter of convenience, as the choice of Euler angles in the parameterization of rotation group elements. Our choice, which optimizes the metric
and  the Laplacian to their simplest possible form, is the following:
\begin{eqnarray}\label{SigmaExp1}
  \Sigma\left[\mathfrak{s}\left(\boldsymbol{\Upsilon} \right)\right]
  & = & \exp\left[\mathfrak{s}_{C}\left(\boldsymbol{\Upsilon} \right)\right]
  \cdot \exp\left[\mathfrak{s}_{long}\left(\boldsymbol{\Upsilon} \right)\right]\cdot\exp\left[\mathfrak{s}_{short}\left(\boldsymbol{\Upsilon}
  \right)\right]\nonumber\\
  \mathfrak{s}_{C}\left(\boldsymbol{\Upsilon} \right) &=& \text{projection onto the
  Cartan subalgebra of $Solv$} \nonumber\\
  \mathfrak{s}_{long}\left(\boldsymbol{\Upsilon} \right) & =& \text{projection onto the
  subspace  spanned by the long root generators}\nonumber\\
  \mathfrak{s}_{short}\left(\boldsymbol{\Upsilon} \right) & =& \text{projection onto the
  subspace  spanned by the short root generators}\nonumber\\
  \mathfrak{s}\left(\boldsymbol{\Upsilon} \right) &=&\mathfrak{s}_{C}\left(\boldsymbol{\Upsilon} \right)\oplus
  \mathfrak{s}_{long}\left(\boldsymbol{\Upsilon} \right) \oplus \mathfrak{s}_{short}\left(\boldsymbol{\Upsilon} \right)
\end{eqnarray}
The image of the $\Sigma$-map will be denoted $\mathbb{L}(\boldsymbol\Upsilon)$,  and for  the reader's convenience we display below, fixing the simplest
value $s=1$, the explicit form of the matrix in the cases $r=1,2$. We have
\begin{equation}\label{Lr1}
  \mathbb{L}(\boldsymbol\Upsilon)^{[1,1]} \, = \, \left(
\begin{array}{cccc}
 e^{\Upsilon _1} & \frac{e^{\Upsilon _1} \Upsilon _{2,1}}{\sqrt{2}} &
   \frac{e^{\Upsilon _1} \Upsilon _{2,2}}{\sqrt{2}} & -\frac{1}{4} e^{\Upsilon _1}
   \left(\Upsilon _{2,1}^2+\Upsilon _{2,2}^2\right) \\
 0 & 1 & 0 & -\frac{\Upsilon _{2,1}}{\sqrt{2}} \\
 0 & 0 & 1 & -\frac{\Upsilon _{2,2}}{\sqrt{2}} \\
 0 & 0 & 0 & e^{-\Upsilon _1} \\
\end{array}
\right)
\end{equation}
for the case $r=1$, while the case $r=2$ is
\begin{multline}\label{exempli2}
   \mathbb{L}(\boldsymbol\Upsilon)^{[2,1]} =  \\[3pt] \left(
\scriptsize \begin{array}{cccccc}
 e^{\Upsilon _1} & \frac{e^{\Upsilon _1} \Upsilon _3}{\sqrt{2}} & \frac{1}{2}
   e^{\Upsilon _1} \left(\sqrt{2} U_1+\Upsilon _3 V_1\right) & \frac{1}{2}
   e^{\Upsilon _1} \left(\sqrt{2} U_2+\Upsilon _3 V_2\right) & -\frac{1}{8}
   e^{\Upsilon _1} \left(4 \mathbf{U}\cdot \mathbf{V}+\sqrt{2} \left(\Upsilon _3 \mathbf{V}^2-4 \Upsilon
   _4\right)\right) & -\frac{1}{4} e^{\Upsilon _1} \left(\mathbf{U}^2+2 \Upsilon _3
   \Upsilon _4\right) \\
 0 & e^{\Upsilon _2} & \frac{e^{\Upsilon _2} V_1}{\sqrt{2}} & \frac{e^{\Upsilon
   _2} V_2}{\sqrt{2}} & -\frac{1}{4} e^{\Upsilon _2} \mathbf{V}^2 & -\frac{e^{\Upsilon _2}
   \Upsilon _4}{\sqrt{2}} \\
 0 & 0 & 1 & 0 & -\frac{V_1}{\sqrt{2}} & -\frac{U_1}{\sqrt{2}} \\
 0 & 0 & 0 & 1 & -\frac{V_2}{\sqrt{2}} & -\frac{U_2}{\sqrt{2}} \\
 0 & 0 & 0 & 0 & e^{-\Upsilon _2} & -\frac{e^{-\Upsilon _2} \Upsilon _3}{\sqrt{2}}
   \\
 0 & 0 & 0 & 0 & 0 & e^{-\Upsilon _1} \\
\end{array}
\right)
\end{multline}

\par
% In eq.~(\ref{exempli2})
Here  we have used the notation $U_i = \Upsilon_{5,i}$ and  $V_i = \Upsilon_{6,i}$, which stresses the existence of two
Paint vectors $\mathbf{U},\mathbf{V}$ in the case $r=2$ and the paint group covariant structure of the solvable group element
$\mathbb{L}(\boldsymbol\Upsilon)^{[2,s]}$. Indeed from eq.~(\ref{exempli2}) it is immediate to deduce the general form of the matrix for any value
of $2s$.
%%%%%%%%%%%%%%%%%%%%%%%
\section{Geodesics in the symmetric spaces $\mathrm{U/H}$}
\label{cosettus} Next, we show that, for all noncompact symmetric spaces $\mathrm{U/H}$ of the type discussed in this paper, the geodesic
equations can be explicitly integrated;   there is indeed an explicit and very simple integration formula that yields the unique geodesic passing
through any two given points $p_{1,2} \in \mathrm{U/H}$.
\subsection{Integrating the geodesic equation}\label{lavatoconperlana}
To begin with, we recall that the Lie algebra $\mathbb{U}$ of $\mathrm{U}$ admits a Cartan involution $\theta$, and   splits  into the  $+1$ and $-1$ eigenspaces
 of $\theta$:\begin{equation}
\label{cartandecompo}
\mathbb{U}=\mathbb{H}\oplus \mathbb{K},
\end{equation}
where   $\mathbb{H}$ is the Lie algebra of $\mathrm{H}$, and $\mathbb{K}$ is its orthogonal complement, namely, the space of coset generators.
In the real matrix representation of $\mathrm{U}$ in the basis provided by  the statement \ref{statamento}, the involution $\theta$ is just
matrix transposition, as we already recalled in eq.~(\ref{thetatranspo}). In this representation   $\mathbb{H}$ consists of   skew-symmetric
matrices while $\mathbb{K}$ of  symmetric matrices. Using this basis  we are relying on the embedding of any of the classical noncompact
simple groups into $\mathrm{SL(N,\mathbb{R})}$ where $\mathrm{N}$ is the dimension of the fundamental representation of $\mathrm{U}$. In our basis
the elements of the group $\mathcal{O}\in\mathrm{SO(r,r+q)}$ and those of the algebra $h\in \so(r,r+q)$ respectively satisfy:
\begin{equation}\label{samarcanda}
  \mathcal{O}^T \, \eta_t \,  \mathcal{O} \, = \, \eta_t \quad ; \quad h^T \, \eta_t \, + \, \eta_t \, h \, = \, 0
\end{equation}
where the symmetric $\eta_t$  matrix is that defined in eq.~(\ref{etatdefi}).
\par
In the $\eta_t$ basis the elements of the $\mathbb{K}$ subspace are both $\eta_t$--skew-symmetric and symmetric tout-court, while the elements of
the subalgebra $\mathbb{H}$ are both $\eta_t$-skew-symmetric and skew-symmetric tout-court. When dealing with elements of $\mathrm{U}$ and of its
Lie algebra, we shall always refer to the chosen matrix representation.
\par
The $\mathbb{K}$ vector space  is isomorphic to the tangent space to $\mathrm{U}/\mathrm{H}$ at the origin $\mathcal{O}$:
\begin{equation}\label{trofarello}
\mathbb{K}\sim T_\mathcal{O}\left(\mathrm{U/H}\right)\,.
\end{equation}
According to the Iwasawa decomposition, a   matrix $\mathfrak{g} \in\mathrm{U}$ can always be written as the product of an upper-triangular
matrix $\mathcal{S}$ and a compact one $h\in \mathrm{H}$: $\mathfrak{g}=\mathcal{S}\,h\,$, $\mathcal{S}$, $h$ being uniquely associated with every
$\mathfrak{g}$. The matrix $\mathcal{S}$, being upper-triangular, is an element of the solvable subgroup
$\exp\left[Solv_{\mathrm{U/H}}\right]\subset \mathrm{U}$.  It follows that within each left coset in $\mathrm{U/H}$ one can find
a \emph{unique} representative $\mathbb{L}\in \exp\left[Solv_{\mathrm{U/H}}\right]$, whose parameters $\Upsilon_A$ appearing in the decomposition
of a generic element of $Solv_{\mathrm{U/H}}$ along a basis of generators $\mathrm{T}^A$, $A=1,\dots, {\rm dim}(\mathrm{U/H})$, that is,
$  \mathfrak{s}\, = \,
 \Upsilon_A \, \mathrm{T}^A
$,
 define the solvable parameterization introduced earlier. The group
 $\mathcal{S}_{\mathrm{U/H}}\equiv
 \exp\left[\Solv_{\mathrm{U/H}}\right]$ has  a transitive
 action on $\mathrm{U/H}$, as we know from previous sections.
 While $\mathrm{T}^A$ provide a basis of generators of $\Solv_{\mathrm{U/H}}$, let us denote  by
 $\mathrm{K}_A$ a basis of generators of the subspace $\mathbb{K}$. As explained in the previous sections, we
 choose the coset representative to be an element $ \mathbb{L}(\boldsymbol{\Upsilon})$.
\par
Utilizing a normalization such that ${\rm Tr}(\mathrm{K}_A\,\mathrm{K}_B)\,=\,\delta_{AB}$ and  ${\rm
Tr}(\mathrm{H}_A\,\mathrm{H}_B)\,=\,\delta_{AB}$, the vielbein of the symmetric space is computed through the left-invariant 1-form of the
solvable group as follows:
\begin{eqnarray}\label{filbone}
\Theta \equiv\mathbb{L}^{-1}d\mathbb{L}\quad ; \quad \mathbb{V} \, = \, V^A \, \mathrm{K}_A \quad ; \quad V^A \, = \,
\text{Tr} \left(\Theta \cdot \mathrm{K}_A\right),
\end{eqnarray}
while the $\mathbb{H}$-connection is provided by
\begin{eqnarray}\label{filbone} \mathbb{Q} \, = \, Q^A \, \mathrm{H}_A \quad ; \quad Q^A \, = \, \text{Tr} \left(\Theta \, \mathrm{H}_A\right),
\end{eqnarray}
and one can write the projections of $\Theta=\mathbb{Q}+\mathbb{V} $. Besides being both $\eta_t$-skew-symmetric, the $1$-form valued matrices
$\mathbb{Q}$ and $\mathbb{V}$ are respectively tout-court skew-symmetric and symmetric, respectively:
\begin{equation}\label{pritaneo}
\mathbb{V}=\mathbb{V}^T=\mathbb{V} \quad ; \quad \mathbb{Q}^T=-\mathbb{Q}
\end{equation}
Because of the Maurer Cartan equation $d\Theta +\Theta \wedge \Theta \, = \, 0$ satisfied by the left-invariant one-form $\Theta$, we have the
identities
\begin{align}
\mathcal{D}\mathbb{V}\equiv d\mathbb{V}+\mathbb{Q}\wedge \mathbb{V}
+\mathbb{V}\wedge \mathbb{Q}=0\,\,,\,\,\,\,\mathfrak{R}[\mathbb{Q}]\equiv d\mathbb{Q}+\mathbb{Q}\wedge \mathbb{Q}=-\mathbb{V}\wedge \mathbb{V}
\end{align}
The $\mathrm{U}$-invariant metric on the symmetric space $\mathrm{U/H}$ has the form
\begin{equation}
ds^2={\rm Tr}(\mathbb{V}^2)=\sum_{A=1}^{dn} \,V^A\otimes V^A
\end{equation}
In order to conveniently describe the geodesics on $\mathrm{U/H}$, we introduce the   matrix
\begin{equation}\label{turlipano}
\mathcal{M}(\Upsilon)\,\equiv\, \mathbb{L}(\Upsilon)\,\mathbb{L}(\Upsilon)^T.
\end{equation}
This   is a function merely of the point on $\mathrm{U/H}$, in that it does not depend on the choice of the coset representative:
\begin{equation}\label{summeo}
 \mathbb{L}(\Upsilon)\rightarrow\,\mathbb{L}(\Upsilon)\,h\,\,\Rightarrow\,\,\,\mathcal{M}\,\rightarrow\,\mathbb{L}(\Upsilon)h\,h^T\,\mathbb{L}(v)^T=
\mathbb{L}(\Upsilon)\mathbb{L}(\Upsilon)^T=\mathcal{M}(\Upsilon)\,,
\end{equation}
where we have used the property that, in the chosen real matrix representation, an element $h\in \mathrm{H}$ is an orthogonal matrix. In terms of
the matrix $\mathcal{M}(\Upsilon)$ we can define the   $\mathbb{U}$-valued 1-form
\begin{equation}\label{perbacco}
\mho\,\equiv \,\mathcal{M}^{-1}\,d\mathcal{M}=2\,\mathbb{L}^{-T}\,\mathbb{V}\,\mathbb{L}^{T}\,.
\end{equation}
The last identity in eq.~(\ref{perbacco}) follows from the short calculation
\begin{alignat}{4}\label{calculatia}
 \mho\ &= &\mathbb{L}^{-T}\,\mathbb{L}^{-1} \,d\mathbb{L} \,\mathbb{L} ^{T}\, + \, \mathbb{L}^{-T}d\mathbb{L}^{T}  &= &\,\mathbb{L}^{-T}\,
 \Theta \,\mathbb{L} ^{T} + \, \mathbb{L}^{-T}\Theta^T\mathbb{L}^{T} \nonumber \\
 \null&=& \mathbb{L}^{-T}\, \left(\Theta^T \,+\, \Theta\right)\mathbb{L}^T
&=& 2  \mathbb{L}^{-T}\, \mathbb{V} \,\mathbb{L}^T,
\end{alignat}
where we have used the short-hand notation $M^{-T}\equiv (M^{-1})^T$ and the fact that the Vielbein 1-form is the projection of the left-invariant
matrix valued one-form  $\Theta$ onto its symmetric part. Thus we rewrite the metric as
\begin{equation}\label{cruccio}
ds^2 \,\propto\, {\rm Tr}(\mho \, \mho)\,.
\end{equation}
Consider now a geodesic described by the functions
\begin{equation}
\Upsilon_A=\Upsilon_A(t;\,\Upsilon^0)\,.
\end{equation}
where  $\Upsilon^0\equiv (\Upsilon^0_A)$ denotes the initial point on the manifold:
\begin{equation}
\Upsilon_A(t=0;\,\Upsilon^0)=\Upsilon^0_A\,.
\end{equation}
The matrix valued vielbein 1-form and $\mathbb{H}$-connection have the following canonical expansion in coordinate differentials:
\begin{equation}\label{pirotello}
  \mathbb{V} \, = \, \mathbb{V}^A(\Upsilon) \, \mathrm{d}\Upsilon_A \quad ; \quad \mathbb{Q} \, = \, \mathbb{Q}^A(\Upsilon) \, \mathrm{d}\Upsilon_A
\end{equation}
The geodesic is a curve in the manifold, namely a map of the real line to the symmetric space
\begin{equation}\label{curvettaman}
  \mu \, : \, \mathbb{R} \, \,\rightarrow \, \frac{\mathrm{U}}{\mathrm{H}}.
\end{equation}
By   pull-back we obtain
\begin{equation}\label{fringuello}
\begin{array}{rclcrcl}
  \mu^\star\left(\mathbb{V}\right) & = & L(t) \,\mathrm{d}t &\quad ; \quad & \mu^\star\left(\mathbb{Q}\right) & = & W(t) \,\mathrm{d}t \\
  L(t) &\equiv &\dot{\Upsilon}_A\,\mathbb{V}^A\left(\Upsilon(t;\,\Upsilon^0)\right) &\quad ; \quad  &
   W(t) & \equiv &\dot{\Upsilon}_A\,\mathbb{Q}^A\left(\Upsilon(t;\,\Upsilon^0)\right)\\
   \end{array}
\end{equation}
where we have denoted
\begin{equation}\label{lampsaco}
  \dot{\Upsilon}_A\equiv \frac{d}{dt}\Upsilon_A(t;\,\Upsilon_0)\,.
\end{equation}
As it was shown in \cite{Chemissany:2009hq},  the geodesic equations can be expressed as  follows:
\begin{statement}
\label{statamento2} The matrix valued $1$-form $\mho$ is constant along geodesics, namely:
\begin{equation}\label{diodorosiculo}
 \frac{\mathrm{d}}{\mathrm{d}t}\mu^\star\left(\mho\right) \, = \, 0
\end{equation}
which implies:
\begin{equation}
\mathcal{M}(\Upsilon(t;\,\Upsilon_0))^{-1}\frac{d}{dt}\mathcal{M}(\Upsilon(t;\,\Upsilon_0))=Q^T={\rm const}\label{geoM}
\end{equation}
The constant matrix $Q\in \mathbb{U}$ encodes the conserved Noether charges of the solution.
\end{statement}
Indeed $Q$ is a constant matrix in the tangent space $T_{\Upsilon_0}(\mathrm{U/H})$ to the Riemannian  manifold $\mathrm{U/H}$ at the initial
point $\Upsilon_0$, and defines the \emph{initial velocity} of the geodesic. The \emph{simple transitive action of the solvable group} on the
symmetric space $\mathrm{U/H}$ implies that $\mathbb{K}$ is isomorphic to the tangent space at the origin $T_\mathbf{e}(\mathrm{U/H})$,
$\mathbf{e}$ being defined by $\Upsilon^A=0$, and that we have:
\begin{equation}
\label{calderondelabarca}
Q\in \mathbb{L}(\Upsilon_0)\,\mathbb{K}\,\mathbb{L}(\Upsilon_0)^{-1}\sim T_{\Upsilon_0}\left(\mathrm{U/H}\right)\,
\end{equation}
Let us show that, from eq.~ (\ref{geoM}), a Lax equation follows:
\begin{align}
0&=\frac{d}{dt}\left(\mathcal{M}(\Upsilon(t;\,\Upsilon_0))^{-1}\frac{d}{dt}\mathcal{M}(\Upsilon(t;\,\Upsilon_0))\right)=
\mathbb{L}^{-T}\,(\dot{L}+[L,\Omega^T])\,
\mathbb{L}^{T}\nonumber\\
&\Downarrow \nonumber\\
0&=\dot{L}+[L,\,\Omega^T]=\dot{L}+[L,\,W^T]=\dot{L}+[W,\,L]\,
\end{align}
where we have used $W^T=-W$ and in the above derivation the coset representative $\mathbb{L}$ is evaluated on $\Upsilon^A(t;\,\Upsilon_0)$. In the
last equation $\dot{L}+[W,\,L]=0$, $W$ and $L$ are, respectively, the skew-symmetric and symmetric projections of the same matrix
$\dot{\Upsilon}_A\,\Omega_A(\Upsilon(t;\,\Upsilon^0))$, and thus:
\begin{equation}
\label{mevoici}
W=L_{>}-L_{<}\,
\end{equation}
eq.~~\eqref{mevoici} provides the precise link between the approach to the integration of geodesic equations originally pursued in
\cite{FreSorin2006,Fre:2009zz,dualborel} and based on Kodama's integration algorithm \cite{kodama1}, and the approach based on Noether charges,  which
directly provides the final integral. Indeed  the former approach was based on putting the geodesic equation for the geodesic velocity in the form
of a Lax equation. After integration of Lax equation a second integration was necessary, but already reduced to quadratures, in order to obtain
the final form of the functions $\Upsilon_A(t)$. Instead in  the Noether approach we directly get $\Upsilon_A(t)$, as we show below.
\par  Given the
initial data
\begin{equation}
\Upsilon^0\equiv (\Upsilon^0_A)\,\,,\,\,\,Q\in T_{\Upsilon^0}(\mathrm{U/H})\,,
\end{equation}
the solution to (\ref{geoM}) is given by the functions $\Upsilon_A(t)$ satisfying the   matrix equation
\begin{equation}
\mathcal{M}(\Upsilon(t;\,\Upsilon^0))=\mathcal{M}(\Upsilon^0)\,\exp\left[Q^T\,t\right]\,.\label{geoMsol}
\end{equation}
We can change the initial point from $\Upsilon^0$ to $\Upsilon^{0\prime}$ by acting on $\mathcal{M}$ by means of the unique element
$\mathfrak{g}\in \mathcal{S}_{\mathrm{U/H}} \subset \mathrm{U}$ that maps $\Upsilon^0$ to $\Upsilon^{0\prime}$ (simple transitivity), and solving
the matrix equation
\begin{equation}\label{gerbone1}
\mathcal{M}(\Upsilon(t;\,\Upsilon^{0\prime}))=\mathfrak{g}\,\mathcal{M}(\Upsilon(t;\,\Upsilon^0))\,\mathfrak{g}^T\,.
\end{equation}
The initial velocity $Q$ (Noether charge) transforms consequently:
\begin{equation}\label{kusnetzov}
 Q\rightarrow Q'=\mathfrak{g}\,Q\,\mathfrak{g}^{-1}\,\in \,  T_{\Upsilon^{0\prime}}(\mathrm{U/H})
 \end{equation}
We can choose $\Upsilon^0=0$. In this case we denote by $\Upsilon(t)\equiv v(t;\,\mathbf{e})$, and equation (\ref{geoMsol}) reads:
\begin{equation}
\mathcal{M}(\Upsilon(t))=\exp\left[{Q^T\,t}\right]=\exp\left[{Q\,t}\right]
\end{equation}
since $Q=Q^T\in \mathbb{K}$ and $\mathcal{M}(\mathbf{e})={\bf 1}$. The exponential on the right-hand side is readily computed  diagonalizing $Q$
by means of a constant orthogonal matrix $\mathcal{O}$:
\begin{equation}
Q=\mathcal{O}\,Q_D\,\mathcal{O}^T\,\,,\,\,\,Q_D={\rm diag}(\lambda_1,\dots, \lambda_N)
\end{equation}
In general, through a transformation $\mathcal{O}\in \mathrm{H}$, $Q$ can be always rotated into the Cartan subalgebra of the coset
$\mathrm{U/H}$.
Thus the number of independent eigenvalues $\lambda$ of a generic $Q$ equals the rank of the coset. The geodesic matrix equation
then becomes
\begin{equation}\label{ciurlaccone}
\mathcal{M}(\Upsilon(t))=\mathcal{O}\,\exp\left[{Q_D\,t}\right]\,\mathcal{O}^T \, = \,
\exp\left[\underbrace{\mathcal{O}\cdot Q_D\cdot \mathcal{O}^T}_{-2\, L_0} \, t\right]
\end{equation}
Indeed what appears in brackets in \eqref{ciurlaccone}  is just proportional by a factor $-2$ to the initial Lax operator $L_0$, which
constitutes the initial condition for the time evolving Lax $L(t)$.  The precise correspondence between the two algorithms is provided by
the identification
\begin{equation}\label{curruespundo2}
  Q \, = \, - \, 2 \, L_0
\end{equation}
As one sees, we have a complete integration  algorithm for the geodesics,
provided that we are able to extract the upper triangular coset representative
from the evolving matrix $\mathcal{M}(t)$ obtained by means of
a simple matrix exponential, that  is, if  we are able to get  $\mathbb{L}(t) \, = \,
\mathbb{L}\left(\Upsilon(t)\right)$ such that:
\begin{equation}\label{ciulifischio}
\mathcal{M}(t)=\mathbb{L}(t)\mathbb{L}(t)^T.
\end{equation}
This goal is achieved through the so-called Cholesky decomposition
(see for instance \cite{Horn-Johnson}).
 In particular one can apply the Cholesky-Banachiewicz and Cholesky-Crout
algorithms. Here are the iterative formulas:
\begin{align}\label{ciulacavolo}
\mathbb{L}_{jj}&=\sqrt{\mathcal{M}_{jj}-\sum_{k=1}^{j-1} \mathbb{L}_{jk}^2}\nonumber\\
\mathbb{L}_{ij}&=\frac{1}{\mathbb{L}_{jj}}\,\left(\mathcal{M}_{ij}-\sum_{k=1}^{j-1}\,\mathbb{L}_{ik}\mathbb{L}_{jk}\right)\,\,,\,\,\,i>j\,,
\end{align}
where $\mathcal{M}=\mathcal{M}(t)$.
%%%%%%%%%%%%%%%%%%%%%%%%%%%%%%%%%%%%%%
%%%%%%%%%%%%%%%%%%%%%%%%%%%
\subsection{The Geodesic Distance} \label{distante} As it  has been a few    times emphasized above, the main reason while the symmetric spaces under
investigation prove relevant for DeepL modeling  is that, similarly to flat space, they admit a global and canonical definition of a distance, the
geodesic one. This is also the prerequisite to study geodesic grids and diffusion equations. In this section we are going to elaborate this aspect in full,
relying on the results of section \ref{lavatoconperlana}. It follows that on these negatively curved manifolds, there is   a unique geodesic
 joinining any two points. We begin by recalling eq.~ (\ref{geoMsol}), which implies that given two points of the
manifold $p_1,p_2$ whose solvable coordinates we collectively denote $\Upsilon^{1,2}$, there is a unique geodesic that joins them defined as the
following curve:
\begin{eqnarray}\label{orgiadiM}
  \mathcal{M}\left(\Upsilon(t)\right) & =& \mathcal{M}\left(\Upsilon^1\right) \, \exp\left[Q^T\, t\right] \nonumber\\
  Q^T & =& \log\left[ \mathcal{M}^{-1}\left(\Upsilon^1\right) \,\mathcal{M}\left(\Upsilon^2\right) \right] \,.
\end{eqnarray}
In particular,   redefining the \emph{time} variable $t$, \emph{i.e.} the affine parameter along the geodesic, we can normalize so that
\begin{equation}\label{cordino}
  \Upsilon^2 \, = \, \Upsilon(1) \quad ; \quad  \Upsilon^1 \, = \, \Upsilon(0).
\end{equation}
As it was observed in section \ref{lavatoconperlana}, the matrix of the Noether charges $Q^T$  is not symmetric, unless the origin point $p_1$  of
the geodesic going toward $p_2$ is  $p_1 \, = \, \mathbf{e}$,   the \emph{neutral element} of the solvable group:
\begin{equation}\label{solvgroupN}
\mathbf{e} \, \in \,   \mathcal{S}_{\mathrm{U/H}} \equiv \exp\left[\mathbb{B}\mathrm{(N,\mathbb{R})}\right] \quad ; \quad \mathbf{e} \,
\equiv \, \exp[\mathbf{0}]
\end{equation}
In solvable coordinates the  \emph{neutral element} is identified by $\Upsilon_A \, = \, 0$ ($A=1,\dots,\dim {\mathrm{U/H}}$). The
invariance under the solvable group $\mathcal{S}_{\mathrm{U/H}}$ allows however to reduce the calculation of any geodesic to the case where the
initial point $p_1 \, = \, \mathbf{e}$ is the origin.  Having clarified the structure of the solvable group, recalling eq.~ (\ref{gerbone1}) and
coming back to eq.~(\ref{orgiadiM}) we define:
\begin{equation}\label{trenuncolo}
\mathcal{S}_{\mathrm{U/H}} \, \ni \, \mathfrak{g}\, = \, \mathbb{L}^{-1}\left(\Upsilon^1\right) \quad ;
\quad \mathbb{L}\left(\tilde{\Upsilon}(t,0)\right) \, = \, \mathfrak{g} \, \mathbb{L}\left(\Upsilon(t,\Upsilon^1)\right)
\end{equation}
so that the trajectory $\tilde{\Upsilon}(t,0)$ satisfies the following conditions:
\begin{equation}\label{conditiones}
 \tilde{\Upsilon}(0,0) \, = \, 0 \, = \, \mathfrak{g}(\Upsilon^1) \quad ; \quad \tilde{\Upsilon}(1,0) \, = \,
 \tilde{\Upsilon}^2 \, \equiv \, \mathfrak{g} \left(\Upsilon^2\right)
\end{equation}
In the above equations (\ref{trenuncolo}) the coset representative $\mathbb{L}^{-1}\left(\Upsilon^1\right)$ is looked at as a realization of the
unique solvable group element $\mathfrak{g}$ which maps the point $p_1$, labeled by coordinates $\Upsilon^1$, into the neutral element
$\mathbf{e}$, labeled by coordinates $\Upsilon =0$. The image by $\mathfrak{g}$ of the second point $p_2$ has, by definition, the new
coordinates $\tilde{\Upsilon}^2$.  In this way we obtain:
\begin{equation}\label{crisalide}
 \mathcal{M}\left(\tilde{\Upsilon}(t,0)\right) \, = \, \exp \left[t \, Q_0\right] \, = \, \mathfrak{g} \, \exp \left[t \, Q^T\right] \,
 \mathfrak{g}^{-1} \quad \Leftrightarrow \quad  Q_0 = \mathfrak{g} \, Q^T \, \mathfrak{g}^{-1}
\end{equation}
and the Noether charge $Q_0$ shifted to the origin is now a symmetric matrix, whatever the points $p_{1,2}$ were,  and such is
$\mathcal{M}\left(\tilde{\Upsilon}(t,0)\right)$.
\par
Let us now come back to equation (\ref{perbacco}), which yields the following compact expression for the line element in the symmetric space:
\begin{equation}\label{fagiolino}
  ds^2 \, = \, Tr\left(\mathcal{M}^{-1}(\Upsilon)d\mathcal{M}(\Upsilon) \, \cdot \,\mathcal{M}^{-1}(\Upsilon)d\mathcal{M}(\Upsilon)\right );
\end{equation}
restricting the line element to the geodesic we have:
\begin{align}\label{borlotto}
  ds^2\mid_{geod.} & = \, Tr\left(\mathcal{M}^{-1}(\Upsilon(t)d\mathcal{M}(\Upsilon(t)) \, \cdot \,\mathcal{M}^{-1}(\Upsilon(t))d\mathcal{M}(\Upsilon(t))\right ) \, = \, Tr\left(Q^T \cdot Q^T\right) \, dt^2 \nonumber\\[3pt]
  ds & = \, \sqrt{Tr\left(Q_0 \cdot Q_0\right)} \, dt
\end{align}
In the second line of equation (\ref{borlotto}) we have used the properties of the trace to rewrite $Tr\left(Q^T \cdot Q^T\right) $ in terms of
the symmetric matrix $Q_0$. We can do even better, and use the \emph{Schur Decomposition}, which  exists for any non degenerate symmetric matrix:
\begin{eqnarray}\label{ladivago}
 Q_0 &= & - 2 \,\mathcal{O} \, L_{\mathrm{diag}} \mathcal{O}^T \quad \text{where}\quad \mathrm{O} \in \mathrm{SO(N)} \quad \text{i.e.} \;\;
 \mathcal{O}\cdot \mathcal{O}^T \, = \, \mathbf{1} \nonumber\\
  L_{\mathrm{diag}} & = &\left(
                           \begin{array}{ccccccc}
                             -\sum_{i=2}^{N}\mu_{i}& 0 & 0 & \dots & \dots & \dots & 0 \\
                             0 & \mu_2 & 0 & 0 & \dots & \dots & 0 \\
                             0 & 0 & \mu_2 & 0 & \dots & \dots & 0 \\
                            \dots & \dots & \dots & \dots & \dots & \dots & \dots \\
                            \dots & \dots & \dots & \dots & \dots & \dots & \dots \\
                             0 & \dots & \dots & \dots & 0 & \mu_{N-1} & 0 \\
                             0 & \dots  & \dots & \dots & \dots  & 0 & \mu_N \\
                           \end{array}
                         \right)
\end{eqnarray}
having denoted by $\mu_{2,\dots,N}$ the $N-1$ independent eigenvalues of the traceless initial Lax operator $L_0$, where conventionally we set
\begin{equation}\label{convention}
  \mu_1 \, = \,  -\sum_{i=2}^{N}\mu_{i} .
\end{equation}
In this way the arc length of   the unique geodesic joining the point $p_1$ with the point $p_2$ is given by:
\begin{equation}\label{edessa}
  \ell\left(p_1,p_2\right) \, = \, \ell\left(\mathbf{e}, \tilde{p}_2\right) \, = \, 4 \,\sqrt{Tr\left(L_{\mathrm{diag}}\cdot L_{\mathrm{diag}} \right) }
  \, \int_0^1 dt \, = \, 4 \,\sqrt{\sum_{i=1}^N \mu_i^2}
\end{equation}
The final appearance of the geodesic distance is intriguing, as it looks like a Euclidean distance. However, it only looks like that, but it is quite
different, as  the $\mu_i$ are not the differences of the coordinates of the points $p_1,p_2$,  but they rather   are the eigenvalues of the Noether
charge matrix $Q_0$. The latter is uniquely identified by the choice of the two points $p_1,p_2$ yet the solvable group is non abelian and the
coordinates $\Upsilon^1 $, $\Upsilon^2 $  identify the coordinates of the \emph{difference point} $\tilde{\Upsilon}$ in the complicated
non-linear way encoded in the above described procedure.
\par
Furthermore, given $p_1,p_2$ we have the explicit form, not only of the eigenvalues, but also of the entire $Q_0$-matrix, and we can calculate the
full  geodesic by means of the integration method described in section \ref{lavatoconperlana} and the Cholevsky-Crout algorithm, to
obtain the time evolving upper triangular representative $\mathbb{L}(t) \, = \,\mathbb{L}\left(\tilde{\Upsilon}(t,0)\right)$. This being given, we
set
\begin{equation}\label{gigione}
  \mathbb{L}\left({\Upsilon}(t,\Upsilon^1)\right)\, = \, \mathfrak{ g}^{-1} \, \mathbb{L}(t)
\end{equation}
and using the inverse of the exponential map $\Sigma$ which is an iterative   algorithm,  we work out the functions
\begin{equation}\label{tuttageo}
  \Upsilon_A(t,\Upsilon^1)  \quad ; \quad \Upsilon_A(0,\Upsilon^1) \, = \Upsilon_A^1, \quad \Upsilon_A(1,\Upsilon^1) \, = \Upsilon_A^2
\end{equation}
%%%%%%%%%%%%%%%%%%%%%%%%%%%%
\subsubsection{Geodesic distance in view of the Tits-Satake projection}\label{lontanolontano} Let us now  reconsider the geodesic problem in view
of the Tits-Satake projection, focusing in particular on the manifolds
\begin{equation}\label{larrybar}
  \mathcal{M}^{[r,s]}_\mathfrak{d} \, \equiv \, \frac{\mathrm{SO(r,r+2s)}}{\mathrm{SO(r) \times SO(r+2s)}}.
\end{equation}
\par
The first step, as we did in eq.~~\eqref{turlipano}, is to introduce the symmetric matrix
\begin{equation}\label{turlipano2}
\mathcal{M}\equiv \mathbb{L}\,\mathbb{L}^T\,,
\end{equation}
and recall eq.~~\eqref{fagiolino}, which expresses the metric in terms of this   invariant object.
\par
Next we recall eq.~~\eqref{orgiadiM}, which expresses the full geodesic in terms of a constant symmetric matrix $Q^T$, called the matrix of Noether
charges, and eq.~(\ref{borlotto}) that yields the final answer for the distance in terms of $Q^T$:
\begin{equation}\label{tracciandosegni}
  d(\mathbf{u},\mathbf{v}) \, = \, \sqrt{Tr(Q^T\cdot Q^T)},
\end{equation}
where for the geodesic described by  the matrix equation (\ref{borlotto}), which  we rewrite in the current notation as
\begin{equation}
\mathcal{M}[\boldsymbol{\Upsilon}(\lambda)]=\mathcal{M}[\boldsymbol{\Upsilon}(0)]
\,e^{Q^T\,t}\,,
\end{equation}
we have denoted as follows the initial and final points:
\begin{equation}
\boldsymbol{\Upsilon}(t=0)={\bf v}\,\,,\,\,\,\,\boldsymbol{\Upsilon}(t=1)={\bf u}
\end{equation}
For consistency, if we choose $\mathbf{u},\mathbf{v}$ arbitrarily, the corresponding Noether charge matrix is defined as
\begin{equation}
Q^T[{\bf u},\,{\bf v}]\equiv\log\left(\mathcal{M}[{\bf u}]^{-1}\,\mathcal{M}[{\bf v}]\right)\,,\label{Quv}
\end{equation}
and the distance $d({\bf u},\,{\bf v})$ is calculated by (\ref{tracciandosegni})   using    the Noether charge $Q^T =Q^T[{\bf u},\,{\bf v}]$.
\subsection{A computational strategy}
The closed analytic formula for the derivatives of the distance is essential in order to construct analytic loss functions depending on the
weight parameters. On the other hand the calculation of the logarithm of a large   matrix is a difficult task, unless it is made  numerically.
However, further elaborations can lead to a simpler and doable task. To this effect let us consider the Cartan orthogonal decomposition in eq.~
(\ref{cartandecompo}) of the Lie algebra $\mathbb{U}$. Instead of the solvable parameterization, one might use the Cartan parameterization in
terms of coordinates $\boldsymbol{\xi}=\{\xi^A\}$ such that the coset representative is chosen to be:
\begin{equation}
\mathbb{L}_K[\boldsymbol{\xi}]=e^{K[\boldsymbol{\xi}]}\,\,,\,\,\,\,K[\boldsymbol{\xi}]\equiv {\xi^A\,\mathrm{K}_A}=K[\boldsymbol{\xi}]^T\,,
\end{equation}
where $\{\mathrm{K}_A\}$ denotes a basis of $\mathbb{K}$. Projecting the left-invariant 1-form $\Omega \, \equiv \,\mathbb{L}^{-1}d\mathbb{L} $ on
the generators $\{K_A\}$, which are always appropriately normalized (\emph{i.e.} $Tr(\mathrm{K}_A\, \mathrm{K}_B)\, =\, \delta_{AB}$), we always
obtain an equivalent definition of  the vielbein $V^A\, = \, Tr\left(\mathrm{K}_A \,\Omega\right)$, independently from the specific choice of coset
representative $\mathbb{L}$. Different choices of $\mathbb{L}$ simply amount to different choices of coordinates on the manifold. We consider and
compare two choices of coordinates, the solvable ones $\boldsymbol{\Upsilon}$ and the Cartan ones $\boldsymbol{\xi}$. The relation between the
two coordinate systems can be  written as
\begin{equation}\label{prigotov}
\mathbb{L}_K[\boldsymbol{\xi}(\boldsymbol{\Upsilon})]=\mathbb{L}_s[\boldsymbol{\Upsilon}]
\,h[\boldsymbol{\Upsilon}]\,\,,\,\,
\,h[\boldsymbol{\Upsilon}]
\in \mathrm{H}\,,
\end{equation}
where $h[\boldsymbol{\Upsilon}]$ is some suitable $\mathrm{H_{c}} $ compensator. We can avoid the involved calculation of
$h[\boldsymbol{\Upsilon}]$, considering that eq.~(\ref{prigotov}) is equivalent to the   matrix equation
\begin{equation}
\mathcal{M}_K[\boldsymbol{\xi}(\boldsymbol{\Upsilon})]\equiv \mathbb{L}_K[\boldsymbol{\xi}(\boldsymbol{\Upsilon})]\,\mathbb{L}_K[\boldsymbol{\xi}(\boldsymbol{\Upsilon})]^T=
\mathbb{L}_s[\boldsymbol{\Upsilon}]\,\mathbb{L}_s[\boldsymbol{\Upsilon}]^T=
\mathcal{M}_s[\boldsymbol{\Upsilon}]\label{MkMs}.
\end{equation}
Indeed the matrix $\mathcal{M}$ is a  well-defined function on the symmetric space and   is independent from the coset representative utilized in its
definition \eqref{turlipano2}\footnote{Actually the set of matrices $\mathcal{M}$  possessing the following two properties,   being
symmetric $\mathcal{M}^T \,=\,\mathcal{M}$,  and  being an element  of the group $\mathrm{U} \,= \,\mathrm{G_R}$, i.e.,
$\mathcal{M}\eta\mathcal{M}\, =\, \eta$, can be regarded as the intrinsic global definition of points of the symmetric space, as there is a
one-to-one correspondence between points of the manifold $\mathrm{U}/\mathrm{H_c}$ and these matrices.}. Let us then rewrite the geodesic distance
in a  different way, using the formal relation between the Cartan and the solvable coordinates. To this end we notice that, using for the two
points ${\bf u},\,{\bf v}$ the solvable parameterization in which they are respectively  described by the coordinates $\Upsilon^A_{{\bf u}}$ and
$\Upsilon^A_{{\bf v}}$, both $\mathbb{L}[{\bf u}]^{-1}=\mathbb{L}_s[\boldsymbol{\Upsilon}_{{\bf u}}]^{-1}$ and  $\mathbb{L}[{\bf
v}]=\mathbb{L}_s[\boldsymbol{\Upsilon}_{{\bf v}}]$ belong to the same solvable Lie group $\mathcal{S} \equiv e^{Solv}$. Their product is therefore
in the same solvable group:
\begin{equation}
\mathbb{L}_s[\boldsymbol{\Upsilon}({\bf u},\,{\bf v})]\equiv \mathbb{L}_s[\boldsymbol{\Upsilon}_{{\bf u}}]^{-1}\,\mathbb{L}_s[\boldsymbol{\Upsilon}_{{\bf v}}]\,\in \,\mathcal{S}\,\label{Lsuv}
\end{equation}
The metric equivalence of the symmetric space $\mathrm{U/H_c}$ with the corresponding solvable group $\mathcal{S}$ implies that every manifold
point uniquely identifies a corresponding element of the abstract solvable group, so that we can write
\begin{equation}\label{talismano}
  \mathbb{L}_s\left[\boldsymbol{\Upsilon}({\bf u},\,{\bf v})\right]\, = \, \mathbb{L}_s\left[\boldsymbol{\Upsilon}_{\mathbf{w}}\right]
\end{equation}
where we have denoted
\begin{equation}\label{sciamano}
  \mathbf{w} \, = \, \mathbf{u}^{-1} \cdot \mathbf{v}
\end{equation}
the abstract element of the solvable group $\mathcal{S}$ obtained as the product of the inverse of the abstract element $\mathbf{u}$ by the
abstract element $\mathbf{v}$.  Clearly, when ${\bf u}={\bf v}$, the abstract element $\mathbf{w}$ is just the identity and
$\mathbb{L}_s[\boldsymbol{\Upsilon}_{\mathbf{w}}]={\rm Id}$. Therefore the generator of the geodesic connecting the two points, given in eq.~~\eqref{Quv}, is expressed solely in terms of the matrix $\mathcal{M}\left(\boldsymbol{\Upsilon}_{\mathbf{w}}\right)$:
\begin{equation}\label{flipper}
Q[{\bf u},\,{\bf v}]\, = \, \log\left[\mathcal{M}\left(\boldsymbol{\Upsilon}_{\mathbf{w}}\right)\right] \quad ; \quad
\mathcal{M}\left(\boldsymbol{\Upsilon}_{\mathbf{w}}\right)\, = \,  \mathbb{L}_s[\boldsymbol{\Upsilon}_{\mathbf{w}}]\cdot
\mathbb{L}_s[\boldsymbol{\Upsilon}_{\mathbf{w}}]^T.
\end{equation}
This leads to the following expression for the distance function:
\begin{equation}\label{maseghepensu}
d^2({\bf u},{\bf v})=\,{\rm Tr}\left[\log\left[\mathcal{M}\left(\boldsymbol{\Upsilon}_{\mathbf{w}}\right)\right]^2\right]\,
\quad ; \quad \mathbf{w}\equiv\mathbf{u}^{-1}\cdot \mathbf{v}
\end{equation}
Let us now use eq.~~\eqref{MkMs}, which states that the matrix $\mathcal{M}(\mathbf{w})$ depends only on the manifold point (i.e., the abstract
element $\mathbf{w}\in \mathcal{S} $ of the solvable group) and not on the coordinates utilized to describe it. Hence in formula
(\ref{maseghepensu}) we can replace the solvable parameterization of $\mathcal{M}(\mathbf{w})$ with  the Cartan one:
\begin{equation}\label{perbacco}
\mathcal{M}_s(\mathbf{w}) =\mathcal{M}_K(\mathbf{w})
 =\mathbb{L}_K[\boldsymbol{\xi}_{\mathbf{w}})]^2=\exp\left(2\,K[\boldsymbol{\xi}_{\mathbf{w}})]\right).
\end{equation}
The second part of equation (\ref{perbacco}) follows from the fact that the coset generators $\mathrm{K}_A$ are symmetric, and therefore  so is
their exponential. Hence the geodesic distance can be written as
\begin{equation}
\mathrm{d}^2({\bf u},{\bf v})=\,{\rm Tr}\left[4\,\,K[\boldsymbol{\xi}_{\mathbf{w}}]^2\right]=2\,\xi^A_{{\bf w}}\, \xi^B_{{\bf w}} \,\delta_{AB}\, \quad ; \quad \mathbf{w}\equiv\mathbf{u}^{-1}\cdot \mathbf{v}. \label{nduv}
\end{equation}
\par
To sum up, the computational steps required by the above procedure are :
\begin{itemize}
\item evaluation of the general relation between the Cartan coordinates $\boldsymbol{\xi}$ and the solvable ones $\boldsymbol{\Upsilon}$:
    ${\xi}^A={\xi}^A(\boldsymbol{\Upsilon})$ from eq.~ (\ref{MkMs});
\item evaluation of the product law $\boldsymbol{\Upsilon}({\bf u},\,{\bf v})$ from eq.~ (\ref{Lsuv});
\item evaluation of the geodesic distance from eq.~ (\ref{nduv}).
\end{itemize}
\par
Next, in view of eq.~(\ref{nduv}), we need to establish the relation between the Cartan and solvable coordinate.
\paragraph{\sc Cartan to solvable coordinates.}
In the fundamental representation of ${\rm SO}(r,\,r+q)$ (where $q=2s$ or $q=2s+1$, depending on the $\mathfrak{d}_\ell$ or $\mathfrak{b}_\ell$
case), using the $\eta_b$-basis, which exists for both algebras, we can write
\begin{equation}
\label{cormaiore}
K_b[\boldsymbol{\xi}]=\left(\begin{array}{c|c}{\bf 0}_{(r+q)\times (r+q)} & \boldsymbol{\xi}\\
\hline
\boldsymbol{\xi}^T & {\bf 0}_{r\times r}\\
\end{array}\right)
\end{equation}
where now $\boldsymbol{\xi}$ denotes the $r\times (r+q)$ matrix  $\boldsymbol{\xi}=(\xi_i{}^m)$, $i=1,\dots, r,\,m=1,\dots,r+q $. The version
$K_t[\boldsymbol{\xi}]$ of the same matrix (\ref{cormaiore}) in the $\eta_t$ basis can be obtained by means of the inverse similarity
transformation with the $\Omega^{-1} = \Omega^T$ matrix.
\par
Equation (\ref{nduv}) can be written, in terms of the rectangular matrix $\boldsymbol{\xi}$, as follows:
\begin{equation}
d^2({\bf u},{\bf v})={\rm Tr}\left[\boldsymbol{\xi}_\mathbf{w}^T \,\boldsymbol{\xi}_\mathbf{w}\right] \quad ;
\quad \mathbf{w}\equiv\mathbf{u}^{-1}\cdot \mathbf{v}\label{nnduv}
\end{equation}
By matrix exponentiation the matrix $\mathcal{M}_K[\boldsymbol{\xi}]$ has the following explicit form:
\begin{equation}
\mathcal{M}_K[\boldsymbol{\xi}]=\exp\left(2K[\boldsymbol{\xi}]\right)=\left(\begin{matrix}
\cosh\left(2(\boldsymbol{\xi}\boldsymbol{\xi}^T)^{\frac{1}{2}}\right) & \sinh\left(2(\boldsymbol{\xi}\boldsymbol{\xi}^T)^{\frac{1}{2}}\right)
\,(\boldsymbol{\xi}\boldsymbol{\xi}^T)^{-\frac{1}{2}}\,\boldsymbol{\xi}\cr \boldsymbol{\xi}^T \sinh\left(2(\boldsymbol{\xi}\boldsymbol{\xi}^T)^{\frac{1}{2}}\right)\,(\boldsymbol{\xi}\boldsymbol{\xi}^T)^{-\frac{1}{2}}& \cosh\left(2(\boldsymbol{\xi}^T\boldsymbol{\xi})^{\frac{1}{2}}\right) \end{matrix}\right)\,.
\end{equation}
The above matrix can be written in the more compact form
\begin{equation}
\label{baltimora}
\mathcal{M}_K[\boldsymbol{\xi}]=\left(\begin{matrix}
\left({\bf 1}+{\bf X}{\bf X}^T\right)^{\frac{1}{2}} & {\bf X}\cr {\bf X}^T& \left({\bf 1}+{\bf X}^T{\bf X}\right)^{\frac{1}{2}} \end{matrix}\right)\,,
\end{equation}
where the rectangular matrix ${\bf X}={\bf X}(\boldsymbol{\xi})$ is related to $\boldsymbol{\xi}$ by the equation
\begin{equation}
\boldsymbol{\xi}=\frac{1}{2}\,{\rm arcsinh}\left[\left({\bf X}{\bf X}^T\right)^{\frac{1}{2}}\right]\,\left({\bf X}{\bf X}^T\right)^{-\frac{1}{2}}\,{\bf X}\,,
\end{equation}
and ${\bf X}$, as a function of the solvable coordinates $\Upsilon^A$, is read off $\mathcal{M}_K=\mathcal{M}_s$.
 we can also write:
\begin{equation}
\boldsymbol{\xi}^T\boldsymbol{\xi}=\frac{1}{4}\,{\rm arccosh}\left[\left({\bf 1}+{\bf X}^T{\bf X}\right)^{\frac{1}{2}} \right]\,.
\end{equation}
There is no need to find the transformation law ${\xi}^A={\xi}^A(\boldsymbol{\Upsilon})$, which is a laboriuous and difficult calculation. It
suffices to note that the formula for the distance (\ref{nnduv}) can then be rewritten in the   form
\begin{equation}
\mathrm{d}^2({\bf u},{\bf v})=\frac{1}{4}\,{\rm Tr}\left[{\rm arccosh}^2\left({\bf 1}+{\bf X}_\mathbf{w}^T{\bf X}_\mathbf{w}\right)^{\frac{1}{2}}\right], \qquad \mathbf{w}\equiv\mathbf{u}^{-1}\cdot \mathbf{v},
\end{equation}
and going one step back to note that the matrix $\left({\bf 1}+{\bf X}_\mathbf{w}^T{\bf X}_\mathbf{w}\right)^{\frac{1}{2}}$, according with
eq.~(\ref{perbacco}) is  identical to the  right lower blower block of the matrix $\mathcal{M}_s(\boldsymbol{\Upsilon}_\mathbf{w})$ in
the $\eta_b$ basis.
\subsection{Summary for  the distance formula} \label{distanzabartolo} Summarizing what we have learned from the above
discussion, we can say that, given the two points $\mathbf{u}$ and $\mathbf{v}$ described in solvable coordinates, the essential calculation is
provided by the formula for the \emph{product law}.  Once we have $\boldsymbol{\Upsilon}[\mathbf{u},\mathbf{b}]$ we construct the corresponding
coset representative $\mathbb{L}(\boldsymbol{\Upsilon}[\mathbf{u},\mathbf{b}])$ and the corresponding $\mathcal{M} = \mathbb{L}\mathbb{L}^T$
matrix. We transform it to the $\eta_b$-basis and we obtain a matrix with the block structure:
\begin{equation}\label{cagnesco}
  \mathcal{M}_{b}^{solv} \, =\,\left(\begin{array}{c|c}
                                 * & * \\
                                 \hline
                                 * & \mathfrak{m} \left( \boldsymbol{\Upsilon}_\mathbf{w}\right)
                               \end{array}\right).
\end{equation}
Then the distance between the two points is obtained by setting
\begin{equation}\label{coriandolidipasta}
\mathrm{d}^2({\bf u},{\bf v})= \mathrm{N^2}\left(\mathbf{w}\right) \, \equiv \, \frac{1}{4}\,{\rm Tr}\left[{\rm arccosh}^2\left(\mathfrak{m}
\left( \boldsymbol{\Upsilon}_\mathbf{w}\right)\right)\right] \quad ; \quad \mathbf{w}\equiv\mathbf{u}^{-1}\cdot \mathbf{v}.
\end{equation}
Note that the size of the matrix $\mathfrak{m}$ is just $r\times r$ which is the noncompact rank, typically a relatively small number, so
that the calculation of the function $\cosh^1$, which requires the diagonalization of the matrix, can be done using only  an $\mathrm{SO(r)}$
transformation, and   is typically doable also symbolically without too much effort.
\par
Furthermore from a conceptual point of view it is important to stress that the distance formula is actually a consequence of two facts:
\begin{enumerate}
  \item the manifold is a group (solvable in this case) $\mathcal{S}$ ,so that to any pair of points $\mathbf{u}$ and $\mathbf{v}$ we can
      associate a third one $\mathbf{w}\equiv \mathbf{u}^{-1}\cdot \mathbf{v}$
  \item on the solvable $\mathcal{S}$ group there exists a map
  \begin{equation}\label{normsq}
    \mathrm{N^2} \colon  \mathcal{S}  \to  \mathbb{R}_+
  \end{equation}
  which is invariant under     the full group of automorphisms of the solvable group, that is,  the group $\mathrm{U}$ of isometries of the
  manifold $\mathrm{U/H_c}$. The value $\mathrm{N^2}(\mathbf{w})$  for any element  $\mathbf{w} \in \mathcal{S}$ of the solvable group can be
  called its norm: $\|\mathbf{w}\|^2 = \mathrm{N^2}(\mathbf{w})$ and the distance of two points $\mathbf{u}$ and $\mathbf{v}$ is just
  $\|\mathbf{u}^{-1}\cdot \mathbf{v}\|^2$, which formally is the same definition of distance as in flat $\mathbb{R}^N$ where the group is not
  only solvable but also abelian, and for the group product one uses the additive notation.
\end{enumerate}
All learning algorithms are based on the extremalization of functionals constructed in terms of such invariant functions. If there exist other
positive functions on the solvable group that are invariant with respect to the full group $\mathrm{U}$ or important subgroups, the
learning algorithms can be refined.
\section{paint group structure of the square norm and of the distance function}
\label{normapitturata} After the general discussion of the previous section \ref{lontanolontano}, we illustrate the explicit form of the square
norm function on the solvable group $\mathcal{S}_{\mathrm{U/H}}$ in two relevant cases $r=1$ and $r=2$, emphasizing their paint group or subpaint
Group covariant structure. The square norm (from which the distance function follows) is the main ingredient in all Deep Learning algorithms; the
newly discovered analytic structure of this ingredient,  in view of the Tits-Satake projection and of the Paint and subpaint groups,  is what opens
  attractive perspectives of applications in the context of Data Science of the PGTS theory of hyperbolic symmetric spaces. Actually the
\emph{Poincar\'e Balls} introduced  in \cite{francesi1,francesi2,francesi3}  are just the tip on iceberg whose submersed base enlarges by
increasing the noncompact rank $r$. Furthermore, even in the case of the Poincar\'e balls, the richer structure due to the Tits-Satake projection
and subpaint spheres was not noticed   by the authors of \cite{francesi1,francesi2,francesi3}.
\subsection{The general distance formula applied to the $r=1$ case}\label{r1lontananza}
In order to calculate the geodesic distance between points of the full manifold $\frac{\mathrm{SO(1,1+2s)}}{\mathrm{SO(1+2s)}}$ we have to recall
eq.~(\ref{Lr1}), displaying the $r=1$ form of the solvable group element, that is,  the explicit form of the $\Sigma$-map, as defined in
eq.~(\ref{sigmaexp}-\ref{SigmaExp1}),  and calculate its inverse:
\begin{equation}\label{sigmamenouno}
  \Sigma^{-1} \, : \quad \mathcal{S}_{\mathrm{U/H}} \, \longrightarrow \, Solv_{\mathrm{U/H}}
\end{equation}
In this case the inversion is very simple and we have
\begin{equation}\label{pirlone}
\Sigma^{-1}\, : \quad \boldsymbol{\Upsilon}_{1} \, = \, \log\left[\mathbb{L}_{1,1} \right]  \quad ; \quad  \boldsymbol{\Upsilon}_{2,i} \, = \, - \,\sqrt{2} \, \mathbb{L}_{1 + i, 2 + 2s} \quad i=1,\dots,2\, s
\end{equation}
Furthermore we have also to recall the explicit product law between two elements in a  given solvable group. Llet us call
$\boldsymbol{U}$ and $\boldsymbol{V}$ the solvable coordinates of two points $\mathbf{u}$ and $\mathbf{v}$: in the $r=1$ case they have the
following form:
      \begin{eqnarray}\label{pointuv1}
        \boldsymbol{U} &=&\left\{u_1, {u}_{2,i}\right\}\nonumber\\
 \boldsymbol{V} &=&\left\{v_1, {v}_{2,i}\right\}
      \end{eqnarray}
  According with eq.~ (\ref{Lsuv}) the next step consists of constructing the product:
 \begin{equation}\label{perbaccus1}
   \boldsymbol{\mathcal{S}}_{[1,s]}\, \ni\, \mathbb{L}^{-1}\left(\mathbf{U}\right)\,\mathbb{L}\left(\mathbf{V}\right) \, \equiv \,
   \mathbb{L}\left(\mathbf{W}\right)
 \end{equation}
 The invertibility of the map $\Sigma$ described abstractly and in general in eqs.~ (\ref{sigmaexp}-\ref{SigmaExp1}) is made explicit
 for the case $r=1$  in eq.~(\ref{pirlone}): it allows  to work out the coordinate vector $\boldsymbol{W}$ from the coordinate vectors
 $\boldsymbol{U},\boldsymbol{V}$. We have:
 \begin{equation}\label{prunello1}
  \begin{array}{rcl}
 w_1& = & v_1\, -\, u_1 \\
 w_{2,i}& = & v_{2,i}\, -\, \exp\left[u_1 \, - \, v_1\right]\, u_{2,i} \quad ; \quad i=1\, \dots, 2s\\
\end{array}
\end{equation}
According with eq.~(\ref{baltimora}) the next step is that of constructing the fundamental harmonic symmetric matrix:
\begin{equation}\label{fiordavalle1}
\mathcal{M}_t(\mathbf{W}) \, \equiv \,\mathbb{ L}(\mathbf{W}) \, \mathbb{L}(\mathbf{W})^T
\end{equation}
and transform it to the block diagonal basis utilizing the case $r=1$ of the matrix $\Omega$ defined for all $r,s$  in eq.~(\ref{omegagrande}):
\begin{eqnarray}\label{stelovalle1}
\mathcal{M}_b(\mathbf{W}) & = & \Omega \, \mathcal{M}_t(\mathbf{W}) \, \Omega^T  \nonumber\\
&=& \hat{\mathbb{ L}}(\mathbf{W}) \, \hat{\mathbb{ L}}(\mathbf{W})^T
\quad \text{where} \quad \hat{\mathbb{ L}}(\mathbf{W}) \, \equiv \, \Omega \, \mathbb{ L}(\mathbf{W}) \, \Omega^T
\end{eqnarray}
According with eq.~(\ref{cagnesco}) we can express the distance in terms of the eigenvalues of the $1\times 1$ matrix $\mathfrak{m}$ which can be
described as:
\begin{equation}\label{matrucolata1}
\mathfrak{m} \, = \, \mathcal{L}\cdot \mathcal{L}
\end{equation}
where we have named $\mathcal{L}$ the last  row  of the matrix $\hat{\mathbb{ L}}(\mathbf{W}) $:
\begin{equation}\label{nendartallo1}
\mathcal{L}_I \, = \, \hat{\mathbb{ L}}(\mathbf{W})_{2s+2,I}  \quad ; \quad I = 1,\, \dots,\, 2+2s
\end{equation}
The very much useful point is that the vector $\mathcal{L}$ can be explicitly extracted from the construction of the base-transformed element
of the solvable group and written in a paint group covariant form for a generic value of the parameter $s$. Its explicit expression is displayed
below and it is manifestly invariant since it depends only from the solvable coordinate $w_1$ that is a $\mathrm{G_{paint}}$ singlet and from
invariant/covariant combinations of the $\mathrm{G_{paint}}$-vector $w_{2,i}$:
\begin{eqnarray}\label{caprino}
 \mathcal{L}(\mathbf{w})& = &  \left\{\frac{1}{2} e^{w_1} \underbrace{w_{2,i}}_{i=1,\dots ,2s},-\frac{1}{8} e^{w_1}
   \left(\|\mathbf{w}_2\|_{\mathrm{Gpaint}}^2-4\right)-\frac{1}{2}\,e^{-w_1},\frac{1}{8} e^{w_1}
   \left(\|\mathbf{w}_2\|_{\mathrm{Gpaint}}^2+4\right)+\frac{1}{2}\,e^{-w_1}\right\}\nonumber \\
   \|\mathbf{w}_2\|_{\mathrm{Gpaint}}^2&=& \sum_{i}^{2s} w_{2,i} \, w_{2,i}
\end{eqnarray}
Finally we can express the distance formula of the two point $\mathbf{U},\mathbf{V}$ as:
\begin{eqnarray}\label{creonte1}
   d^2(\mathbf{u},\mathbf{v}) & = & \mathrm{N^2}(\mathbf{w}) \, = \,
    \left({\rm arccosh}[ \mathcal{L(\mathbf{w})}\cdot \mathcal{L(\mathbf{w})} ]\right)^2\,
   \quad ; \quad \mathbf{w} = \mathbf{u}^{-1} \cdot \mathbf{v}
\end{eqnarray}
Explicitly we get:
\begin{equation}\label{bistolfo}
{\mathrm{N}^2}(\mathbf{w}) \, = \,  \left( {\rm arccosh}\left[\frac{1}{32} \left(16 e^{-2 w_1}+e^{2 w_1}+8 \|\mathbf{w}_2\|_{Gpaint}^2+
   e^{2 w_1} \, \left(\|\mathbf{w}_2\|_{Gpaint}^2+4\right){}^2\right)\right]\right)^2
\end{equation}
%%%%%%%%%%%%%%%%%%%%%%%%%%%%%%%%%%%%%%%%%%%%%%%%%%%%%
\subsubsection{Paint invariant distance formula $r=1$}
Referring to eq.~(\ref{pointuv1}) we find it convenient to rename
\begin{eqnarray}\label{pointuvNew}
\boldsymbol{U} &=&\left\{\lambda, {u}_{i}\right\},\nonumber\\
\boldsymbol{V} &=&\left\{\mu, {v}_{i}\right\}.
\end{eqnarray}
Applying the procedure described in the previous section up to eq.~(\ref{creonte1}), we easily arrive at the following result:
\begin{eqnarray}
\label{paintinvadista1}
   d^2(\mathbf{U},\mathbf{V})&=& \left({\rm arccosh}[ \mathfrak{A}_{U,V} ]\right)^2,
\end{eqnarray}
where, defining
\begin{equation}\label{maggiorana}
\mathbb{X}_i \, = \, e^{\frac{\lambda -\mu }{2}}\, u_i \quad ; \quad \mathbb{Y}_i \, = \, e^{\frac{-\lambda +\mu }{2}}\, v_i
\end{equation}
we obtain
\begin{equation}\label{cumino}
\mathfrak{A}_{U,V} \, =\,   \cosh \left[2 (\lambda -\mu )\right]+\frac{1}{2} \cosh \left[\lambda -\mu\right] \|
   \mathbb{X}-\mathbb{Y}\|^2+\frac{1}{32} \|\mathbb{X}-\mathbb{Y}\| ^4
\end{equation}
\subsubsection{The squared norm and Paint invariance for $r=1$}\label{sparko}
The square-norm ${\mathrm{N}^2}(\mathbf{w})$ has a manifest $\mathrm{G_{paint}}$ invariant structure as it depends only on the two
$\mathrm{G_{paint}}$ invariants $w_1$ and $\|\mathbf{w}_2\|_{Gpaint}^2$. This is very interesting; it means that if we consider the sublocus
$\mathfrak{B}_{\mu,R}\subset \mathcal{M}^{[1,s]}$ of the full manifold characterized as
\begin{equation}\label{sottoluogoP}
  \mathfrak{B}_{\mu,R} \, = \, \left\{\mathbf{w} \in \mathcal{M}^{[1,s]} \,  | \, w_1=\mu, \, \, \|\mathbf{w}_2\|_{Gpaint}^2 \, = \, R^2 \right\},
\end{equation}
then  all points in $\mathfrak{B}_{\mu,R}$ have the same norm, and really they make up something similar to a ball since the
 condition:
 \begin{equation}\label{sferella}
   \|\mathbf{w}_2\|_{Gpaint}^2 \, = \, R^2
 \end{equation}
 defines an $\mathbb{S}^{2s-1}$ sphere. Next consider that if we use not only the Paint group, but  rather
 the full maximal compact subgroup $H_c$, we can always rotate any element of the solvable Lie algebra $Solv(\mathrm{U/H})$ into the Cartan
 subalgebra, which means that we can rotate $\mathbf{w}$ into a $\tilde{\mathbf{w}}$ such that $\tilde{w}_1 \, = \, \lambda \neq \mu$
 and $\|\tilde{\mathbf{w}}_2\|_{Gpaint}^2 =0$. Obviously the $H_c$ rotations are a symmetry of the squared norm, so that we have
 \begin{equation}\label{rabaton}
   {\mathrm{N}^2}(\mathbf{w}) \, = \, {\mathrm{N}^2}(\tilde{\mathbf{w}}) \, = \, \left( {\rm arccosh} \left[ \cosh( 2\, \lambda )\right]\right)^2
   = 4 \lambda^2.
 \end{equation}
This means that in a shell of fixed squared norm ${\mathrm{N}^2}(\mathbf{w})=4\lambda^2$, the   Paint  group balls $\|\mathbf{w}_2\|_{Gpaint}^2 \, =
\, R^2$ have a maximum radius that can be calculated by eq.~ (\ref{bistolfo}).
\subsubsection{The squared norm and subpaint invariance for $r=1$}\label{okraps}
Next let us analyze the squared norm formula from the point of view of the $\mathrm{G_{subpaint}}$ group. That is easily done by splitting the
$G_{\mathrm{paint}}$  vector into its first component, which is a \emph{subpaint singlet}, plus the following $2s-1$ components that form
instead a $\mathrm{G_{subpaint}}$ multiplet:
\begin{eqnarray}\label{fornarello}
 w_{2,1} &=&w_2 \nonumber\\
 w_{2,1+\alpha} & = & \hat{w}_{\alpha} \quad ; \quad \alpha = 1,\dots, \, 2s-1 \nonumber \\
   \|\mathbf{w}_2\|_{\mathrm{Gpaint}}^2 &= & w_2^2 + \|\hat{\mathbf{w}}_2\|_{\mathrm{subpaint}}^2.
\end{eqnarray}
Inserting equation (\ref{fornarello}) into eq.~ (\ref{bistolfo}) we obtain
\begin{eqnarray}\label{panellenica}
{\mathrm{N}^2}(\mathbf{w}) & = &  \left( {\rm arccosh}\left[ \mathcal{A}_{TS} \left(\mathbf{w}_{TS}\right) \, + \,
\Delta\left(\mathbf{w}_{TS},\|\hat{\mathbf{w}}\|^2\right)\right]\right)^2\nonumber\\
\mathcal{A}_{TS} \left(\mathbf{w}_{TS}\right) &=& \frac{1}{32} \left(16 e^{-2 w_1}+8 w_2^2+e^{2 w_1} \left(w_2^2+4\right){}^2\right)\nonumber\\
\Delta\left(\mathbf{w}_{TS},\|\hat{\mathbf{w}}\|^2\right)&=& \frac{1}{32}\|\hat{\mathbf{w}}\|^2 \left(e^{2 w_1} \left(\|\hat{\mathbf{w}}\|^2+8\right)
+2 e^{2 w_1} w_2^2+8\right),
\end{eqnarray}
 where $\mathbf{w}_{TS}=\{w_1,w_2,0,\dots,0\} $ is the solvable coordinate vector of a point in the Tits-Satake submanifold (in this case the
 upper complex plane or the unit disk, as one prefers), while $\|\hat{\mathbf{w}}_2\|_{subpaint}^2 \, = \, R^2$ is a sphere $S^{2s-2}$ of points
 that project onto the same Tits-Satake point. Obviously if we fix the norm $\mathrm{N^2}(\mathbf{w})$ and  a Tits-Satake point
 $\mathbf{w}_{TS}$,
 whose norm must necessarily be smaller $\mathrm{N^2}(\mathbf{w}_{TS})\leq \mathrm{N^2}(\mathbf{w})$,
 then the radius $R$   of  $S^{2s-2}$ sphere  is fixed by equation eq.~ (\ref{panellenica}).
\subsection{The general distance formula applied to the $r=2$ case}\label{r2lontananza}
In the case $r=2$ the paint group covariant form utilizes  four solvable scalar coordinates $\upsilon_{1,2,3,4}$ (the coordinates on the two
hyperbolic planes forming the \emph{sub Tits-Satake submanifold}), and two solvable vector coordinates $\pmb{\Upsilon}_{5,6}$ each of which has
$2s$ components. Indeed, according to the general formula (\ref{kyrillov}),  the generic form of a solvable coordinate vector with $r(r+2s)$
components is the
\begin{equation}\label{kyrillov2}
  \boldsymbol{\Upsilon} \, =\,  \left\{ \underbrace{\upsilon_1, \upsilon_2}_{\text{n.c. Cartan}},
  \underbrace{ \upsilon_3, \upsilon_4}_{\text{long roots}}\,  \underbrace{\pmb{\Upsilon}_{5}}_{2s \, \,\text{preimages  of the short root 5}},
  \underbrace{\pmb{\Upsilon}_6}_{2s \, \,\text{preimages  of the short root 6}} \right\}
\end{equation}
Two points $p_u,p_v \in \mathcal{M}^{[2,s]}_\mathfrak{d}$ are specified by two solvable coordinate vectors  with the structure mentioned in
eq.~(\ref{kyrillov2}):
\begin{eqnarray}\label{pointuv}
\boldsymbol{U} &=&\left\{u_1, \,u_2,\, u_3, \, u_4, \,  \boldsymbol{U}_{5}, \, \boldsymbol{U}_{6} \right\}\nonumber\\
\boldsymbol{V} &=&\left\{v_1, \,v_2,\, v_3, \, v_4, \,  \boldsymbol{V}_{5}, \, \boldsymbol{V}_{6} \right\}
\end{eqnarray}
According to eq.~ (\ref{Lsuv}), the next step consists in constructing the product
\begin{equation}\label{perbaccus}
\boldsymbol{\mathcal{S}}_{[2,s]}\, \ni\, \mathbb{L}^{-1}\left(\mathbf{U}\right)\,\mathbb{L}\left(\mathbf{V}\right) \,
\equiv \, \mathbb{L}\left(\mathbf{W}\right)
\end{equation}
The invertibility of the map $\Sigma$ described abstractly and in general in eqs.~(\ref{sigmaexp},\ref{SigmaExp1}) and made explicit for the case
$r=2$  in Table \ref{activr2}, allows one  to work out the coordinate vector $\mathbf{W}$ from
\begin{equation}\label{succabrep}
   \boldsymbol{W} \, = \, \Sigma^{-1} \left( \mathbb{L}^{-1}\left(\mathbf{U}\right)\,\mathbb{L}\left(\mathbf{V}\right)\right).
\end{equation}
The explicit form in the case $r=2$ of the coordinate vector $\boldsymbol{W} $ in terms of the coordinate vectors
 $\boldsymbol{U},\boldsymbol{V}$ is
 \begin{equation}\label{prunello}
  \begin{array}{rcl}
 w_1& = & v_1\, -\, u_1 \\
 w_2& = & v_2\, -\, u_2 \\
 w_3& = & v_3\, -\, e^{u_1-u_2-v_1+v_2}\, u_3 \, \\
 w_4& = & v_4 \, + \, \frac{1}{\sqrt{2}}\,{e^{u_2-v_1-v_2} \left(e^{u_1}
   \, \mathbf{U}_{5} \cdot \mathbf{U}_{6}\, -e^{v_1} \, \mathbf{U}_{6}\cdot
   \mathbf{V}_{5}\right)}-\frac{1}{4} \, e^{2
   \left(u_2-v_2\right)} \, v_3 \, \mid \mathbf{U}_{6}\mid^2 \, \\
   \null && \quad \,\, +\,\frac{1}{4}
   e^{u_1+u_2-v_1-v_2} \left(u_3
   \, \mid \mathbf{U}_{6}\mid^2 \,-\,4\,
   u_4\right) \\
 \mathbf{W}_5& = & \mathbf{V}_{5}  \, + \,\frac{1}{\sqrt{2}} \, e^{u_2-v_2}\, v_3 \, \mathbf{U}_6
   \, -\, \frac{1}{2} e^{u_1-v_1}
   \left(\sqrt{2} \, u_3 \, \mathbf{U}_6\,+ \, 2\,
   \mathbf{U}_5\right)\\
 \mathbf{W}_6& = & \mathbf{V}_6\, -\, e^{u_2-v_2} \, \mathbf{U}_6  \\
\end{array}
\end{equation}
According to eq.~~\eqref{baltimora}, the next step is   to  construct the fundamental harmonic symmetric matrix
\begin{equation}\label{fiordavalle}
\mathcal{M}_t(\mathbf{W}) \, \equiv \,\mathbb{ L}(\mathbf{W}) \, \mathbb{L}(\mathbf{W})^T
\end{equation}
and transform it to the block diagonal basis utilizing the matrix $\Omega$ defined for all $r,s$  in eq.~(\ref{omegagrande}): :
\begin{eqnarray}\label{stelovalle}
   \mathcal{M}_b(\mathbf{W}) & = & \Omega \, \mathcal{M}_t(\mathbf{W}) \, \Omega^T  \nonumber\\
    &=& \hat{\mathbb{ L}}(\mathbf{W}) \, \hat{\mathbb{ L}}(\mathbf{W})^T
   \quad \text{where} \quad \hat{\mathbb{ L}}(\mathbf{W}) \, \equiv \, \Omega \, \mathbb{ L}(\mathbf{W}) \, \Omega^T
\end{eqnarray}
According to eq.~ (\ref{cagnesco}) we can express the distance in terms of the eigenvalues of the $2\times 2$ matrix $\mathfrak{m}$, which can be
described as
\begin{equation}\label{matrucolata}
    \mathfrak{m} \, = \, \left(
                           \begin{array}{c|c}
                             \mathcal{L}^{[1]}\cdot \mathcal{L}^{[1]} & \mathcal{L}^{[1]}\cdot \mathcal{L}^{[2]} \\
                             \hline
                             \mathcal{L}^{[1]}\cdot \mathcal{L}^{[2]} & \mathcal{L}^{[2]}\cdot \mathcal{L}^{[2]} \\
                           \end{array}
                         \right),
\end{equation}
where we   called $\mathcal{L}^{[1,2]}$ the last two rows of the matrix $\hat{\mathbb{ L}}(\mathbf{W}) $:
\begin{equation}\label{nendartallo}
   \mathcal{L}^{[1]}_i \, = \, \hat{\mathbb{ L}}(\mathbf{W})_{2s+3,i} \quad ; \quad  \mathcal{L}^{[2]}_i \, = \, \hat{\mathbb{ L}}(\mathbf{W})_{2s+4,i}.
\end{equation}
  The very   useful thing is that the vectors $\mathcal{L}^{[1,2]}$ can be explicitly extracted from the construction of the base-transformed
  element  of the solvable group and written in a Paint group covariant form for a generic value of the parameter $s$. Their explicit form is
  displayed in Table \ref{vettorini}. Finally, we can express the distance   between the two point $\mathbf{U},\mathbf{V}$ in terms of
  the eigenvalues of the matrix in eq.~(\ref{matrucolata}). We have
  \begin{equation}\label{gromiko}
       \Lambda_\pm \, =\, \mathcal{L}^{[1]}\cdot\mathcal{L}^{[1]}+\mathcal{L}^{[2]}\cdot\mathcal{L}^{[2]}
       \pm\sqrt{\left(\mathcal{L}^{[1]}\cdot\mathcal{L}^{[1]}-\mathcal{L}^{[2]}\cdot\mathcal{L}^{[2]}\right)^2
       +4\left(\mathcal{L}^{[1]}\cdot\mathcal{L}^{[2]}\right)^2}
      \end{equation}
 and the squared distance  is   \begin{equation}\label{creonte}
   d^2(\mathbf{U},\mathbf{V}) \, = \, \left({\rm arccosh}[\Lambda_+]\right)^2\, + \, \left({\rm arccosh}[\Lambda_-]\right)^2.
 \end{equation}

%%%%%%%%%%%%%%%%%%%%%%%%%%%%%%%%%%%%%%%
\begin{table}[!tb]
\begin{center}
$
\begin{array}{|l|c|l|}
\hline
\mathbb{ L}_{1,1} &=&e^{\Upsilon _1} \\ \hline
\mathbb{ L}_{1,2}  &=& \frac{e^{\Upsilon _1} \Upsilon _3}{\sqrt{2}} \\ \hline
\pmb{\mathbb{ L}}_{1,2+i} &=& \frac{1}{2} e^{\Upsilon _1} \left(\sqrt{2}
   \boldsymbol{\Upsilon} _{5,i}\, + \, \Upsilon _3 \,\boldsymbol{\Upsilon}
   _{6,i}\right) \\ \hline
\mathbb{ L}_{1,3+2s} &=&e^{\Upsilon _1} \left(-\,\frac{1}{2} \,\boldsymbol{\Upsilon}
   _{5}\cdot \boldsymbol{\Upsilon}_{6}\, -\, \frac{\Upsilon _3
   \mid\boldsymbol{\Upsilon} _{6}\mid^2}{4 \sqrt{2}}+\frac{\Upsilon
   _4}{\sqrt{2}}\right) \\ \hline
\mathbb{ L}_{1,4+2s} &=& -\frac{1}{2} e^{\Upsilon _1} \left(\frac{1}{2}
   \,\mid \boldsymbol{\Upsilon}
    _{5}\mid ^2 \, +\, \Upsilon _3 \Upsilon _4\right) \\ \hline
\mathbb{ L}_{2,2}&=& e^{\Upsilon _2}\\ \hline
\pmb{\mathbb{ L}}_{2,2+i} &=&\frac{e^{\Upsilon _2} }{\sqrt{2}}\boldsymbol{\Upsilon} _{6,i} \\ \hline
\mathbb{ L}_{2,3+2s} &=& -\frac{1}{4} e^{\Upsilon _2} \, \mid\boldsymbol{\Upsilon}
   _{6}\mid^2 \\ \hline
\mathbb{ L}_{2,4+2s} &=& -\frac{e^{\Upsilon _2} }{\sqrt{2}} \, \Upsilon _4 \\ \hline
\pmb{\mathbb{ L}}_{2+i,2+j}  &=& \delta_{i,j} \quad ; \quad (j\geq i)  \\ \hline
\mathbb{ L}_{3+2s,3+2s} &=& e^{-\Upsilon _2}  \\ \hline
\mathbb{ L}_{3+2s,4+2s}  &=& -\frac{e^{-\Upsilon _2} }{\sqrt{2}} \, \Upsilon _3 \\ \hline
\pmb{\mathbb{ L}}_{2+i,3+2s} &=& -\frac{1}{\sqrt{2}} \, \boldsymbol{\Upsilon }_{6,i}\\ \hline
\pmb{\mathbb{ L}}_{2+i,4+2s}  &=&  -\frac{1}{\sqrt{2}} \, \boldsymbol{\Upsilon }_{5,i} \\ \hline
 \mathbb{ L}_{4+2s,4+2s} &=& e^{-\Upsilon _1} \\ \hline
 \end{array}
 $
 \end{center}
 \caption{The explicit form in the case $r=2$, $s={generic}$ of the \emph{exponential map} defined
 in eqs.~(\ref{sigmaexp}-\ref{SigmaExp1}) from the Solvable Lie algebra $Solv$ to the Solvable Group Manifold $\mathcal{S}_{[2,s]} $.
 This map which is defined in a precise mathematical way and it is explicitly invertible is the substitute of the point-wise activation functions utilized in neural network constructions. The above table is written in a paint group covariant way: the boldface symbols are $2s$-component vector objects that transform in the fundamental defining representation of $\mathrm{SO(2s)}$.\label{activr2}}
\end{table}

%%%%%%%

\begin{table}[!tb]
\begin{center}
$
\begin{array}{|lcl|}
\hline
\pmb{\mathcal{L}}^{[1]} & = & \frac{1}{2\sqrt{2}} \, e^{w_1} \, \left(\sqrt{2} \, \mathbf{W}_5 \, + \, w_3 \, \mathbf{W}_6 \right)  \\
\mathcal{L}^{[1]}_{2s+1}  & = & -\frac{1}{2}\, e^{-w_1} \, - \, \frac{1}{8} \, e^{w_1} \, \left(-4\, + \, 2\, w_3 \, w_4 \,
+ \, \mid \mathbf{W}_5\mid^2 \right) \\
\mathcal{L}^{[1]}_{2s+2}& = & -\frac{1}{16} \, e^{w_1} \, \left( -4\sqrt{2} \, w_4 \, + \,  4 \, \mathbf{W}_5 \cdot \mathbf{W}_6
\, + \, \sqrt{2} \, w_3 \, \left(-4 \, + \, \mid \mathbf{W}_6\mid^2 \right)\right)\\
\mathcal{L}^{[1]}_{2s+3} & = & \frac{1}{8} \, e^{w_1}  \, \left( 4 \, + \, e^{2 \, w_1} \, \left(4\, + \, 2 \, w_3 \, w_4 \,+ \, \mid \mathbf{W}_5 \mid^2 \right) \right)\\
\mathcal{L}^{[1]}_{2s+4} & = & \frac{1}{16} \, e^{w_1} \,\left( -4\sqrt{2} \, w_4 \, + \,  4 \, \mathbf{W}_5 \cdot \mathbf{W}_6
\, + \, \sqrt{2} \, w_3 \, \left(4 \, + \, \mid \mathbf{W}_6\mid^2 \right)\right)\\
\hline
\pmb{\mathcal{L}}^{[2]} & = & \frac{1}{2} \, e^{w_2} \, \mathbf{W}_6\\
\mathcal{L}^{[2]}_{2s+1}  & = & \frac{1}{2\sqrt{2}} \, e^{-w_2}\, \left(w_3 \, - \, e^{2 \, w_2} \, w_4 \right) \\
\mathcal{L}^{[2]}_{2s+2}& = &  -\frac{1}{2}\, e^{-w_2} \, - \, \frac{1}{8} \, e^{w_2} \, \left(-4 \, + \, \mid \mathbf{W}_6 \mid^2 \right) \\
\mathcal{L}^{[2]}_{2s+3} & = & \frac{1}{2\sqrt{2}} \, e^{-w_2}\, \left(- \, w_3 \, + \, e^{2\, w_2} \, w_4\right)\\
\mathcal{L}^{[2]}_{2s+4} & = & \frac{1}{2} \, e^{-\,w_2} \, \left( 4\, + \, e^{2\, w_2} \, \left( 4\,  +\, \mid \mathbf{W}_6 \mid^2\right) \right) \\
\hline
\end{array}
 $
 \end{center}
 \caption{The explicit form of the two vectors $\mathcal{L}^{[1,2]}$ necessary for the calculation of the
 eigenvalues of the matrix $\mathfrak{m}$ which define the final geodesic distance function among any two points in the symmetric
 manifold $\mathcal{M}^{[2,s]}_\mathfrak{d}$.  The notation that we have adopted here show the complete covariance with respect
 to the paint group. We have used bold faced symbols to denote $2s$-vectors like $\mathbf{W}_5$ and $\mathbf{W}_6$ whose scalar product
 is denoted with a dot.\label{vettorini}}
\end{table}

%%%%%%%%%%

\subsubsection{The squared norm in the $r=2$ case} As it was done for the case $r=1$, also here the Paint invariant description of the squared
distance is   the definition of the squared norm of a solvable group element. Denoting
%%%%%%%%%%%%%%%%%%%%%%%%%%%
\begin{equation}\label{upponer2}
  \Upsilon_{\mathbf{w}} \, =\, \left\{w_1,w_2,w_3,w_4,\mathbb{U}_i,\mathbb{V}_j \right\} \quad ;  \quad i,j \, = \, 1,\dots, 2s
\end{equation}
the solvable coordinate vector that identifies the generic solvable group element via the $\Sigma$ exponential map, we find that the squared norm
of the solvable group element is  \begin{eqnarray}
\label{xorosho}
 \mathrm{ N^2}\left(\mathbf{w} \right) &=& \frac{1}{2} \left({\rm arccosh}\left[\frac{1}{2}
   \left(\mathfrak{A}_\mathbf{w}-\sqrt{\mathfrak{B}_\mathbf{w}}\right)\right]\right)^2+\frac{1}{2} \left({\rm arccosh}\left[\frac{1}{2}
   \left(\mathfrak{A}_\mathbf{w}+\sqrt{\mathfrak{B}_\mathbf{w}}\right)\right]\right)^2
   \end{eqnarray}
   where
\begin{eqnarray}
\label{Ainv}
\mathfrak{A}_\mathbf{w} &=& \frac{1}{128} \left(32 \mathbb{V}^2+4 e^{2 w_1} \left(\mathbb{U}^2+2 w_3
   w_4+4\right){}^2+32 \left(\mathbb{U}^2+e^{-2 w_2} w_3^2+e^{2 w_2} w_4^2\right)\right.\nonumber\\
   &&\left.+e^{2
   w_1} \left(4 \mathbb{U}\cdot\mathbb{V}+\sqrt{2} \left(\left(\mathbb{V}^2+4\right) w_3-4
   w_4\right)\right){}^2+2 \left(\mathbb{V}^2 e^{w_2}+8 \sinh
   \left(w_2\right)\right){}^2\right.\nonumber\\
   &&\left.+2 \left(\mathbb{V}^2 \sinh
   \left(w_2\right)+\left(\mathbb{V}^2+8\right) \cosh \left(w_2\right)\right){}^2+64
   e^{-2 w_1}\right)
\end{eqnarray}
and
\begin{eqnarray}
\label{Binv}
 \mathfrak{B}_\mathbf{w}&=& \frac{1}{16384} \,\left( 16 e^{-4 w_2} \left(e^{2 \left(w_1+w_2\right)} \left(\mathbb{U}^2+2 w_3 w_4+4\right)
   \left(4 \mathbb{U}\cdot\mathbb{V}+\sqrt{2} \left(\left(\mathbb{V}^2+4\right) w_3-4
   w_4\right)\right)\right.\right.\nonumber\\
   &&\left.\left.+4 \left(4 \mathbb{U}\cdot\mathbb{V} e^{2 w_2}+\sqrt{2}
   \left(\mathbb{V}^2+4\right) e^{4 w_2} w_4-4 \sqrt{2} w_3\right)\right){}^2+\left(32
   \mathbb{V}^2-4 e^{2 w_1} \left(\mathbb{U}^2+2 w_3 w_4+4\right){}^2 \right.\right.\nonumber\\
   &&\left.\left.-32
   \left(\mathbb{U}^2+e^{-2 w_2} w_3^2+e^{2 w_2} w_4^2\right)+e^{2 w_1} \left(4
   \mathbb{U}\cdot\mathbb{V}+\sqrt{2} \left(\left(\mathbb{V}^2+4\right) w_3-4
   w_4\right)\right){}^2+2 \left(\mathbb{V}^2 e^{w_2}+8 \sinh
   \left(w_2\right)\right){}^2\right.\right.\nonumber\\
   &&\left.\left.+2 \left(\mathbb{V}^2 \sinh
   \left(w_2\right)+\left(\mathbb{V}^2+8\right) \cosh \left(w_2\right)\right){}^2-64
   e^{-2 w_1}\right){}^2 \right)
\end{eqnarray}
We have thus shown the explicit expression  of the solvable group squared norm in an   paint group invariant form. The number of
paint group invariants from which it depends is 7 rather than 6. Indeed, besides the four singlets $w_1,w_2,w_3,w_4$, there are also the Paint
invariant squared norm $\mathbb{U}^2$, $\mathbb{V}^2$ and the scalar product $\mathbb{U}\mathbb{V}$. When we are on the Tits-Satake submanifold
the vectors $\mathbb{U}$ and $\mathbb{V}$ are one dimensional and the invariant $\mathbb{U}\mathbb{V}$ is not independent from the first two. In
the case of solvable group elements that do not entirely on the Tits-Satake submanifold, the third invariant is independent since it includes the
angle between the two vectors, and it is significant. The analysis of the squared norm formula with respect to the subpaint  group can be done in
exactly the same way as in the $r=1$ case and we do not repeat it here.
%%%%%%%%%%%%%%%%%%%%%%%%%%%%%%%%%%%%%%%%%%%%%%%%%%%%%%%%%%%%%%%%%%%%%%%%%%%%%%%%%%%%%%%%%%
\section{Systematics of low-rank $\mathcal{M}^{[r,s]}$ manifolds, clustering and Grassmannian leaves}
\label{lollobrigida}
%%%%%%%%%%%%%%%%%%%%%%%%%%%%%%%%%%%%%%%%%%%%%%%%%%%%%%%%%%%%%%%%%%%%%%%%%%%%%%%%%%%%%%%%%%
In this section we show that the manifolds $\mathcal{M}_\mathfrak{d}^{[r,s]}$ and their Tits-Satake projections (see eq.~(\ref{paniepesci})), for
small values of the noncompact rank $r=1,2,3,4$, for which  the distance formula  (\ref{coriandolidipasta}) can be
written explicitly in terms of radicals, fall into a surprisingly systematic scheme.
\begin{table}[htb]
\begin{center}
{
\scriptsize{
\begin{tabular}{||c||ccl||c||ccl|c|c||}
  \hline
  \hline
 \null & \emph{Full Manif.} & \null & \emph{Type} & \emph{Paint}  &  \emph{Tits-Satake} & \null & \emph{Type} & \emph{c}& \null\\
 $r_{nc}$ & $\mathcal{M}$ & \null & $\mathrm{Rie}_m, \mathcal{K}_n$&\emph{Group}& $\mathcal{M}_{TS}$ &\null &
 $\mathrm{Rie}_m, \mathcal{K}_n$ & \emph{map} & \null \\
 \null & $\frac{\mathrm{U}}{\mathrm{H_c}}$ & \null & $\mathcal{SK}_n,\mathcal{QM}_{4m}$&$\mathrm{G_{Paint}}$& $\frac{\mathrm{U_{TS}}}{\mathrm{H_{TS}}}$
 &\null &
 $\mathcal{SK}_n,\mathcal{QM}_{4m}$ & \null & \null \\
 \null & \null &\null & \null & \null &\null & \null &\null & \null & \null \\
 \hline
 \hline
 \null & \null &\null & \null & \null &\null & \null &\null & \null & \null \\
 $1$& $\frac{\mathrm{SO(1,1+2s)}}{\mathrm{SO(1+2 s)}}$ & $ =$ &$ \mathrm{Rie}_{1+2s}$& $\mathrm{SO(2s)}$ &
 $\frac{\mathrm{SO(1,2)}}{\mathrm{SO(2)}}$ &$\sim$&$\mathcal{SK}_1 \equiv \frac{\mathrm{SU(1,1)}}{\mathrm{U(1)}}$&
 $\stackrel{c}{\rightarrow}$ & $ \underbrace{\frac{\mathrm{G_{2(2)}}}{\mathrm{SU(2)\times SU(2)}}}_{\mathcal{QM}_{8}}  $ \\
 \null & \null &\null & \null & \null &\null & \null &\null & \null & \null\\
\hline
 \hline
 \null & \null &\null & \null & \null &\null & \null &\null & \null & \null \\
 \null& $\frac{\mathrm{SO(2,2+2s)}}{\mathrm{SO(2)\times SO(2+2 s)}}$ & $ =$ &$ \mathcal{K}_{2+s}$& \null &
 $\frac{\mathrm{SO(2,3)}}{\mathrm{SO(2)\times SO(3)}}$ &$=$&$\mathcal{K}_3 = \frac{\mathrm{Sp(4,\mathbb{R})}}{\mathrm{U(2)}}$&
\null & \null \\
$2$ & $\downarrow$ &\null & $\downarrow$ & $\mathrm{SO(2s)}$&$\downarrow$ & \null&$\downarrow$ & \null & \null \\
\null &$ \frac{\mathrm{SL(2,\mathbb{R})}}{\mathrm{SO(2)}} \times \frac{\mathrm{SO(2,2+2s)}}{\mathrm{SO(2)\times SO(2+2 s)}}$
 &$=$& $\mathcal{SK}_{3+2s}$&\null & $\frac{\mathrm{SL(2,\mathbb{R})}}{\mathrm{SO(2)}} \times \frac{\mathrm{SO(2,3)}}{\mathrm{SO(2)\times SO(3)}}$
  & $=$ &$\mathcal{SK}_{3}$ & $\stackrel{c}{\rightarrow}$
  & $\underbrace{\frac{\mathrm{SO(4,5)}}{\mathrm{SO(4)\times SO(5)}}}_{\mathcal{QM}_{20}}  $ \\
 \hline\hline
 %%%%%%%%%%%%%%%%%%%%%%%%%%%%%%%%%%%%%%%%%%%%%%%%%%%
 \null & \null &\null & \null & \null &\null & \null &\null & \null & \null \\
 $3$ & $\frac{\mathrm{SO(3,3+2s)}}{\mathrm{SO(3)\times SO(3+2 s)}}$ & $ =$ &$ \mathrm{Rie}_{9+6s}$& $\mathrm{SO(2s)}$ &
 $\frac{\mathrm{SO(3,4)}}{\mathrm{SO(3)\times SO(4)}}$ &$=$&$\mathcal{Q}_{12} $&
$\stackrel{\mathrm{Id}}{\leftrightarrow}$ & $\underbrace{\frac{\mathrm{SO(3,4)}}{\mathrm{SO(3)\times SO(4)}}}_{\mathcal{QM}_{12}}  $ \\
\null& \null&\null & \null & \null&$ \frac{\mathrm{SL(2,\mathbb{R})}}{\mathrm{SO(2)}} \times \frac{\mathrm{SL(2,\mathbb{R})}}{\mathrm{SO(2)}}
 $ & $=$ &$\mathcal{SK}_2 $& $\stackrel{c^{-1}}{\longleftarrow}$& $\hookleftarrow$\\
 \null & \null &\null & \null & \null &\null & \null &\null & \null & \null \\
 \hline\hline
  %%%%%%%%%%%%%%%%%%%%%%%%%%%%%%%%%%%%%%%%%%%%%%%%%%%
 \null & \null &\null & \null & \null &\null & \null &\null & \null & \null \\
 $4$ & $\frac{\mathrm{SO(4,4+2s)}}{\mathrm{SO(4)\times SO(4+2 s)}}$ & $ =$ &$ \mathcal{Q}_{16+8s}$& $\mathrm{SO(2s)}$ &
 $\frac{\mathrm{SO(4,5)}}{\mathrm{SO(4)\times SO(5)}}$ &$=$&$\mathcal{Q}_{20} $&
$\stackrel{\mathrm{Id}}{\leftrightarrow}$ & $\underbrace{\frac{\mathrm{SO(4,5)}}{\mathrm{SO(4)\times SO(5)}}}_{\mathcal{QM}_{20}}  $ \\
\null& \null&\null & \null & \null&$ \frac{\mathrm{SL(2,\mathbb{R})}}{\mathrm{SO(2)}} \times \frac{\mathrm{Sp(4,\mathbb{R})}}{\mathrm{U(2)}} $
 & $=$ &$\mathcal{SK}_3 $& $\stackrel{c^{-1}}{\longleftarrow}$& $\hookleftarrow$\\
 \null & \null &\null & \null & \null &\null & \null &\null & \null & \null \\
 \hline\hline
\end{tabular}
}}
\end{center}
\caption{Table of the noncompact symmetric manifolds of type $\mathrm{SO(r,r+2s)/SO(r)\times SO(r+2s)}$ and of
their Tits-Satake submanifolds, with emphasis on their relation via $c$-map. \label{rmeno4}}
\end{table}
\par
Indeed, apart from the practical reasons of convenience due to the algebraic solvability in terms of radicals, the limitation to noncompact rank
$r_{nc} \, \leq \, 4$ has deeper mathematical reasons, rooted in Special K\"ahler Geometry, special Quaternionic Geometry and the $c$-map
that relates them. One should  better not   ignore such structural connections, even if  one's final aim is the application of the geometrical theory
to Data Science. This becomes quite clear if we combine the results displayed in Tables  \ref{skTS} and \ref{homomodelisti} as we have done  in  Table
\ref{rmeno4}. Inspecting this table we realize that, at the level of the Tits-Satake submanifolds, which is what really matters if we rely on
\emph{paint group} and \emph{sub paint group} \emph{invariance/covariance}, as indeed we plan to do, and if we utilize both conceptual
resources at our disposal, namely the \emph{Tits-Satake projection}  and the \emph{$c$-map}, then there are only three Special K\"ahler Geometries that
codify all possible geometrical questions involved in the mapping of Data to symmetric manifolds of type  $\mathrm{SO(r,r+2s)/SO(r)\times
SO(r+2s)}$ with $r \leq 4$, namely:
\begin{equation}\label{quattrobello}
  \mathcal{SK}_1 \equiv \frac{\mathrm{SL(2,\mathbb{R})}}{\mathrm{SO(2)}},\quad\mathcal{SK}_2 \equiv \frac{\mathrm{SL(2,\mathbb{R})}}{\mathrm{SO(2)}}
  \times \frac{\mathrm{SL(2,\mathbb{R})}}{\mathrm{SO(2)}},
  \quad\mathcal{SK}_3 \equiv \frac{\mathrm{SL(2,\mathbb{R})}}{\mathrm{SO(2)}}\times \frac{\mathrm{Sp(4,\mathbb{R})}}{\mathrm{U(2)}},
\end{equation}
From the purely metric point of view the only manifold ingredients are just two, namely the hyperbolic Lobachevsky-Poincar\'e plane
\begin{equation}\label{lobapoinca}
  \mathbb{H}_2 \, \equiv \, \frac{\mathrm{SL(2,\mathbb{R})}}{\mathrm{SO(2)}}
\end{equation}
and the Siegel upper half plane of rank (or genus) $g=2$ (see below)
\begin{equation}\label{siegeluppo}
  \mathbb{SH}_2 \, \equiv \, \frac{\mathrm{Sp(4,\mathbb{R})}}{\mathrm{U(2)}}
\end{equation}
The Special K\"ahler structure is more than metric, as it involves the determination of \emph{the flat symplectic bundle} and of the matrix
$\mathcal{N}_{\Lambda\Sigma}$ (see sections \ref{scrittaN}), yet for Special K\"ahler manifolds that are homogeneous symmetric spaces like those
in eq.~\eqref{quattrobello} this amounts only to the determination of the symplectic $\mathbb{W}$-representation. Once this  substantially easy step is completed, all the rest follows automatically, via  the $c$-map construction, at least  for the manifolds that are in the image of the $c$-map. From Table \ref{rmeno4} we see that the Tits-Satake submanifolds of the cases $r=3,4$
are   of that type, and   their Special K\"ahler progenitors are in the list (\ref{quattrobello}), hence the previous observation
applies. Given the relevance of the solvable group--metric manifold equivalence, two aspects are particularly inspiring: first, the insight into the algebraic structure of the quaternionic K\"ahler manifolds provided by the $c$-map, and second, the $c$-map construction of the vielbein (eq.~ (\ref{filibaine})) which satisfies the Maurer-Cartan equations of the full solvable group (see eq.~~\eqref{MCSolv}). Specifically relevant is the \emph{Golden Split} (see section 1.7
of \cite{advancio} eq.~(1.7.3)) of the quaternionic Lie algebra of the Quaternionic K\"ahler manifolds:
\begin{equation}
\mbox{adj}(\mathbb{U}_{\mathcal{Q}}) =\mbox{adj}(\mathrm{SL(2,\mathbb{R})_E})\oplus
\mbox{adj}(\mathbb{U}_{\mathcal{SK}})\oplus
\mathfrak{W}_{(2,\mathbf{W})} \label{goldendecompo}
\end{equation}
where $\mathrm{SL(2,\mathbb{R})_E}$ is the Ehlers group, $\mathbb{U}_{\mathcal{SK}}$ is the Special K\"ahler Geometry Lie Group and $\mathbf{W}$
is the symplectic representation of $\mathbb{U}_{\mathcal{SK}}$ required by Special K\"ahler Geometry. In our case the $\mathbf{W}$-representations
are $(\mathbf{2},\mathbf{3})$ and $(\mathbf{2},\mathbf{5})$  for  $\mathrm{SO(4,3)}$ and $\mathrm{SO(4,5)}$.
\par
The golden decomposition and its restriction to the solvable subgroup is    precious for the issue of discretizations which, as we already
stressed, cannot avoid the consideration of discrete subgroups of the  solvable group and of its discrete automorphism group. Indeed the golden
decomposition is a guide for the construction of discrete subgroups of the larger group utilizing, as  bricks, the already determined discrete
subgroups of the building blocks. This  observation applies to discrete subgroups of all types, both those with elliptic-like rotations and
parabolic-like translations or those with elliptic like rotations and Fuchsian-like translation (see \cite{tassellandum} for the
latter case). In the first case (parabolic like translations), besides the $\mathrm{Sp(4,\mathbb{Z})}$ example discussed in section
\ref{plantageneti} we have the general construction, valid for all $\mathrm{SO(r,r+q,\mathbb{Z})}$ groups presented in section
\ref{generaloneSOpq}.
\par
Having outlined the perspective, let us turn  to the details.
%%%%%%%%%%%%%%%%%%%%%%%%%%%%%%%%%%%%%%%%%%%%%%%
\subsection{The manifolds $\mathrm{SO(r,r+2s)}/\mathrm{SO(r)\times SO(r+2s)}$  for small $r$} Let us start from the following strategic
observation. Referring to eq.~ (\ref{cagnesco}) we note that the time dependent matrix $\mathcal{M}(t)$ has the   structure
\begin{equation}\label{trucciolo}
  \mathcal{M}(t) \, = \, \exp\left[-2 \, t \, L \right]\, = \, \exp\left[-2 \, t \, \mathcal{O}\,L_{diag}\, \mathcal{O}^T \right] \, = \,
  \mathcal{O}\, \exp\left[-2 \, t\,   L_{diag}  \right]\,\mathcal{O}^T,
\end{equation}
where $L$ is the symmetric Lax operator at the origin, whose diagonalized form, obtained by an orthogonal transformation $\mathcal{O}\in
\mathrm{SO(2r+2s)}$,  is
\begin{equation}\label{diagostrut}
  \mathrm{diag}\left(L_{diag}\right) \, = \, \left\{\mu_1,\mu_2,\dots,\mu_r,0,\dots,0,-\mu_r,-\mu_{r-1},\dots,-\mu_1\right\}.
\end{equation}
Any geodesic,     specified by   an initial point and an initial Lax operator, i.e., a generic constant matrix $L$ satisfying  the
  two conditions
\begin{equation}\label{laxone}
  L=L^T \quad ; \quad L\,\eta_t + \eta_t \, L \, = \, 0,
\end{equation}
can be worked out  in an algoritmic way if we are able to explicitly  solve  the secular equation for the eigenvalues $\lambda_i$
($i=1,\dots,2r+2s$) and construct the corresponding ortho-normalized eigenvectors $\mathbf{v}_i$, which constitute the columns of the orthogonal
matrix $\mathcal{O}$. In view of eq.~\eqref{diagostrut}, we see that the secular equation of the Lax operator has the   general form
\begin{equation}\label{capitano}
  \mathfrak{P}(\lambda) \, = \, \lambda^{2s} \, \prod_{i=1}^r\left(\lambda^2 - \mu_i^2\right).
\end{equation}
Hence,   the eigenvalues $\mu_i$ admit an explicit expression in terms of radicals of the Lax operator components for all
values of the noncompact rank $r=1,2,3,4$. In other words, the geodesic problem is algebraically solvable for all such values of the noncompact
rank. For large values of the parameter $s$ difficulties might arise in finding an appropriate basis of orthonormal vectors in the null-eigenspace
(the eigenspace corresponding to the eigenvalue $\lambda=0$), but certainly, for $r=1,2,3,4$, we can work out explicitly  all
  geodesics in the corresponding Tits-Satake projected manifolds
\begin{equation}\label{colonnello}
  \mathcal{M}^{[r,r+1]}_{TS} \equiv \, \frac{\mathrm{SO(r,r+1)}}{\mathrm{SO(r)} \times \mathrm{SO(r)}}.
\end{equation}
In the previous sections we   already illustrated the structure of the square norm function and of the geodesic construction for the cases $r=1,2$,
which   provide all the building blocks for all the cases $r=1,2,3,4$. A further study
of those  building blocks will follow in section \ref{plantageneti}.  We observe that a similar consideration applies, in view of
eq.~\eqref{coriandolidipasta}, to the determination of the geodesic distance between any two points (both in the full manifold
$\mathcal{M}_\mathfrak{d}^{[r,s]}$ and in its Tits-Satake projection $\mathcal{M}^{[r,r+1]}_{TS}$). The argument is very simple: the distance
$\mathrm{d}^2({\bf u},{\bf v})$ is given purely in terms of the eigenvalues of the matrix
\begin{equation}\label{baldacchino}
 \mathfrak{m} \, = \, \mathfrak{m} \left( \boldsymbol{\Upsilon}[\mathbf{u},\mathbf{b}]\right)
\end{equation}
defined in eq.~(\ref{cagnesco}), in that, writing
\begin{equation}\label{carriola}
  \mathfrak{m} \, = \,\mathfrak{O}\,
  \left(
  \begin{array}{cccccc}
   \lambda_{1} & 0 & \ldots & \ldots & 0 & 0 \\
   0 & \lambda_{2} & 0 & \ldots & \ldots & 0 \\
   \vdots & \ldots & \ddots & \ldots & \ldots & \vdots \\
   \vdots & \ldots & \ldots & \ddots & \ldots & \vdots \\
   0 & \ldots & \ldots & 0  & \lambda_{r-1} & 0 \\
   0 & \ldots & \ldots & \ldots & 0 & \lambda_r \\
  \end{array}
  \right) \mathfrak{O}^T
\end{equation}
where $\lambda_i=\lambda_i[\mathbf{u},\mathbf{v}]$ are the eigenvalues of $\mathfrak{m}$ (not to be confused with the eigenvalues $\mu_i$ of the
Lax operator) and $\mathfrak{O}\in \mathrm{SO(}r\mathrm{)}$ is the diagonalizing orthogonal matrix (not to be confused with $\mathcal{O}\in
\mathrm{SO(2r+2s)}$),   the distance formula (\ref{coriandolidipasta}) becomes
\begin{equation}\label{pasticcione}
  \mathrm{d}^2({\bf u},{\bf v})=\frac{1}{4}\,\sum_{i=1}^r\,{\rm arccosh}^2\left(\lambda_i[\mathbf{u},\mathbf{v}]\right).
\end{equation}
Once again, if $r\leq 4$ the secular polynomial is of maximum degree $r=4$ and the eigenvalues can be explicitly written in terms of radicals, so
that the dependence on the solvable coordinates $\boldsymbol{\Upsilon}_{\mathbf{u}}$ and $\boldsymbol{\Upsilon}_{\mathbf{v}}$ of the two points is
explicit and written in terms of algebraic functions and elementary transcendentals.
\subsection{Clustering}
The first main issue in \emph{Data Science Applications} of \emph{Symmetric Spaces with a non trivial Tits-Satake Projection} is the
\emph{grouping of data}.  The fact that any data set can be embedded into a manifold $\mathrm{SL(N,\mathbb{R})/SO(N)}$ has been discussed at
length in the previous sections. The idea is that whatever     differentiable manifold \(\mathcal{M}_{\text{data}}\) might be the
hosting environment for the set of data, certainly that manifold can be embedded into $\mathrm{SL(N,\mathbb{R})/SO(N)}$:
\begin{equation}\label{isocrate}
\mathcal{M}_{\text{data}} \overset{\Phi }{\longrightarrow } \text{SL}(\mathrm{N},\mathbb{R})/\text{SO}(\mathrm{N})
\end{equation}
for a suitable value of $\mathrm N$. Here
  $\Phi $ denotes a suitable set of functions, described by their developments into harmonics, whose vanishing locus is
\(\text{\emph{$\mathcal{M}_{\text{\rm data}}$}}\).
Hence the embedding into some $\mathrm{SL(N,\mathbb{R})/SO(N)}$ cannot be easily determined by \emph{unsupervised neural networks}, as the
set of possibilities is too large and actually infinite: the uncertainty on the choice of $\mathrm{N}$, the infinite number of coefficients in
each function $\Phi $ and the a priori undetermined number of them. Such a task may be the target only of \emph{supervised learning} in
presence of suitable training sets.
\subsubsection{A nontrivial hypothesis}\label{nontrihyp}
Assuming that the manifold hosting the data may be embedded into a non-maximally split symmetric space with a nontrivial Tits-Satake
projection such as $\frac{\mathrm{SO}(r,r+2s)}{\mathrm{SO}(r)\times \mathrm{SO}(r+2s)}$, is on the contrary a strong \emph{hypothesis}, which implies a clustering of the data into \emph{Tits-Satake fibers}. Conceptually this is so because
$\frac{\mathrm{SO}(r,r+2s)}{\mathrm{SO}(r)\times \mathrm{SO}(r+2s)}$ is naturally embedded into    $\mathrm{SL(N,\mathbb{R})/SO(N)}$ for
$\mathrm{N}=2(r+s)$ by means of the triangular embedding $\boldsymbol{\Phi_{\mathrm{triang}}}$ advocated in Statement \ref{statamento} and
explicitly encoded in the construction of the matrix  $\eta_t$ (see eq.~(\ref{etatdefi})):
\begin{equation}\label{certosino}
\mathcal{M}_{\text{data}}\quad\overset{\Phi _{\text{unkwown}}}{\longrightarrow } \quad \frac{\text{SO}(r,r+2s)}{\text{SO}(r)\times \text{SO}(r+2s)}\quad
\overset{\Phi _{\mathrm{triang}}}{\longrightarrow} \,\,\frac{\text{SL}(2r+2s,\mathbb{R})}{\text{SO}(2r+2s)}.
\end{equation}
Looking at this equation we see that embedding the Data into a specific manifold of type $\frac{\mathrm{SO}(r,r+2s)}{\mathrm{SO}(r)\times
\mathrm{SO}(r+2s)}$ with a specific value of $r$, in particular $r=1,2,3,4$,  is a hypothesis, since it corresponds to \emph{choosing a priori} a
part of the embedding functions. The triangular   embedding  $\boldsymbol{\Phi_{\mathrm{triang}}}$  is determined by setting
to zero a suitable    set of harmonic functions. Actually these functions are just matrix elements of the fundamental harmonic \(\mathcal{M}^{\text{AB}}\) that are set to zero, or
related to each other by quadratic equations which amounts to a fully intrinsic and coordinate independent description of the immersion of
$\mathrm{SO}(r,r+2s)/\mathrm{SO}(r)\times \mathrm{SO}(r+2s)$ into \(\text{SL}(2r+2s,\mathbb{R})/\mathrm{SO}(2r+2s)\). Indeed the matrix
\(\mathcal{M}^{\text{AB}}\) is required to be an element of the group $\mathrm{SO}(r,r+2s)$, so it satisfies the equations
\begin{equation}\label{domenicano}
\mathcal{M} \, \eta _t \, \mathcal{M} =\eta _t \quad ; \quad \mathcal{M}\, = \, \mathcal{M}^T.
\end{equation}
It follows that assuming the first embedding in eq.~(\ref{certosino}) is a hypothesis and, as such, it might be the target of an
\emph{unsupervised neural network learning}. What is the indicator that this hypothesis is viable?
\subsubsection{Tits-Satake grouping}
\label{eulerando} If the data admit the embedding into $\mathrm{SO}(r,r+2s)/\mathrm{SO}(r)\times \mathrm{SO}(r+2s)$, it  follows that  data
that have the same Tits-Satake projection may not coincide, yet they must be grouped  in some way, forming similarity classes. How this grouping
is performed is a geometrical issue that can be analyzed in different ways adopting different viewpoints. It is a fundamental question that is
left for the future while trying to construct new neural networks based on the mathematical tokens provided by the PGTS theory presented in this
paper; we just remind the reader that the constructions of references \cite{francesi1,francesi2,francesi3}, have implicitly but unconsciously
assumed the above mentioned hypothesis and this in its more restrictive declination, i.e.,  the choice $r=1$. We proceed to sum up the viewpoint
based on the Euler parameterization of the manifold $\mathcal{M}_{\mathfrak{d}}^{[r,s]}$ and on Grassmanians.

\subsection{The structure of Tits-Satake fibers from the Euler viewpoint}
In the general case \(\text{SL}(\mathrm{N},\mathbb{R})\)/SO(N) we have the alternative Euler parameterization of points of the manifold that is
provided by the pair $\left\{\exp[\boldsymbol{\mu}],\mathcal{O}\right\}$ formed by an element of the maximal torus $\exp[\boldsymbol{\mu}]$,
namely, the exponential of any element of the Cartan subalgebra, $\boldsymbol{\mu}\in\mathcal{H}\subset \slal(\mathrm{N},\mathbb{R})$, and any
element $\mathcal{O}\in \mathrm{SO(N)}$.  Let us briefly outline the logic. A point
in the manifold is uniquely determined by a symmetric matrix $\mathcal{M}^{AB}$ with the required properties, i.e., just $\mathrm{det} \mathcal{M} \, =
\, 1$ in the general \(\text{SL}(\mathrm{N},\mathbb{R})\)/SO(N) case, condition (\ref{domenicano}) in the \(
\mathrm{SO}(r,r+2s)/\mathrm{SO}(r)\times \mathrm{SO}(r+2s)\) case. Such matrix is uniquely determined by its spectrum of eigenvalues $\mu_i$ and
by the orthogonal matrix $\mathcal{O}$ which diagonalizes it). On the other hand, by means of the Cholevsky-Crout algorithm, every symmetric
matrix can be written as $\mathcal{M}\, =\,\mathbb{L}\cdot\mathbb{L}^T$, where $\mathbb{L}$ is a triangular matrix namely an element of the
solvable group $\mathcal{S}=\exp[Solv]$. Furthermore by means of the inverse map $\Sigma^{-1}$ from $\mathbb{L}$ we retrieve the corresponding
element of the Solvable Lie algebra, that is, the corresponding solvable coordinates $\boldsymbol{\Upsilon}$. If we use the simple exponential for
$\Sigma$ we see that
\begin{eqnarray}\label{famigerato}
  \mathcal{O}^T \, \mathcal{M} \, \mathcal{O} & = & \exp\left[ \mu_i \mathcal{H}^i \right] \quad ; \quad \mathcal{H}^i =
  \text{Cartan generators} \nonumber\\
  & \Downarrow & \nonumber\\
 \exp\left[ \mu_i \mathcal{H}^i \right]\,  & = &  \mathcal{O}^T \, \mathbb{L} \, \mathcal{O} \, \cdot \, \mathcal{O}^T \,
  \mathbb{L}^T \, \mathcal{O} \quad \Rightarrow \quad  \mathcal{O}^T \, \mathbb{L} \, \mathcal{O} \, = \,
  \exp\left[ \ft 12 \, \mu_i \mathcal{H}^i \right] \nonumber\\
 & \Downarrow & \nonumber\\
 \mathcal{O}^T \, \exp[\mathfrak{s}] \, \mathcal{O} & = & \exp\left[\mathcal{O}^T \, \mathfrak{s} \, \mathcal{O}\right]
 \, = \, \exp\left[ \mu_i \mathcal{H}^i \right], \quad \text{where $\mathfrak{s} \in Solv$} .
  \end{eqnarray}
Implicitly, the above argument shows that we have a linear action of the orthogonal group on the solvable Lie algebra which rotates a generic
element into the Cartan subalgebra. Hence the solvable Lie algebra element $\mathfrak{s}$ that maps to a point in the symmetric manifold can be
parameterized by a Cartan subalgebra element $\boldsymbol{\mu} \cdot \boldsymbol{\mathcal{H}}$ and a finite orthogonal group element $\mathcal{O}$
as we stated; this is the Euler parameterization.
\subsubsection{Euler parameterization in the non-maximally split cases}
In the non maximally split cases the Euler parameterization is more subtle than in the maximally split case, since the element
$\boldsymbol{\mu}_{n.c.} \in \mathcal{H}_{n.c.}$ must belong to the noncompact Cartan subalgebra and the element $\mathcal{O}\in \mathrm{H}_c $
in the maximal compact subgroup is defined only up to the stabilizer of $\boldsymbol{\mu}_{n.c.}$ which is the paint group $\mathrm{G_{Paint}}$.
In other words the Euler parameterization is given by a pair $\left\{\exp[\boldsymbol{\mu}_{n.c.}], \mathfrak{p} \right\} $, where
\begin{equation}\label{panamerican}
  \mathfrak{p} \,\in \,\frac{\mathrm{H}_c}{\mathrm{G_{Paint}}}.
\end{equation}
We can easily check the dimensions. The dimension of the complete manifold is
\begin{equation}\label{dimMrs}
 \mathrm{ dim}\, \mathcal{M}^{[r,s]}_\mathfrak{d} \, = \, \mathrm{ dim} \, \frac{\mathrm{SO}(r,r+2s)}{\mathrm{SO}(r)\times\mathrm{ SO}(r+2s)}
 \, = \, r \, (r+2\, s),
\end{equation}
while
\begin{equation}\label{borlengo}
\mathrm{ dim}\, \frac{\mathrm{H}_c}{\mathrm{G_{Paint}}} \, = \, \mathrm{ dim}\, \frac{\mathrm{SO}(r)\times\mathrm{ SO}(r+2s)}{\mathrm{SO}(2s)}
\, = \, r^2 \, + \, r \, (2\, s \, -\, 1),
\end{equation}
so  that we  can verify that
\begin{equation}\label{cromatografo}
  \mathrm{ dim}\, \mathcal{M}^{[r,s]}_\mathfrak{d}  \, = \, \underbrace{\mathrm{ dim} \, \mathcal{H}_{n.c}}_{= \, r} \, + \, \mathrm{ dim}\,
  \frac{\mathrm{H}_c}{\mathrm{G_{Paint}}},
\end{equation}
where $\mathcal{H}_{n.c}$ denotes the noncompact Cartan subalgebra.
\par
Next we want to decompose the manifolds in the considered series with respect to the Tits-Satake submanifold
 \begin{equation}\label{tittosatacchio}
   \mathcal{M}^{[r,s]}_\mathfrak{d} \, \stackrel{\pi_{TS}}{\longrightarrow} \, \mathcal{M}^{[r,r+1]}_{TS} \,
   \sim \, \frac{\mathrm{SO}(r,r+1)}{\mathrm{SO}(r)\times\mathrm{ SO}(r+1)} \, \subset \, \mathcal{M}^{[r,s]}_\mathfrak{d},
 \end{equation}
 where $\pi_{TS}$ denotes the Tits-Satake projection. We have
 \begin{equation}\label{dimmatittasatakka}
 \mathrm{ dim} \, \mathcal{M}^{[r,r+1]}_{TS} \, =  \, \mathrm{ dim}\, \frac{\mathrm{SO}(r,r+1)}{\mathrm{SO}(r)\times\mathrm{ SO}(r+1)}  \, = \, r^2 +r.
 \end{equation}
It follows that for each point $p \in  \mathcal{M}^{[r,r+1]}_{TS}$ the fiber over it, namely $\pi_{TS}^{-1}(p)$, must have     dimension
\begin{equation}\label{dimmafibba}
  \mathrm{ dim} \, \pi_{TS}^{-1}(p) \, = \, r \, (2s-1).
\end{equation}
How do we interpret this result? Let us discuss it algebraically.
The Lie algebra $\mathbb{H}_{\mathrm{c}}$ of the maximally compact subgroup $\mathrm{H_c}$ admits the following decomposition as a vector space:
\begin{eqnarray}\label{vettadecompa}
  \mathbb{H}_{\mathrm{c}} & =& \mathbb{G}_{\mathrm{paint}}\, \oplus \, \mathbb{H}_{\mathrm{TS}} \, \oplus \, \mathbb{F}_{{\mathrm{TS}}}\nonumber\\
  \mathbb{G}_{\mathrm{paint}}& = & \so(2s)  \quad \text{ Lie algebra of the paint group} \nonumber\\
   \mathbb{G}_{\mathrm{TS}}\supset \mathbb{H}_{\mathrm{TS}}& = & \so(r) \oplus \so (r+1)  \quad \text{maximal compact subalgebra of $\mathbb{G}_{\mathrm{TS}} = \so(r,r+1)$} \nonumber\\
   \mathbb{F}_{\mathrm{TS}} & = & \text{vector subspace which does not close a Lie algebra} \quad \mathrm{dim} \, \mathbb{F}_{\mathrm{TS}} \, = \, r(2s-1)
\end{eqnarray}
Considering now the vector subspace $\mathbb{F}_{\mathrm{TS}} $, that we call    space of \emph{TS fiber generators}, we find
\begin{equation}\label{FTScommu}
  \left[\mathbb{F}_{\mathrm{TS}}\, , \, \mathbb{F}_{\mathrm{TS}}\right]   \,\equiv \, \mathbb{H}_{\mathrm{F}} \, = \, \so(r) \oplus \mathbb{G}_{\mathrm{subpaint}} \, = \, \so(r) \oplus \so(2s-1)
\end{equation}
Thus,  the commutator of $\mathbb{F}_{\mathrm{TS}}$ with itself generates a compact Lie algebra $\mathbb{H}_{\mathrm{F}}$ which is identified in
equation (\ref{FTScommu}). This suggests that the direct sum of vector spaces
\begin{equation}\label{corsarobianco}
 \mathbb{G}_{\mathrm{F}} \,\equiv \, \mathbb{H}_{\mathrm{F}}\oplus \mathbb{F}_{\mathrm{TS}} \sim \so(r+2s-1)
\end{equation}
 should be a  subalgebra of the full compact algebra $\mathbb{G}_{\mathrm{F}} \subset \mathbb{H}_c$, that is,    we should have
 \begin{equation}\label{corsaroverde}
\left[ \mathbb{G}_{\mathrm{F}} \,,\, \mathbb{G}_{\mathrm{F}} \right] \subset   \mathbb{G}_{\mathrm{F}},
\end{equation}
which indeed turns out to be true. The explicit structure of the Lie algebra $\mathbb{G}_{\mathrm{F}}$ is easily identified:
\begin{eqnarray}\label{corsaronero}
   \mathbb{G}_{\mathrm{F}} &=& \mathbb{H}_{\mathrm{F}}\oplus \mathbb{F}_{\mathrm{TS}} \sim \so(r+2s-1) \subset \so(r+2s)\nonumber \\
  \left[ \mathbb{H}_{\mathrm{F}}\, , \, \mathbb{F}_{\mathrm{TS}}\right]  &=& \mathbb{F}_{\mathrm{TS}}\nonumber \\
  \mathbb{G}_{\mathrm{F}} &\subset& \mathbb{H}_{\mathrm{c}} \, = \, \so(r) \oplus \so(r+2s)
\end{eqnarray}
Lifting the above results from the Lie algebra to the Lie group level we realize that the \emph{TS fiber generators} span in its origin the
tangent space to the following Grassmannian manifold:
\begin{equation}\label{crinolina}
  \mathcal{F}_{{TS}} \, = \,  \frac{\mathrm{SO}(r+2s-1)}{\mathrm{ SO}(r) \times \mathrm{SO}(2s-1)}
\end{equation}
This analysis suggests a new parameterization of the full manifold that is midway between the solvable   and the Euler
parameterization. The rationale is the following. In the Euler parameterization of $\mathcal{M}^{[r,s]}_\mathfrak{d}$ we can proceed in two steps.
First in eq.~(\ref{panamerican}) we choose
\begin{equation}\label{panamericanTS}
  \mathfrak{p}_{\mathrm{TS}} \,\in \,\mathrm{H}_{\mathrm{TS}}.
\end{equation}
The pair $\left\{\exp[\boldsymbol{\mu}_{n.c.}],\mathfrak{p}_{\mathrm{TS}} \right\} $ singles out a point in the Tits-Satake submanifold, and every
point $p\in\mathcal{M}_{TS}$ admits such a representation for a suitable $\boldsymbol{\mu}_{n.c.}(p)$ and a suitable
$\mathfrak{p}_{\mathrm{TS}}(p)$.  Next we choose a further generic element
\begin{equation}\label{panamericanTS}
  \mathfrak{p}_{\mathrm{F}} \,\in \, \frac{\mathrm{G}_{\mathrm{F}}}{\mathrm{H}_{\mathrm{F}}}
\end{equation}
and as $\mathfrak{p}_{\mathrm{F}}$ varies in $\mathcal{F}_{{TS}} $ and $\boldsymbol{\mu}_{n.c.}(p), \mathfrak{p}_{\mathrm{TS}}(p)$ vary with
$p\in\mathcal{M}_{TS}$,  we obtain all the points of $\mathcal{M}^{[r,s]}_\mathfrak{d}$. In other words the action of $\mathfrak{p}_{\mathrm{F}}$
on the Tits-Satake submanifold generates a homeomorphic copy of the latter in the full manifold. Similary for each point $p\in\mathcal{M}_{TS}$
considering its orbit under all elements  $\mathfrak{p}_{\mathrm{F}} \,\in \, \frac{\mathrm{G}_{\mathrm{F}}}{\mathrm{H}_{\mathrm{F}}}$  we obtain
a copy of the Grasmannian $\mathcal{F}_{{TS}}$, which we call a \emph{Grasmannian leaf}. Although each leaf  is associated with an initial
point $p\in \mathcal{M}_{TS}$, the leaf is not the Tits-Satake fiber $\pi^{-1}_{TS}(p)$. Indeed the
Grasmmanian leaf  projects in $\mathcal{M}_{TS}$ onto a domain $\mathcal{D}(p)\subset \mathcal{M}_{TS}$ of dimension
$\mathrm{dim}\left[\mathcal{D}(p)\right]\,=\,r$.
\par
Instead a Tits-Satake fiber $\pi^{-1}_{TS}(p)$  is not
compact;  rather, it can be defined in terms of a suitable \emph{normal subgroup of the solvable group}
$\mathcal{S}_{[r,s]}$, metrically equivalent to $\mathcal{M}^{[r,s]}_\mathfrak{d}$  as we show in the next section.
\subsection{The Tits-Satake projection revisited in a different perspective}\label{TSnormalsub}
We illustrate the idea mentioned in the last lines of the previous section by focusing on an explicit example namely on the case $r=2$, where for
simplicity we also choose $s=1$. The general form of solvable group element in terms of the solvable coordinates was provided in
eq.~(\ref{exempli2}). We consider that expression and for our present convenience we rename the coordinates $\Upsilon_{i}\to w_{i}$ and
$\Upsilon_{5,6,i }\to w_{5,6,i}$. The Tits-Satake projection corresponds to setting $w_{5,i}=w_{6,i}=0$ for ($i\geq 2$). This leads to the
solvable group element   \begin{eqnarray}\label{cornelioprisco}
 \mathbb{L}_{TS}(\mathbf{w}_{TS}) & = &  \left(
\begin{array}{cccccc}
 e^{w_1} & \frac{e^{w_1} w_3}{\sqrt{2}} & \frac{1}{2} e^{w_1} \left(\sqrt{2} w_5+w_3
   w_6\right) & 0 & e^{w_1} \left(-\frac{w_3 w_6^2}{4 \sqrt{2}}-\frac{w_5
   w_6}{2}+\frac{w_4}{\sqrt{2}}\right) & -\frac{1}{2} e^{w_1} \left(\frac{w_5^2}{2}+w_3
   w_4\right) \\
 0 & e^{w_2} & \frac{e^{w_2} w_6}{\sqrt{2}} & 0 & -\frac{1}{4} e^{w_2} w_6^2 &
   -\frac{e^{w_2} w_4}{\sqrt{2}} \\
 0 & 0 & 1 & 0 & -\frac{w_6}{\sqrt{2}} & -\frac{w_5}{\sqrt{2}} \\
 0 & 0 & 0 & 1 & 0 & 0 \\
 0 & 0 & 0 & 0 & e^{-w_2} & -\frac{e^{-w_2} w_3}{\sqrt{2}} \\
 0 & 0 & 0 & 0 & 0 & e^{-w_1} \\
\end{array}
\right) \nonumber\\
\mathbf{w}_{TS} & = & \left\{w_1,w_2,w_3,w_4,w_5,w_6\right\}
\end{eqnarray}
As we see, the embedding of the TS solvable group into the full   solvable group is described by  a $6\times 6$ matrix which has one-column (the third) and
one row (the fourth) with all zeros except identity at place $(4,3)$. In the general case we have $2s-1$ columns and $(2s-1)$ rows of that type.
\subsubsection{A normal subgroup of the full solvable group contains the TS fibers}
At variance with the generation of the Tits-Satake fibers by means compact rotations associated with the Grassmanian, as we did in section
\ref{eulerando}, we want here to  associate them with a suitable normal subgroup $\mathcal{NS}_{fib} \subset \mathcal{S}_{U/H}$ of the full
solvable group.
\par
As a candidate we consider the following  three parameter ansatz for the generic subgroup element (the parameters $\rho,p,q$ being respectively
associated with $w_4, w_{5,2},w_{6,2}$ and the generalization to higher values of $s$ being also clear and straightforward):
\begin{equation}\label{crodolino}
 \mathfrak{L}(\rho,p,q) \, \equiv \, \left(
\begin{array}{cccccc}
 1 & 0 & 0 & 2 p & \rho -4 p q & -2
   p^2 \\
 0 & 1 & 0 & 2 q & -2 q^2 & -\rho
   \\
 0 & 0 & 1 & 0 & 0 & 0 \\
 0 & 0 & 0 & 1 & -2 q & -2 p \\
 0 & 0 & 0 & 0 & 1 & 0 \\
 0 & 0 & 0 & 0 & 0 & 1 \\
\end{array}
\right).
\end{equation}
These matrices form a group as
\begin{eqnarray}\label{cervicale}
   \mathfrak{L}(\rho_1,p_1,q_1)\cdot \mathfrak{L}(\rho_2,p_2,q_2) &= & \mathfrak{L}(4p_1\,q_2 +\rho_1+\rho_2,p_1+p_2,q_1+q_2)
   \nonumber\\
   \mathfrak{L}^{-1} (\rho_1,p,q)& = & \mathfrak{L}(4p\,q -\rho,-p,-q).
 \end{eqnarray}
 This subgroup  is normal since
 \begin{eqnarray}
 \label{normalpazzo}
   &&\mathbb{L}_{\mathcal{M}}\cdot \mathfrak{L}(\rho,p,q)\cdot \mathbb{L}_{\mathcal{M}}^{-1} \,=\, \mathfrak{L}(\lambda,h,k) \nonumber\\
   &&  \lambda \, = \,  e^{w_1+w_2} \left(\sqrt{2} p w_{6,2}-\sqrt{2} q w_{5,2}+\sqrt{2} q^2
   w_3+\rho \right),\quad h\, = \,  \frac{1}{2} e^{w_1} \left(2 p+\sqrt{2} q w_3\right),\quad k\, = \,  q
   e^{w_2}. \nonumber\\
 \end{eqnarray}
 It follows that the full solvable group can be written in terms of lateral classes of the normal subgroup as
 \begin{equation}\label{crepuscolo}
  \mathcal{S}_{\mathrm{U/H}} \ni \mathbb{L}_{\mathcal{M}} \, = \, \mathbb{L}_{TS}(\mathbf{w}_{TS}) \cdot \mathfrak{L}(\rho,p,q)
 \end{equation}
 There seems to be a double counting since the parameter $w_4$ is contained both in the Tits-Satake group element and
 in the normal solvable group element. That is not the case  because of the following very useful and significant fact.
 In both cases --- the solvable Tits-Satake group $\mathcal{S}_{TS}$ and the normal subgroup $\mathcal{NS}_{fib} $ ---  the group elements
\begin{equation}\label{subnormale}
   \mathfrak{N}(\rho) \, = \, \left(
\begin{array}{cccccc}
 1 & 0 & 0 & 0 & \rho  & 0 \\
 0 & 1 & 0 & 0 & 0 & -\rho  \\
 0 & 0 & 1 & 0 & 0 & 0 \\
 0 & 0 & 0 & 1 & 0 & 0 \\
 0 & 0 & 0 & 0 & 1 & 0 \\
 0 & 0 & 0 & 0 & 0 & 1 \\
\end{array}
\right)
\end{equation}
 make up a normal subgroup since one has
\begin{eqnarray}\label{carpeneto}
   \mathfrak{L}(\rho,p,q)\cdot \mathfrak{N}(\rho)\cdot \mathfrak{L}^{-1}(\rho,p,q) &=& \mathfrak{N}(\rho) \\
   \mathbb{L}_{TS}(\mathbf{w}_{TS})\cdot \mathfrak{N}(\rho)\cdot \mathbb{L}^{-1}_{TS}(\mathbf{w}_{TS}) &=& \mathfrak{N}\left(\rho  \, e^{w_1+w_2}\right)
\end{eqnarray}
In other words we have two  chains of subgroups and normal subgroups:
\begin{eqnarray}\label{golovabolit}
   \mathcal{S}_{\mathrm{U/H}}& \rhd & \mathcal{NS}_{fib}  \, \rhd \, \mathfrak{N}\\
   \mathcal{S}_{\mathrm{U/H}} &\supset & \mathcal{S}_{\mathrm{TS}} \, \rhd \, \mathfrak{N}
\end{eqnarray}
Correspondingly, we can consider the Tits-Satake solvable group as \emph{the quotient of full solvable group with respect to the quotient of
normal subgroup $\mathcal{NS}_{fib}$ with respect to the group $\mathfrak{N}$}. In practice this means that we can write the lateral classes in
eq.~(\ref{crepuscolo}) as:
\begin{equation}\label{crepuscolodue}
  S_{\mathcal{M}} \ni \mathbb{L}_{\mathcal{M}} \, = \, \mathbb{L}_{TS}(\mathbf{w}_{TS})|_{w_4 = 0} \cdot \mathfrak{N}(\rho)\cdot \mathfrak{L}(0,p,q)
\end{equation}
Indeed since $\mathfrak{N}$ is a normal subgroup both for $\mathcal{S}_{TS}$ and for $\mathcal{NS}_{fib} $ the left and right lateral classes are
identical and this allows to put $\mathfrak{N}(\rho)$ in the middle. We stress that this point of view is quite instrumental in the problem  of
determining discrete parabolic subgroups of the groups $\mathrm{SO(r,r+s,\mathbb{Z})}$ that is addressed and solved in section
\ref{generaloneSOpq}.
\par
\subsubsection{Three complementary view points on the structure of the full manifold $\mathrm{U/H}$ in terms of its Tits-Satake submanifold}
\label{balubbus} Summarizing, in the present section we have illustrated how the full manifold can be viewed in three different complementary ways
in relation with  its Satake submanifold
\begin{equation}\label{pastrengo}
  \mathrm{\frac{U_{TS}}{H_{TS}} \, \subset \, \frac{U}{H}}.
\end{equation}
Given the Tits-Satake projection
\begin{equation}\label{tiziosata}
  \pi_{\mathrm{TS}} \, = \, \mathrm{\frac{U}{H}}\, \longrightarrow \, \mathrm{\frac{U_{TS}}{H_{TS}}}
\end{equation}
\begin{description}
  \item[1)] in the first picture we view $\mathrm{\frac{U}{H}}$ as the total space of a \emph{vector bundle} over the base manifold
      $\mathcal{M}_{\mathrm{base}}=\mathrm{\frac{U_{TS}}{H_{TS}}}$ with structure group $\mathrm{G_{struc}}\, =\, \mathrm{SO(r)}\times
      \mathrm{G}_{\mathrm{subpaint}}$. Indeed given a point $p\in\mathrm{\frac{U_{TS}}{H_{TS}}}$ we have:
      \begin{equation}\label{vettorispazi}
        \pi_{\mathrm{TS}}(p) \, \sim \, \mathbb{V}_{1}\oplus \mathbb{V}_{2} \oplus \dots \oplus\mathbb{V}_{r}
      \end{equation}
      where  $\mathbb{V}_{i}$ are $r$ copies of a $(2s-1)$-dimensional vector space supporting the fundamental representation of the subpaint
  group $\mathrm{G_{subpaint}}= \mathrm{SO(2s-1)}$. The group $\mathrm{SO(r)}$ included in the isotropy group of the base manifold acts on the
  $r$ vector spaces $\mathbb{V}_{i}$ rotating one into the other. The first picture is entirely located in the solvable
  parameterization of the full manifold $\mathrm{\frac{U}{H}}$.
 \item[2)] The second picture, as we have explained above, relies on a hybrid parameterization of $\mathrm{U/H}$, partially solvable, partially
     Euler. Given any point $p\in \mathrm{\frac{U_{TS}}{H_{TS}}}$, we consider its orbit with respect to all elements  $\mathfrak{p}_{\mathrm{F}}
     \,\in \, \mathcal{F}_{\mathrm{TS}}\, \equiv \,\frac{\mathrm{G}_{\mathrm{F}}}{\mathrm{H}_{\mathrm{F}}}$ that are compact rotations. The
     orbit $\mathcal{F}_{\mathrm{TS}}(p)\subset\mathrm{\frac{U}{H}}$ is a submanifold of the full symmetric space, diffeomorphic to the
     Grassmanian  $\frac{\mathrm{G}_{\mathrm{F}}}{\mathrm{H}_{\mathrm{F}}}$. The full manifold is
     \begin{equation}\label{paleocappo}
       \mathrm{\frac{U}{H}} \, = \, \bigcup_{p\in \mathrm{\frac{U_{TS}}{H_{TS}}}} \,\mathcal{F}_{\mathrm{TS}}(p).
     \end{equation}
  \item[3)] Since both the full manifold $\mathrm{\frac{U}{H}}$ and the Tits-Satake submanifold $\mathrm{\frac{U_{TS}}{H_{TS}}}$ are
      diffeomorphic, as differentiable manifolds, to their respective solvable groups $\mathcal{S}_{\mathrm{U/H}}$ and
      $\mathcal{S}_{\mathrm{\frac{U_{TS}}{H_{TS}}}}$, we can inquire what is the relation between the former and the latter solvable groups. We
      find that the latter is the quotient of the former with respect to a normal subgroup $\mathcal{NS}_{fib}$, as we have explained in eq.~\eqref{crepuscolodue}.
\end{description}
The way of looking at the Tits-Satake fibers as  elements of a normal subgroup (picture 3) is not contradictory with, rather it is complementary
to, the other representation of the same in the $\mathrm{G_{subpaint}}$ language (picture 1). The last viewpoint is   handier in constructing
functionals to be extremized, in order to map data points to the manifold. Conversely the vie point in terms of normal solvable subgroups is
appropriate while discussing discretization of the space. Since the space is a group (the solvable group $\mathcal{S}$), the only meaningful way to
perform a discretization is by means of discrete subgroups of the solvable group (see section \ref{generaloneSOpq}).
\par Furthermore, let us add that while the  present article was paper for publication, research on the application of the PGTS theory
to the formulation of the new paradigm in the construction of neural networks advanced and the second view point (picture 2) happened to be very
much handier  as a basic ingrediens of the new paradigm. For this we refer the reader to the forthcoming paper \cite{naviga}.
\section{The two K\"ahlerian building blocks, $\mathbb{H}_2$ and $\mathbb{SH}_2$}
In this section we study in detail the two K\"ahlerian building blocks of the various constructions outlined in section \ref{lollobrigida}.
\subsection{The Poincar\'e-Lobachevsky Hyperbolic plane $\mathbb{H}_2$}
In this subsection we collect all the relevant items and necessary formulas relative to the Tits-Satake submanifold of the  $r=1$ series of
symmetric spaces $\mathcal{M}_\mathfrak{d}^{[r,s]}$, namely,
\begin{equation}\label{rubinetto}
 \mathbb{H}_2 \, \equiv \,  \mathcal{M}_{\mathrm{TS}}^{[1,s]} \, = \, \frac{\mathrm{SO}(1,2)}{\mathrm{SO}(2)}.
\end{equation}
\subsubsection{The various descriptions of the $r=1$ Tits-Satake submanifold} The main special feature of the $r=1$ case is related to the
fact that the group  $\mathrm{ U_{TS}}$ may be looked at as one of the isomorphic groups $\mathrm{SL(2,\mathbb{R})}$,
$\mathrm{Sp(2, \mathbb{R})}$    or $\mathrm{SU(1,1)}$, while the group $\mathrm{SO(1,2)}$ is locally isomorphic to them.
The presentation in terms of $\mathrm{SL(2,\mathbb{R})}$ versus  the one in  terms of $\mathrm{SO(1,2)}$ corresponds to the use of the $2$-dimensional
spinor representation versus the $3$-dimensional vector representation \footnote{As it is thoroughly discussed in 
\cite{tassellandum} the use of the spinor versus the vector representation of the manifold has relevant consequences at the level of the harmonic
expansion of functions on the manifold (see section \ref{micompiacio}),  hence of  the neural networks algorithms utilized to describe for instance
pictures or other types of data mapped to hypersurfaces inside the hyperbolic symmetric space.}. On the other hand
the $\mathrm{SL(2,\mathbb{R})}$   and $\mathrm{SU(1,1)}$
presentation of
the Tits-Satake group are related by a Cayley transformation and the Tits-Satake submanifold can be alternatively viewed in either one of the two
standard geometric models:
\begin{eqnarray}\label{gracile}
\mathbb{H}_2 & = &  \frac{\mathrm{SL(2,\mathbb{R})}}{\mathrm{SO(2)}} \quad = \quad
\left\{ z \in \mathbb{C}\,|\, \mathrm{Im} z >0 \right\} \nonumber\\
\mathbb{H}_2 & = &  \frac{\mathrm{SU(1,1)}}{\mathrm{U(1)}} \quad
=  \quad \left\{ \omega \in \mathbb{C}\,|\, |\omega|^2 < 1 \right\}
\end{eqnarray}
In order to emphasize the general features of the Tits-Satake projection we   start from the $\mathrm{SO(1,2)}$ description of the
Tits-Satake subgroup $\mathrm{ G_{TS}}$.
\subsubsection{Special algebraic features of the $r=1$ case} algebraically the $r=1$ case has the
following distinctive simplifying features:
\begin{description}
  \item[A)] The maximal compact  subgroup $\mathrm{H_c}$  is   not only  semisimple, but simple:
  \begin{equation}\label{boldrini}
    \mathrm{H_c}\, = \, \mathrm{SO}(2+2\,s)  \qquad  \text{(the  $\mathrm{SO}(r)$  factor is missing for $r=1$)}.
  \end{equation}
  \item[B)] In the construction of the Lie algebra $\so(1,1+2s)$, the long roots are absent and we have only one paint group $2s$-multiplet of
      short roots.
  \item[C)] In virtue of the above feature, the solvable group $\mathcal{S} = e^{Solv}$, metrically equivalent to the noncompact symmetric
      manifold $\mathcal{M}_\mathfrak{d}^{[1,s]}$ has a very simple structure: its Lie algebra $Solv$ contains an abelian ideal of dimension
      $2s$ and one Cartan generator.  The generators are
      \begin{equation}\label{solvablealgebra}
      \mathrm{T} \, = \, \left\{  C, S_i \right\} \quad ; \quad \left[C\, , \, S_i\right] \, = \, S_i \quad ;
      \quad \left[ S_i \, , \, S_j \right] \, = \,0
      \end{equation}
  Correspondingly, the solvable group element parameterizing the symmetric space is the one given in eq.~\eqref{Lr1}.
  \item[D)] The Grassmannian leaf of eq.~(\ref{crinolina}) reduces to  the  $2s-1$ dimensional sphere:
 \begin{equation}\label{sferotta}
   \mathcal{F}_{TS}\, = \, \frac{\mathrm{SO}(2s)}{\mathrm{SO}(2s-1)} \, \sim \,\mathbb{S}^{2s-1}
   \end{equation}
   This geometric characterization of the Grasmmanian leaves is what justifies the terminology    ``Poincar\'e Balls.''
\end{description}
\subsubsection{The Tits-Satake projection and the spinor representation} The spinor representation of the full group $\mathrm{SO}(1,1+2s)$
certainly exists for any value of $s$, its dimension being much larger than that of the ($2+2s$)-dimensional defining  vector representation.
The Tits-Satake subgroup $\mathrm{U}_{\mathrm{TS}} \, = \,  \mathrm{SO(1,2)}$ benefits instead from the exceptional low-rank Lie algebra automorphisms
and its spinor representation has dimension $2$. Indeed, as already recalled, the spinor group which provides a double covering of the identity component of
$\mathrm{SO(1,2)}$ is the maximally split simple group $\mathrm{SL(2,\mathbb{R})}$. It is very simple and quite important to recall the precise
form of this double covering map.
\par
As usual we start from the construction of the $\mathrm{SO(1,2)}$ gamma-matrices in a base adapted to the chosen basis of the mother group
$\mathrm{SO}(1,1+2s)$ where the solvable subgroup elements are upper triangular as displayed in eq.~(\ref{Lr1}). The solvable coordinate vector is
 very simple in the $r=1$ case and the Tits-Satake projection is equally simple:
\begin{equation}\label{carnevale}
  \pi_{TS}\left(\boldsymbol{\Upsilon}\right) \, = \, \left\{\Upsilon_1, \Upsilon_{2,1},0,\dots ,0  \right\}
\end{equation}
Correspondingly, starting from equation (\ref{Lr1}), we have:
\begin{equation}\label{carnegomme}
  \forall \, \mathbb{L}(\boldsymbol{\Upsilon}) \in e^{Solv_{\so(1,1+2s)}} \quad : \quad \pi_{TS}\left(\mathbb{L}(\boldsymbol{\Upsilon})\right) \, = \,
  \, \left(
  \begin{array}{c|c|c}
  e^{\Upsilon_1} & e^{\Upsilon_1} \, \frac{1}{\sqrt{2}}\, \Upsilon_{2,1} & -\, e^{\Upsilon_1} \, \frac{1}{4}\,  \Upsilon_{2,1}^2\\
  \hline
  0 & 1 & - \frac{1}{2}\, \Upsilon_{2,1} \\
  \hline
  0& 0 & e^{-\,\Upsilon_1} \\
  \end{array}
  \right) \, \equiv \, \mathbb{L}(\widehat{\boldsymbol{\Upsilon}}) \in \, e^{Solv_{\so(1,2)}}
\end{equation}
where we have denoted by $Solv_{\so(1,1+2s)}$ the solvable Lie algebra associated with the non-maximally split full manifold and $Solv_{\so(1,2)}$
the solvable Lie algebra associated with the Tits-Satake maximally split submanifold (\ref{rubinetto}). Furthermore we have denoted
\begin{equation}\label{attouppo}
  \widehat{\boldsymbol{\Upsilon}} \, =\, \left\{\Upsilon_{1}, \, \Upsilon_{2,1} \right\}
\end{equation}
In the upper triangular basis inherited from the mother group $\mathrm{SO}(1,1+2s)$, the  invariant metric of the Tits-Satake subgroup
$\mathrm{SO(1,2)}$ is
\begin{equation}\label{etat12}
  \eta_t \, = \,\left(
                  \begin{array}{ccc}
                    0 & 0 & 1 \\
                    0 & 1 & 0 \\
                    1 & 0  & 0 \\
                  \end{array}.
                \right)
\end{equation}
A triplet of  gamma matrices   satisfies the   Clifford algebra relations
\begin{equation}\label{gammimatri}
  \left\{\gamma_{i} \, , \, \gamma_{j}\right\} \, = \, \eta_{t|ij} \, \mathbf{1}_{2\times 2}\quad ; \quad i,j\, = \, 1,2,3.
\end{equation}
An explicit representation is
\begin{equation}\label{gammuti}
  \gamma_1 \, = \, \left(
                     \begin{array}{cc}
                       0 & 0 \\
                       1 & 0 \\
                     \end{array}
                   \right) \quad ; \quad \gamma_2 \, = \, \frac{1}{\sqrt{2}} \, \left(
                     \begin{array}{cc}
                       -1 & 0 \\
                       0 & 1 \\
                     \end{array}
                   \right) \quad ; \quad \gamma_3 \, = \,\left(
                     \begin{array}{cc}
                       0 & 1 \\
                       0 & 0 \\
                     \end{array}
                   \right).
\end{equation}
Then the map from the group $\mathrm{SL(2,\mathbb{R})}$, whose elements are $2\times 2$ real valued matrices of determinant one:
\begin{equation}\label{grunwald}
 \mathrm{SL(2,\mathbb{R})}\, \ni\, \mathfrak{m} \, = \, \left(
 \begin{array}{cc}
 \alpha & \beta \\
 \gamma & \delta \\
 \end{array}
 \right) \quad ; \quad \alpha \, \delta \, - \, \gamma \, \beta \, = \, 1
\end{equation}
to elements of the group $\mathrm{SO(1,2)}$, namely, to $3\times 3 $ real valued matrices $\mathcal{O}$ satisfying the constraint
$\mathcal{O}^T\,\eta_t \, \mathcal{O} \, = \, \eta_t$ is provided by the following quadratic relation:
\begin{equation}\label{coriaceo}
\mathcal{O}_{i}^{\phantom{i}j}\left[\mathfrak{m}\right] \, = \, \text{Tr}\left(\mathfrak{m}^{-1} \, \gamma_{i} \, \mathfrak{m} \,
\gamma_{k} \right) \, \eta^{kj}
\end{equation}
Inside the group $\mathrm{SL(2,\mathbb{R})}$ we have a 2-dimensional solvable subgroup (a Borel subgroup) $\mathfrak{T}\subset\mathrm{SL(2,\mathbb{R})}$ formed by
the upper triangular matrices whose image in $\mathrm{SO(1,2)}$ is also upper triangular
\begin{equation}\label{Tsubbo} \mathfrak{T}\, \ni \, \mathfrak{t} \, = \, \left( \begin{array}{cc} \alpha & \beta \\ 0 & \alpha^{-1} \\
\end{array} \right) \quad ; \quad \mathcal{O}_{i}^{\phantom{i}j}\left[\mathfrak{t}\right] \, = \, \text{upper triangular}
\end{equation}
Actually we find that:
\begin{equation}\label{curiente}
  \mathbb{L}(\widehat{\boldsymbol{\Upsilon}}) \, = \,  \mathcal{O}\left[\mathfrak{t}\right] \quad \text{if and only if }
  \quad \alpha = \pm \,e^{\frac{\Upsilon_1}{2}} \,\, , \,\, \beta \, = \,
  \pm \frac{1}{2} \, e^{\frac{\Upsilon_1}{2}} \,\Upsilon_{2,1} .
\end{equation}
Recalling that the upper complex plane (the first choice in eq.~(\ref{gracile})) is invariant under the action of $\mathrm{SL(2,\mathbb{R})}$, we
see that the solvable subgroup $\mathfrak{T}$ has a simple transitive action on it. Indeed under the fractional linear transformation
\begin{equation}\label{perdapiedi}
 \forall z \in \mathbb{C} \, ,\, \forall \mathfrak{g} = \left(
 \begin{array}{cc}
  \alpha & \beta \\
 \gamma & \delta \\
 \end{array}
 \right)\in \mathrm{PSL(2,\mathbb{R})} \quad : \quad \mathfrak{g}(z) \, = \, \frac{\alpha \, z + \beta}{\gamma\, z + \delta}
\end{equation}
 the upper
complex plane is just the orbit under the solvable subgroup  $\mathcal{S} \subset \mathrm{PSL(2,\mathbb{R})} $ of the special point $i$. Hence we
can set
\begin{equation}\label{carotide}
  z \, = \, i \, \alpha^2 \, + \, \alpha \, \beta \, = \, \frac{1}{2} e^{\Upsilon_1} \, \Upsilon_{2,1}\, +\, i \, e^{\Upsilon_1} \,
  \in \, \left\{z\in \mathbb{C} \, |\, \text{Im} z >0\right \} .
\end{equation}
The representation of the Tits-Satake submanifold as the interior of the unit circle is obtained by means of the standard Cayley transformation
\begin{equation}\label{caullo}
  \omega \, = \, \frac{z\, - \, i}{z \, + \, i} \, = \, \frac{e^{2 \Upsilon _1} \left(\Upsilon _{2,1}^2+4\right)-4}{e^{2 \Upsilon _1}
   \left(\Upsilon _{2,1}^2+4\right)+8 e^{\Upsilon _1}+4}\, -\, i \,\frac{4  e^{\Upsilon
   _1} \Upsilon _{2,1}}{e^{2 \Upsilon _1} \left(\Upsilon _{2,1}^2+4\right)+8
   e^{\Upsilon _1}+4}.
\end{equation}
Eq.~\eqref{caullo} is the explicit final form of the Tits-Satake projection that maps any point $p\in \mathcal{M}_\mathfrak{d}^{[1,s]}$ to a point
of the Tits-Satake submanifold $\mathcal{M}_{TS}^{[1,2]}$ represented as the interior of the unit-circle in the complex plane $\omega$. The point
of the full manifold is identified by the vector of solvable coordinates
\begin{equation}\label{ipsilon}
  \boldsymbol{\Upsilon} \, = \, \{\Upsilon_1, \Upsilon_{2,i} \} \quad ; \quad i=1,\dots, 2s.
\end{equation}
 One disregards all the coordinates
$\Upsilon_{2,i>1}$ and out of $\Upsilon _1$ and $\Upsilon _2$ constructs the complex number $\omega$. This can be done for a discrete number of
points $p$ or for entire submanifolds or lines. In particular the Tits-Satake projection can be done for a geodesic. Suppose that $\boldsymbol{
\Upsilon}_{geo}(t)$ represents a geodesic of the full manifold; by means of the  Tits-Satake projection we can see the image of the ambient space
geodesic inside the unit-circle:
\begin{equation}\label{ptsgeo}
  \omega_g(t) \, = \, \pi_{TS} \left(\boldsymbol{
\Upsilon}_{geo}(t)\right).
\end{equation}
The question arises whether the image of a geodesic of the full space inside the Tits-Satake submanifold is a geodesic of the latter. The answer
is no! What  is true is rather the following:
\begin{description}
  \item[A)] Any geodesic of the Tits-Satake submanifold $\mathcal{M}_{TS}^{[1,2]}$ \emph{is also a geodesic} with respect to the metric of the
      full manifold $\mathcal{M}_{\mathfrak{d}}^{[1,1+2s]}$.
  \item[B)] Generically the Tits-Satake projection of a geodesic of the full manifold $\mathcal{M}_{\mathfrak{d}}^{[1,1+2s]}$ \emph{is not a
      geodesic} with respect to the metric of $\mathcal{M}_{TS}^{[1,2]}$, unless it is fully contained in $\mathcal{M}_{TS}^{[1,2]}$.
\end{description}
In the next subsection we  illustrate the above facts by means of the simplest example $r=s=1$.
\subsubsection{Illustration in the simplest case $r=s=1$}\label{lustrascarpe11}
We begin by considering the explicit form of the solvable group element providing the solvable coset representative of the manifold
$\frac{\mathrm{SO(1,1+2)}}{\mathrm{SO(3)}}$:
\begin{equation}\label{r1s1Lfat}
 \mathbb{L}_{\mathfrak{d}|1,s}(\boldsymbol{\Upsilon}) \,=\,\left(
\begin{array}{cccc}
 e^{\Upsilon_1} & \frac{e^{\Upsilon_1 } \Upsilon_ {2,1}}{\sqrt{2}} &
   \frac{e^{\Upsilon_1} \Upsilon_{2,2}}{\sqrt{2}} & -\frac{1}{4} e^{\Upsilon
   (1)} \left(\Upsilon_{2,1}^2+\Upsilon_{2,2}^2\right) \\
 0 & 1 & 0 & -\frac{\Upsilon_{2,1}}{\sqrt{2}} \\
 0 & 0 & 1 & -\frac{\Upsilon_{2,2}}{\sqrt{2}} \\
 0 & 0 & 0 & e^{-\Upsilon_1} \\
\end{array}
\right),
\end{equation}
which corresponds to the generic solvable coordinate vector
\begin{equation}\label{upvecgen11}
  \boldsymbol{\Upsilon}\, = \, \left\{\Upsilon_1,\, \Upsilon_{2,1},\, \Upsilon_{2,2} \right\}
\end{equation}
The Tits-Satake projection on the Tits-Satake submanifold defined  is obtained by setting $\Upsilon_{2,2}=0$. For convenience we can label a point
in  $ \mathcal{M}_{TS}^{[1,2]}$ by $\Upsilon_1=w_1,\Upsilon_{2,1}=w_2$, yielding the point (\ref{carotide}) in the upper complex plane. In the
intrinsic parameterization of the full manifold $\mathcal{M}^{[1,1]}_\mathfrak{d}$ by means of the symmetric matrix
$\mathcal{M}(\boldsymbol{\Upsilon}) \equiv \mathbb{L}(\boldsymbol{\Upsilon}) \,\mathbb{L}^T(\boldsymbol{\Upsilon})$, the points of the Tits-Satake
submanifold are in one-to-one correspondence with the matrices
\begin{equation}\label{corlanino}
  \mathcal{M}_{TS}(w_1,w_2) \, = \, \left(
\begin{array}{cccc}
 \frac{1}{16} e^{2 w_1} \left(w_2^2+4\right){}^2 & \frac{e^{w_1} w_2
   \left(w_2^2+4\right)}{4 \sqrt{2}} & 0 & -\frac{w_2^2}{4} \\
 \frac{e^{w_1} w_2 \left(w_2^2+4\right)}{4 \sqrt{2}} & \frac{1}{2}
   \left(w_2^2+2\right) & 0 & -\frac{e^{-w_1} w_2}{\sqrt{2}} \\
 0 & 0 & 1 & 0 \\
 -\frac{w_2^2}{4} & -\frac{e^{-w_1} w_2}{\sqrt{2}} & 0 & e^{-2 w_1} \\
\end{array}
\right).
\end{equation}
The geometry of the Tits-Satake projection and the foliation of the full manifold into Grassmannian leaves can now be understood by recalling the
basic principle of the Euler parameterization. The entire space can be organized into the union of the Grassmaniann leaves  associated with  the
points of the Tits-Satake submanifold,  that is,  the Hyperbolic plane.  Each Grassmannian leaf is the orbit of a point $p\in \mathcal{M}_{TS}$  by
 the action of  $\exp\left[\mathbb{F}_{\mathrm{TS}}\right]$ (see eqs.~(\ref{corsarobianco}-\ref{corsaronero})). In the $r=1,s=1$ case the
subspace $\mathbb{F}_{\mathrm{TS}}$ is one-dimensional and   is spanned by the compact generator
\begin{equation}\label{fJ}
  J_{F} \, = \, \left(
\begin{array}{cccc}
 0 & 0 & \frac{1}{\sqrt{2}} & 0 \\
 0 & 0 & 0 & 0 \\
 -\frac{1}{\sqrt{2}} & 0 & 0 & -\frac{1}{\sqrt{2}} \\
 0 & 0 & \frac{1}{\sqrt{2}} & 0 \\
\end{array}
\right).
\end{equation}
Hence we have
\begin{equation}\label{Lattus}
  \mathcal{O}_F(\phi) \, = \, \exp\left[\phi \, J_F \right]\, = \, \left(
\begin{array}{cccc}
 \cos ^2\left[\frac{\phi }{2}\right] & 0 & \frac{\sin [\phi ]}{\sqrt{2}} &
   \frac{1}{2} [\cos [\phi ]-1] \\
 0 & 1 & 0 & 0 \\
 -\frac{\sin [\phi ]}{\sqrt{2}} & 0 & \cos [\phi ] & -\frac{\sin [\phi
   ]}{\sqrt{2}} \\
 \frac{1}{2} [\cos [\phi ]-1] & 0 & \frac{\sin [\phi ]}{\sqrt{2}} & \cos
   ^2\left[\frac{\phi }{2}\right] \\
\end{array}
\right)
\end{equation}
and the points of the full manifold  $\mathcal{M}^{[1,1]}_\mathfrak{d}$ are in  a one-to-one correspondence with the symmetric matrices
\begin{equation}\label{carnegietrust}
  \mathcal{M}\left[w_1,w_2,\phi\right] \, = \, \mathcal{O}_F(\phi) \, \mathcal{M}_{TS}(w_1,w_2) \, \mathcal{O}^T_F(\phi).
\end{equation}
In this way one explicitly realizes the foliation structure of the full manifold. With  each point $p\in \mathcal{M}_{\mathrm{TS}}^{[1,2]}$ of the
Tits-Satake submanifold --- singled out by the coordinates of a point in the upper complex plane or in the unit circle, eventually $w_{1,2}$ ---  one
associates an entire sphere $\mathbb{S}^{2s-1}$, in the case $s=1$ a circle $\mathbb{S}^{1}$. Explicitly utilizing the Cholevsky-Crout algorithm
in order to express the symmetric matrix $\mathcal{M}\left[w_1,w_2,\phi\right]$ in terms of an upper triangular matrix
$\mathbb{L}(\boldsymbol{\Upsilon})$ and then the inverse of the $\Sigma$-map to work out $\boldsymbol{\Upsilon}$ one finds
\begin{eqnarray}\label{cromacino}
  \boldsymbol{\Upsilon}\left(w_1,w_2,\phi\right)& = & \left\{\frac{1}{2} \log \left(\frac{1}{\left(e^{2 w_1} \left(w_2^2+4\right) (\cos (\phi
   )-1)-4 (\cos (\phi )+1)\right){}^2}\right)+w_1+\log (8),w_2,\right.\nonumber\\
   &&\left.-\frac{1}{4}
   \left(e^{w_1} w_2^2+8 \sinh \left(w_1\right)\right) \sin (\phi )\right\}.
\end{eqnarray}
Calling $x,y$ the real and imaginary parts of the upper complex plane $z$-point in the map (\ref{carotide}), we finally get the following
parameterization of the  points of the full-manifold:
\begin{eqnarray}\label{puntacchia}
 && \left(x,y,w_{2,2}\right) \, = \, \nonumber\\
 &&\left\{\frac{2 \sigma }{-\left(\rho ^2+\sigma ^2-1\right) \cos (\phi )+\rho
   ^2+\sigma ^2+1},\frac{2 \rho }{-\left(\rho ^2+\sigma ^2-1\right) \cos (\phi
   )+\rho ^2+\sigma ^2+1},-\frac{\left(\rho ^2+\sigma ^2-1\right) \sin (\phi
   )}{\rho }\right\},\nonumber
\end{eqnarray}
where we have redefined
\begin{equation}\label{cristallus}
  w_1 \, = \, \log[\rho] \quad ; \quad w_2 \, = \, 2 \, \frac{\sigma}{\rho}
\end{equation}
\begin{figure}[!hbt]
\begin{center}
\iffigs
\includegraphics[width=65mm]{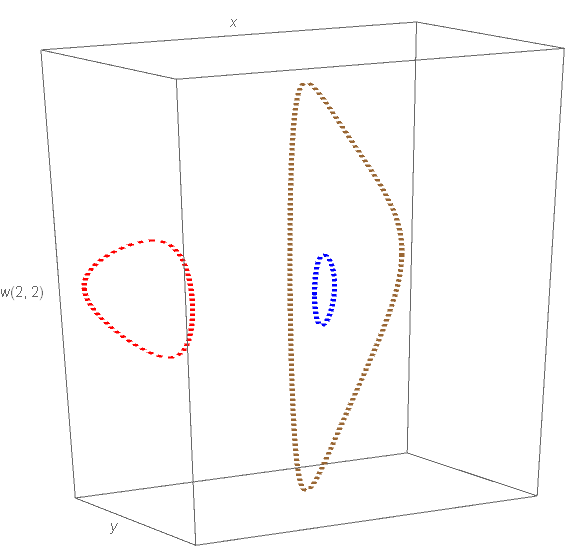}
 \includegraphics[width=65mm]{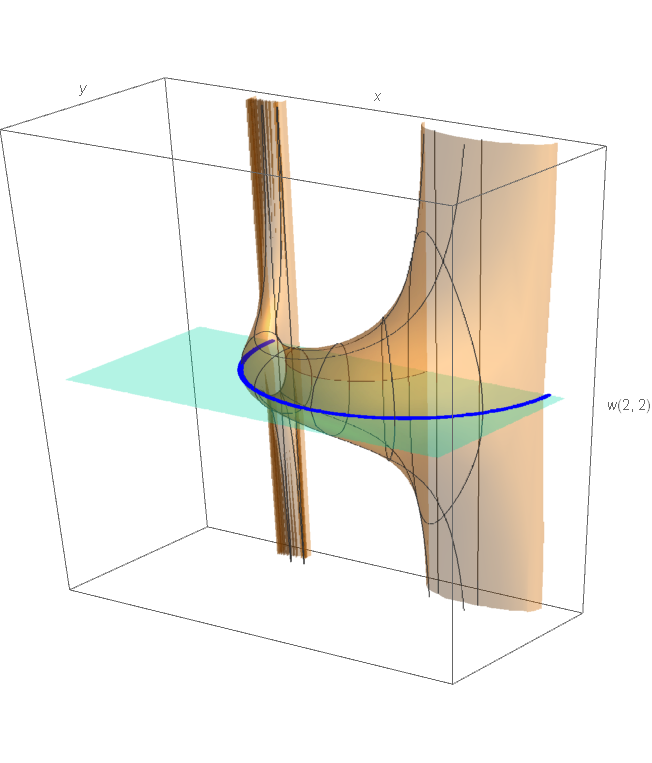}
\else
\end{center}
 \fi
\caption{The foliation structure encoded in eq.~(\ref{puntacchia}).
On the left there are  three circles respectively associated with three different points of the upper complex plain,
 namely $z=-\ft 12+ i \ft 32$, $z= \ft 23 + i \ft 12$ and
$z= 2+ i\ft 32$. On the right we display the surface generated by taking
all the circles associated with each of the points along a line in the upper complex plane
that is the bluish horizontal plane cutting the surface. The chosen line is   a geodesic of the Tits-Satake submanifold,
  a circle with center on the real line. This is
 the thick blue line described by the formula $z[\psi]\, = \, -\,\ft 13 \, + \, \ft  32 \, e^{i \psi}$.
\label{circolino1}}
 \iffigs
 \hskip 1.5cm \unitlength=1.1mm
 \end{center}
  \fi
\end{figure}

In Figure \ref{circolino1} we illustrate the foliation of the considered manifold in two ways: with three examples of topological circles associated
with three different points of the upper complex plane model of the hyperbolic space, and with the surface generated by taking the circles
associated with each of the points in a line (actually a geodesic) of  the same upper complex plane.
\par
\begin{figure}[!hbt]
\begin{center}
\iffigs
\includegraphics[width=55mm]{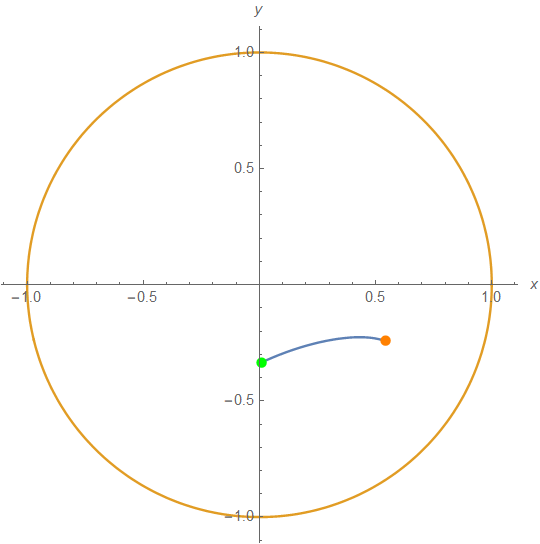}
\includegraphics[width=55mm]{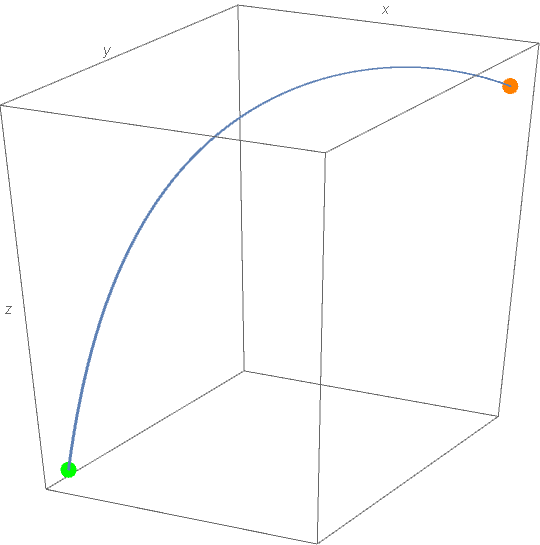}
\includegraphics[width=55mm]{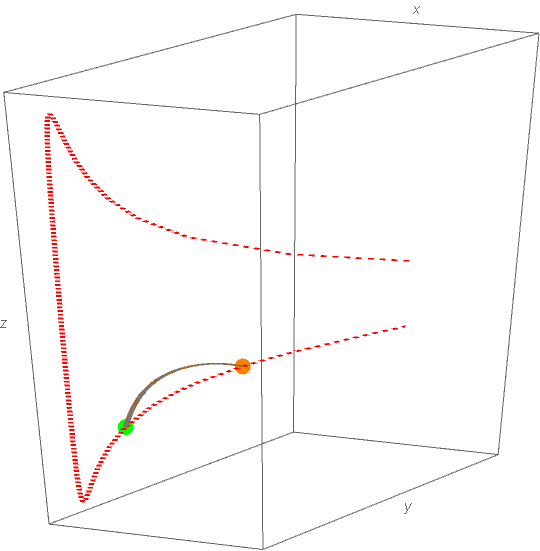}
\else
\end{center}
 \fi
\caption{The circular leaf associated with the point $w_1=2, \, w_2 =3/2$ of the Tits-Satake submanifold.
The leaf is shown in the three dimensional space $\{x,y,z=w_{2,2}\}$ where $x,y$ are the real and imaginary part of the point
$\zeta(w_1,w_2)$ in the upper complex plane.
It corresponds to the close curve, in dashed red   in the   picture on the right.
Two points on the leaf, respectively corresponding to $\phi=\pi/10$ and $\phi=\pi/5$,
are joined by a geodesic, in solid brown, that, as one sees, does not lie on the leaf.
The Tits-Satake projection of such geodesic is displayed in the first picture,
while its entire three dimensional shape is shown in the second picture.
\label{r1s1esempiofibra}}
 \iffigs
 \hskip 1.5cm \unitlength=1.1mm
\end{center}
  \fi
\end{figure}
\begin{figure}[!hbt]
\begin{center}
\iffigs
\includegraphics[width=55mm]{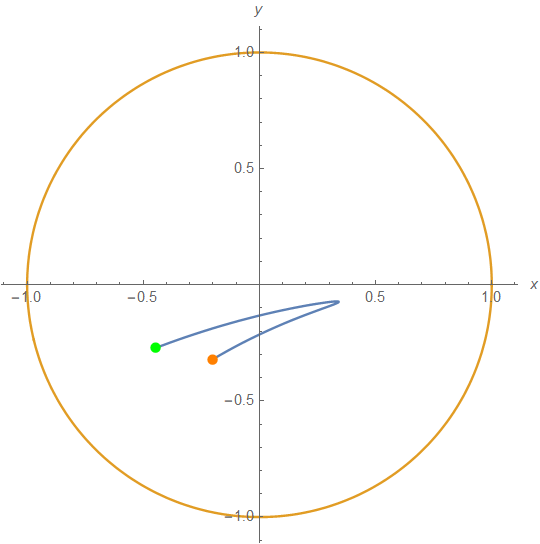}
\includegraphics[width=55mm]{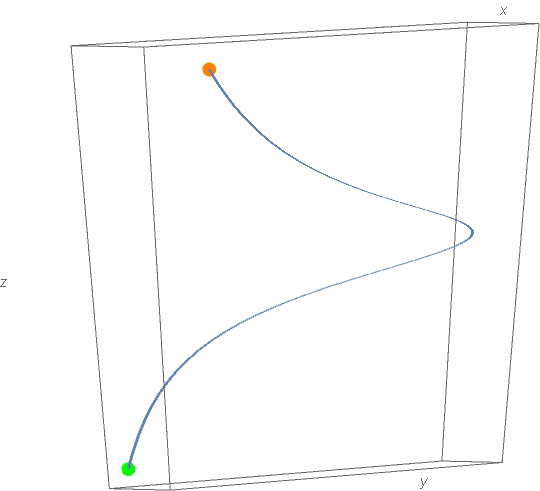}
\includegraphics[width=55mm]{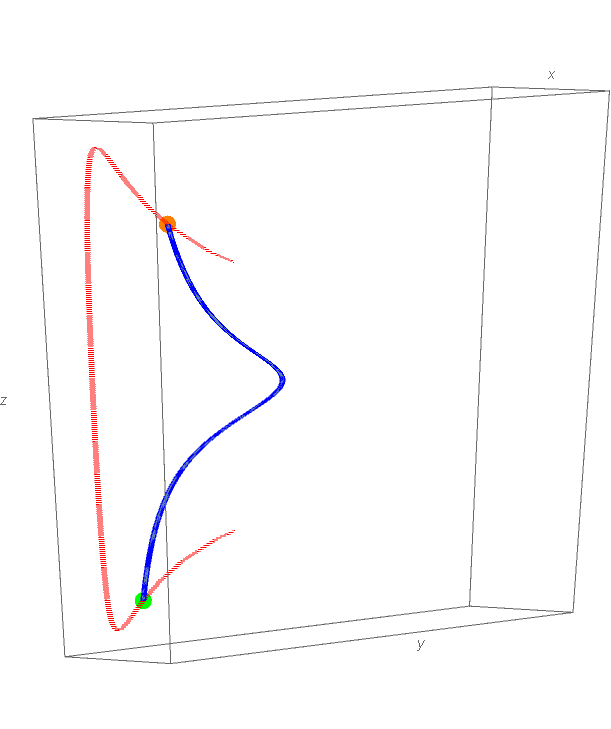}
\else
\end{center}
 \fi
\caption{The circular leaf associated with the point $w_1=2, \, w_2 =3/2$ of the Tits-Satake submanifold.
The leaf is shown in the three dimensional space $\{x,y,z=w_{2,2}\}$ where $x,y$ are the real and imaginary part of the point
$\zeta(w_1,w_2)$ in the upper complex plane.
It corresponds to the close curve,  in dashed   red color in the third picture.
Two points on the leaf, respectively corresponding to $\phi=-\pi/4$ and $\phi=\pi/3$, are joined by a geodesic, in  solid blue,
that, as one sees, does not lie on the leaf. The Tits-Satake projection (in its unit disk model) of
this  geodesic is displayed in the first picture, while its entire three dimensional shape is shown in the second picture.
\label{r1s1esempiofib2}}
\iffigs
\hskip 1.5cm \unitlength=1.1mm
\end{center}
\fi
\end{figure}

Next let us consider Figure \ref{r1s1esempiofibra}: it displays the geodesic that joins two points belonging to the same leaf,   the one
singled out by the parameters $w_1=2, \, w_2 =3/2$.
The important point to be stressed is that the geodesic joining two points of the same leaf does not lie on the leaf (see the caption of
Figure \ref{r1s1esempiofibra} for more details). The same phenomenon is  emphasized even more strongly in  Figure \ref{r1s1esempiofib2}.
Furthermore the example in the quoted figure emphasizes another fundamental aspect. As it is visually clear from the picture, the projection of a
geodesic of the ambient space onto the Tits-Satake submanifold, namely the interior of the unit disk, is not a geodesic of the TS manifold
(i.e., an arc of a circle with center at the boundary of the unit disk). Indeed  the projected geodesic can be a much more complicated
curve, even having nodes,  as it happens in the present example.
\par
\begin{figure}[!hbt]
\begin{center}
\iffigs
\includegraphics[width=55mm]{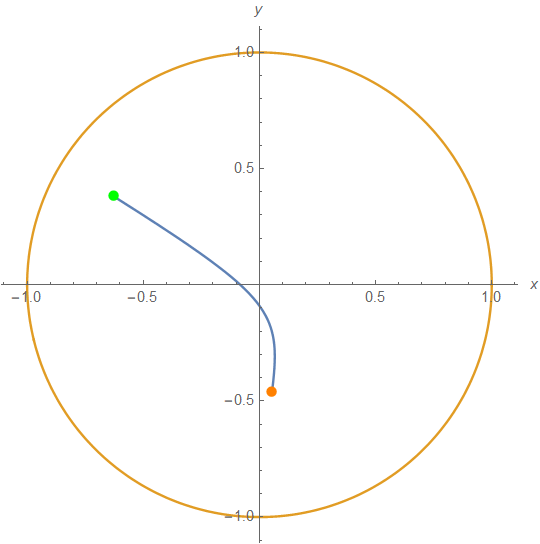}
\includegraphics[width=55mm]{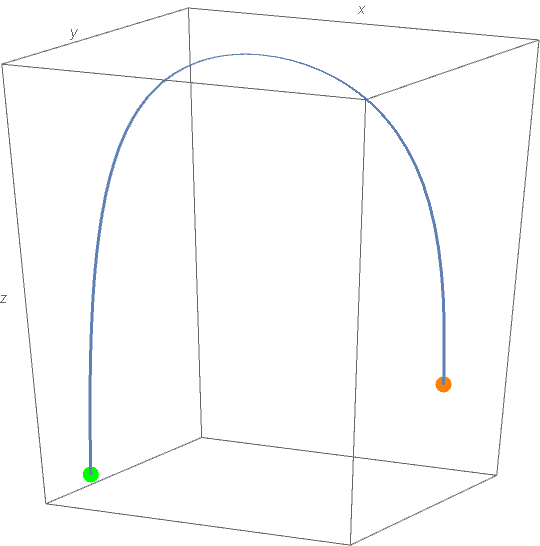}
\includegraphics[width=55mm]{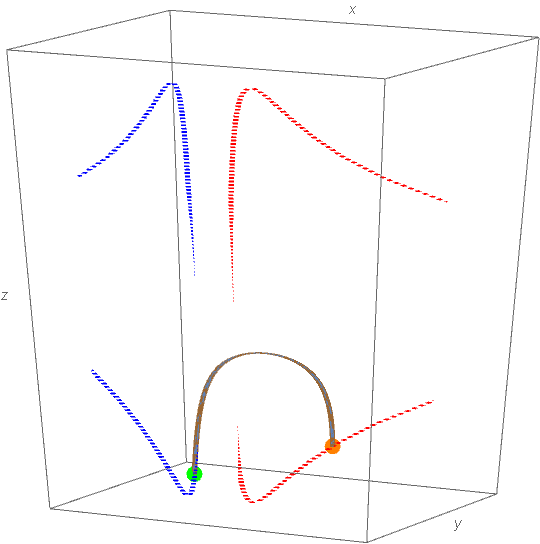}
\else
\end{center}
 \fi
\caption{A geodesic that joins two points located on different leaves.
The first circular leaf is associated with the point $w_1=6/5, \, w_2 =7/3$ of the Tits-Satake submanifold,
while the second leaf is associated with the  point $w_1=3/4, \, w_2 =-10/3$ of the same manifold.
The two leaves are shown in the third picture  in the three dimensional space $\{x,y,z=w_{2,2}\}$, where $x,y$ are the real and imaginary part of
the point $\zeta(w_1,w_2)$ in the upper complex plane, and   correspond to the close curves,
in  dashed  red and blue, respectively. The point on the first leaf, corresponding to $\phi=\pi/4$,
is joined to the point corresponding to $\phi=2\pi/3$ on the second leaf by a geodesic, in solid brown.
The projection of such a geodesic in the Tits-Satake submanifold is displayed in the first picture,
while the second picture displays its development in the full three dimensional space. \label{saltafoglia}}
\iffigs
\hskip 1.5cm \unitlength=1.1mm
\end{center}
\fi
\end{figure}

In Figure \ref{saltafoglia} we display instead the case of a geodesic joining two points on two different leaves.
\subsubsection{Squared norm and Grassmannian leaves}
Last but not least we reconsider, in view of the Grassmannian leaf picture, the squared norm on the solvable group. To this effect we focus on the
\emph{Grassmannian leaf} associated with a point  of the Tits-Satake submanifold. Utilizing the $s=1$ case in order to obtain simple results that
we later generalize by paint invariance, we consider the orbit under the $\mathrm{U(1)}$ group generated by the unique fiber generator that was
already  presented in eq.~\eqref{cromacino}. We compute the squared norm   $\mathrm{N^2} (\boldsymbol{\Upsilon}[w1,w2,\phi])$,
finding
 \begin{equation}\label{grassouomonormale}
  \mathrm{N^2} (\boldsymbol{\Upsilon}[w1,w2,\phi])\, = \,
  \left( {\rm arccosh}\left[  \frac{1}{32} \left(16 e^{-2 w_1}+8 w_2^2+e^{2 w_1} \left(w_2^2+4\right){}^2\right)\right]\right)^2,
\end{equation}
which is independent from the angle $\phi$. So all points in the Grassmann leaf have the same norm as their unique representative point lying
on the Tits-Satake submanifold.
This  is a completely general result that holds true for $r=1$ and all values of $s$, but it is also true for higher values of the
noncompact rank $r$ and for all corresponding various of $s$.
\par
This displays the true meaning of the description of the full manifold as provided by eq.~(\ref{paleocappo}). Each point $p\in
\mathrm{\frac{U_{TS}}{H_{TS}}}$ of the Tits-Satake submanifold has a well defined and finite squared norm $N^2(p) \, = \,R^2$. There is however in
the full manifold  $\mathrm{\frac{U}{H}}$ and entire compact submanifold, diffeomorphic to the Grassmannian (\ref{crinolina}), whose points have
the same squared norm.
\par
The use of this property in devising optimized procedures of mapping data to hyperbolic manifolds is one of the several open perspectives for
applications of   PGTS theory to Data Science.
%%%%%%%%%%%%%%%%%%%%%%%%%%%%%%%%%%%%%%%%%%%%%%%%%%%%%%%%%%%%%%%%%%%%%%%%%%%
\subsection{The building block $\mathbb{SH}_2$}
\label{pargnocco} In this section we focus on the case $r=2$ of the series of symmetric manifolds $\frac
{\mathrm{SO(r,r+2\,s)}}{\mathrm{SO(r)\times SO(r+2s)}}$ in a completely synoptic setup with respect to our previous treatment of the case $r=1$.
The motivation for this synopsis is twofold:
\begin{enumerate}
  \item On  the one hand we want to stress that the $r=1$ case is completely aligned with all the subsequent $r>1$ ones, and that the Tits-Satake
      projection, which went unnoticed and unexploited by the authors in \cite{francesi1,francesi2,francesi3}, is actually the conceptual
      backbone for all the members of the considered series of manifolds.
  \item On the other hand we want  to emphasize that the $r=2$ and $r=1$ cases are twins inside the entire series, as their respective Tits-Satake submanifolds are  just the first and the second instance of a \emph{Siegel upper complex plane}, which is the appropriate
      generalization of the Lobachevsky-Poincar\'e hyperbolic plane. Instead, for values $r>2$, the Tits-Satake submanifold, which, by definition,
      is always  a maximally split symmetric space, is not a further instance of a Siegel upper complex plane. Indeed the appearance of the
      first two Siegel planes is strictly linked with the low-rank sporadic isomorphisms of simple Lie algebras.
  \item In spite of  its exceptionality, the Siegel plane twin relation between the $r=2$ and $r=1$ case has nontrivial consequences for the
      generalization to the higher noncompact rank $r=3,4$ via the $c$-map relations that we have emphasized in section \ref{lollobrigida}.  In
      particular, as we stressed there, the standard hyperbolic plane and the Siegel upper plane of order 2 are the K\"ahler
      building blocks for all cases $r=1,2,3,4$. This is particularly relevant at the level of the search for discrete subgroups; the discrete
      subgroups of the two quaternionic cases $3$ and $4$ will be indeed built  starting from tensor products of the discrete subgroups of
      $\mathrm{SL(2,\mathbb{R})}$ and $\mathrm{Sp(4,\mathbb{R})}$, as in a Chinese box game.
\end{enumerate}
In view of what we said above, we propose to call the $r=2$ case \emph{Siegel lumps} opposed to the \emph{Poincar\'e balls} of the case $r=1$.
The word \emph{lump} opposed to \emph{ball} recalls the fact that, in the $r=2$ case, Tits-Satake fibers are not related to spherical leaves, but
rather with leaves that have the shape of a more general Grassmannian.
\par
In subsection \ref{siegelpiatto} we discuss in depth the structure of the Tits-Satake submanifold
\begin{equation}\label{guardacaccia}
  \mathbb{SH}_2 \, \equiv \, \frac{\mathrm{SO(2,3)}}{\mathrm{SO(2)\times SO(3)}} \, \simeq \,
  \frac{\mathrm{Sp(4,\mathbb{R})}}{\mathrm{U(1)\times SU(2)}}
\end{equation}
which is indeed the natural generalization of the $r=1$ Tits-Satake manifold, \emph{i.e.} of the hyperbolic plane, being the next instance of a
Siegel upper plane of which the hyperbolic plane is the first instance.
\subsubsection{The Siegel upper plane and the Tits-Satake projection in
the \pmb{$ r=2$} case}\label{siegelpiatto} According with the statement made at the end of section \ref{pargnocco}, we come now to the complete study of
the Tits-Satake manifold (\ref{guardacaccia}). The equivalence shown in eq.~(\ref{guardacaccia}) is due to the equivalence of the complex Lie
algebra $\mathfrak{b}_2$ with the complex Lie algebra $\mathfrak{c}_2$, which is one of the few low-rank sporadic isomorphisms:
\begin{equation}\label{critollo}
  \so(5,\mathbb{C}) \, \simeq \, \spalg (4,\mathbb{C}).
\end{equation}
At the level of maximal compact real section eq.~(\ref{critollo}) implies
\begin{equation}\label{critollocom}
  \so(5) \, \simeq \, \usp (4),
\end{equation}
while at the level of maximally split real section it yields
\begin{equation}\label{critollosplit}
  \so(2,3) \, \simeq \, \spalg (4,\mathbb{R}).
\end{equation}
Eq.~\eqref{critollosplit} is the one relevant to our purposes. The isomorphism is explicitly realized by means of the spinor representation of the
group $\mathrm{SO(2,3)}$, so that the first step in its construction goes through the set up of a convenient set of $4\times 4$ gamma matrices well
adapted to the chosen $\eta_t$:
\begin{equation}\label{eta5tria}
\eta_t \,=\, \left(
\begin{array}{ccccc}
 0 & 0 & 0 & 0 & 1 \\
 0 & 0 & 0 & 1 & 0 \\
 0 & 0 & 1 & 0 & 0 \\
 0 & 1 & 0 & 0 & 0 \\
 1 & 0 & 0 & 0 & 0 \\
\end{array}
\right).
\end{equation}
\subsubsection{Gamma matrices and the spinor algebra of $\mathrm{SO(2,3)}$}
We have determined a set of five  $4\times 4$ matrices $\Gamma_i$ that possess the following
properties:
\begin{enumerate}
  \item they satisfy the required Clifford algebra relations
  \begin{equation}\label{clifford5}
    \left[\, \Gamma_i \, , \, \Gamma_j \right] \, = \, 2 \, \eta_{t|ij} \, 1_{4\times 4} ;
  \end{equation}
  \item   their matrix entries are real and they are all traceless:
  \begin{equation}\label{straccione}
    \text{Tr} \, \Gamma_i \, = \, 0 \quad \quad i=1\, \dots, \, 5;
  \end{equation}
  \item they admit a  skew-symmetric charge conjugation matrix
  \begin{equation}\label{chargeconjug}
   \mathrm{C_s} \, = \, \left(
\begin{array}{cc|rr}
 \phantom{-}0 & \phantom{-}0 \phantom{-}  & \phantom{-} 1 & \phantom{-} 0 \\
 \phantom{-}0 & \phantom{-}0 \phantom{-}  & \phantom{-} 0 & \phantom{-} 1 \\
 \hline
 -1 & \phantom{-} 0 \phantom{-}  & \phantom{-} 0 & \phantom{-} 0 \\
 \phantom{-}0 & -1\phantom{-}  & \phantom{-} 0 & \phantom{-}0 \\
\end{array}
\right)
  \end{equation}
  such that
  \begin{equation}\label{antisymrep}
    \mathrm{C_s}\Gamma_i \, = \,- \, \left( \mathrm{C_s}\Gamma_i\right)^T.
  \end{equation}
  namely, the matrices $\mathrm{C_s}\Gamma_i$ are all skew-symmetric.
  \item The spinor generators of $\mathrm{SO(2,3)}$
  \begin{equation}\label{spingenso23}
    J_{ij} \, \equiv \, \ft 14 \, \left[\Gamma_i \, , \, \Gamma_j \right]
  \end{equation}
  multiplied by the charge conjugation matrix are instead symmetric, that is, they satisfy the condition
  \begin{equation}\label{symplcondo}
    \mathrm{C_s} \, J_{ij} \, + \, J_{ij}^T \, \mathrm{C_s} \, = \, \mathbf{0}_{4\times 4}.
  \end{equation}
\end{enumerate}
Eq.~\eqref{symplcondo} is the key   to  the isomorphism. Indeed it states that the $10$ generators of $\mathrm{SO(2,3)}$ are $10$ linearly
independent matrices satisfying the condition to be elements of the $\spalg (4,\mathbb{R})$ Lie algebra with symplectic invariant matrix
$\mathrm{C_s}$ as given in equation (\ref{chargeconjug}). On the other hand, the two conditions (\ref{straccione},\ref{antisymrep}) tell us that
the vector representation is obtained from the symplectic traceless skew-symmetric bispinor representation. As we show below, this is the key
ingredient to realize the explicit map from $5\times 5$ group elements of $\mathrm{SO(2,3)}$ and $4\times 4 $ group elements of $\Sp
(4,\mathbb{R})$.
\par Having clarified the required properties and their group-theoretical interpretation, we write the explicit gamma matrix basis we have found:
\begin{equation}
\begin{array}{|ccc|ccc|}
  \hline
     \Gamma_1 & = & \left(
\begin{array}{cccc}
 0 & 0 & 0 & \sqrt{2} \\
 0 & 0 & -\sqrt{2} & 0 \\
 0 & 0 & 0 & 0 \\
 0 & 0 & 0 & 0 \\
\end{array}
\right) & \Gamma_2 & = & \left(
\begin{array}{cccc}
 0 & \sqrt{2} & 0 & 0 \\
 0 & 0 & 0 & 0 \\
 0 & 0 & 0 & 0 \\
 0 & 0 & \sqrt{2} & 0 \\
\end{array}
\right) \\
     \hline
     \Gamma_3 & = & \left(
\begin{array}{cccc}
 1 & 0 & 0 & 0 \\
 0 & -1 & 0 & 0 \\
 0 & 0 & 1 & 0 \\
 0 & 0 & 0 & -1 \\
\end{array}
\right) & \Gamma_4 & = & \left(
\begin{array}{cccc}
 0 & 0 & 0 & 0 \\
 \sqrt{2} & 0 & 0 & 0 \\
 0 & 0 & 0 & \sqrt{2} \\
 0 & 0 & 0 & 0 \\
\end{array}
\right) \\
     \hline
     \Gamma_5 & = & \left(
\begin{array}{cccc}
 0 & 0 & 0 & 0 \\
 0 & 0 & 0 & 0 \\
 0 & -\sqrt{2} & 0 & 0 \\
 \sqrt{2} & 0 & 0 & 0 \\
\end{array}
\right) & \null & \null & \null\\
\hline
   \end{array}
\end{equation}
\subsubsection{The map from the spinor to the vector representation}
Having   established the relation between the $\so(2,3)$ and the $\spalg (4,\mathbb{R})$ Lie algebras,  we can construct the   map
between spinor symplectic group elements and pseudo-orthogonal ones.
 Let $\mathcal{S} \in \Sp(4,\mathbb{R})$ be a symplectic group element, namely a $4\times 4 $ real matrix such that the following condition holds true:
\begin{equation}\label{crinetto}
 \mathcal{S}^T \,  \mathrm{C_s} \,  \mathcal{S} \, = \, \mathrm{C_s}.
\end{equation}
The corresponding group element $\mathcal{O}[\mathcal{S}] \in \mathrm{SO(2,3)}$ is determined by the relation
\begin{equation}\label{pernacchia}
  \mathcal{S}^{-1} \, \Gamma_j \, \mathcal{S} \, = \, \mathcal{O}^i_{\phantom{i}j}[\mathcal{S}] \, \Gamma_i,
\end{equation}
which implies
\begin{equation}\label{cogrutto}
  \mathcal{O}^i_{\phantom{i}j}[\mathcal{S}] \, = \, \ft 14 \, \text{Tr} \, \left( \Gamma^T_{i} \, \mathcal{S}^{-1} \, \Gamma_j \, \mathcal{S}\right),
\end{equation}
as the gamma-matrices are normalized so  that
\begin{equation}\label{curlano}
  \ft{1}{4}\, \text{Tr} \left( \Gamma^T_{i} \, \Gamma_p \right) \, = \, \delta^i_p
\end{equation}
Note that eq.~\eqref{cogrutto} is a perfect $r=2$ counterpart to eq.~\eqref{coriaceo} which  applies to the Tits-Satake group of the $r=1$ case.
\subsubsection{The pre-image in the symplectic group $\Sp(4,\mathbb{R})$ of the solvable subgroup of $\mathrm{SO(2,3)}$}
A question that naturally arises at this point is to o f identify the pre-image in the spinor symplectic group of the solvable group which is
metrically equivalent to the Tits-Satake symmetric space $\mathrm{SO(2,3)}/\mathrm{SO(2)}\times\mathrm{SO(3)}$. At the level of Lie algebra we
find that the 6-dimensional solvable Lie algebra of the Tits-Satake submanifold admits, in the symplectic 4-dimensional representation, the
following generators, where we do not care about the precise correspondence with the generators in the vector $5$-dimensional representation:
\begin{equation}\label{sympgensolv}
  \begin{array}{|ccc|ccc|}
  \hline
     T^{[Sp]}_1 & = & \left(
\begin{array}{cccc}
 -\frac{1}{2} & 0 & 0 & 0 \\
 0 & -\frac{1}{2} & 0 & 0 \\
 0 & 0 & \frac{1}{2} & 0 \\
 0 & 0 & 0 & \frac{1}{2} \\
\end{array}
\right) &  T^{[Sp]}_2 & = & \left(
\begin{array}{cccc}
 -\frac{1}{2} & 0 & 0 & 0 \\
 0 & \frac{1}{2} & 0 & 0 \\
 0 & 0 & \frac{1}{2} & 0 \\
 0 & 0 & 0 & -\frac{1}{2} \\
\end{array}
\right) \\
     \hline
     T^{[Sp]}_3 & = & \left(
\begin{array}{cccc}
 0 & 0 & 0 & 0 \\
 0 & 0 & 0 & \sqrt{2} \\
 0 & 0 & 0 & 0 \\
 0 & 0 & 0 & 0 \\
\end{array}
\right) &  T^{[Sp]}_4 & = & \left(
\begin{array}{cccc}
 0 & 0 & -\sqrt{2} & 0 \\
 0 & 0 & 0 & 0 \\
 0 & 0 & 0 & 0 \\
 0 & 0 & 0 & 0 \\
\end{array}
\right) \\
     \hline
     T^{[Sp]}_5 & = & \left(
\begin{array}{cccc}
 0 & 0 & 0 & 1 \\
 0 & 0 & 1 & 0 \\
 0 & 0 & 0 & 0 \\
 0 & 0 & 0 & 0 \\
\end{array}
\right) &  T^{[Sp]}_6 & = & \left(
\begin{array}{cccc}
 0 & 1 & 0 & 0 \\
 0 & 0 & 0 & 0 \\
 0 & 0 & 0 & 0 \\
 0 & 0 & -1 & 0 \\
\end{array}
\right) \\
     \hline
   \end{array}
\end{equation}
As one sees, considering the block structure of the symplectic group elements
\begin{equation}\label{pantalone}
  \mathcal{S} \, = \, \left( \begin{array}{c|c}
                               A & B \\
                               \hline
                               C & D
                             \end{array}
  \right),
\end{equation}
where $A,B,C,D$ are $2\times 2$ matrices, we find that the elements of the solvable subgroup are characterized by the vanishing of the $C$-block,
both at the level of Lie algebra and of the full group. The additional constraints stemming  from the request that the group element should be in
the exponential map of the Lie algebra generated by the $6$ generators in eq.~(\ref{sympgensolv}) leads, after suitable redefinitions, to the
following general form of the $4\times 4$ matrix depending on 6-parameters $\mathbf{a}=a_i$:
\begin{equation}\label{sympletsolv}
  \mathcal{S}_{s}[\mathbf{a}]\, = \, \left(
\begin{array}{cccc}
 \frac{1}{\sqrt{a_1} \sqrt{a_2}} & \frac{a_6}{\sqrt{a_1} \sqrt{a_2}} &
   -\frac{\sqrt{2} a_4}{\sqrt{a_1} \sqrt{a_2}}-\frac{a_5 a_6}{\sqrt{a_1}
   \sqrt{a_2}} & \frac{a_5}{\sqrt{a_1} \sqrt{a_2}} \\
 0 & \frac{\sqrt{a_2}}{\sqrt{a_1}} & \frac{\sqrt{a_2}
   a_5}{\sqrt{a_1}}-\frac{\sqrt{2} \sqrt{a_2} a_3 a_6}{\sqrt{a_1}} &
   \frac{\sqrt{2} \sqrt{a_2} a_3}{\sqrt{a_1}} \\
 0 & 0 & \sqrt{a_1} \sqrt{a_2} & 0 \\
 0 & 0 & -\frac{\sqrt{a_1} a_6}{\sqrt{a_2}} & \frac{\sqrt{a_1}}{\sqrt{a_2}} \\
\end{array}
\right).
\end{equation}
As one can see,  the block $C$ vanishes,  while $D=\left(A^{-1}\right)^T$. The fact that the matrices of type (\ref{sympletsolv}) constitute a subgroup can be
explicitly verified, as  we have
\begin{equation}\label{leggegruppa}
  \mathcal{S}_{s}[\mathbf{a}]\cdot\mathcal{S}_{s}[\mathbf{b}]\, = \, \mathcal{S}_{s}[\mathbf{c}]
\end{equation}
where
\begin{equation}\label{gruppalaw}
\begin{array}{rcl}
 c_1&=& a_1 b_1 \\
 c_2&=& a_2 b_2 \\
 c_3&=& \frac{a_3 b_1}{b_2}+b_3 \\
 c_4&=& -a_6^2 b_3 b_2^2+a_4 b_1 b_2-\sqrt{2} a_6 b_5 b_2+b_4 \\
 c_5&=& a_5 b_1+\sqrt{2} a_6 b_2 b_3+b_5 \\
 c_6&=& a_6 b_2+b_6.
\end{array}
\end{equation}
Utilizing the spinor map of eq.~(\ref{cogrutto}) we obtain
\begin{equation}\label{pennabianca}
  \mathcal{O}\left[\mathcal{S}_{s}[\mathbf{a}] \right]\, = \, \left(
\begin{array}{ccccc}
 a_1 & \sqrt{2} a_2 a_3 & \sqrt{2} a_5 & \frac{\sqrt{2} a_4}{a_2} &
   -\frac{a_5^2+2 a_3 a_4}{a_1} \\
 0 & a_2 & \sqrt{2} a_6 & -\frac{a_6^2}{a_2} & \frac{a_6 \left(\sqrt{2} a_3
   a_6-2 a_5\right)-\sqrt{2} a_4}{a_1} \\
 0 & 0 & 1 & -\frac{\sqrt{2} a_6}{a_2} & \frac{2 a_3 a_6-\sqrt{2} a_5}{a_1} \\
 0 & 0 & 0 & \frac{1}{a_2} & -\frac{\sqrt{2} a_3}{a_1} \\
 0 & 0 & 0 & 0 & \frac{1}{a_1}
\end{array}.
\right)
\end{equation}
The corresponding solvable group element built  with the $\Sigma$-exponential map in terms of the solvable coordinates
\begin{equation}\label{craniato}
  \boldsymbol{W}=\left\{w_1,w_2,w_3,w_4,w_5,\underbrace{0,\dots,0}_{(2s-1)},w_6,\underbrace{0,\dots,0}_{(2s-1)}\right\},
\end{equation}
where the zero  correspond to the Tits-Satake projection, have the   explicit form
\begin{equation}\label{LBmat}
  \mathbb{L}\left(\mathbf{W}\right)\, =\, \left(
\begin{array}{ccccc}
 e^{ {w_1}} & \frac{e^{ {w_1}}  {w_3}}{\sqrt{2}} & \frac{1}{2}
   e^{ {w_1}} \left( {w_3}  {w_6}+\sqrt{2}  {w_5}\right) &
   \frac{1}{8} e^{ {w_1}} \left(-\sqrt{2}  {w_3}  {w_6}^2+4 \sqrt{2}
    {w_4}-4  {w_5}  {w_6}\right) & -\frac{1}{4} e^{ {w_1}} \left(2
    {w_3}  {w_4}+ {w_5}^2\right) \\
 0 & e^{ {w_2}} & \frac{e^{ {w_2}}  {w_6}}{\sqrt{2}} & -\frac{1}{4}
   e^{ {w_2}}  {w_6}^2 & -\frac{e^{ {w_2}}  {w_4}}{\sqrt{2}} \\
 0 & 0 & 1 & -\frac{ {w_6}}{\sqrt{2}} & -\frac{ {w_5}}{\sqrt{2}} \\
 0 & 0 & 0 & e^{- {w_2}} & -\frac{e^{- {w_2}}  {w_3}}{\sqrt{2}} \\
 0 & 0 & 0 & 0 & e^{- {w_1}}
\end{array}
\right).
\end{equation}
Imposing the identification
\begin{equation}\label{lanerossi}
  \mathcal{O}\left[\mathcal{S}_{s}[\mathbf{a}] \right] \, = \, \mathbb{L}\left(\mathbf{W}\right)
\end{equation}
one finds a unique solution for the $a_i$ parameters  in terms of the $w_i$ parameters. Replacing that solution into eq.~(\ref{pennabianca}) one
obtains the spinor image of the Tits-Satake projected solvable group element:
\begin{equation}\label{cartoon}
  \mathcal{S}_{s}[\mathbf{W}] \, =\,\left(
\begin{array}{cccc}
 e^{\frac{1}{2} \left(-w_1-w_2\right)} & \frac{1}{2} e^{\frac{1}{2}
   \left(w_2-w_1\right)} w_6 & \frac{1}{4} e^{\frac{1}{2} \left(w_1+w_2\right)}
   \left(w_5 w_6-2 \sqrt{2} w_4\right) & \frac{1}{4} e^{\frac{1}{2}
   \left(w_1-w_2\right)} \left(2 w_5+\sqrt{2} w_3 w_6\right) \\
 0 & e^{\frac{1}{2} \left(w_2-w_1\right)} & \frac{1}{2} e^{\frac{1}{2}
   \left(w_1+w_2\right)} w_5 & \frac{e^{\frac{1}{2} \left(w_1-w_2\right)}
   w_3}{\sqrt{2}} \\
 0 & 0 & e^{\frac{1}{2} \left(w_1+w_2\right)} & 0 \\
 0 & 0 & -\frac{1}{2} e^{\frac{1}{2} \left(w_1+w_2\right)} w_6 & e^{\frac{1}{2}
   \left(w_1-w_2\right)} \\
\end{array}
\right)
\end{equation}
Note that eq.~(\ref{cartoon}) is the precise counterpart of eqs.~ (\ref{Tsubbo},\ref{curiente}) that appear in the discussion of the $r=1$ case.
\subsubsection{The Siegel upper plane}
The \emph{Siegel upper complex plane} of degree (or genus) $g$ is the generalization to higher dimensions of the Lobachevsky-Poincar\'e
hyperbolic plane.
Just as the standard hyperbolic plane with Poincar\'e metric is a complex analytic realization of a maximally split
symmetric space, namely $\mathrm{SL(2,\mathbb{R})/SO(2)}$, in the same way, the upper Siegel plane of degree $g$ is the complex analytic
realization of the symmetric space
\begin{equation}\label{siegelsymspac}
  \mathcal{M}_{Siegel} \, = \, \frac{\Sp\mathrm{(2\,g, \mathbb{R})}}{\mathrm{S[U(1)\times U(g)]}}.
\end{equation}
The key observation is that, just in the same way as the \emph{fractional linear transformation}
\begin{equation}\label{fraclin}
  z \, \rightarrow \, \tilde{z} \, \equiv \, \frac{a \,z +b }{c\, z + d } \quad ; \quad \left(
  \begin{array}{cc}
     a & b \\
     c & d \\
     \end{array}
 \right) \, \in \, \mathrm{PSL(2,\mathbb{R})}
\end{equation}
maps complex numbers $z$ with strictly positive imaginary part into complex numbers $\tilde{z}$ with the same property, \emph{the fractional
linear matrix transformation}
\begin{equation}\label{fraclinmat}
  Z_{g\times g} \, \rightarrow \, \tilde{Z}_{g\times g} \, \equiv \, (A_{g\times g} \,Z_{g\times g}
  +B_{g\times g} )\cdot(C_{g\times g}\,Z_{g\times g} + D_{g\times g})^{-1} \quad ; \quad \left(
  \begin{array}{c|c}
   A_{g\times g} & B_{g\times g} \\
   \hline
   C_{g\times g} &  D_{g\times g} \\
   \end{array}
\right) \, \in \, \Sp(2\,g,\mathbb{R})
\end{equation}
maps \emph{complex symmetric matrices}
\begin{equation}\label{Zetone}
  Z_{g\times g} \, = \,  X_{g\times g} \, + \, \mathit{i} \, Y_{g\times g} \quad ; \quad  Z_{g\times g}^T \, = \,  Z_{g\times g},
\end{equation}
whose imaginary part $Y_{g\times g}$ is positive definite (namely has strictly positive eigenvalues),  into \emph{complex symmetric matrices}
$\tilde{Z}_{g\times g}$ with the same property. The relations among the $g \times g$ blocks
\begin{equation}\label{blocrela}
  A^T \, C \, = \, C^T \, A \quad ; \quad B^T \, D \, = \, D^T \, B \quad ; \quad A^T D\, - \, C^T \, B \, = \, \mathbf{1 },
\end{equation}
following from the very definition of the $\Sp(2\,g,\mathbb{R})$ group, are instrumental in the lengthy yet straightforward proof of what was
stated above.
\par
The number of real components of $Z_{g\times g}$ exactly matches the dimension of the \emph{maximally split symmetric space} defined in
eq.~(\ref{siegelsymspac}), so that the upper Siegel plane is its holomorphic realization. Furthermore the choice of the Borel solvable
subgroup inside $\Sp(2\,g,\mathbb{R})$ provides a convenient parameterization of the matrix $Z_{g\times g}$. Indeed this   is the orbit under
the fractional linear action of the Borel subgroup of the special matrix $Z_0 \, = \,\mathit{i} \, \mathbf{1}_{g\times g}$.
\par
Applying this idea to our case $g=r=2$  and utilizing the parameterization of the Borel solvable subgroup provided in
eq.~(\ref{cartoon}), we obtain
\begin{eqnarray}\label{muradilegno}
  Z & = & \, X \, + \, \mathit{i} \, Y \nonumber\\
  X & = & \left(
\begin{array}{cc}
 \frac{1}{8} \left(w_6 \left(4 w_5+\sqrt{2} w_3 w_6\right)-4 \sqrt{2} w_4\right)
   & \frac{1}{4} \left(2 w_5+\sqrt{2} w_3 w_6\right) \\
 \frac{1}{4} \left(2 w_5+\sqrt{2} w_3 w_6\right) & \frac{w_3}{\sqrt{2}} \\
\end{array}
\right) \nonumber\\
  Y & = & \left(
\begin{array}{cc}
 \frac{1}{4} e^{-w_1-w_2} \left(e^{2 w_2} w_6^2+4\right) & \frac{1}{2}
   e^{w_2-w_1} w_6 \\
 \frac{1}{2} e^{w_2-w_1} w_6 & e^{w_2-w_1} \\
\end{array}
\right).
\end{eqnarray}
Eqs.\eqref{muradilegno} are the precise counterpart in the $r=2$ case of eqs.~(\ref{puntacchia}-\ref{cristallus}). It is also  appropriate to stress
once more what was already recalled above, namely, that in the series of non-maximally split manifolds (\ref{paniepesci}) the appearance of the
Siegel upper plane as realization of the Tits-Satake sub-manifold $\mathcal{M}_{TS} \subset \mathcal{M}^{[r,r+q]}$ where $q=2s$, or $q=2s+1$ is
limited to $r=1,2$ and it is strictly due to sporadic low-rank Lie algebra isomorphisms. The Tits-Satake submanifold exists in all cases but it is
not equivalent to a Siegel symmetric space as defined in eq.~(\ref{siegelsymspac}). For instance in the case $r=3$, the Tits-Satake submanifold is
\begin{equation}\label{mangiagente}
  \mathcal{M}^{[3,4]}_{TS} \, = \,\frac{\mathrm{SO(3,4)}}{\mathrm{SO(3)\times SO(4)}},
\end{equation}
where the numerator   $\mathrm{SO(3,4)}$ corresponds to the maximally split real section $\so(3,4)$ of the complex $\so(7,C)$ Lie algebra. The
spinor representation of $\mathrm{SO(3,4)}$ is therefore an appropriate $\mathrm{Spin(3,4)}$ subgroup of $\mathrm{SO(4,4)}$  which is made of
pseudo-orthogonal rather than symplectic matrices. The upper Siegel plane with $g>2$  can instead appear as Tits-Satake submanifold in the Tits
Satake projection of other non maximally split manifolds originating from real sections of the $\spalg(2n,\mathbb{C})$ Lie algebra, which we do
not study in the present article.  Yet in view of the analysis presented in the introductory section \ref{lollobrigida} the upper Siegel plane of
degree $2$ appears as building block also in the $c$-map construction of cases $3,4$.
\par
As we stressed several times, in applications to Data Science a crucial issue is that of \emph{discretization schemes} of the symmetric spaces
$\mathrm{U/H}$ to which we map the data. In view of the fundamental theorem \ref{teoUHsolvrep} such \emph{discretization schemes} necessarily
involve the consideration of discrete subgroups of the full isometry group that should intersect the solvable group by means of suitable infinite
order normal subgroups. Such observation provides a precious key in the search for discrete subgroups suitable for discretization schemes. In
section \ref{plantageneti}  we illustrate such an issue with the determination of a discrete subgroup of $\mathrm{Sp(4,\mathbb{Z})}$ of the requested
type, which up to our knowledge was so far unknown.

\subsection{$c$-map construction of the Tits-Satake manifold $\mathrm{SO(3,4)}/\mathrm{SO(3)\times SO(4)}$}
Following the last comments reported in the above subsection, we conclude with the explicit demonstration of the $c$-map construction of the Tits
Satake manifold for the $r=3$ universality class starting from the special K\"ahler manifold
\begin{equation}\label{formaggio}
  \mathcal{SK}_2 \, \equiv \, \frac{\mathrm{SL(2,\mathbb{R})}}{\mathrm{SO(2)}} \times \frac{\mathrm{SO(1,2)}}{\mathrm{SO(2)}}
\end{equation}
The holomorphic symplectic section of eq.~(\ref{ololo}) is in this case given by
\begin{eqnarray}\label{holsymsec}
  \Omega_h & = & \{X^\Lambda, F_\Sigma\}\nonumber\\
   & = &  \left\{\frac{z^2}{2 \sqrt{2}},\frac{z}{2},-\frac{1}{2 \sqrt{2}},-\frac{S}{2
   \sqrt{2}},\frac{S z}{2},\frac{S z^2}{2 \sqrt{2}}\right\},
\end{eqnarray}
where $S$ is the complex coordinate spanning the first hyperbolic plane $\frac{\mathrm{SL(2,\mathbb{R})}}{\mathrm{SO(2)}}$, while $z$ is the
complex coordinate spanning the second hyperbolic plane $\frac{\mathrm{SO(1,2)}}{\mathrm{SO(2)}}$;  in the special structure the two planes  have distinct and
not interchangeable roles. Let us also call
\begin{equation}\label{buonasera}
  x \, = \, \mathrm{Re} z \quad ; y \, = \, \mathrm{Im} z \quad ; \quad P \, = \, \mathrm{Re }S \quad ; Q \, = \, \mathrm{Im }S.
\end{equation}
From its definition given in eqs.~(\ref{etamedia}-\ref{intriscripen}) we can calculate the $ 3\times 3$ complex period matrix
$\mathcal{N}_{\Lambda\Sigma}$,  which has the   form
\begin{equation}\label{NLS}
  \mathcal{N}_{\Lambda\Sigma}\, =\, \left(
\begin{array}{ccc}
 \frac{i Q \left(x^2+y^2\right)^2}{y^2} & \frac{i \sqrt{2} Q x
   \left(x^2+y^2\right)}{y^2} & P-\frac{i Q x^2}{y^2} \\
 \frac{i \sqrt{2} Q x \left(x^2+y^2\right)}{y^2} & P+i Q \left(\frac{2
   x^2}{y^2}+1\right) & -\frac{i \sqrt{2} Q x}{y^2} \\
 P-\frac{i Q x^2}{y^2} & -\frac{i \sqrt{2} Q x}{y^2} & \frac{i Q}{y^2} \\
\end{array}
\right).
\end{equation}
We can easily extract from eq.~\eqref{NLS}  the real and imaginary parts of $\mathcal{N}_{\Lambda\Sigma}$ and     compute  the $6\times 6$ matrix
$\mathcal{M}_4^{-1}$ defined in eq.~(\ref{inversem4}):
\begin{equation}\label{pirolino}
\mathcal{M}_4^{-1} \, = \, \left(
\begin{array}{cccccc}
 \frac{1}{Q y^2} & -\frac{\sqrt{2} x}{Q y^2} & -\frac{x^2}{Q y^2} & -\frac{P x^2}{Q
   y^2} & -\frac{\sqrt{2} P x}{Q y^2} & \frac{P}{Q y^2} \\
 -\frac{\sqrt{2} x}{Q y^2} & \frac{\frac{2 x^2}{y^2}+1}{Q} & \frac{\sqrt{2} x
   \left(x^2+y^2\right)}{Q y^2} & \frac{\sqrt{2} P x \left(x^2+y^2\right)}{Q y^2} &
   \frac{P \left(2 x^2+y^2\right)}{Q y^2} & -\frac{\sqrt{2} P x}{Q y^2} \\
 -\frac{x^2}{Q y^2} & \frac{\sqrt{2} x \left(x^2+y^2\right)}{Q y^2} &
   \frac{\left(x^2+y^2\right)^2}{Q y^2} & \frac{P \left(x^2+y^2\right)^2}{Q y^2} &
   \frac{\sqrt{2} P x \left(x^2+y^2\right)}{Q y^2} & -\frac{P x^2}{Q y^2} \\
 -\frac{P x^2}{Q y^2} & \frac{\sqrt{2} P x \left(x^2+y^2\right)}{Q y^2} & \frac{P
   \left(x^2+y^2\right)^2}{Q y^2} & \frac{\left(P^2+Q^2\right)
   \left(x^2+y^2\right)^2}{Q y^2} & \frac{\sqrt{2} x \left(P^2+Q^2\right)
   \left(x^2+y^2\right)}{Q y^2} & -\frac{x^2 \left(P^2+Q^2\right)}{Q y^2} \\
 -\frac{\sqrt{2} P x}{Q y^2} & \frac{P \left(2 x^2+y^2\right)}{Q y^2} &
   \frac{\sqrt{2} P x \left(x^2+y^2\right)}{Q y^2} & \frac{\sqrt{2} x
   \left(P^2+Q^2\right) \left(x^2+y^2\right)}{Q y^2} & \frac{\left(P^2+Q^2\right)
   \left(2 x^2+y^2\right)}{Q y^2} & -\frac{\sqrt{2} x \left(P^2+Q^2\right)}{Q y^2} \\
 \frac{P}{Q y^2} & -\frac{\sqrt{2} P x}{Q y^2} & -\frac{P x^2}{Q y^2} & -\frac{x^2
   \left(P^2+Q^2\right)}{Q y^2} & -\frac{\sqrt{2} x \left(P^2+Q^2\right)}{Q y^2} &
   \frac{P^2+Q^2}{Q y^2} \\
\end{array}
\right)
\end{equation}
At this point we have, ready in explicit form, all the items necessary to write the Quaternionic K\"ahler metric in the image of the $c$-map as
given in eq.~(\ref{geodaction}). However it is better to pause and recall that the considered Special K\"ahler model is a homogeneous symmetric
space, so that also its quaternionic K\"ahler image has to be a symmetric space, according to what is explained in section \ref{omosymmetro}.
Furthermore, according with eqs.~(\ref{qcosetto}-\ref{turnaconto}), the symplectic image $\Lambda[\mathfrak{g}^{-1} \in \mathrm{Sp(6,\mathbb{R})}$
of the $\mathrm{SL(2,\mathbb{R}) \times SO(1,2)}$   element which takes  $z,S$ to their respective origins $z_0 \, = \,i$ and $S_0 \, = \, i$
should take the matrix $\mathcal{M}_4^{-1}$ to a constant simple form. In our notation the mentioned symplectic group element is
\begin{equation}\label{lambone}
  \Lambda(\mathfrak{g}^{-1})\, = \, \left(
\begin{array}{cccccc}
 \frac{1}{\sqrt{Q} y} & -\frac{\sqrt{2} x}{\sqrt{Q} y} & -\frac{x^2}{\sqrt{Q} y} &
   -\frac{P x^2}{\sqrt{Q} y} & -\frac{\sqrt{2} P x}{\sqrt{Q} y} & \frac{P}{\sqrt{Q}
   y} \\
 0 & \frac{1}{\sqrt{Q}} & \frac{\sqrt{2} x}{\sqrt{Q}} & \frac{\sqrt{2} P x}{\sqrt{Q}}
   & \frac{P}{\sqrt{Q}} & 0 \\
 0 & 0 & \frac{y}{\sqrt{Q}} & \frac{P y}{\sqrt{Q}} & 0 & 0 \\
 0 & 0 & 0 & \sqrt{Q} y & 0 & 0 \\
 0 & 0 & 0 & \sqrt{2} \sqrt{Q} x & \sqrt{Q} & 0 \\
 0 & 0 & 0 & -\frac{\sqrt{Q} x^2}{y} & -\frac{\sqrt{2} \sqrt{Q} x}{y} &
   \frac{\sqrt{Q}}{y} \\
\end{array}
\right),
\end{equation}
and indeed we have
\begin{equation}\label{cunegondax}
  \mathcal{M}_4^{-1} \, = \, \Lambda(\mathfrak{g}^{-1})^T\cdot\Lambda(\mathfrak{g}^{-1}).
\end{equation}
This enables us to right directly the 12 vielbein $1$-forms of the symmetric space (\ref{mangiagente}) according with eq.~(\ref{filibaine}):
\begin{eqnarray}\label{confiscato}
  E^1 &=& dU, \,\,\,\, E^2 =\frac{dP}{2\, Q},\,\,\,\, E^3 =\frac{dQ}{2\, Q}, \,\,\,\,E^4 =\frac{dx}{2\, y},\,\,\,\, E^5 =\frac{dy}{2\, y} \nonumber\\
  E^6 &=& \exp[-U] \, \left(da + \frac{\alpha}{2} \mathbf{Z}^T \cdot \mathbb{C} d\mathbf{Z}\right),
  \,\,\,\, E^{6+i} =\exp\left[\frac{U}{2}\right] \,\left(\Lambda(\mathfrak{g}^{-1})\cdot \mathbf{Z}\right)^i \quad i=1,\dots,6,
\end{eqnarray}
where $\mathbf{Z}=\{Z^1,Z^2,Z^3,Z^4,Z^5,Z^6\}$, $\alpha$ is some constant to be determined  and
\begin{equation}\label{cmatrazza}
  \mathbb{C} \, = \,  \left(
                              \begin{array}{cc}
                                \mathbf{0}_{3\times 3} & \mathbf{1}_{3\times 3} \\
                                -\mathbf{1}_{3\times 3} & \mathbf{0}_{3\times 3} \\
                              \end{array}
                            \right).
\end{equation}
The $1$-forms $E^I$ satisfy the   Maurer Cartan equations by construction:
\begin{equation}\label{MCcm}
  \begin{array}{rcl}
0&=&  d{E^{1}} \\
0&=&  d{E^{2}}-2 E^{2}\wedge E^{3} \\
0&=&  d{E^{3}} \\
0&=&  d{E^{4}}-2 E^{4}\wedge E^{5} \\
0&=&   d{E^{5}} \\
0&=&   d{E^{6}}-\alpha  E^{7}\wedge E^{10}-\alpha  E^{8}\wedge E^{11}-\alpha  E^{9}\wedge
      E^{12}+E^{1}\wedge E^{6} \\
0&=&   d{E^{7}}+\frac{1}{2} E^{1}\wedge E^{7}-2 E^{2}\wedge E^{12}+E^{3}\wedge E^{7}+2 \sqrt{2}
   E^{4}\wedge E^{8}+2 E^{5}\wedge E^{7} \\
0&=&   d{E^{8}}+\frac{1}{2} E^{1}\wedge E^{8}-2 E^{2}\wedge E^{11}+E^{3}\wedge E^{8}-2 \sqrt{2}
   E^{4}\wedge E^{9} \\
0&=&   d{E^{9}}+\frac{1}{2} E^{1}\wedge E^{9}-2 E^{2}\wedge E^{10}+E^{3}\wedge E^{9}-2
   E^{5}\wedge E^{9} \\
0&=&   d{E^{10}}+\frac{1}{2} E^{1}\wedge E^{10}-E^{3}\wedge E^{10}-2 E^{5}\wedge E^{10} \\
0&=&   d{E^{11}}+\frac{1}{2} E^{1}\wedge E^{11}-E^{3}\wedge E^{11}-2 \sqrt{2} E^{4}\wedge E^{10}
   \\
0&=&   d{E^{12}}+\frac{1}{2} E^{1}\wedge E^{12}-E^{3}\wedge E^{12}+2 \sqrt{2} E^{4}\wedge
   E^{11}+2 E^{5}\wedge E^{12} \\
\end{array}
\end{equation}
On the other hand, if we consider the standard parameterization of the solvable group $S_{U/H}$ for the maximally split $\mathrm{U/H}$ manifold of
eq.~(\ref{mangiagente}) according with eqs.~(\ref{generalini}-\ref{kyrillov}), we can  construct the left-invariant one-form $\Theta$ defined in
eq.~(\ref{filbone}), and   project it onto the generators of the solvable group $T_A$ rather than onto the generators $K_A$ of the coset manifold:
\begin{equation}\label{quaccherini}
  \Theta \, = \, \sum_{A=1}^{12} \mathfrak{e}^A \, T_A;
\end{equation}
thus we obtain another set of $1$-form $\mathfrak{e}^A$ satisfying the   Maurer Cartan equations:
\begin{equation}\label{MCsol}
  \begin{array}{rcl}
0&=& d{\mathfrak{e}^{1}} \\
0&=& d{\mathfrak{e}^{2}} \\
0&=& d{\mathfrak{e}^{3}} \\
0&=& d{\mathfrak{e}^{4}}+\mathfrak{e}^{1}\wedge \mathfrak{e}^{4}-\mathfrak{e}^{2}\wedge \mathfrak{e}^{4} \\
0&=& d{\mathfrak{e}^{5}}+\mathfrak{e}^{1}\wedge \mathfrak{e}^{5}-\mathfrak{e}^{3}\wedge \mathfrak{e}^{5}+\frac{\mathfrak{e}^{4}\wedge
 \mathfrak{e}^{6}}{\sqrt{2}} \\
0&=& d{\mathfrak{e}^{6}}+\mathfrak{e}^{2}\wedge \mathfrak{e}^{6}-\mathfrak{e}^{3}\wedge \mathfrak{e}^{6} \\
0&=& d{\mathfrak{e}^{7}}+\mathfrak{e}^{1}\wedge \mathfrak{e}^{7}+\mathfrak{e}^{2}\wedge \mathfrak{e}^{7}-\frac{\mathfrak{e}^{5}\wedge
   \mathfrak{e}^{9}}{\sqrt{2}}+\frac{\mathfrak{e}^{6}\wedge \mathfrak{e}^{8}}{\sqrt{2}}-\frac{\mathfrak{e}^{10}\wedge
   \mathfrak{e}^{11}}{\sqrt{2}} \\
0&=& d{\mathfrak{e}^{8}}+\mathfrak{e}^{1}\wedge \mathfrak{e}^{8}+\mathfrak{e}^{3}\wedge \mathfrak{e}^{8}+\frac{\mathfrak{e}^{4}\wedge
   \mathfrak{e}^{9}}{\sqrt{2}}-\frac{\mathfrak{e}^{10}\wedge \mathfrak{e}^{12}}{\sqrt{2}} \\
0&=& d{\mathfrak{e}^{9}}+\mathfrak{e}^{2}\wedge \mathfrak{e}^{9}+\mathfrak{e}^{3}\wedge \mathfrak{e}^{9}-\frac{\mathfrak{e}^{11}\wedge
 \mathfrak{e}^{12}}{\sqrt{2}} \\
0&=& d{\mathfrak{e}^{10}}+\mathfrak{e}^{1}\wedge \mathfrak{e}^{10}+\frac{\mathfrak{e}^{4}\wedge \mathfrak{e}^{11}}{\sqrt{2}}+\frac{\mathfrak{e}^{5}\wedge
   \mathfrak{e}^{12}}{\sqrt{2}} \\
0&=& d{\mathfrak{e}^{11}}+\mathfrak{e}^{2}\wedge \mathfrak{e}^{11}+\frac{\mathfrak{e}^{6}\wedge \mathfrak{e}^{12}}{\sqrt{2}} \\
0&=& d{\mathfrak{e}^{12}}+\mathfrak{e}^{3}\wedge \mathfrak{e}^{12} \\
\end{array}
\end{equation}
The two sets of Maurer Cartan equations coincide if we set $\alpha \, = \, - \frac{1}{64}$ and if we make the following identification:
%%%%%%%%%%%%%%%%%%%%%%%%%%%%%%%%%%%%%%%%%%%%%%%%%%%%%%%%%%%%%%%%%%%%%%%%%%%%%%%%%%%%%%%%%%%%%%%%%%%%%%%%%%%%%%%%%%%%%%%%%%%%%%%%%%%%%%%
\begin{equation}\label{agatasco}
  \begin{array}{rcl}
E^1 & = & \mathfrak{e}^{1}+\mathfrak{e}^{2} \\
E^2 & = & \mathfrak{e}^{4} \\
E^3 & = & \frac{1}{2} (\mathfrak{e}^{1}-\mathfrak{e}^{2}) \\
E^4 & = & \mathfrak{e}^{12} \\
E^5 & = & \frac{\mathfrak{e}^{3}}{2} \\
E^6 & = & \mathfrak{e}^{7} \\
E^7 & = & -32 \sqrt{2} \mathfrak{e}^{8} \\
E^8 & = & -8 \sqrt{2} \mathfrak{e}^{10} \\
E^9 & = & -2 \sqrt{2} \mathfrak{e}^{5} \\
E^{10} & = & \mathfrak{e}^{6} \\
E^{11} & = & 4 \mathfrak{e}^{11} \\
E^{12} & = & 16 \mathfrak{e}^{9} \\
\end{array}
\end{equation}
We do not show the explicit form of the $1$-form $\mathfrak{e}^i$ because is too big. What we have shown proves that the $c$-map construction
provides a much simpler form of the metric. The invariant Einstein metric can be obtained from the $1$ forms $E^i$ with a few simple relative
rescalings that we also omit for brevity.
\section{A grid generating discrete subgroup of $\mathrm{Sp(4,\mathbb{Z})}$}\label{plantageneti}
The perspective applications of  PGTS theory to Data Science will always meet the problem of discretization, in one way or another, since data,
even if very big, are anyhow discrete sets: hence their mapping to symmetric spaces $\mathrm{U/H}$  necessarily introduces a finite set of points.
%or of finite subsets into the latter.
The  ways to look at such discrete set of   points are various, yet  eventually they always
result into some kind of tessellation of the hosting manifold that can be more or less regular, depending on the case.
So, if distribute $m$ points $P=\{p_1,\dots,p_m\}$ in a hosting manifold $\mathcal M$, we have a tesselation $\{C_1,\dots,C_m\}$, where
the cell $C_i$ is formed by the points of $\mathcal M$ whose closest point of $P$ is $p_i$.
Hence tessellation schemes and
their associated graphs and groups are an unavoidable chapter of the tale. In  \cite{tassellandum} it was  discussed the general
setup of tessellations of the Hyperbolic Plane $\mathbb{H}_2$ within the framework of Coxeter groups, introducing also some new, or less studied
ones. In that context  the distinction between polygonal and apeirogon tessellations (an apeirogon is a polygon with an infinite
number of edges) was emphasized, discussing  the relative merit of both in the context of Data mapping. Obviously tessellations  of the Siegel halfspace
$\mathbb{SH}_2$  have also to be constructed and this requires the study and classification of discrete subgroups of $\mathrm{Sp(4,\mathbb{R})}$,
which is a much wider and more complex task. In line with the main conceptual stream of the present paper, we want to show that the representation
of the entire symmetric space as the group manifold of a solvable Lie group provides new insights and suggestions for the construction of infinite
discrete subgroups of $\mathrm{Sp(4,\mathbb{R})}$, certainly related with tessellation schemes of $\mathbb{SH}_2$. In particular, in the present
section we show an example of a discrete subgroup of $\mathrm{Sp(4,\mathbb{Z})}$ which appears to be  a generalization for $\mathbb{SH}_2$ of what
the full modular group $\mathrm{PSL(2,\mathbb{Z})}$ is for the Hyperbolic Plane $\mathbb{H}_2$, namely, the generator of an apeirogon tessellation.
The main point of interest in the construction we are going to present in this section is the very construction algorithm. It is based on a
discretization of the solvable subgroup of $\mathrm{Sp(4,\mathbb{R})}$.
\par
In line with what we stated above, we concentrated on the following idea. The imaginable discretizations of the full space require  a  two
steps process:
\begin{enumerate}
   \item in a first step one singles out a discrete subgroup $\Delta\subset\mathrm{SO(r,r+1)}$ of the Tits-Satake isometry group with an ensuing
       tessellation of the Tits-Satake submanifold $\mathcal{M}_{TS} \, =\, \frac{\mathrm{SO(r,r+1)}}{\mathrm{SO(r) \times SO(r+1)}}$.
   \item In a second step one might try to extend the group $\Delta\subset\mathrm{SO(r,r+1)}$ with new fiber generators that necessarily belong
       to the normal subgroup discussed in previous section \ref{TSnormalsub}. Actually in section \ref{generaloneSOpq} we show how we can
       directly determine the parabolic subgroups of all $\mathrm{SO(r,r+p,\mathbb{Z})}$ for all $r,p\in\mathbb{N}$.
\end{enumerate}
\par
The requirements for step number one are rather clear and the result has been inspiring in order to find the general solutions discussed in sect.
\ref{generaloneSOpq}. For \emph{apeirogon-like tessellation} the group should be of the type $\Delta_{\infty,3,2}$ which, in the
$\mathrm{PSL(2,\mathbb{R})}$ case is just the modular group $\mathrm{PSL(2,\mathbb{R})}$ itself. This means that the group $\Delta$
should have several  generators $T_i$ ($i=1,..,p$) of infinite order, typically obtained as finite parabolic elements of the solvable group or
of its subgroups. It is however necessary that such generators close to a group, which, in particular means that the inverse of each of them must be
a word composed with the generators as letters. These infinite dimensional subgroup is the analogous of what in crystallography is  the
\emph{lattice}. In addition to the lattice, we need the analogue of the \emph{Point Group}, which should map the lattice into itself. If we
call $\mathcal{I}$  the infinite group generated by the $T_i$, the point group $\mathfrak{P}$ is the external automorphism group of $\mathcal{I}$,
and it must be a finite subgroup of the compact isotropy group of the manifold, namely, we must have $\mathfrak{P} \subset \mathrm{SO(r) \times
SO(r+1)}$.
\subsection{The case $r=2$ and the discrete group $\Delta_{32,8}$}
In the case $r=2$ the Tits-Satake isometry group $\mathrm{SO(2,3)}$ is locally isomorphic to $\mathrm{Sp(4,\mathbb{R})}$, whose fundamental
defining representation is the spinor representation of $\mathrm{SO(2,3)}$. Hence what we are looking for is a discrete subgroup
$\Delta_{2m,n}\subset \mathrm{Sp(4,\mathbb{R})}$ which will be the double covering of a subgroup $\hat{\Delta}_{m,n}\subset \mathrm{SO(2,3)}$. The
notation that we utilize is the following:
\begin{eqnarray}\label{patergloriosus}
  2m & = & |\mathfrak{P}_{2m}| \quad \text{where}  \quad \mathfrak{P}_{2m}\subset \mathrm{SU(2)\times U(1)} \quad \text{ = Point Group} ;\nonumber\\
  n & = & \# \,\,\text{ of generators $T_i$}.
\end{eqnarray}
Utilizing the guiding lines sketched out  above and requesting that the discrete subgroup  $\Delta_{2m,n}$ should contain two copies of the modular
group $\mathrm{PSL(2,\mathbb{Z})}$, we have been able to determine a group
\begin{equation}\label{m32n8}
  \Delta_{32,8} \subset \mathrm{Sp(4,\mathbb{Z})}
\end{equation}
which is actually a subgroup of the integer-valued symplectic group. We presently describe it.
\par
Prior to that let us explain the reason behind the request of the two copies of  $\mathrm{PSL(2,\mathbb{Z})}$, which proved to be   decisive.
The reason is that, in accordance with equation \ref{fraclinmat}, the action of any element $g\in \mathrm{Sp(4,\mathbb{R})}$ on the coset manifold is
the  fractional linear transformation of the symmetric complex matrix:
\begin{equation}\label{zetone}
  Z \, = \,\left(
\begin{array}{cc}
 z & \omega  \\
 \omega  & \zeta  \\
\end{array}
\right) \quad ; \quad z,\omega,\zeta \in \mathbb{C},
\end{equation}
which represents  the entire manifold.  When $Z$ is diagonal, i.e.~when $\omega = 0$, the two remaining complex entries $z,\zeta$ represent the
coordinates of two hyperbolic upper planes. The subgroup $\Gamma \subset \mathrm{Sp(4,\mathbb{R})}$ which respects diagonality, namely the
condition $\omega =0$, is $\Gamma \, = \, \mathrm{SL(2,\mathbb{R})} \times \mathrm{SL(2,\mathbb{R}) }$. Hence it is reasonable to ask that the
intersection
 $\Delta\bigcap \Gamma$ of   the sought for  discrete group with the stability subgroup $\Gamma$ of the submanifold $\omega =0$ should be
  \begin{equation}\label{pignattoso}
   \Delta\bigcap \Gamma \, = \, \mathrm{PSL(2,\mathbb{Z})_{1}} \times\mathrm{ PSL(2,\mathbb{Z})_{2}}.
 \end{equation}
In this way the two disconnected hyperbolic planes $z$ and $\zeta$  might be covered with apeirogon tiles. Subsequentely, introducing additional
generators of the group $\Delta$, we generate a tiling  also of the third plane that interacts with the first two creating new faces, edges and
vertices.
\par
In the above outlined perspective it is convenient to recall eq.~(\ref{muradilegno}), which provides the parameterization of the of the complex
matrix (\ref{zetone}) in terms of the solvable coordinates $w_1,w_2,w_3, w_4,w_5,w_6$. This    can be summarized by stating
\begin{eqnarray}\label{lisciva}
  z &=& \frac{1}{8} \left(-4 \sqrt{2} w_4+w_6 \left(4 w_5+\sqrt{2} w_3 w_6\right)+2 i
   e^{-w_1-w_2} \left(e^{2 w_2} w_6^2+4\right)\right) \nonumber\\
  \zeta &=& \frac{w_3}{\sqrt{2}}+i e^{w_2-w_1} \nonumber\\
 \omega &=& \frac{1}{4} \left(2 w_5+2 i e^{w_2-w_1} w_6+\sqrt{2} w_3 w_6\right).
\end{eqnarray}
As one sees from eq.~\eqref{lisciva},  the diagonalization condition of the matrix $Z$ corresponds to setting $w_5 = w_6 =0$,  which implies that the
other two solvable coordinates $w_3,w_4$ obtain the interpretation of real parts of the complex coordinates $\zeta$ and $z$, respectively. This is
a formidably helpful information. The parabolic generators of the sought for discrete group $\Delta$ have to be found as finite elements of the
solvable group associated  with $w_3,w_4$ and $w_5, w_6$, respectively. There is a conceptual distinction between the two pairs. In the solvable group
which is metrically equivalent to the Siegel half-plane, $w_3,w_4$ correspond to the exponentials of the long root generators that are paint group
singlets, while $w_5,w_6$ correspond to exponentials of the short generators that form paint group multiplets. This observation might be very much
relevant  while addressing the issue of extending the tessellation group from the Tits-Satake submanifold to the entire manifold.
\par
Inspired by our previous analysis of the normal subgroup of the solvable group associated with the Tits-Satake fibers, and relying on paint group
invariance that states that  in the Tits-Satake projection it is irrelevant  which of the $2s$ copies of the coordinates $ w_{5,i}, w_{6,i}$ we
keep, setting to zero all   other copies, we consider the $4\times 4$ image $\mathcal{S}_{s}[\mathbf{W}] $ of the solvable Tits-Satake group
inside $\mathrm{Sp(4,\mathbb{R})}$ which was presented in eq.~(\ref{cartoon}). In that matrix we make the substitution
\begin{equation}\label{carnoso}
  w_{1,2,3} \to 0 \quad ; \quad w_4\to \sqrt{2} \rho\quad ; \quad w_5\to 2 \sqrt{2} p \quad ; \quad w_6\to 2 \sqrt{2} q
\end{equation}
and    obtain the symplectic image of the normal subgroup associated with the short roots $w_{5,6}$ and the highest long root $w_{4}$:
\begin{equation}\label{bassotto}
  \mathfrak{L}\left(\rho,p,q\right) \, =\,\left(
\begin{array}{cccc}
 1 & \sqrt{2} q & \frac{1}{4} (8 p q-4 \rho ) & \sqrt{2} p \\
 0 & 1 & \sqrt{2} p & 0 \\
 0 & 0 & 1 & 0 \\
 0 & 0 & -\sqrt{2} q & 1 \\
\end{array}
\right).
\end{equation}
It is from the normal subgroup $\mathfrak{L}\left(\rho,p,q\right)$ that we  selected the integer-valued generators able to translate the solvable
coordinates $w_5,w_6$,  and hence generate translations of the complex number $\omega$ appearing in the parameterization of the Siegel upper plane.
\par
On the other hand we consider the subgroup obtained by setting in $\mathcal{S}_{s}[\mathbf{W}] $ the values
\begin{equation}\label{carnoso}
  w_{1,2} \to 0 \quad ; \quad w_3\to  \sqrt{2} h \quad ; \quad w_{4} \to \sqrt{2} k  \quad ; \quad w_{5,6} \to 0
\end{equation}
and we obtain the 2-parameter subgroup made up by the  matrices
\begin{equation}\label{longhetto}
  \mathfrak{T}\left(h,k\right) \, =\,\left(
\begin{array}{cccc}
 1 & 0 & -k & 0 \\
 0 & 1 & 0 & h \\
 0 & 0 & 1 & 0 \\
 0 & 0 & 0 & 1 \\
\end{array}
\right)
\end{equation}
It is from the above subgroup that we have taken the integer-valued generators able to translate the solvable coordinates $w_3,w_4$. Our staring
point has been the following.
\paragraph{\sc Choice of 4 infinite order translation generators.}
We have selected:
\begin{eqnarray}
  T_1 \,=&  \mathfrak{T}\left(0,-1\right) =& \left(
\begin{array}{cccc}
 1 & 0 & 1 & 0 \\
 0 & 1 & 0 & 0 \\
 0 & 0 & 1 & 0 \\
 0 & 0 & 0 & 1 \\
\end{array}
\right) \\
  T_2 \,=&  \mathfrak{T}\left(1,0\right)  = & \left(
\begin{array}{cccc}
 1 & 0 & 0 & 0 \\
 0 & 1 & 0 & 1 \\
 0 & 0 & 1 & 0 \\
 0 & 0 & 0 & 1 \\
\end{array}
\right)\\
  T_3 \,=&  \mathfrak{L}\left(0,\frac{1}{\sqrt{2}},0\right) = &  \left(
\begin{array}{cccc}
 1 & 0 & 0 & 1 \\
 0 & 1 & 1 & 0 \\
 0 & 0 & 1 & 0 \\
 0 & 0 & 0 & 1 \\
\end{array}
\right)\\
  T_4 \,=& \mathfrak{L}\left(0,0,\frac{1}{\sqrt{2}}\right) = & \left(
\begin{array}{cccc}
 1 & 1 & 0 & 0 \\
 0 & 1 & 0 & 0 \\
 0 & 0 & 1 & 0 \\
 0 & 0 & -1 & 1 \\
\end{array}
\right)
\end{eqnarray}
Defining a symbol
\begin{equation}\label{frattale}
  \mathbf{FL}\left[\mathfrak{g},Z\right] \, \equiv \,  \, (A_{\mathfrak{g}} \,Z +B_{\mathfrak{g}} )\cdot(C_{\mathfrak{g}}\,Z + D_{\mathfrak{g}})^{-1}  \quad ; \quad \mathfrak{g} = \left(\begin{array}{c|c}
                                                                                            A_{\mathfrak{g}} & B_{\mathfrak{g}} \\
                                                                                            \hline
                                                                                            C_{\mathfrak{g}} &  D_{\mathfrak{g}} \\
                                                                                          \end{array}
                                                                                        \right),
\end{equation}
which denotes the fractional linear action of an $\mathrm{Sp(4,\mathbb{R})}$ matrix on a complex symmetric matrix $Z$, we find
\begin{alignat}{3}\label{gromiko}
  \mathbf{FL}\left[T_1,Z\right] & =\, \left(
\begin{array}{c|c}
 z+1 & \omega  \\
 \hline
 \omega  & \zeta  \\
\end{array}
\right) \quad  &; \quad \mathbf{FL}\left[T_2,Z\right] & = \, \left(
\begin{array}{c|c}
 z & \omega  \nonumber\\
 \hline
 \omega  & \zeta +1 \\
\end{array}
\right) \quad\\
  \mathbf{FL}\left[T_3,Z\right] & =\,  \left(
\begin{array}{c|c}
 z & \omega +1 \\
 \hline
 \omega +1 & \zeta  \\
\end{array}
\right)\quad  &; \quad \mathbf{FL}\left[T_4,Z\right] & = \, \left(
\begin{array}{c|c}
 \zeta +2 \omega +z & \zeta +\omega  \\
 \hline
 \zeta +\omega  & \zeta  \\
\end{array}
\right) \quad
\end{alignat}
The above defined generators are all of infinite order  and create a regular grid of points in the three complex planes. Since we are interested in
creating two copies of the modular group, one acting on the plane $z$, and  the other acting on the plane $\zeta$, and since the modular group is
  generated by a generator $T$ of infinite order and by another $S$ of order 2 such that $(ST)^3= \mathrm{Id}$, we had to find two extra
generators $S_1,S_2$ of our discrete  group, to be extracted from the compact subgroup $\mathrm{SU(2) \times U(1) \subset Sp(4,\mathbb{R})}$, such
that
\begin{equation}\label{coronatasso}
  \left[S_1\, , \,S_2\right] \, = \, 0 \quad S_1^2 \, = \, S_2^2\, = \, \mathrm{Id} \quad ; \quad (S_1\cdot T_1)^2 \,= \,(S_2\cdot T_2)^2 \, = \, \mathrm{Id}.
\end{equation}
%%%%%%%%%%%%%%%%%%%%%%%%%%%%%%%%%%%%%%%%%
In the above equation one should be careful that the relation $\dots \, = \, \mathrm{Id} $ must be understood in the projective sense and not
necessarily matrixwise. In other words what matters is that the matrix product should act as the identity in the fractional linear
transformation. With some ingenuity we found:
\begin{equation}\label{caterino}
  S_1 \, = \, \left(
\begin{array}{cccc}
 0 & 0 & 1 & 0 \\
 0 & 1 & 0 & 0 \\
 -1 & 0 & 0 & 0 \\
 0 & 0 & 0 & 1 \\
\end{array}
\right) \quad ; \quad S_2 \, = \, \left(
\begin{array}{cccc}
 1 & 0 & 0 & 0 \\
 0 & 0 & 0 & 1 \\
 0 & 0 & 1 & 0 \\
 0 & -1 & 0 & 0 \\
\end{array}
\right)
\end{equation}
In this way we obtained the goal of introducing two copies of $\mathrm{SL(2,\mathbb{Z})}$ in our candidate group, yet if we wanted to include also
the generators $T_3,T_4$ we had to find additional elliptic generators, namely included in $\mathrm{SU(2)\times U(1)}$ able to create the inverse
of $T_3,T_4$.
\par
After some attempts and restricting our attention to integer-valued matrices we found the following two generators:
\begin{equation}\label{panegirico}
  Q_1 \, = \,  \left(
\begin{array}{cccc}
 0 & -1 & 0 & 0 \\
 1 & 0 & 0 & 0 \\
 0 & 0 & 0 & -1 \\
 0 & 0 & 1 & 0 \\
\end{array}
\right)         \quad  ;  \quad Q_2 \, = \, \left(
\begin{array}{cccc}
 0 & 0 & 0 & 1 \\
 0 & 0 & 1 & 0 \\
 0 & -1 & 0 & 0 \\
 -1 & 0 & 0 & 0 \\
\end{array}
\right)
\end{equation}
which, in a proper matrix sense, satisfy the   relations
\begin{equation}\label{corego}
  Q_1^4 \, = \, Q_2^4 \, = \,  Q_1\, T_3 \, Q_1^3\, = \,  Q_2\, T_4 \, Q_2^3\, = \, \mathbf{1}_{4\times4}.
\end{equation}
Eq.~\eqref{corego} demonstrates that by themselves $\{Q_1,Q_2, T_3, T_4\}$ generate a group, because the inverse of any word can be found among the
words themselves.
\par
The transformations associated with the four elliptic generators so far introduced  are the following:
\begin{alignat}{3}\label{inversiones}
  \mathbf{FL}\left[S_1,Z\right] & =  \left(
\begin{array}{c|c}
 -\frac{1}{z} & -\frac{\omega }{z} \\
 \hline
 -\frac{\omega }{z} & \zeta -\frac{\omega ^2}{z} \\
\end{array}
\right) \quad & \text{order of the transformation} & = 4\, \nonumber\\
  \mathbf{FL}\left[S_2,Z\right] & = \left(
\begin{array}{c|c}
 z-\frac{\omega ^2}{\zeta } & -\frac{\omega }{\zeta } \\
 \hline
 -\frac{\omega }{\zeta } & -\frac{1}{\zeta } \\
\end{array}
\right) \quad & \text{order of the transformation} & = 4\, \nonumber\\
  \mathbf{FL}\left[Q_1,Z\right] & =  \left(
\begin{array}{c|c}
 \zeta  & -\omega  \\
 \hline
 -\omega  & z \\
\end{array}
\right) \quad & \text{order of the transformation} & = 2\, \nonumber\\
  \mathbf{FL}\left[Q_2,Z\right] & =  \left(
\begin{array}{c|c}
 \frac{z}{\omega ^2-\zeta  z} & \frac{\omega }{\zeta  z-\omega ^2} \\
 \hline
 \frac{\omega }{\zeta  z-\omega ^2} & \frac{\zeta }{\omega ^2-\zeta  z} \\
\end{array}
\right)\quad & \text{order of the transformation} & = 2\,
\end{alignat}
\paragraph{\sc The finite elliptic subgroup.}
Having gone so far, the next interesting question  is what is the finite group generated by the four  elliptic generators introduced. The answer
was obtained utilizing the appropriate routines included in the background {\sc mathematica}  NoteBook \emph{NonComSymSpacNeuNet.nb} devoted to the
generation of a finite group from a set of elements, and to the subsequent structural analysis of the same by its reorganization into conjugation
classes.
\par
The result is a  finite subgroup  $\mathfrak{P}_{32} \subset \mathrm{Sp(4,Z)}$ containing $32$ elements, as the given name indicates. The order of
the generators is, depending on the conjugacy class\footnote{There are 14 conjugacy classes} , one of the numbers $2,4,8$. With some effort we
were able to identify two generators respectively, called $\mathcal{A},\mathcal{B}$, which  are the following:
\begin{alignat}{3}
  \mathcal{A} &=& Q_1\, S_1 & = \left(
\begin{array}{cccc}
 0 & -1 & 0 & 0 \\
 0 & 0 & 1 & 0 \\
 0 & 0 & 0 & -1 \\
 -1 & 0 & 0 & 0 \\
\end{array}
\right) \\
  \mathcal{B} &=& Q_2^2 \, Q_2 & = \left(
\begin{array}{cccc}
 0 & 0 & 0 & -1 \\
 0 & 0 & -1 & 0 \\
 0 & 1 & 0 & 0 \\
 1 & 0 & 0 & 0 \\
\end{array}
\right).
\end{alignat}
 They satisfy the following relations as $\mathrm{Sp(4,\mathbb{Z})}$ matrices:
 \begin{equation}\label{parlatus}
   \mathcal{A}^8 \, = \, \mathcal{B}^4 \, = \, \left(\mathcal{B\, A}\right)^4 \, = \, \mathbf{1}_{4\times 4}.
 \end{equation}
 The explicit action of the elliptic subgroup generators on the Siegel plane is instead the following one:
 \begin{alignat}{3}\label{Abactiones}
  \mathbf{FL}\left[\mathcal{A},Z\right] & =  \left(
\begin{array}{c|c}
 \zeta -\frac{\omega ^2}{z} & \frac{\omega }{z} \\
 \hline
 \frac{\omega }{z} & -\frac{1}{z} \\
\end{array}
\right)\quad & \text{order of the transformation} & = 4\,; \nonumber\\
  \mathbf{FL}\left[\mathcal{B},Z\right] & = \left(
\begin{array}{c|c}
 \frac{z}{\omega ^2-\zeta  z} & \frac{\omega }{\zeta  z-\omega ^2} \\
 \frac{\omega }{\zeta  z-\omega ^2} & \frac{\zeta }{\omega ^2-\zeta  z} \\
\end{array}
\right) \quad & \text{order of the transformation} & = 2\,; \nonumber\\
  \mathbf{FL}\left[\mathcal{B} \,\mathcal{A},Z\right] & = \left(
\begin{array}{c|c}
 z-\frac{\omega ^2}{\zeta } & -\frac{\omega }{\zeta } \\
 \hline
 -\frac{\omega }{\zeta } & -\frac{1}{\zeta } \\
\end{array}
\right) \quad & \text{order of the transformation} & = 2\, .
\end{alignat}
A complete list of the 32 group elements of $\mathfrak{P}_{32}$   is provided in Table  \ref{32elem}.
%%%%%%%%%%%%%%%%%%%%%%%%%%%%%%%%%%%%%%%%%%%%%%%%%%%%%%%%%%%%%%%%%%
\begin{table}
\begin{center}
{\tiny{
$$
\begin{array}{||c|rclcc||c|rclcc||}
\hline\hline
1& \mathcal{A} &=& \mathcal{Q}_1.\mathcal{S}_1 &  = & \left(
\begin{array}{cccc}
 0 & -1 & 0 & 0 \\
 0 & 0 & 1 & 0 \\
 0 & 0 & 0 & -1 \\
 -1 & 0 & 0 & 0 \\
\end{array}
\right) & 17 &\mathcal{A}.\mathcal{B}.\mathcal{B}.\mathcal{B} &=&
   \mathcal{Q}_1.\mathcal{Q}_1.\mathcal{S}_1 &  = & \left(
\begin{array}{cccc}
 0 & 0 & -1 & 0 \\
 0 & -1 & 0 & 0 \\
 1 & 0 & 0 & 0 \\
 0 & 0 & 0 & -1 \\
\end{array}
\right) \\
\hline
 2 &\mathcal{B} &=& \mathcal{Q}_1.\mathcal{Q}_1.\mathcal{Q}_2 &
    = & \left(
\begin{array}{cccc}
 0 & 0 & 0 & -1 \\
 0 & 0 & -1 & 0 \\
 0 & 1 & 0 & 0 \\
 1 & 0 & 0 & 0 \\
\end{array}
\right) & 18 &\mathcal{B}.\mathcal{A}.\mathcal{B}.\mathcal{A} &=&
   \mathcal{S}_2.\mathcal{S}_2 &  = & \left(
\begin{array}{cccc}
 1 & 0 & 0 & 0 \\
 0 & -1 & 0 & 0 \\
 0 & 0 & 1 & 0 \\
 0 & 0 & 0 & -1 \\
\end{array}
\right) \\
\hline
3 & \mathcal{A}.\mathcal{A} &=&
   \mathcal{Q}_1.\mathcal{Q}_2.\mathcal{S}_2.\mathcal{S}_2 &  = & \left(
\begin{array}{cccc}
 0 & 0 & -1 & 0 \\
 0 & 0 & 0 & -1 \\
 1 & 0 & 0 & 0 \\
 0 & 1 & 0 & 0 \\
\end{array}
\right) & 19 &\mathcal{B}.\mathcal{A}.\mathcal{B}.\mathcal{B} &=&
   \mathcal{Q}_1.\mathcal{Q}_1.\mathcal{S}_2 &  = & \left(
\begin{array}{cccc}
 -1 & 0 & 0 & 0 \\
 0 & 0 & 0 & -1 \\
 0 & 0 & -1 & 0 \\
 0 & 1 & 0 & 0 \\
\end{array}
\right) \\
\hline
4 & \mathcal{A}.\mathcal{B} &=& \mathcal{S}_1 &  = & \left(
\begin{array}{cccc}
 0 & 0 & 1 & 0 \\
 0 & 1 & 0 & 0 \\
 -1 & 0 & 0 & 0 \\
 0 & 0 & 0 & 1 \\
\end{array}
\right) & 20 &\mathcal{B}.\mathcal{B}.\mathcal{B}.\mathcal{B} &=&
   \mathcal{Q}_1.\mathcal{Q}_2.\mathcal{Q}_2.\mathcal{Q}_1 &  = & \left(
\begin{array}{cccc}
 1 & 0 & 0 & 0 \\
 0 & 1 & 0 & 0 \\
 0 & 0 & 1 & 0 \\
 0 & 0 & 0 & 1 \\
\end{array}
\right) \\
\hline
5 & \mathcal{B}.\mathcal{A} &=& \mathcal{S}_2 &  = & \left(
\begin{array}{cccc}
 1 & 0 & 0 & 0 \\
 0 & 0 & 0 & 1 \\
 0 & 0 & 1 & 0 \\
 0 & -1 & 0 & 0 \\
\end{array}
\right) & 21 &\mathcal{A}.\mathcal{A}.\mathcal{A}.\mathcal{B}.\mathcal{A} &=&
   \mathcal{Q}_1 &  = & \left(
\begin{array}{cccc}
 0 & -1 & 0 & 0 \\
 1 & 0 & 0 & 0 \\
 0 & 0 & 0 & -1 \\
 0 & 0 & 1 & 0 \\
\end{array}
\right) \\
\hline
 6 &\mathcal{B}.\mathcal{B} &=& \mathcal{Q}_2.\mathcal{Q}_2 &  = & \left(
\begin{array}{cccc}
 -1 & 0 & 0 & 0 \\
 0 & -1 & 0 & 0 \\
 0 & 0 & -1 & 0 \\
 0 & 0 & 0 & -1 \\
\end{array}
\right) & 22 & \mathcal{A}.\mathcal{A}.\mathcal{A}.\mathcal{B}.\mathcal{B} &=&
   \mathcal{Q}_1.\mathcal{S}_1.\mathcal{S}_1.\mathcal{S}_2 &  = & \left(
\begin{array}{cccc}
 0 & 0 & 0 & -1 \\
 -1 & 0 & 0 & 0 \\
 0 & 1 & 0 & 0 \\
 0 & 0 & -1 & 0 \\
\end{array}
\right) \\
\hline
7 & \mathcal{A}.\mathcal{A}.\mathcal{A} &=&
   \mathcal{Q}_1.\mathcal{S}_2.\mathcal{S}_2.\mathcal{S}_2
   &  = & \left(
\begin{array}{cccc}
 0 & 0 & 0 & 1 \\
 1 & 0 & 0 & 0 \\
 0 & -1 & 0 & 0 \\
 0 & 0 & 1 & 0 \\
\end{array}
\right) & 23 & \mathcal{A}.\mathcal{A}.\mathcal{B}.\mathcal{A}.\mathcal{B} &=&
   \mathcal{Q}_1.\mathcal{S}_1.\mathcal{S}_1.\mathcal{S}_1
  &  = & \left(
\begin{array}{cccc}
 0 & -1 & 0 & 0 \\
 0 & 0 & -1 & 0 \\
 0 & 0 & 0 & -1 \\
 1 & 0 & 0 & 0 \\
\end{array}
\right) \\
\hline
8 & \mathcal{A}.\mathcal{A}.\mathcal{B} &=&
   \mathcal{Q}_1.\mathcal{S}_1.\mathcal{S}_1 &  = & \left(
\begin{array}{cccc}
 0 & -1 & 0 & 0 \\
 -1 & 0 & 0 & 0 \\
 0 & 0 & 0 & -1 \\
 0 & 0 & -1 & 0 \\
\end{array}
\right) & 24 &\mathcal{A}.\mathcal{A}.\mathcal{B}.\mathcal{B}.\mathcal{B} &=&
   \mathcal{Q}_1.\mathcal{S}_2.\mathcal{S}_2 &  = & \left(
\begin{array}{cccc}
 0 & 1 & 0 & 0 \\
 1 & 0 & 0 & 0 \\
 0 & 0 & 0 & 1 \\
 0 & 0 & 1 & 0 \\
\end{array}
\right) \\
\hline
9 & \mathcal{A}.\mathcal{B}.\mathcal{A} &=&
   \mathcal{Q}_1.\mathcal{S}_1.\mathcal{S}_2 &  = & \left(
\begin{array}{cccc}
 0 & 0 & 0 & -1 \\
 0 & 0 & 1 & 0 \\
 0 & 1 & 0 & 0 \\
 -1 & 0 & 0 & 0 \\
\end{array}
\right) & 25 &\mathcal{A}.\mathcal{B}.\mathcal{A}.\mathcal{B}.\mathcal{A} &=&
   \mathcal{Q}_1.\mathcal{S}_1.\mathcal{S}_2.\text{$\mathcal{S}$
   2} &  = & \left(
\begin{array}{cccc}
 0 & 1 & 0 & 0 \\
 0 & 0 & 1 & 0 \\
 0 & 0 & 0 & 1 \\
 -1 & 0 & 0 & 0 \\
\end{array}
\right) \\
\hline
10 & \mathcal{A}.\mathcal{B}.\mathcal{B} &=&
   \mathcal{Q}_1.\mathcal{Q}_2.\mathcal{Q}_2.\mathcal{S}_1 &  = & \left(
\begin{array}{cccc}
 0 & 1 & 0 & 0 \\
 0 & 0 & -1 & 0 \\
 0 & 0 & 0 & 1 \\
 1 & 0 & 0 & 0 \\
\end{array}
\right) & 26 &\mathcal{A}.\mathcal{B}.\mathcal{A}.\mathcal{B}.\mathcal{B} &=&
   \mathcal{Q}_2.\mathcal{S}_1.\mathcal{S}_1 &  = & \left(
\begin{array}{cccc}
 0 & 0 & 0 & 1 \\
 0 & 0 & -1 & 0 \\
 0 & -1 & 0 & 0 \\
 1 & 0 & 0 & 0 \\
\end{array}
\right) \\
\hline
11 &  \mathcal{B}.\mathcal{A}.\mathcal{B} &=& \mathcal{Q}_1.\mathcal{S}_2 &
    = & \left(
\begin{array}{cccc}
 0 & 0 & 0 & -1 \\
 1 & 0 & 0 & 0 \\
 0 & 1 & 0 & 0 \\
 0 & 0 & 1 & 0 \\
\end{array}
\right) & 27 &\mathcal{B}.\mathcal{A}.\mathcal{B}.\mathcal{B}.\mathcal{B} &=&
   \mathcal{Q}_2.\mathcal{S}_1 &  = & \left(
\begin{array}{cccc}
 0 & 0 & 0 & 1 \\
 -1 & 0 & 0 & 0 \\
 0 & -1 & 0 & 0 \\
 0 & 0 & -1 & 0 \\
\end{array}
\right) \\
\hline
12 & \mathcal{B}.\mathcal{B}.\mathcal{B} &=& \mathcal{Q}_2 &  = & \left(
\begin{array}{cccc}
 0 & 0 & 0 & 1 \\
 0 & 0 & 1 & 0 \\
 0 & -1 & 0 & 0 \\
 -1 & 0 & 0 & 0 \\
\end{array}
\right) & 28 &\mathcal{A}.\mathcal{A}.\mathcal{A}.\mathcal{B}.\mathcal{A}.\mathcal{B} &=&
   \mathcal{Q}_2.\mathcal{Q}_1 &  = & \left(
\begin{array}{cccc}
 0 & 0 & 1 & 0 \\
 0 & 0 & 0 & -1 \\
 -1 & 0 & 0 & 0 \\
 0 & 1 & 0 & 0 \\
\end{array}
\right) \\
\hline
13 & \mathcal{A}.\mathcal{A}.\mathcal{A}.\mathcal{B} &=&
   \mathcal{S}_2.\mathcal{S}_2.\mathcal{S}_2 &  = & \left(
\begin{array}{cccc}
 1 & 0 & 0 & 0 \\
 0 & 0 & 0 & -1 \\
 0 & 0 & 1 & 0 \\
 0 & 1 & 0 & 0 \\
\end{array}
\right) & 29 &\mathcal{A}.\mathcal{A}.\mathcal{A}.\mathcal{B}.\mathcal{B}.\mathcal{B} &=&
   \mathcal{S}_1.\mathcal{S}_1.\mathcal{S}_2 &  = & \left(
\begin{array}{cccc}
 -1 & 0 & 0 & 0 \\
 0 & 0 & 0 & 1 \\
 0 & 0 & -1 & 0 \\
 0 & -1 & 0 & 0 \\
\end{array}
\right) \\
\hline
14 & \mathcal{A}.\mathcal{A}.\mathcal{B}.\mathcal{A} &=&
   \mathcal{S}_1.\mathcal{S}_1.\mathcal{S}_1 &  = & \left(
\begin{array}{cccc}
 0 & 0 & -1 & 0 \\
 0 & 1 & 0 & 0 \\
 1 & 0 & 0 & 0 \\
 0 & 0 & 0 & 1 \\
\end{array}
\right) & 30 & \mathcal{A}.\mathcal{A}.\mathcal{B}.\mathcal{A}.\mathcal{B}.\mathcal{A} &=&
   \mathcal{Q}_1.\mathcal{Q}_2 &  = & \left(
\begin{array}{cccc}
 0 & 0 & -1 & 0 \\
 0 & 0 & 0 & 1 \\
 1 & 0 & 0 & 0 \\
 0 & -1 & 0 & 0 \\
\end{array}
\right) \\
\hline
 15 &\mathcal{A}.\mathcal{A}.\mathcal{B}.\mathcal{B} &=&
   \mathcal{S}_1.\mathcal{S}_2 &  = & \left(
\begin{array}{cccc}
 0 & 0 & 1 & 0 \\
 0 & 0 & 0 & 1 \\
 -1 & 0 & 0 & 0 \\
 0 & -1 & 0 & 0 \\
\end{array}
\right) & 31 &\mathcal{A}.\mathcal{A}.\mathcal{B}.\mathcal{A}.\mathcal{B}.\mathcal{B} &=&
   \mathcal{S}_1.\mathcal{S}_2.\mathcal{S}_2 &  = & \left(
\begin{array}{cccc}
 0 & 0 & 1 & 0 \\
 0 & -1 & 0 & 0 \\
 -1 & 0 & 0 & 0 \\
 0 & 0 & 0 & -1 \\
\end{array}
\right) \\
\hline
16 &  \mathcal{A}.\mathcal{B}.\mathcal{A}.\mathcal{B} &=&
   \mathcal{S}_1.\mathcal{S}_1 &  = & \left(
\begin{array}{cccc}
 -1 & 0 & 0 & 0 \\
 0 & 1 & 0 & 0 \\
 0 & 0 & -1 & 0 \\
 0 & 0 & 0 & 1 \\
\end{array}
\right) & 32 &
   \mathcal{A}.\mathcal{A}.\mathcal{A}.\mathcal{B}.\mathcal{A}.\mathcal{B}.\mathcal{B} &=&
   \mathcal{Q}_1.\mathcal{Q}_2.\mathcal{Q}_2 &  = & \left(
\begin{array}{cccc}
 0 & 1 & 0 & 0 \\
 -1 & 0 & 0 & 0 \\
 0 & 0 & 0 & 1 \\
 0 & 0 & -1 & 0 \\
\end{array}
\right) \\
\hline\hline
\end{array}
$$}}
\end{center}
\caption{\label{32elem} The complete list of the 32 elements of the elliptic subgroup $\mathfrak{P}_{32}$ of
the discrete subgroup $\Delta_{[32,8]} \subset \mathrm{Sp(4,\mathbb{Z})}$. In the column after the progressive number
one finds the identification of the group element as a word in terms of the generators $\mathcal{A},\mathcal{B}$,
while in the next column there is the identification of the same element as a word in the generators $\mathcal{S}_{1,2},\mathcal{Q}_{1,2}$.
The last column provides the explicit expression of the element as an integer-valued $4\times 4 $ symplectic matrix, namely
as an element of $\mathrm{Sp(4,\mathbb{Z})}$}
\end{table}
\paragraph{\sc The normal parabolic subgroup.}
Having constructed the elliptic point group $\mathfrak{P}_{32}$, we turn our attention to the parabolic translation generators that we already
introduced in such a way that special elements of the point group can construct their inverse.  We take the set $\mathfrak{S}$ of $8$ elements
composed by $T_{1,2,3,4}$ plus their inverses $T^{-1}_{1,2,3,4}$  and we create the orbit of that set with respect to the conjugation with respect
to all elements of $\mathfrak{P}_{32}$. We discover that such orbit is a set $\mathfrak{U}$ containing 16 elements and that the inverse of each of
the elements  in $\mathfrak{U}$   is contained in $\mathfrak{U}$ as well. Hence we conclude that we can single out a subset of  $\mathfrak{G}
\subset \mathfrak{U}$  containing 8 elements and no inverse of any of them. Clearly the complement of $\mathfrak{G}$ in $\mathfrak{U}$ contains
the $8$ inverses.
\par
In Table \ref{8genni} the list of the $8$ parabolic generators of order $\infty$ forming the set $\mathfrak{G}$  is  presented in explicit matrix
form as elements of $\mathrm{Sp(4,\mathbb{Z})}$.
\par
Next it is very important to emphasize another series of relations satisfied by the generators $T_i$ not involving the elliptic generators of the
finite subgroup $\mathfrak{P}_{32}$.
%%%%%%%%%%%%%%%%%%%%%%%%%%%%%%%%%%%%%%%%%%%%%%
% Altra tavola
%%%%%%%%%%%%%%%%%%%%%%%%%%%%%%%%%%%%%%%%%%%%%%
\begin{table}
\begin{center}
$$
\begin{array}{||ccc||ccc||}
\hline\hline
 T_1 & = & \left(
\begin{array}{cccc}
 1 & 0 & 1 & 0 \\
 0 & 1 & 0 & 0 \\
 0 & 0 & 1 & 0 \\
 0 & 0 & 0 & 1 \\
\end{array}
\right) & T_2 & = & \left(
\begin{array}{cccc}
 1 & 0 & 0 & 0 \\
 0 & 1 & 0 & 1 \\
 0 & 0 & 1 & 0 \\
 0 & 0 & 0 & 1 \\
\end{array}
\right)\\
\hline
 T_3 & = & \left(
\begin{array}{cccc}
 1 & 0 & 0 & 1 \\
 0 & 1 & 1 & 0 \\
 0 & 0 & 1 & 0 \\
 0 & 0 & 0 & 1 \\
\end{array}
\right) & T_4 & = & \left(
\begin{array}{cccc}
 1 & 1 & 0 & 0 \\
 0 & 1 & 0 & 0 \\
 0 & 0 & 1 & 0 \\
 0 & 0 & -1 & 1 \\
\end{array}
\right) \\
\hline
 T_5 & = & \left(
\begin{array}{cccc}
 1 & 0 & 0 & 0 \\
 -1 & 1 & 0 & 0 \\
 0 & 0 & 1 & 1 \\
 0 & 0 & 0 & 1 \\
\end{array}
\right) & T_6 & = & \left(
\begin{array}{cccc}
 1 & 0 & 0 & 0 \\
 0 & 1 & 0 & 0 \\
 -1 & 0 & 1 & 0 \\
 0 & 0 & 0 & 1 \\
\end{array}
\right)\\
\hline
 T_7 & = & \left(
\begin{array}{cccc}
 1 & 0 & 0 & 0 \\
 0 & 1 & 0 & 0 \\
 0 & -1 & 1 & 0 \\
 -1 & 0 & 0 & 1 \\
\end{array}
\right)& T_8 & = & \left(
\begin{array}{cccc}
 1 & 0 & 0 & 0 \\
 0 & 1 & 0 & 0 \\
 0 & 0 & 1 & 0 \\
 0 & -1 & 0 & 1 \\
\end{array}
\right) \\
\hline\hline
\end{array}
$$
\end{center}
\caption{The 8 generators  of the normal subgroup  $\mathfrak{T}_{8}$ of the discrete subgroup $\Delta_{[32,8]} \subset \mathrm{Sp(4,\mathbb{Z})}$. \label{8genni}}
\end{table}
Such relations are the following ones:
\begin{alignat}{8}\label{granduca}
 (T_1\, T_6)^6 & = &  (T_2\, T_8)^6 & = &  (T_3\, T_7)^6 & =  & (T_4\, T_5)^6 & =  \mathbf{1}_{4\times 4}
\end{alignat}
and they are extremely important. They suffice to show that the inverse of each of the 8 generators can be written as a word in the same letters
and from that it follows also that any word in $T_i$ has an inverse that can be expressed as a word in the same symbols. The infinite set of words
written in terms of the eight symbols $T_i$ forms a group $\mathfrak{T}$. The finite group $\mathfrak{P}_{32}$ is by construction an automorphism
group of $\mathfrak{T}$, the action of  $\mathfrak{P}_{32}$ on $\mathfrak{T}$ being by conjugation:
\begin{equation}\label{construzza}
  \forall \gamma \in \mathfrak{P}_{32} \quad : \quad C_\gamma \, : \, \mathfrak{T} \longrightarrow \mathfrak{T} \quad ; \quad \forall t\in
  \mathfrak{T}  \, : \,  C_\gamma(t) \, = \, \gamma\cdot t \cdot \gamma^{-1} \in \mathfrak{T}
\end{equation}
Because of the above property the set of all elements of $\mathrm{Sp(4,\mathbb{Z})}$ of the form $\gamma \cdot t$ where $\gamma \in
\mathfrak{P}_{32}$ and $t \in \mathrm{T}$ is a group since also their product and their inverse can be written in the same way. We name such a
group $\Delta_{32,8}$.
\par
The group $\mathfrak{T}$  is a normal subgroup in $\Delta_{32,8}$.
\section{A family of subgroups of ${\rm SO}(r,r+q,\,\mathbb{Z})$}
\label{generaloneSOpq}
We wish to generalize the construction of the subgroup $\Delta_{32,8}$ of ${\rm Sp}(4,\mathbb{Z})$ discussed in the previous section  to a discrete subgroup $\hat{\Delta}^{[r,q]}$ of ${\rm SO}(r,r+q)$.
Let us introduce the following notation for the generators of $Solv$ in $\mathfrak{so}(r,r+q)$:
\begin{align}
  &  H_i=H_{\boldsymbol{\epsilon}_i},\,E_{\boldsymbol{\epsilon}_i\pm \boldsymbol{\epsilon}_j},\,E_{\boldsymbol{\epsilon}_i^I}\,\,;\,\,\,\, i=1,\dots, r\,\,,\,\,\,\,I=1,\dots q\,,
\end{align}
where $\boldsymbol{\epsilon}_i$ is an orthonormal basis and $\boldsymbol{\epsilon}_i\pm \boldsymbol{\epsilon}_j$  are the long positive roots $\boldsymbol{\alpha}$, while $\boldsymbol{\epsilon}_i^I$ are the restricted short positive roots $\boldsymbol{\sigma}^I$; the index $I$, spanning the fundamental representation of the paint group ${\rm SO}(q)$, labels their multiplicity. These generators are normalized to satisfy the   conditions
\begin{align}
    [H_i,\,E_{\boldsymbol{\alpha}}]]&= \alpha_i\,E_{\boldsymbol{\alpha}}\,,\,\,\,[H_i,\,E_{\boldsymbol{\sigma}^I}]] =\sigma^I_i\,E_{\boldsymbol{\sigma}^I}\,,\nonumber\\
    [E_{\boldsymbol{\alpha}},\,E_{-\boldsymbol{\alpha}}]&=4\,\alpha_iH_i\,,\,\,\,[E_{\boldsymbol{\sigma}^I},\,E_{-\boldsymbol{\sigma}^I}]=2\,\sigma^I_iH_i\,,\nonumber\\
    [E_{\boldsymbol{\epsilon}_i^I},\,E_{\boldsymbol{\epsilon}_j^J}]&=-\delta_{IJ}\,E_{\boldsymbol{\epsilon}_i+\boldsymbol{\epsilon}_j}\,,\,\,\,[E_{\boldsymbol{\epsilon}_i^I},\,E_{-\boldsymbol{\epsilon}_j^J}]=-\delta_{IJ}\,E_{\boldsymbol{\epsilon}_i-\boldsymbol{\epsilon}_j}\,,
\end{align}
where we take $E_{-\boldsymbol{\alpha}}=E_{\boldsymbol{\alpha}}^T$ and  $E_{-\boldsymbol{\sigma}^I}=E_{\boldsymbol{\sigma}^I}^T$.
Next, we consider, in the fundamental representation of ${\rm SO}(r,r+q)$, the following integer generators with integer generators
\begin{align}
    J_{ij}^{\pm}&\equiv \exp\left(\frac{\pi}{4}(E_{\boldsymbol{\epsilon}_i\pm\boldsymbol{\epsilon}_j}-E_{\boldsymbol{\epsilon}_i\pm\boldsymbol{\epsilon}_j}^T)\right)\,,\nonumber\\
    J_i^I&=\exp\left(\frac{\pi}{2}(E_{\boldsymbol{\epsilon}_i^I}-E_{\boldsymbol{\epsilon}_i^I}^T)\right)\,.
\end{align}
These matrices are all elliptic and generate the group $\hat{\mathfrak{P}}^{[r,q]}$,  which we shall call \emph{generalized Weyl group} \cite{Fre:2011ns}.
They satisfy the conditions
$$(J_i^I)^2=(J_{ij}^{\pm})^4={\bf 1}\,.$$
We also define the following  parabolic generators with integer entries
\begin{equation}
    T_{ij}^{\pm}\equiv \exp\left(2E_{\boldsymbol{\epsilon}_i\pm\boldsymbol{\epsilon}_j}\right)\,\,;\,\,\,T_i^I\equiv \exp\left(2E_{\boldsymbol{\epsilon}_i^I}\right)\,.
\end{equation}
One can verify that the following relations hold:
\begin{align}
(J_{ij}^{\pm})^{-1}\cdot T_{ij}^{\pm}\cdot J_{ij}^{\pm}=(T_{ij}^{\pm})^{-T}\,\,,\,\,\,\,(J_{i}^{I})^{-1}\cdot T_{i}^{I}\cdot J_{i}^{I}=(T_{}^{I})^{-T}\,.
\end{align}
In general, one can verify that:
\begin{equation}
(T_i^I)^{-1},\,(T_i^I)^{T},\,(T_i^I)^{-T}\,\in \,\hat{\mathfrak{P}}^{[r,q]\,-1}\cdot T_i^I\cdot \hat{\mathfrak{P}}^{[r,q]}\,.
\end{equation}
In order to generate the inverses of the $T_{ij}^{\pm}$ we can use the following relations:
\begin{align}
\left(T_i^I\cdot T_j^I\cdot (T_i^I)^{-1}\cdot(T_j^I)^{-1} \right)\cdot T_{ij}^{+}&=(T_{ij}^{+})^{-1}\,,\nonumber\\
\left(T_i^I\cdot (T_j^I)^T\cdot (T_i^I)^{-1}\cdot(T_j^I)^{-T} \right)\cdot T_{ij}^{-}&=(T_{ij}^{-})^{-1}\,. \label{newrelations}
\end{align}
The group $\hat{\Delta}^{[r,q]}$ is generated by $J_{ij}^\pm,\,J_i^I,\,T_{ij}^{\pm},\,T_i^I$. In Appendix \ref{peritonite} we give the explicit
matrix form of the generators for the cases $r=2,q=1$ and $r=3,q=1$. Here we wish to show how the generators of $\Delta_{32,8}$, discussed in the
previous section, fit  general construction of the present section. One can verify that
\begin{align}
T_1&=\exp\left(E^{[{\rm Sp}]}_{\boldsymbol{\epsilon}_1+\boldsymbol{\epsilon}_2}\right)\,\,,\,\,\,T_2=\exp\left(E^{[{\rm Sp}]}_{\boldsymbol{\epsilon}_1-\boldsymbol{\epsilon}_2}\right)\,\,,\,\,\,T_3=\exp\left(E^{[{\rm Sp}]}_{\boldsymbol{\epsilon}_1}\right)\,\,,\,\,\,T_4=\exp\left(E^{[{\rm Sp}]}_{\boldsymbol{\epsilon}_2}\right)\,,\nonumber\\
T_5&=T_4^{-T}\,\,,\,\,T_6=T_1^{-T}\,\,,\,\,\,T_7=T_3^{-T}\,\,,\,\,\,T_8=T_2^{-T},
\end{align}
where $E^{[{\rm Sp}]}_{\boldsymbol{\alpha}}$ denotes the generator $E_{\boldsymbol{\alpha}}$ in the fundamental 4-dimensional representation of ${\rm Sp}(4,\mathbb{Z})$.
The Cartan generators in the same representation read
$$H_1^{[{\rm Sp}]}=\frac{1}{2}\,\left(\begin{matrix}1 & 0 & 0 & 0\cr 0 & 1 & 0 & 0\cr 0 & 0 & -1 & 0\cr 0 & 0 & 0 & -1\end{matrix}\right)\,\,,\,\,\,H_2^{[{\rm Sp}]}=\frac{1}{2}\,\left(\begin{matrix}1 & 0 & 0 & 0\cr 0 & -1 & 0 & 0\cr 0 & 0 & -1 & 0\cr 0 & 0 & 0 & 1\end{matrix}\right)\,.$$
As for the elliptic generators of $\mathfrak{P}_{32}$ we have:
\begin{align}
    Q_1&=\exp\left(\frac{\pi}{2}\,(E^{[{\rm Sp}]}_{\boldsymbol{\epsilon}_2}-E^{[{\rm Sp}]}_{-\boldsymbol{\epsilon}_2})\right)\,\,,\,\,\,Q_2=\exp\left(\frac{\pi}{2}\,(E^{[{\rm Sp}]}_{\boldsymbol{\epsilon}_1}-E^{[{\rm Sp}]}_{-\boldsymbol{\epsilon}_1})\right)\,,\nonumber\\
    S_1&=\exp\left(\frac{\pi}{2}\,(E^{[{\rm Sp}]}_{\boldsymbol{\epsilon}_1+\boldsymbol{\epsilon}_2}-E^{[{\rm Sp}]}_{-\boldsymbol{\epsilon}_1-\boldsymbol{\epsilon}_2})\right)\,\,,\,\,\, S_2=\exp\left(\frac{\pi}{2}\,(E^{[{\rm Sp}]}_{\boldsymbol{\epsilon}_1-\boldsymbol{\epsilon}_2}-E^{[{\rm Sp}]}_{-\boldsymbol{\epsilon}_1+\boldsymbol{\epsilon}_2})\right)
\end{align}
Note the different powers in the definition of the generators in terms of exponentials. This is related to ${\rm Sp}(4,\mathbb{R})$ being the double cover of ${\rm SO}(2,3)$. Aside from this, $T_1$ corresponds to $T_{12}^+$, $T_2$ to $T_{12}^-$, $T_3$ to $T_{1}^{I=1}$ and $T_4$ to $T_{2}^{I=1}$, $Q_{1}$ to $J_2^{I=1}$, $Q_{2}$ to $J_1^{I=1}$, $S_{1}$ to $J_{12}^{+}$, $S_{2}$ to $J_{12}^{-}$. Note that for the parabolic generators of $\Delta_{32,8}$, we find relations which are analogous to \eqref{newrelations}:
\begin{equation}
    T_3\cdot T_4\cdot T_3^{-1}\cdot T_4^{-1}=T_1^{-2}\,,\,\,\, T_3\cdot T_4^T\cdot T_3^{-1}\cdot T_4^{-T}=T_2^{-2}\,.
\end{equation}
We observe that the order of $\hat{\mathfrak{P}}^{[2,1]}\subset {\rm SO}(2,3)$ is 16, half the order of ${\mathfrak{P}}_{32}$.\par
In general, we find:
\begin{align}
|\hat{\mathfrak{P}}^{[1,1]}|=2\,,\,\,|\hat{\mathfrak{P}}^{[2,1]}|=16\,,\,\, |\hat{\mathfrak{P}}^{[3,1]}|=192\,,\,\, |\hat{\mathfrak{P}}^{[4,1]}|=3072\,.
\end{align}
These numbers for $|\hat{\mathfrak{P}}^{[r,1]}|$ are reproduced by the expression $2^{2r-1}\,r!$. As for certain non-Tits-Satake manifolds, we also find:
\begin{align}
|\hat{\mathfrak{P}}^{[2,2]}|=32\,,\,\,|\hat{\mathfrak{P}}^{[2,3]}|=64\,,\,\,|\hat{\mathfrak{P}}^{[2,4]}|=128\,,\,\,|\hat{\mathfrak{P}}^{[3,2]}|=384\,,\,\, |\hat{\mathfrak{P}}^{[3,3]}|=768\,.
\end{align}

\section{The Laplacian Operator and the Harmonics on $\mathrm{U/H}$} \label{micompiacio}
One of the main targets for the use of the PGTS theory of noncompact symmetric spaces in Data Science problems is related to the following
fact. Once data have been mapped to the chosen $\mathrm{U/H}$ in a prescribed optimal way, any functional relation among the data is encoded in
functions $\Phi$ defined over the hosting manifold $\mathrm{U/H}$. All functions on a symmetric space (compact or noncompact) admit a canonical
representation as a series development
\begin{equation}\label{larvaditopo}
\Phi(\Upsilon) \, = \, \sum_{n=1}^\infty \, \,\sum_{\widehat{\pmb{\lambda}}_{2n}\in \widehat{\mathcal{I}}_{2n}}
\, C^\ell_{ \widehat{\pmb{\lambda}}_{2n}}\, \pmb{\mathfrak{harm}}^{\widehat{\pmb{\lambda}}_{2n}}_\ell(\Upsilon) \,
\end{equation}
in harmonics  $\pmb{\mathfrak{harm}}^{\widehat{\pmb{\lambda}}_{2n}}_\ell(\Upsilon)$, which are eigenfunctions of the Laplace-Beltrami operator and
admit a unique algebraic group-theoretical construction. Furthermore, the spectrum of the Laplacian is also determined a priori in terms of
group theory. This means that any kind of learning algorithm should focus only on the determination of the series coefficients $ C^\ell_{
\widehat{\pmb{\lambda}}_{2n}}$ obtaining in this way an intrinsic, coordinate independent representation of the functional relations satisfied by
the data.
\par
For all the noncompact symmetric spaces considered in the PGTS theory, the general construction of harmonic analysis on $\mathrm{U/H}$
specializes in such a way that all the harmonics $\pmb{\mathfrak{harm}}^{\widehat{\pmb{\lambda}}_{2n}}_\ell(\Upsilon)$ are actually constructed as
suitable tensor products of the symmetric matrix $\mathcal{M}(\Upsilon)$ introduced in eq.~\eqref{turlipano}.
\par
In the present section we just provide a concise summary of the harmonic analysis set up in the case of the maximally split symmetric manifold
$\frac{\mathrm{SL(N,\mathbb{R})}}{\mathrm{SO(N)}}$, illustrating the  eigenvalue  and the eigenfunction spectrum. Due to the always existing
triangular embedding (see statement \ref{statamento})
 \begin{equation}\label{perdignolo}
   \mathrm{\frac{U}{H} } \, \stackrel{\Phi_{triang}}{\Longrightarrow} \,\frac{\mathrm{SL(N,\mathbb{R})}}{\mathrm{SO(N)}},
 \end{equation}
the harmonics of the various $\mathrm{U/H}$ can be obtained from those of the mother manifold $\frac{\mathrm{SL(N,\mathbb{R})}}{\mathrm{SO(N)}}$
through the implementation of the additional algebraic constraints satisfied by the matrix $\mathcal{M}$ when it is restricted to the subgroup
$\mathbb{U}\subset \mathrm{SL(N,\mathbb{R})}$.
\par
The present section is somewhat sketchy and, as we already said, it is limited to the maximal case $\mathrm{SL(N,\mathbb{R})}$, yet it is an
essential conceptual thread in the texture of the PGTS theory in view of Data Science applications.
\subsection{Laplace-Beltrami operators on a symmetric space} A fundamental concept relative to a symmetric space $\mathcal{M} \,=
\,\mathrm{U/H}$ is that of rank (not to be confused with the noncompact rank $r_{n.c.}$ previously discussed). The coset-rank can be defined in
equivalent ways that are rooted either in the algebraic structure of the isometry group $\mathrm{U}$ and of its subroup $\mathrm{H}$ or in the
spectrum of invariant differential operators.
\par Quoting from \cite{castdauriafre},
we recall that on a coset manifold $\mathcal{M}\, = \, \mathrm{U/H}$, if we call $T_A$ a set of generators of the $\mathbb{U}$ Lie algebra,
$\mathbf{k}_A(\Upsilon)\, =\, {k}_A^\mu(\Upsilon)\,\partial_\mu$ the corresponding Killing vector fields generating the isometries of the
$\mathrm{U}$-invariant metric, a differential operator $\pmb{\triangle}$ of the second  or  higher  order is called an \emph{invariant operator} if
and only if
\begin{equation}\label{lappobeltro} \left [\pmb{\triangle} \, , \, \mathcal{L}_A \right] \, = \, 0 \qquad \forall  \mathcal{L}_A\, \in \, \imath^\star\left(\mathbb{U}\right),
\end{equation}
where
\begin{equation}\label{gromello}
  \mathcal{L}_A \, \equiv \, \mathbf{i}_{\mathbf{k}_A} \circ \mathrm{d} \, + \, \mathrm{d}\circ \mathbf{i}_{\mathbf{k}_A}
\end{equation}
is the Lie derivative along the Killing vector $\mathbf{k}_A$ which realizes a map from the abstract Lie algebra $\mathbb{U}$ to its image inside
the infinite dimensional Lie algebra of vector fields, i.e., into the space of sections of the tangent bundle to $\mathcal{M}$:
\begin{equation}\label{ranuncolo}
  \imath^\star \quad : \quad \mathbb{U} \, \longrightarrow \, \imath^\star\left(\mathbb{U}\right) \, \subset \,
  \Gamma \left[\mathcal{ T} \left(\mathrm{\frac{U}{H}}\right)\right ].
\end{equation}
Indeed the Lie derivatives satisfy the commutation relations of the $\mathbb{U}$ Lie algebra with the same structure constants:
\begin{equation}\label{isomorco}
  \left[ \mathcal{L}_A \, , \, \mathcal{L}_B \right] \, = \, f^{\phantom{AB}C}_{AB} \, \mathcal{L}_C \quad \Leftrightarrow \quad
    \left[ T_A \, , \, T_B \right] \, = \, f^{\phantom{AB}C}_{AB} \, T_C.
\end{equation}
Now we note that:
\begin{enumerate}
  \item the invariant operators are diagonal on \emph{harmonics}, which can be indeed defined as their eigenfunctions. The \emph{Laplace
      Beltrami operators} are a complete set of invariant operators and every $\pmb{\triangle}$ satisfying eq.~\eqref{lappobeltro}  is a function
      of them.
  \item There are $\mathfrak{r}=$ \emph{coset rank} independent Laplace-Beltrami operators on each coset manifold $\mathrm{{U}/{H}}$. The
      coset rank of $\mathrm{{U}/{H}}$ is defined as the number of irreducible representations of the subalgebra $\mathbb{H}$ in which its
      orthogonal complement $\mathbb{K}$ splits in the decomposition
      \begin{eqnarray}\label{cornodafrica}
        \mathbb{U}& = &\mathbb{H}\, \oplus \, \mathbb{K} \quad ; \quad \mathbb{K} \, \bigoplus_{\alpha=1}^{r} \, \mathbb{K}_\alpha \nonumber\\
        && \left[ \mathbb{H} \, , \, \mathbb{K}_\alpha\right] \, \subset \, \mathbb{K}_\alpha \quad \text{no invariant subspace in each $\mathbb{K}_\alpha$}.
      \end{eqnarray}
  \item In a coset manifold of coset rank $\mathfrak{r}$ the first lower degree $\mathfrak{r}$ Laplace-Beltrami operators can be chosen as a
      complete basis of invariant operators. The remaining higher order Laplace-Beltrami  operators are functionally dependent on the first $\mathfrak{r}$
      operators.
\end{enumerate}
It should be stressed that what we said above applies to the invariant differential operators $\pmb{\triangle}_{s}$ acting  on the space of
functions in each open chart of $\mathfrak{U}\subset\mathrm{{U}/{H}}$:
\begin{equation}\label{scalarLap}
  \pmb{\triangle}_s \, : \, \quad  \mathcal{C}^\infty \left(\mathfrak{U} \right) \, \longrightarrow \,  \mathcal{C}^\infty \left(\mathfrak{U} \right).
\end{equation}
Given the coset rank, there are just $\mathfrak{r}$ independent Laplace-Beltrami operators $\pmb{\triangle}_s$. If we consider the space of
sections of each associated bundle to the principal  bundle
\begin{equation}\label{lavapiatti}
  \mathrm{U}\, \stackrel{\pi}{\longrightarrow} \, \mathrm{U/H} \quad : \quad \forall p\in \mathrm{U/H} \quad \pi^{-1}(p) \simeq \mathrm{H}
\end{equation}
then for each of such bundles
\begin{equation}\label{maleperloro}
   \mathcal{E}_\alpha \stackrel{\pi_\alpha}{\longrightarrow} \, \mathrm{\frac{U}{H}}
   \quad : \quad \forall p\in \mathrm{\frac{U}{H}}   \quad \pi^{-1}(p) \simeq \mathbb{V}_\alpha,
\end{equation}
where $\mathbb{V}_\alpha$ is the carrier vector space of a linear representation $\mathrm{D}_\alpha$ of the structure group $\mathrm{H}$, there
are  typically $\mathfrak{r}$  Laplace Beltrami operators $\pmb{\triangle}_\alpha$   acting on sections of the bundle:
\begin{equation}\label{vectorLap}
  \pmb{\triangle}_\alpha \, : \, \quad  \Gamma\left(\mathcal{E}_\alpha, \mathrm{\frac{U}{H}} \right) \, \longrightarrow \,
  \Gamma\left(\mathcal{E}_\alpha, \mathrm{\frac{U}{H}} \right)
\end{equation}
\par
In the case of interest to us the coset rank of the coset $\mathrm{\frac{SL(N,\mathbb{R})}{SO(N)}}$  is always $\mathfrak{r}=1$ since the
$\mathbb{K}$ subspace is always in the irreducible representation:
\begin{equation}\label{roveda}
  \mathrm{D_{2sym}(SO(N))} \, =\, \mathbf{fund} \bigotimes_{symm} \mathbf{fund} \, - \, \text{Trace}
\end{equation}
where $\mathbf{fund}$ denotes the N--dimensional defining representation of $\mathrm{SO(N)}$ and the dimension of  the representation in
eq.~(\ref{roveda}) is
\begin{equation}\label{cristore}
  \mathrm{dim}_\mathbb{R} \, \mathrm{D_{2sym}(SO(N))} \, = \, \frac{N(N+1)}{2} -1.
\end{equation}
In the case of the symmetric spaces $\mathcal{M}^{[r,s]}_\mathfrak{d} \,= \mathrm{SO(r,r+2s)}/\mathrm{SO(r) }\times \mathrm{SO(r+2s)}$ the coset
rank is also $\mathfrak{r}=1$ since the $\mathbb{K}$-space is uniformly in the bifundamental irreducible represention of $\mathrm{SO(r)}\times
\mathrm{SO(r+2s)}$.
\par
We note just in passing that the flat Euclidian space used in most neural network algorithms and deep learning constructions is also a coset
manifold of rank $r=1$, precisely given by
\begin{equation}\label{flattomanifoldo}
  \text{Flat Space} \, = \, \mathrm{\frac{ISO(N)}{SO(N)}},
\end{equation}
with the difference that $\mathbb{U}=\mathbb{ISO}(N)$ is not a semi-simple Lie algebra  and the representation to which $\mathbb{K}$ is assigned
is the fundamental N-dimensional one rather than the traceless symmetric of rank 2. The translation group is abelian and coincides with the
manifold itself. The spectrum of the Laplace-Beltrami operator is continuous rather than discrete, as it is the case with the symmetric manifolds
here under consideration, as we are going to illustrate in the sequel. In other words all the structures that introduce a sequence of parameters
suitable to be \emph{target of learning} are lost.
\subsection{Harmonics on coset spaces} \label{cimentoarmonia} Let us consider as a first step
 the Lie group $\mathrm{U}$. A complete functional basis on $\mathrm{U}$ is given by the matrix elements of the irreducible Representations of $\mathrm{U}$. Indeed any function $\Phi$ on $\mathrm{U}$
can be expanded as
\begin{equation} \Phi\left(g\right)=\sum_{\left(\mu \right)}\sum_{m,n=1}^{{\rm
dim}\left(\mu\right)}{c_{mn}^{\left(\mu\right)}
D_{mn}^{\left(\mu\right)}\left(g\right)},
\end{equation}
where $\left(\mu\right)$ are the UIRs of $\mathrm{U}$, $m,n$ run in these representations, and $D$ are the elements of these representations. In
fact, the $D_{mn}^{\left(\mu\right)}\left(g\right)$ satisfy  orthogonality and completeness relations that for continuous groups are the obvious
generalization of the same orthogonality relations that hold true for finite groups.
\begin{eqnarray}
\int_\mathrm{U}
\,dg\,D_{mn}^{\left(\mu\right)}\left(g\right)D_{sr}^{\left(\nu\right)}\left(g^{-1}\right)={{\rm
vol}\left(\mathrm{U}\right)\over{\rm vol}\left(\mu\right)}\delta_{mr}
\delta_{ns}\delta^{\left(\mu\right)\left(\nu\right)}\nonumber \\
\sum_{\left(\mu\right)}D_{mn}^{\left(\mu\right)}(g)D_{nm}^{\left(\mu\right)}({g'}^{-1}){\rm
dim}\left(\mu\right)= \delta\left(g-g'\right){\rm
vol}\left(\mathrm{U}\right).
\end{eqnarray}
If $\Phi\left(g\right)$ transforms according to an irreducible representation $\left(\mu\right)$ of $\mathrm{U}$, for example under left multiplication,
namely
\begin{equation}
\Phi^{\left(\mu\right)}_m\left(g'g\right)=D^{\left(\mu\right)}_{mn}\left(g'\right)\Phi^{\left(\mu\right)}_n\left(g\right),
\end{equation}
then only a subset of the complete functional basis is present,  i.e., the $D$'s that transform in the same way, that is,
\begin{equation}
D^{\left(\mu\right)}_{mn}~~~~~\mu,m~{\rm fixed} \label{balengus}
\end{equation} So in
this case the expansion is shorter:
\begin{equation}
\Phi^{\left(\mu\right)}_m\left(g\right)=\sum_nc^{\left(\mu\right)}_nD^{\left(\mu\right)}_{mn}\left(g\right).
\end{equation}
\par
Let us now consider the functions $\Phi\left(y\right)$ on a coset manifold $\mathrm{U/H}$. The matrix elements
$D^{\left(\mu\right)}_{mn}\left(\mathbb{L}\left(y\right)\right)$,
where $\mathbb{L}\left(y\right)$ is a coset representative of the coset space\index{coset space}, and $y$ are the coset manifold coordinates,
constitute a complete functional basis on $\mathrm{U/H}$:
\begin{equation} \Phi\left(y\right)=\sum_{\left(\mu
\right)}\sum_{m,n=1}^{{\rm
dim}\left(\mu\right)}{c_{mn}^{\left(\mu\right)}
D_{mn}^{\left(\mu\right)}\left(\mathbb{L}\left(y\right)\right)}
\end{equation}
satisfying
\begin{eqnarray}
\int_{\mathrm{U/H}}d\mathfrak{m}\!\left(y\right)D_{mn}^{\left(\mu\right)}\left(\mathbb{L}\left(y\right)\right)D_{sr}^{\left(\nu\right)}
\left(\mathbb{L}\left(y\right)^{-1}\right)={{\rm
vol}\left(\mathrm{\mathrm{U/H}}\right)\over{\rm
vol}\left(\mu\right)}\delta_{mr}
\delta_{ns}\delta^{\left(\mu\right)\left(\nu\right)}\nonumber \\
\sum_{\left(\mu\right)}D_{mn}^{\left(\mu\right)}\left(\mathbb{L}\left(y\right)\right)
D_{nm}^{\left(\mu\right)}(\mathbb{L}\left(y^\prime\right)^{-1}){\rm
dim}\left(\mu\right)= \delta\left(y-y'\right){\rm
vol}\left(\mathrm{U/H}\right)
\end{eqnarray}
where $d\mathfrak{m}\!\left(y\right)$ is the invariant measure on $\mathrm{U/H}$.
\par
We are interested in functions $\Phi\left(\mathbb{L}\left(y\right)\right)$ on which a linear action of the subgroup $\mathrm{H}\subset \mathrm{U}$
is well defined, namely that transform within an irreducible linear representation $\left(\rho\right)$ of $\mathrm{H}$:
\begin{equation}
\label{turcomanno} h\cdot\Phi^{\left(\rho\right)}_i\left(\mathbb{L}\left(y\right)\right)\,
 \equiv \,
D_{ij}\left(h\right)^{\left(\rho\right)}\Phi^{\left(\rho\right)}_j\left(\mathbb{L}\left(y\right)\right) \qquad  \forall h \in \mathrm{H}
\end{equation}
the index $i$ running in $\left(\rho\right)$.
\par
Recalling, as it was done in eq.~(\ref{lavapiatti}), that a coset manifold $\mathrm{U/H}$ can be interpreted as the base manifold of a principal
fiber bundle having the subgroup $\mathrm{H}$ as structure group, the functions $\Phi^{\left(\rho\right)}_i\left(\mathbb{L}\left(y\right)\right)$
  in eq.~(\ref{turcomanno}) are actually sections of an associated vector bundle in the representation $(\rho)$ of the principal bundle
(\ref{lavapiatti}). These objects, which are relevant both in classical and quantum field theories and in particular in Kaluza Klein
compactifications of higher dimensional theories , are also considered the main ingredients of recent \emph{Geometrical Deep Learning
constructions}, as we already remarked.
So the issue is  how do we determine such vector bundle sections constructed in terms of  representation matrix elements
\begin{equation} D_{mn}^{\left(\mu\right)}~~~m,n~{\rm
running~in~}\left(\mu\right)~{\rm of~}\mathrm{U}.
\end{equation}
The answer is simple. We have to choose  matrix elements of the   type
\begin{equation}
\label{harmcondition} D_{in}^{\left(\mu\right)}~~~i~{\rm
running~in~}\left(\rho\right)~{\rm of~}H,~n~{\rm
running~in~}\left(\mu\right)~{\rm of~}\mathrm{U},
\end{equation}
hence we have to choose only those irreps $(\mu)$ that satisfy the following condition: the decomposition of $\left(\mu\right)$ with respect to
$\mathrm{H\subset \mathrm{U}}$ must contain the $\mathrm{H}$ irreducible representation $\left(\rho\right)$:
\begin{equation}
\label{hdec}
\left(\mu\right)\stackrel{H}{\longrightarrow}\cdots\oplus\left(\rho\right)\oplus\cdots.
\end{equation}
Only in this case $D_{mn}^{\left(\mu\right)}$ decomposes in $\left(\dots,D_{in}^{\left(\mu\right)},D_{i'n}^{\left(\mu\right)},\dots\right)$ and
the $D_{in}^{\left(\mu\right)}$ actually exists. The functions constructed with matrix elements of the representations satisfying (\ref{hdec}) are
called {\bf $\mathrm{H}-$harmonics} on $\mathrm{U/H}$, and make up a complete basis for the sections
$\Phi^{\left(\rho\right)}_i\left(\mathbb{L}\left(y\right)\right)$ of vector bundle on $\mathrm{U/H}$ associated with the representation
 $\rho$-vector bundle on $\mathrm{U/H}$.
The  expansion of such a section reads
\begin{equation} \label{salamexpansion}
\Phi^{\left(\rho\right)}_i\left(y
\right)={\sum_{\left(\mu\right)}}^\prime\sum_n
c_n^{\left(\mu\right)}D_{in}^{\left(\mu\right)}\left(\mathbb{L}\left(y\right)\right),
\end{equation}
where $\sum^\prime$ means a sum only on the representations $\left(\mu\right)$ satisfying the property (\ref{hdec}). Notice that the
$\mathrm{H}-$harmonics have both an index running in an irreducible representation of $\mathrm{U}$ (on the right) and an index running in an
irreducible representation of $\mathrm{H}$ (on the left). The coefficients of the expansion $c_n^{\left(\mu\right)}$ have an index of a
representation of $\mathrm{U}$ present in the expansion $\sum'$.
\par
Notice that when we deal with functions on the manifold rather than with sections of vector bundles the relevant representation $\rho$ of
$\mathrm{H}$ is simply the singlet representation, so that the sum $\sum^\prime$ is extended only to those representations of $U$ that when
decomposed with respect to $\mathrm{H}$ they contain the singlet.
%%%%%%%%%%%%%%%%%%%%%%%%%%%
\subsection{Differential operators on $\mathrm{H}$ harmonics}
%%%%%%%%%%%%%%%%%%%%%%%%%%%%%%%%%
$\mathrm{H}$--harmonics have a very powerful property:  the action  of differential operators on them  can be expressed in a completely algebraic
way. The $\mathrm{H}$-covariant external differential
\begin{equation}\label{Hcovaction}
    \mathcal{D}^{\mathrm{H}}\mathbb{L}^{-1}(y)\, \equiv \,\left(d+\omega^iT_i\right)\mathbb{L}^{-1}(y)
\end{equation}
becomes
\begin{equation}
\mathcal{D}^{\mathrm{H}}\mathbb{L}^{-1}(y)\, =\, -V^aT_a \,
\mathbb{L}^{-1}(y),
\end{equation}
or, expanding on the vielbein\index{vielbein} $\mathcal{D}^\mathrm{H}\, =\, {E}^a\, \mathcal{D}^\mathrm{H}_a$,
\begin{equation}
\label{rescHcovaction} \mathcal{D}_a^\mathrm{H}\,
\mathbb{L}^{-1}(y)=-\, T_a\, \mathbb{L}^{-1}(y)
\end{equation}
where  $T_a$  are the generators of the subspace $\IK$ defined by the orthogonal decomposition $\mathbb{U}=\IK\oplus\IH$, in the representation in which the
inverse coset representative is expressed.
\par
Since the harmonics {\it are} the inverse coset representatives in the representation $\left(\mu\right)$ we have
\begin{equation}
D^{\left(\mu\right)}\left(\mathbb{L}^{-1}(y)\right)^i_{~n}={\mathbb{L}^{-1}(y)}^i_{~n}.
\end{equation}
More precisely, the harmonic in the $\left(\rho,\mu\right)$, representations ${{D^{\left(\mu\right)}}^i_{~n}}$ is obtained by making  the
decomposition (\ref{hdec}) of the first index of ${{D^{\left(\mu\right)}}^m_{~n}}={\mathbb{L}^{-1}}^m_{~n}$ and taking the $\left(\rho\right)$
term. We consider the {\it inverse} coset representative because for simplicity of notation, we want $\mathrm{H}$ to act on their left, while it
acts on the right of coset representatives.
\par
Hence the action of $\mathrm{H}$-covariant derivative on the harmonic $D^i_{~n}(y) \equiv D^i_{~n}\left(\mathbb{L}(y)\right)$ is
\begin{equation}
\mathcal{D}^\mathrm{H}_a
\,D^i_{~n}\left(y\right)=-\left(T_a
D(y)\right)^i_{~n}=-\left(T_a\right)^i_{m}D^m_{~n},
\label{actionTa}
\end{equation}
where $\left(T_a\right)^i_{~m}$ is defined as the $\left(\rho\right)$ term in the decomposition (\ref{hdec}) of the index $n$ in
$\left(T_a\right)^n_{~m}$, namely,
\begin{equation}
\left(T_a\right)^n_{~m}=\left\{\dots,\left(T_a\right)^i_{~m},\dots\right\}.
\end{equation}
\par
Utilizing
\begin{equation}\label{fittofatto}
    \omega^{ab} \, \eta_{bc} \, \equiv \, \omega^{a}_{\phantom{a}c} \, = \, \frac{1}{2}
    \, C^a_{\phantom{a}cd}\, V^d \, + \, C^a_{\phantom{a}ci}\,\omega^i
\end{equation}
 which expresses the spin
$\so(n)$-connection $\omega^a_{\phantom{a}b}$ of the $n$-dimensional coset manifold $\mathrm{U/H}$ in terms of the $\mathrm{U}$-group structure
constants and of the $\mathrm{H}$-connection $\omega^i$ , it is evident that once the $\mathrm{H}$-covariant differential is reduced to an
algebraic action on the harmonics, the same will be true of the standard $\so(n)$-covariant differential:
\begin{equation}\label{brallo}
   \mathcal{D} \, \equiv \, \mathcal{D}^{\so(n)} \, = \, \mathrm{d} \, + \,
   \omega^{ab} \, J_{ab}
\end{equation}
where $J_{ab}$ are the $\so(n)$ generators.
\par
It follows that the Laplace-Beltrami invariant operators,  which are all  $\so(n)$ invariant operators quadratic in the
$\so(n)$ derivatives, will have an algebraic action  on harmonics, and will have on them a numerical eigenvalue that can be calculated a priori in
terms of the generators of the $\mathbb{U}$ Lie algebra and depends only on the $\mathrm{U}$ irreducible representation $(\mu)$. In conclusion, when
we consider sections of the various associated vector bundles to the principal bundle (\ref{lavapiatti}), the harmonics make up a complete
basis for the functional space over $\mathrm{U/H}$ and all differential operations in this functional space are reduced to algebraic ones.
\par
Thanks to group theory, analysis is reduced to algebra.
\subsection{The eigenvalues of the Laplacian in general} Further elaborating the above
basic facts,  in eqs.~(V.3.28) and (V.3.29) of the second volume of \cite{castdauriafre} it was proved that the scalar Laplace-Beltrami operator
$\pmb{\triangle}_s$ acting on the scalar harmonic in the irreducible representation $(\mu)$
\begin{equation}\label{basharm}
  \mathfrak{Harm}^{(\mu)}_\ell \, \equiv \, D^{(\mu)}_{0,\ell}\left(\mathbb{L}^{-1}(y)\right) \qquad
  \text{($\ell$ runs in $(\mu)$ and  $0$ means singlet)}
\end{equation}
can be represented as
\begin{equation}\label{ReaECrono}
  \pmb{\triangle}_s \, \mathfrak{Harm}^{(\mu)}_\ell \, = \, \left(g^{AB} \, \mathcal{L}_A \, \mathcal{L}_\mathcal{B}\right) \, \mathfrak{Harm}^{(\mu)}_\ell \, = \,
  \, \underbrace{\left(g^{AB} \, T_A \,T_B \right)^{(\mu)}_{\ell j} }_{\text{Casimir operator}}\, \mathfrak{Harm}^{(\mu)}_j,
\end{equation}
where $g^{AB}$ is the Killing metric on the $\mathbb{U}$ Lie algebra, $ \mathcal{L}_{A,B}$ denote, as above, the covariant Lie derivatives along
the Killing vectors, and $\left(g^{AB} \, T_A \,T_B \right)$ is the unique quadratic Casimir operator on the same Lie algebra. Since in each
irreducible representation the Casimir operator is proportional to the identity matrix with an overall constant factor $\mathcal{C}^{(\mu)}$,
the harmonics $\mathfrak{Harm}^{(\mu)}_\ell$ are eigenfunctions of the Laplace operator with  the value of the Casimir as eigenvalue:
\begin{equation}\label{Laplacvalli}
  \pmb{\triangle}_s \, \mathfrak{Harm}^{(\mu)}_\ell  \, = \, \mathcal{C}^{(\mu)} \, \mathfrak{Harm}^{(\mu)}_\ell
\end{equation}
Finally, due to Weyl and Freudenthal  there is a beautiful and concise formula which expresses the value of the Casimir operator in any irreducible
representation of the simple Lie algebra. Let $\Lambda^{(\mu)}$ be the highest weight vector of the irreducible representation $(\mu)$. It has the
form:
\begin{equation}\label{sapore}
  \Lambda^{(\mu)} \, = \, \sum_{i=1}^r \, n^{(\mu)}_i \, \mathbf{w}^i  \quad ; \quad n^{(\mu)}_i\,\in \, \mathbb{N},
\end{equation}
where $\mathbf{w}^i $ are the simple weights of the Lie algebra and the positive integers $n^{(\mu)}_i $ are the so-named \emph{Dynkin labels}
of the irrep $(\mu)$. Then the value of the quadratic Casimir in the $(\mu)$ irrep is given by
\begin{equation}\label{carallo}
\mathcal{C}^{(\mu)} \, = \, \langle \Lambda^{(\mu)} + 2 \rho \, , \, \Lambda^{(\mu)}\rangle \quad ; \quad \rho \, \equiv \, \ft 12 \, \sum_{\alpha >0} \, \alpha.
\end{equation}
According to our notations $\rho$ is the half sum of all the positive roots of the Lie algebra.
%%%%%%%%%%%%%%%%%%%%%%%%%%%%%%%%%%%%%%%%%%%%%%
\subsection{The spectrum of the Laplacian on the symmetric space $\mathrm{SL(N,\mathbb{R})/SO(N)}$} \label{spettrale} In the case of the symmetric
space of interest to us, the fundamental defining N-dimensional representation of the group $\mathrm{U}= \mathrm{SL(N,\mathbb{R})}$  is
irreducible with respect to the subgroup $\mathrm{H}=\mathrm{SO(N)}$, so that it cannot contribute to scalar harmonics. Indeed the lowest
representation that contains a singlet when reduced to the $\mathrm{H}$-subgroup is the $2$-symmetric  one with unit determinant obtained as
\begin{equation}\label{ceffus}
  \young(AB)\, \equiv \, \mathbb{L}^A_i \mathbb{L}^B_j \, \delta^{ij} \, = \, \left(\mathbb{L} \, \mathbb{L}^T\right)^{AB} \, \equiv \,.
  \mathcal{M}^{AB}
\end{equation}
Here the matrix $\mathcal{M}(\Upsilon)$ is the same as the one  previously introduced in eq.~\eqref{turlipano}, whose independent
components are exactly as many as the dimension of the coset manifold, namely, $\ft 12 N(N+1) -1$. For this reason the matrix
$\mathcal{M}({\Upsilon})$ can be regarded as a coordinate for the  points of the manifold, and the metric in terms of $\mathcal{M}({\Upsilon})$ was
written in eq.~(\ref{perbacco},\ref{cruccio}). When the matrix is parameterized in terms of the solvable coordinate or of any other set of
independent coordinates, it automatically satisfies the constraint $\text{Det} \mathcal{M }\, =\, 1$. In order to use $\mathcal{M}$ as a coordinate
such a highly non-linear constraint has to be imposed by hand.
\par
\begin{figure}
 \begin{center}
\includegraphics[width=8cm]{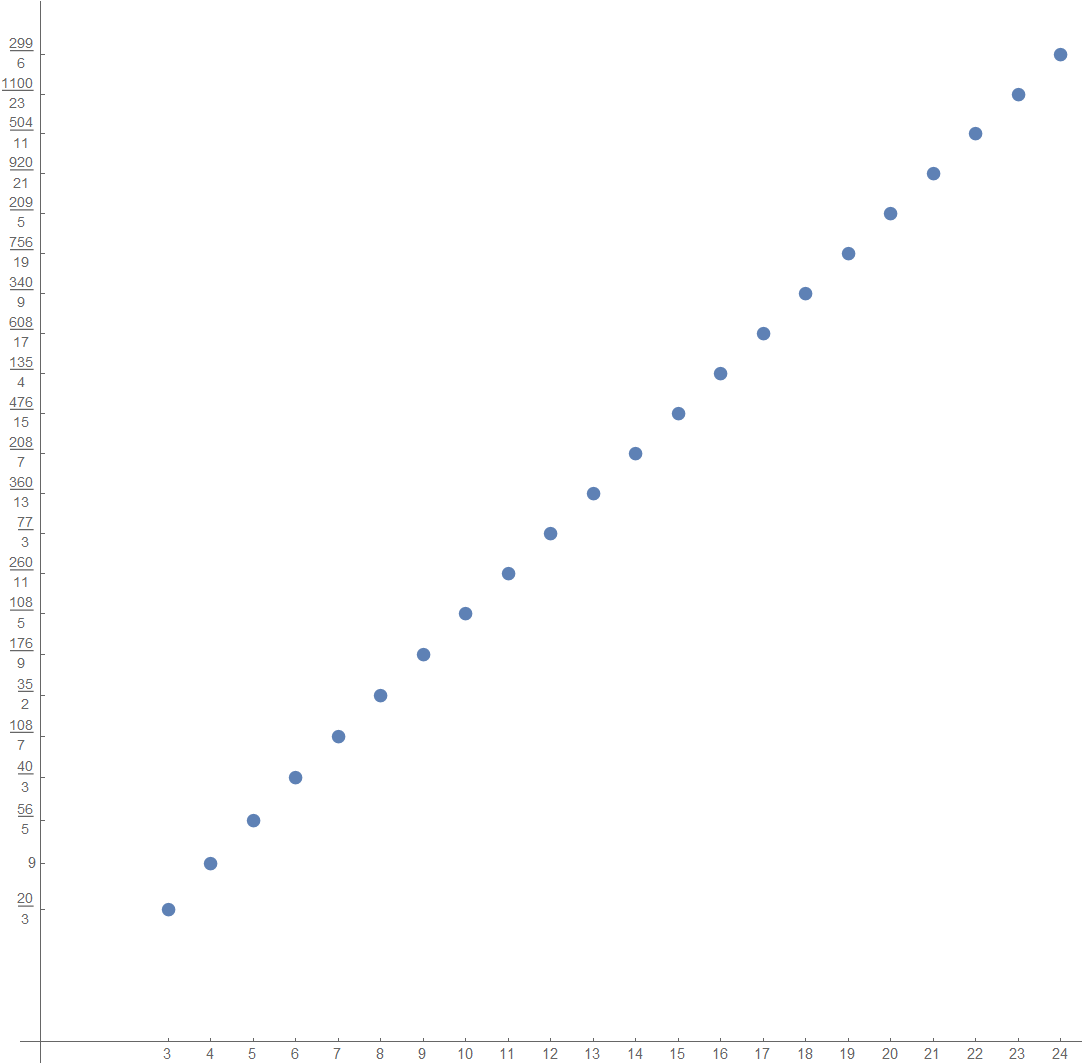}
 \caption{\label{lapatto} Approximately linear behavior of the quadratic Casimir acting  on the fundamental harmonic
 as a function  of  $N$.}
\end{center}
\end{figure}

What is very important is that now we see that the matrix $\mathcal{M}({\Upsilon})$ is also the fundamental building block of harmonics, in terms
of which all the other higher harmonics are constructed. For this reason, we have utilized the Young Tableau notation which is very handy and
economical. We can write:
\begin{equation}\label{limatus}
  \pmb{\triangle}_s \, \young(AB) \, = \, \mathcal{C}^2(N) \, \young(AB) \quad.
\end{equation}
As one can visually see from the ListPlot displayed in Figure \ref{lapatto}, the growth of the eigenvalue of the quadratic Casimir operator on the
fundamental harmonic is approximately linear in $N$, with a slope which starts with  $20/9 \simeq 2.2222$  at $N=3$  and gradually approaches $2$ as
$N\to \infty$.
%%%%%%%%%%%%%%%%%%%%%%%%%%%%%%%%%%%%%%%%%%%%%%%%%%%%%%%%%%%%%%%%%%%%%%%%%%%%%%%
\subsection{The general construction of Harmonics}
\label{cimentoinvenzione} Now the general form of the  harmonics can  be easily constructed by taking into account that for the group
$\mathrm{SL(N,\mathbb{R})}$ the irreducible linear representations are in a  one-to-one correspondence with the Young Tableaux, namely,  with the
partition of an integer number $1 \leq n\in \mathbb{N}$ \emph{i.e.} by  the $r$-tuples $\pmb{\lambda} $ of integer numbers (with $r\leq n$)
such that
\begin{equation}\label{corifeo}
\pmb{\lambda} \, = \, \{\lambda_1, \lambda_2, \dots,\lambda_r\} \quad ; \quad \lambda_1 \geq \lambda_2 \geq\dots \geq \lambda_{r-1}
\geq \lambda_r \quad ; \quad \sum_{i=1}^r \, \lambda_i \, = \, n.
\end{equation}
%%%%%%%%%%%%%%%%%%%%%%%%%%%%%%%%%%%%%%%%%%%%%%%%%%%%%%%%%%%%%%%%%%%%%%%%%%%%%%%
As it is well known, with each partition $\pmb{\lambda} $ one associates a corresponding tableau which is the disposition of $n$ boxes in $r$ lines,
each line not longer than the one above. This is  not only a mnemonic representation of the partition but also a graphical representation of an
operator acting on the indices $A_1, \dots,\, A_n$ of a tensor:
\newcommand{\aone}{\mbox{$A_1$}}
\newcommand{\atwo}{\mbox{$A_2$}}
\newcommand{\alone}{\mbox{$A_{\lambda_1}$}}
\newcommand{\alonep}{\mbox{$A_{\lambda_1+1}$}}
\newcommand{\altwo}{\mbox{$A_{\lambda_2}$}}
\newcommand{\almenol}{\mbox{$A_{n+1-\lambda_r}$}}
\newcommand{\alno}{\mbox{$A_{n}$}}
\begin{equation}\label{tablo}
 \pmb{\lambda}\,\Rightarrow \, {\Yboxdim28pt
 \young(\aone\atwo\bullet\bullet\bullet\bullet\bullet\alone,\alonep\bullet\bullet\bullet\bullet\altwo,\bullet\bullet\bullet\bullet\bullet,\bullet\bullet\alno)}
\end{equation}
Indeed the Young operator acting on the set of $n$ indices associated with the Young tableau is classically defined as follows:
\begin{eqnarray}\label{youngoperator}
  \mathbb{Y}_{\pmb{\lambda}} & = & Q\circ P \nonumber\\
  P & = & \sum_{p} p \, = \, \text{symmetrizer of each line, summing on all permutations of $\lambda_i$ objects}\nonumber\\
  Q&=&\sum_q q \,(-1)^{\delta_p} \, = \,\text{skew-symmetrizer of each column,} \nonumber \\[-12pt] && \parbox{24.2mm}{\hfill}  \text{summing on all permutations of the  columns.}
\end{eqnarray}
In this way, one constructs an irreducible representation of the symmetric group $S_n$, yet for a general theorem the irreducible representations
of the group $\mathrm{SL(N,\mathbb{R})}$ are in one-to-one correspondence with the irreducible $n$-index tensors, namely, those whose indices have
been (skew)symmetrized with one of the available  Young tableaux, which are as many as the partitions of $n$ into integers. Obviously the
\text{irreps} of  $\mathrm{SL(N,\mathbb{R})}$  are infinitely many because the number $n$ of indices of the tensor is unbounded.
\par
The great advantage of working with the group $\mathrm{SL(N,\mathbb{R})}$, which, for suitable $N$, contains, as a subgroup, any other of the
$\mathrm{U}$ groups under consideration, is that the relation between the description of the irreducible representations in terms of Young
Tableaux and in terms of maximal weights (\emph{i.e.}) of \emph{Dynkin labels}, is immediate. The Dynkin labels defined in eq.~(\ref{sapore}) are
immediately related to the Young labels $\pmb{\lambda}$ by means of the   rule
\begin{flalign}\label{chiccobello}
  n_1 &= \lambda_1 -\lambda_2 \nonumber\\
  n_2 &= \lambda_2 -\lambda_3 \nonumber\\
  \dots & = \dots \nonumber\\
  n_r &=\lambda_{r-1} -\lambda_r.
\end{flalign}
%%%%%%%%%%%%%%%%%%%
\subsection{The irreps contributing to harmonic expansions on $\mathrm{SL(N,\mathbb{R})/SO(N)}$} In view of what we discussed before the rule to
select the irreps that contribute to the harmonic expansion of a generic function on the manifold is very simple. The only available Young
tableaux are those where the number of boxes on each line is even since the tensor must be constructed from the tensor product of an integer
number of fundamental harmonics $\young(AB)$. For instance, if we consider the harmonics that can obtained from the square of two fundamental
harmonics we have:
\begin{equation}\label{rattone}
  \young(AB)\otimes \young(CD) \, = \, \young(ABCD) \oplus \young(AB,CD)
\end{equation}
A priori we have:
\begin{equation}\label{micione}
 \yng(2) \otimes \yng(2) \, = \, \yng(4) \oplus \yng(2,2) \oplus \yng(3,1)
\end{equation}
yet applying the Young operator corresponding to the last tableau on the r.h.s. of equation (\ref{micione}) to the product of the two symmetric
matrices we obtain identically zero just because of their symmetry:
\begin{equation}\label{cucciolo}
 \mathbb{ Y}_{\yng(3,1)} \,\left[\mathcal{M}(A,B)\,\mathcal{M}(C,D)\right] \, = \, 0
\end{equation}
Similarly it happens for all higher tensor products of the matrix $\mathcal{M}$ , or equivalently of its inverse:
\begin{equation}\label{inversamente}
  \mathcal{N} \equiv \mathcal{M}^{-1}
\end{equation}
\subsubsection{Examples in the cases $N=4$}
\paragraph{\sc N=4 case.} Utilizing the explicit differential form of the Laplacian operator in the $U/H$ invariant metric written in terms of
solvable coordinates we obtain
\begin{equation}\label{conncticut}
  \pmb{\triangle}_s \, \mathcal{M}(\Upsilon) \, = \, 9 \, \mathcal{M}(\Upsilon),
\end{equation}
in agreement with the corresponding Casimir value calculated with the Weyl-Freudenthal formula. Next,  let us consider the level two harmonics, namely,
those obtained from the product of two fundamental ones as in eq.~\eqref{micione}. For instance we can consider the following object:
\begin{equation}\label{4rep}
  \young(1234) \, = \, \mathcal{M}^{12}\, \mathcal{M}^{34} + \mathcal{M}^{13}\, \mathcal{M}^{24}+ \mathcal{M}^{23}\, \mathcal{M}^{14}
\end{equation}
and applying the Laplace-Beltrami operator we find
\begin{equation}\label{muccarosata}
  \pmb{\triangle}_s \,\young(1234)  \, = \, 24 \, \young(1234).
\end{equation}
The eigenvalue $24$ is precisely the one  expected from the calculation of the Casimir for the representation $\yng(4)$ in the case $N=4$. Indeed it
corresponds to the highest weight vector
\begin{equation}\label{chepeso}
  \Lambda^{\tiny\yng(4)} \,=\, 4 \, \mathbf{w}_1 \, \stackrel{\text{for N=4}}{=} \, \left\{\frac{10}{3},-\frac{2}{3},-\frac{2}{3}\right\},
\end{equation}
which inserted in the Weyl-Freudenthal formula (\ref{carallo}) together with
\begin{equation}\label{cherombo}
  \rho \, \stackrel{\text{for N=4}}{=} \,\left\{2,1,0\right\}
\end{equation}
yields precisely
\begin{equation}\label{checasino}
  \mathcal{C}^{\tiny\yng(4)}(4) \, = \, 24.
\end{equation}
Consider next the other irreducible representation in the same $n=4$ level of the case $N=4$. We have for instance
\begin{equation}\label{picciotto}
  \young(12,34) \, = \, 2 \mathcal{M}^{12}\, \mathcal{M}^{34} \, - \, \mathcal{M}^{32}\, \mathcal{M}^{14} \, - \, \mathcal{M}^{13}\, \mathcal{M}^{24},
\end{equation}
and by explicit evaluation, we find
\begin{equation}\label{crinotto}
 \pmb{\triangle}_s \,\young(12,34)  \, = \, 12 \, \young(12,34)
\end{equation}
which is   consistent with the calculation of the Casimir eigenvalue by means of formula \eqref{carallo}. Indeed we have
\begin{equation}\label{ giorgino}
  \Lambda^{\phantom{i}{\tiny\yng(2,2)}} \, = \, 2 \, \mathbf{w}_2 \, \stackrel{\text{for N=4}}{=} \, \left\{\frac{4}{3},\frac{4}{3},--\frac{2}{3}\right\}
  \quad \Rightarrow \quad \mathcal{C}^{\phantom{i}{\tiny\yng(2,2)}} (4)\, = \, 12.
\end{equation}
\subsection{Summarizing}
\label{summasumma} We conclude by saying that any function on the symmetric space can be developed into harmonics
that are in one-to-one correspondence with the irreducible representations of the isometry group pertaining to the special set clarified above,
for which the parameters $\lambda_i$ are all even.
Npw we want to    introduce  a more compact notation. First of all, we introduce the set of all contributing irreps, which  we call
$\widehat{\mathcal{I}}$. The irreps in $\widehat{\mathcal{I}}$ can be alternatively singled out by the Dynkin labels or by the integer partitions
$\pmb{\lambda}$. They have a degree which is just the sum of the even $\lambda$'s, hence we can write
\begin{equation}\label{sgrunt}
\widehat{\mathcal{I}} \, = \, \bigoplus_{n=1}^\infty \, \widehat{\mathcal{I}}_{2n}.
\end{equation}
By definition $\widehat{\mathcal{I}}_{2n}$ is the set of all partitions of the even number $2n$ into even integers; we call its    elements   $\widehat{\pmb{\lambda}}_{2n}$. The set of $2n$ indices arranged with the appropriate symmetry provided by the
corresponding Young tableau can be regarded as a unique index $\ell$ which takes as many values as there are independent ways of filling the Young
tableau boxes. The number of the latter is the dimension of the irreducible representation. Hence we introduce the symbol
\begin{equation}\label{harmonina}
\pmb{\mathfrak{harm}}^{\widehat{\pmb{\lambda}}_{2n}}_\ell(\Upsilon) \, = \, \text{harmonic of degree $2n$ in $\widehat{\pmb{\lambda}}_{2n}$ irrep},
\end{equation}
and for any function  $\Phi(\Upsilon)$ we  can write the harmonic expansion (\ref{larvaditopo}).
\par
The coefficients $C^\ell_{ \widehat{\pmb{\lambda}}_{2n}}$ are constant numbers and constitute the \emph{generalized Fourier components} of the
chosen function in the harmonic expansion. As we already stated at the beginning of the present section, \emph{generalized Fourier components}
are the \emph{target of supervised learning} in Data Science applications once the Data have been mapped to the $\mathrm{U/H}$ symmetric space
in an optimal way, as the result of some appropriate \emph{unsupervised optimization algorithm}. Obviously when $\mathrm{U\subset SL(N,R)}$ is a
proper subgroup, which is our general case of manifolds with no-trivial Tits-Satake submanifold, one has to decompose the irreducible
representations of the larger group $\mathrm{SL(N,R)}$ with respect to smaller one $\mathrm{U}$. Typically this can be done in an automatic and
algebraic way by diagonalizing the Casimir operator
\begin{equation}\label{criogenico}
  \mathbf{QC} \, = \, \sum_{i=1}^{N-1} \, \mathcal{H}_i \cdot \mathcal{H}_i \, + \, \sum_{\alpha >0} \left(  \mathcal{E}^{\alpha}
  \cdot \mathcal{E}^{-\alpha} + \mathcal{E}^{-\alpha} \cdot \mathcal{E}^{\alpha} \right)
\end{equation}
on the symmetric matrix $\mathcal{M}$ and on its tensor products.
%%%%%%%%%%%%%%%%%%%%%%%%%%%%%%%%%%%%%%%%%%%%%%%%%
%Titto Satako %%%%%%%%%%%%%%%%%%%%%%%%%%%%%%%%%%%%%%%%
%%%%%%%%%%%%%%%%%%%%%%%%%%%%%%%%%%%%%%%%%%%%%%%%%
\section{Conclusions}
\label{zakluchenye}
Since all the main conceptual issues have been repeatedly emphasized both in the introduction
and throughout all the expository sections we do not repeat them; we just reiterate what we said at the beginning: the PGTS theory of noncompact symmetric spaces is a
self-contained and rich theoretical scheme that provides its perspective utilizer in Data Science with an ample spectrum of items to be used in
constructing learning algorithms more structured than the usual ones, merely based on the distance function in flat spaces with the consequent need of \textit{point-wise activation functions} that spoil interpretability and covariance. 
\par
This was true in March 2025 when the first version of this foundational paper was posted on ArXiv and continues to be true at the present time. Yet on the basis of the progresses made in the meantime in the development of \textbf{Cartan  Neural Network}, encoded in the corpus of eight papers of which the present one is only the first \cite{pgtstheory,TSnaviga,naviga,tassellandum,axialbeltra,geotermico,secondtemperature,terzatemperatura}
we are in a  position to illustrate the original motivation with ample facts.
\par
It was challenging in \cite{naviga}, to test,  just in the simplest case of the $r=1$ TS universality class of hyperbolic spaces, that \textbf{Cartan  Neural Network}, unspotted from  \textit{point-wise activation functions}, are at least competitive with classical neural networks. The construction of $r=2$ CaNNs based on Calabi-Vesentini manifolds \cite{cvclassic},\cite{Hua1946}:
\begin{equation}\label{Cvmanigoldi}
  \mathcal{M}_{CV}^{[2,q]}\,  \equiv \, \frac{\mathrm{SO(2,2+q)}}{\mathrm{SO(2)\times SO(2+q)}} 
\end{equation}
each of which, at all time,  must be regarded as one of the two cofactors constituting an item in following infinite series of Special K\"ahler manifolds:
\begin{equation}\label{speckal}
 \mathcal{SK}_{3+q} \, \equiv \, \frac{\mathrm{SL(2,\mathbb{R})}}{\mathrm{SO(2)}}\times
 \frac{\mathrm{SO(2,2+q)}}{\mathrm{SO(2) \times SO(2+q)}}
\end{equation} 
is presently under development \cite{r2paperone}.  Yet much more has happened, as we already mentioned in the introduction, between the posting of the original version of the present paper on ArXiv in March 2025 and the present time when we submit it for publication on a journal. The most relevant development from the point of view of application to ML is that concerned with Information Geometry\cite{raone,cenzone,amarone}, its relation with Geometrical Thermodynamics\cite{gianno1,gianno2,lychaginlecture,Kushner_2020,Lychagin_2020,ludaed,ludaed2,
Ruppeiner_2012b,Ruppeiner_2010,Ruppeiner_2012,Ruppeiner_2013,ruppoRdiag}  and Souriau thermodynamics\cite{souriaub1,souriaub2,barbarpapad4,
marlentropia,caldobarbaresco,barbaresco2,barbaresco3,marlegibbs,barbarpapad1,barbarpapad2,barbarpapad3,nebbo}.
Indeed in a series of three papers \cite{geotermico,secondtemperature,terzatemperatura} two of us, in collaboration with A.S. Sorin, succeeded in constructing a general scheme for the construction of a Thermo-Dynamical K\"ahler metric given a suitable generalized Partition Function and to apply such a scheme to the case of (generalized) Souriau thermodynamics on Calabi Vesentini manifolds, for which the exact partition function was computed by means of  newly discovered abelian structures. These developments provide the best perspectives in order to analyse categorical perception
 \cite{groundwork,neurocoding} both in \textbf{trained supervised CaNNs} and in new \textbf{CaNNs} of the \textbf{unsupervised with reinforcement type} to be developed using the same CV manifolds as layers.
 \par 
 Behind all such developments stands the essential focal point of the present \textbf{foundational paper} namely the \textbf{metric equivalence} of \textbf{non-compact symmetric spaces} with a \textbf{solvable Lie group} and more generally the intuition that the geometrical lore developed within the 50 year long development of Supergravity Theory is not only useful, rather it is essential for a systematic overall geometrical, interpretable and covariant reformulation of AI algorithms. What in March 2025 was a motivated and perspective suggestion
 tends to be in July 2026 a collection of rapidly developing explicit implementations.
 \par  
 On the other hand, from a purely mathematical viewpoint, the reassembling of
noncompact symmetric space geometrical lore in the setup here presented with the inclusion of several new items and results opens new visions and
research directions, the most fertile and attractive being the one associated with the determination of discrete groups of the tessellation type.
\par
As we stated in the introduction the original aim of the present foundational paper was establishing a new PGTS theoretical paradigm for neural network architectures; the latter was finally  christened \textbf{Cartan Neural Networks} in \cite{TSnaviga,naviga}. In the original version of March 2025 we presented a  sketchy summary of the subsequent research lines as follows:
{\textit{
\begin{description}
  \item[A)] Use of the generalized and progressive Tits-Satake projection as a mathematical framework for the interpretation and construction of
      multi-layer neural networks, the various layers belonging to the same Tits-Satake Universality class \cite{naviga}. In this setup a basic
      ingredient is the geometrical interpretation of the separator hypersurfaces in the  organization of data by classes.  The already obtained
      results for $r=1$ have to be generalized to $r=2$ and larger.
  \item[B)] The use of harmonic expansion on $\mathrm{U/H}$ in the construction of neural networks where the data can be regarded as a sequence
      of points in section of a vector bundle such that $\mathrm{U/H}$ is the total bundle space, the base manifold is the Tits-Satake
      submanifold $\mathcal{M}_{TS}$ and the vector fibers are $\pi_{TS}^{-1}(p)$, $\forall p \in \mathcal{M}_{TS}$. Such an issue is addressed
      in a preliminary analysis in the forthcoming paper \cite{tassellandum}.
  \item[C)] Within the scope of the framework to be presented in the paper \cite{naviga},  we plan to study and devise inspection and control additions to
      the neural network architecture that allow to study the evolution of  paint group invariants through epochs.
  \item[D)] Study of the Heat Kernel on the various Tits-Satake universality classes in view of establishing diffusion processes both in the
      continuous and in the discrete case in the forthcoming paper \cite{tassellandum}.
  \item[E)] Embedding of tree graph structures in noncompact hyperbolic manifolds of various noncompact rank.
\end{description}}
By the time the present foundational paper, after updating, is submitted for publication, the \textbf{CaNN} project has been developed to a large extent as we mentioned above. Point A) of the above list is not yet completed, yet it is in fieri \cite{r2paperone}. Point B) and C)  have been partially attained and partially remoduled in view of the new discoveries in Souriau thermodynamics, which were not foreseen in March 2025. The evolution of the distributions of data from the point of view of Paint Group invariants can now be addressed relying on the Casimir Invariants of the Paint Group subalgebras that enter the scheme of spontaneous magnetizations 
\cite{secondtemperature,terzatemperatura} rather than on the Harmonic Expansions. Point D) has been solved for the hyperbolic spaces in \cite{tassellandum} and postponed for other TS classes in view of the different visions provided by the thermodynamical results. Similarly for point E). The present research perspectives after the completion of the following eight papers   \cite{pgtstheory,TSnaviga,naviga,tassellandum,axialbeltra,geotermico,secondtemperature,terzatemperatura}
were outlined as follows in \cite{secondtemperature}.
\par
\textit{
The above results are completely new and open very interesting perspectives for further investigations and study that have an unexpected and potentially deep impact on the overall structure of AI procedures. In order to appreciate such a point let us recall some facts and observations. 
\begin{enumerate}
  \item On one side, in his recent very interesting popularizing book \cite{bennetto}, Max Bennett  remarked the following: \textit{While most AI advancements that occurred in the early 2000s involved applications of supervised-learning models, many of the recent advancements have been applications of generative models. Deepfakes, AI-generated art, and language models like GPT-3 are all examples of generative models at work. Helmholtz suggested that much of human perception is a process of inference-a process of using a generative model to match an inner simulation of the world to the sensory evidence. \dots It turns out that there is, in fact, an abundance of evidence that the neocortical microcircuit is implementing such a generative model.} 
  \item In their very interesting papers \cite{groundwork,neurocoding} on \textbf{Categorical perception} in Biological and Artificial Neural Networks, Bonnasse-Gahot and Nadal consider the internal geometrical representation of the exterior world, in relation with the distribution of given exterior entities into categories. Such geometrical internal representation is presumably realized, to use the terminology of Bennett, in the \textit{microcircuitry of the neocortical columns} in mammal brains, while it is provided by appropriate unsupervised generative algorithms in artificial neural networks. In their study of such an issue, Bonnasse-Gahot and Nadal utilize
      Information Geometry and focus on the analysis of the local geometry around the boundaries between different categories.       
  \item The boundaries between categories, when one pursues a
  \textbf{classification task} by means of a \textbf{supervised Cartan  Neural Network}, are provided by 
  the theory of \textbf{separators}, recently introduced by two of us  in collaboration with other authors, in \cite{tassellandum}. \textbf{Separators are codimension one, homogeneous but not symmetric submanifolds} of the microscopic K\"ahler manifolds constituting each layer of the network. 
  \item According to  Bennett's vision, \textbf{pattern recognition}, alias \textbf{categorical perception}, and 
      \textbf{the generative construction of geometrical representations}, namely \textbf{imagination}, are two different operational modalities of the same neocortical microcircuitry in biological brains or of the same artificial neural network, the \textbf{awake modus operandi} the former, the \textbf{sleeping modus operandi} the latter, namely \textbf{supervised} or \textbf{unsupervised} in the artificial case. Hence the relation between the categorical partition of the same manifold induced by separators and its foliation induced by spontaneous or forced (at non vanishing generalized magnetic fields $\mathbf{h}$) \textbf{magnetizations}, emerges to strategic prominence. 
  \item The extra observables whose mean values are the magnetizations are  associated with  Casimir functions of the Paint-Group and of its subgroups. It follows that the already envisaged similarity grouping of data according to their distribution into Paint fibres might be realized by the spontaneous magnetization mechanism. Like in physical ferromagnets where, below the critical temperature one observes the microscopic \textit{spin up} and \textit{spin down islands}, in the same way, the Souriau like Gibbs distributions might trace out \textit{islands} in the manifold where certain features are dominant and other are depressed and viceversa.  
\end{enumerate}
}
Indeed as it was pointed out in the conclusions to \cite{terzatemperatura}
\textit{
the use of extended Souriau thermodynamics in Machine Learning can occur at two levels that are strictly correlated:
\begin{enumerate}
  \item In a posteriori analysis of the distributions of information data on the various layers of a \textbf{supervised and trained Cartan Neural Network} in order to study the geometry of categorical representations in the \textit{categorical perception} \cite{groundwork,neurocoding} à la Bonnasse-Gahot and Nadal.
  \item In \textbf{unsupervised learning with reinforcement}, namely generative AI. In that capacity (still to be accurately studied) the variation of magnetic fields might be the mathematical realization of the \textbf{agents} that either reinforce or punish the spontaneous evolution of the dynamical system (geodesics in thermodynamical space) 
\end{enumerate}
}
These quotations from the other papers originated from the foundational present one have been done to stress \textit{a posteriori} the correctness of the attribute \textbf{foundational} that we used since its very first postingo to ArXiv in March 2005. The intuition, that the mathematical geometrical structures of Supergravity in particular that of Special K\"ahler Geometry (of the local type or projective type to use Freed's nomenclature
\cite{Freed}) is correct we can further mention the forthcoming paper \cite{toinepietromario} where we are going to show that the so named L(q,P) homogeneous but non symmetric Special K\"ahler manifolds correspond to a different metric with less isometries imposed  on the same solvable Lie group manifolds as those underlying CV manifolds. Hence such spaces can be included as layers in the Cartan Neural Networks based on CV manifolds! 
\section*{Acknowledgements}
The present paper was the first result of an intense research activity pursued since late 2023 spring by the three of us, with some useful contributions from the colleagues Michele Caselle and Guido Sanguinetti and  of the junior members of our collaboration team, the PhD  students Federico Milanesio and Matteo Santoro. The research project, of which the present paper is the first step, has the ambitious aim of providing an innovative mathematical framework for a full fledged geometrical reformulation of ML algorithms and it was  launched because of the vision and generosity of our dear friend Sauro Additati, who incited us to accept the challenge, generously co-financing the  PhD  fellowships of the younger members within the PNRR programme 117,  and the full dedication to this cause of the elder and retired member (P.F.). Since October 2025 the further development of the same project was conducted only by two of us (P.F. and M.T.) with the strategic and essential contribution of our long time collaborator and excellent friend Alexander S. Sorin whom we would like to thank most warmheartedly.
\newpage
\appendix
\section{Explicit matrix representation of the $J_{ij}^\pm,\,J_i^I,\,T_{ij}^\pm,\,T_i^I$ generators.}
\label{peritonite}
In this appendix, we give the explicit matrix representation of the generators $J_{ij}^\pm,\,J_i^I,\,T_{ij}^\pm,\,T_i^I$ for
the discrete subgroups $\hat{\Delta}^{[r,q]}$ of ${\rm SO}(r,r+q)$ discussed in Section \ref{generaloneSOpq}. We work in the basis in which the
invariant matrix is:
\begin{equation}
\diag(\overbrace{+,\dots, +}^{r},\overbrace{-,\dots, -}^{r+q})\,.
\end{equation}
\paragraph{Case $r=2,\,q=1$.}
\begin{align}
    J_{12}^+&=\left(
\begin{array}{ccccc}
 0 & -1 & 0 & 0 & 0 \\
 1 & 0 & 0 & 0 & 0 \\
 0 & 0 & 0 & 1 & 0 \\
 0 & 0 & -1 & 0 & 0 \\
 0 & 0 & 0 & 0 & 1 \\
\end{array}
\right)\,,\,\,  J_{12}^-=\left(
\begin{array}{ccccc}
 0 & -1 & 0 & 0 & 0 \\
 1 & 0 & 0 & 0 & 0 \\
 0 & 0 & 0 & -1 & 0 \\
 0 & 0 & 1 & 0 & 0 \\
 0 & 0 & 0 & 0 & 1 \\
\end{array}
\right)\,,\nonumber\\
J_1^{I=1}&=\left(
\begin{array}{ccccc}
 1 & 0 & 0 & 0 & 0 \\
 0 & 1 & 0 & 0 & 0 \\
 0 & 0 & -1 & 0 & 0 \\
 0 & 0 & 0 & 1 & 0 \\
 0 & 0 & 0 & 0 & -1 \\
\end{array}
\right)\,,\,\,J_2^{I=1}=\left(
\begin{array}{ccccc}
 1 & 0 & 0 & 0 & 0 \\
 0 & 1 & 0 & 0 & 0 \\
 0 & 0 & 1 & 0 & 0 \\
 0 & 0 & 0 & -1 & 0 \\
 0 & 0 & 0 & 0 & -1 \\
\end{array}
\right)\,,\nonumber\\
T_{12}^-&=\left(
\begin{array}{ccccc}
 1 & -2 & 0 & 2 & 0 \\
 2 & 1 & 2 & 0 & 0 \\
 0 & 2 & 1 & -2 & 0 \\
 2 & 0 & 2 & 1 & 0 \\
 0 & 0 & 0 & 0 & 1 \\
\end{array}
\right)\,,\,\,T_{12}^+=\left(
\begin{array}{ccccc}
 1 & -2 & 0 & -2 & 0 \\
 2 & 1 & 2 & 0 & 0 \\
 0 & 2 & 1 & 2 & 0 \\
 -2 & 0 & -2 & 1 & 0 \\
 0 & 0 & 0 & 0 & 1 \\
\end{array}
\right)\,,\nonumber\\
T_1^{I=1}&=\left(
\begin{array}{ccccc}
 3 & 0 & 2 & 0 & 2 \\
 0 & 1 & 0 & 0 & 0 \\
 -2 & 0 & -1 & 0 & -2 \\
 0 & 0 & 0 & 1 & 0 \\
 2 & 0 & 2 & 0 & 1 \\
\end{array}
\right)\,,\,\,T_2^{I=1}=\left(
\begin{array}{ccccc}
 1 & 0 & 0 & 0 & 0 \\
 0 & 3 & 0 & 2 & 2 \\
 0 & 0 & 1 & 0 & 0 \\
 0 & -2 & 0 & -1 & -2 \\
 0 & 2 & 0 & 2 & 1 \\
\end{array}
\right)\,.
\end{align}
\paragraph{Case $r=3,\,q=1$.}
\begin{align}
    J_{12}^+&=\left(
\begin{array}{ccccccc}
 0 & -1 & 0 & 0 & 0 & 0 & 0 \\
 1 & 0 & 0 & 0 & 0 & 0 & 0 \\
 0 & 0 & 1 & 0 & 0 & 0 & 0 \\
 0 & 0 & 0 & 0 & 1 & 0 & 0 \\
 0 & 0 & 0 & -1 & 0 & 0 & 0 \\
 0 & 0 & 0 & 0 & 0 & 1 & 0 \\
 0 & 0 & 0 & 0 & 0 & 0 & 1 \\
\end{array}
\right)\,,\,\,  J_{12}^-=\left(
\begin{array}{ccccccc}
 0 & -1 & 0 & 0 & 0 & 0 & 0 \\
 1 & 0 & 0 & 0 & 0 & 0 & 0 \\
 0 & 0 & 1 & 0 & 0 & 0 & 0 \\
 0 & 0 & 0 & 0 & -1 & 0 & 0 \\
 0 & 0 & 0 & 1 & 0 & 0 & 0 \\
 0 & 0 & 0 & 0 & 0 & 1 & 0 \\
 0 & 0 & 0 & 0 & 0 & 0 & 1 \\
\end{array}
\right)\,,\nonumber\\
J_{13}^+&=\left(
\begin{array}{ccccccc}
 0 & 0 & -1 & 0 & 0 & 0 & 0 \\
 0 & 1 & 0 & 0 & 0 & 0 & 0 \\
 1 & 0 & 0 & 0 & 0 & 0 & 0 \\
 0 & 0 & 0 & 0 & 0 & 1 & 0 \\
 0 & 0 & 0 & 0 & 1 & 0 & 0 \\
 0 & 0 & 0 & -1 & 0 & 0 & 0 \\
 0 & 0 & 0 & 0 & 0 & 0 & 1 \\
\end{array}
\right)\,,\,\,J_{13}^-=\left(
\begin{array}{ccccccc}
 0 & 0 & -1 & 0 & 0 & 0 & 0 \\
 0 & 1 & 0 & 0 & 0 & 0 & 0 \\
 1 & 0 & 0 & 0 & 0 & 0 & 0 \\
 0 & 0 & 0 & 0 & 0 & -1 & 0 \\
 0 & 0 & 0 & 0 & 1 & 0 & 0 \\
 0 & 0 & 0 & 1 & 0 & 0 & 0 \\
 0 & 0 & 0 & 0 & 0 & 0 & 1 \\
\end{array}
\right)\,,\nonumber\\
J_{23}^+&=\left(
\begin{array}{ccccccc}
 1 & 0 & 0 & 0 & 0 & 0 & 0 \\
 0 & 0 & -1 & 0 & 0 & 0 & 0 \\
 0 & 1 & 0 & 0 & 0 & 0 & 0 \\
 0 & 0 & 0 & 1 & 0 & 0 & 0 \\
 0 & 0 & 0 & 0 & 0 & 1 & 0 \\
 0 & 0 & 0 & 0 & -1 & 0 & 0 \\
 0 & 0 & 0 & 0 & 0 & 0 & 1 \\
\end{array}
\right)\,,\,\,J_{23}^-=\left(
\begin{array}{ccccccc}
 1 & 0 & 0 & 0 & 0 & 0 & 0 \\
 0 & 0 & -1 & 0 & 0 & 0 & 0 \\
 0 & 1 & 0 & 0 & 0 & 0 & 0 \\
 0 & 0 & 0 & 1 & 0 & 0 & 0 \\
 0 & 0 & 0 & 0 & 0 & -1 & 0 \\
 0 & 0 & 0 & 0 & 1 & 0 & 0 \\
 0 & 0 & 0 & 0 & 0 & 0 & 1 \\
\end{array}
\right)\,,\nonumber\\
J_1^{I=1}&=\left(
\begin{array}{ccccccc}
 1 & 0 & 0 & 0 & 0 & 0 & 0 \\
 0 & 1 & 0 & 0 & 0 & 0 & 0 \\
 0 & 0 & 1 & 0 & 0 & 0 & 0 \\
 0 & 0 & 0 & -1 & 0 & 0 & 0 \\
 0 & 0 & 0 & 0 & 1 & 0 & 0 \\
 0 & 0 & 0 & 0 & 0 & 1 & 0 \\
 0 & 0 & 0 & 0 & 0 & 0 & -1 \\
\end{array}
\right)\,,\,\,J_2^{I=1}=\left(
\begin{array}{ccccccc}
 1 & 0 & 0 & 0 & 0 & 0 & 0 \\
 0 & 1 & 0 & 0 & 0 & 0 & 0 \\
 0 & 0 & 1 & 0 & 0 & 0 & 0 \\
 0 & 0 & 0 & 1 & 0 & 0 & 0 \\
 0 & 0 & 0 & 0 & -1 & 0 & 0 \\
 0 & 0 & 0 & 0 & 0 & 1 & 0 \\
 0 & 0 & 0 & 0 & 0 & 0 & -1 \\
\end{array}
\right)\,,\nonumber\\
J_3{}^{I=1}&=\left(
\begin{array}{ccccccc}
 1 & 0 & 0 & 0 & 0 & 0 & 0 \\
 0 & 1 & 0 & 0 & 0 & 0 & 0 \\
 0 & 0 & 1 & 0 & 0 & 0 & 0 \\
 0 & 0 & 0 & 1 & 0 & 0 & 0 \\
 0 & 0 & 0 & 0 & 1 & 0 & 0 \\
 0 & 0 & 0 & 0 & 0 & -1 & 0 \\
 0 & 0 & 0 & 0 & 0 & 0 & -1 \\
\end{array}
\right)\,.\nonumber
\end{align}
\begin{align}
T_{12}^-&=\left(
\begin{array}{ccccccc}
 1 & -2 & 0 & 0 & 2 & 0 & 0 \\
 2 & 1 & 0 & 2 & 0 & 0 & 0 \\
 0 & 0 & 1 & 0 & 0 & 0 & 0 \\
 0 & 2 & 0 & 1 & -2 & 0 & 0 \\
 2 & 0 & 0 & 2 & 1 & 0 & 0 \\
 0 & 0 & 0 & 0 & 0 & 1 & 0 \\
 0 & 0 & 0 & 0 & 0 & 0 & 1 \\
\end{array}
\right)\,,\,\,T_{12}^+=\left(
\begin{array}{ccccccc}
 1 & -2 & 0 & 0 & -2 & 0 & 0 \\
 2 & 1 & 0 & 2 & 0 & 0 & 0 \\
 0 & 0 & 1 & 0 & 0 & 0 & 0 \\
 0 & 2 & 0 & 1 & 2 & 0 & 0 \\
 -2 & 0 & 0 & -2 & 1 & 0 & 0 \\
 0 & 0 & 0 & 0 & 0 & 1 & 0 \\
 0 & 0 & 0 & 0 & 0 & 0 & 1 \\
\end{array}
\right)\,\nonumber\\
T_{13}^-&=\left(
\begin{array}{ccccccc}
 1 & 0 & -2 & 0 & 0 & 2 & 0 \\
 0 & 1 & 0 & 0 & 0 & 0 & 0 \\
 2 & 0 & 1 & 2 & 0 & 0 & 0 \\
 0 & 0 & 2 & 1 & 0 & -2 & 0 \\
 0 & 0 & 0 & 0 & 1 & 0 & 0 \\
 2 & 0 & 0 & 2 & 0 & 1 & 0 \\
 0 & 0 & 0 & 0 & 0 & 0 & 1 \\
\end{array}
\right)\,,\,\,T_{13}^+=\left(
\begin{array}{ccccccc}
 1 & 0 & -2 & 0 & 0 & -2 & 0 \\
 0 & 1 & 0 & 0 & 0 & 0 & 0 \\
 2 & 0 & 1 & 2 & 0 & 0 & 0 \\
 0 & 0 & 2 & 1 & 0 & 2 & 0 \\
 0 & 0 & 0 & 0 & 1 & 0 & 0 \\
 -2 & 0 & 0 & -2 & 0 & 1 & 0 \\
 0 & 0 & 0 & 0 & 0 & 0 & 1 \\
\end{array}
\right)\,,\nonumber\\
T_{23}^-&=\left(
\begin{array}{ccccccc}
 1 & 0 & 0 & 0 & 0 & 0 & 0 \\
 0 & 1 & -2 & 0 & 0 & 2 & 0 \\
 0 & 2 & 1 & 0 & 2 & 0 & 0 \\
 0 & 0 & 0 & 1 & 0 & 0 & 0 \\
 0 & 0 & 2 & 0 & 1 & -2 & 0 \\
 0 & 2 & 0 & 0 & 2 & 1 & 0 \\
 0 & 0 & 0 & 0 & 0 & 0 & 1 \\
\end{array}
\right)\,,\,\,T_{23}^+=\left(
\begin{array}{ccccccc}
 1 & 0 & 0 & 0 & 0 & 0 & 0 \\
 0 & 1 & -2 & 0 & 0 & -2 & 0 \\
 0 & 2 & 1 & 0 & 2 & 0 & 0 \\
 0 & 0 & 0 & 1 & 0 & 0 & 0 \\
 0 & 0 & 2 & 0 & 1 & 2 & 0 \\
 0 & -2 & 0 & 0 & -2 & 1 & 0 \\
 0 & 0 & 0 & 0 & 0 & 0 & 1 \\
\end{array}
\right)\,,\nonumber\\
T_1^{I=1}&=\left(
\begin{array}{ccccccc}
 3 & 0 & 0 & 2 & 0 & 0 & 2 \\
 0 & 1 & 0 & 0 & 0 & 0 & 0 \\
 0 & 0 & 1 & 0 & 0 & 0 & 0 \\
 -2 & 0 & 0 & -1 & 0 & 0 & -2 \\
 0 & 0 & 0 & 0 & 1 & 0 & 0 \\
 0 & 0 & 0 & 0 & 0 & 1 & 0 \\
 2 & 0 & 0 & 2 & 0 & 0 & 1 \\
\end{array}
\right)\,,\,\,T_2^{I=1}=\left(
\begin{array}{ccccccc}
 1 & 0 & 0 & 0 & 0 & 0 & 0 \\
 0 & 3 & 0 & 0 & 2 & 0 & 2 \\
 0 & 0 & 1 & 0 & 0 & 0 & 0 \\
 0 & 0 & 0 & 1 & 0 & 0 & 0 \\
 0 & -2 & 0 & 0 & -1 & 0 & -2 \\
 0 & 0 & 0 & 0 & 0 & 1 & 0 \\
 0 & 2 & 0 & 0 & 2 & 0 & 1 \\
\end{array}
\right)\,,\nonumber\\
T_3^{I=1}&=\left(
\begin{array}{ccccccc}
 1 & 0 & 0 & 0 & 0 & 0 & 0 \\
 0 & 1 & 0 & 0 & 0 & 0 & 0 \\
 0 & 0 & 3 & 0 & 0 & 2 & 2 \\
 0 & 0 & 0 & 1 & 0 & 0 & 0 \\
 0 & 0 & 0 & 0 & 1 & 0 & 0 \\
 0 & 0 & -2 & 0 & 0 & -1 & -2 \\
 0 & 0 & 2 & 0 & 0 & 2 & 1 \\
\end{array}
\right)\,.
\end{align}
%%%%%%%%%%%%%%%%%%%%%%%%%%%%%%%%%%%%%%%%%%%%%%%%%%%%%%%%%%%%%%%%%

\end{document}